\documentclass[3p,times]{elsarticle}
\usepackage{epstopdf}
\usepackage{bbm}
\usepackage{bm}
\usepackage{array} 
\usepackage{color}
\usepackage{xcolor}

\usepackage{graphicx}
\usepackage{psfrag}

\usepackage{amssymb}
\usepackage{amsmath}
\usepackage{dsfont}
\usepackage{rotating}
\usepackage{pgf,tikz}
\usetikzlibrary{arrows}

\definecolor{ttzzqq}{rgb}{0.2,0.6,0}

\definecolor{qqttcc}{rgb}{0,0.2,0.8}

\definecolor{qqttzz}{rgb}{0,0.2,0.6}
\definecolor{ffqqqq}{rgb}{1,0,0}
\definecolor{qqwuqq}{rgb}{0,0.39,0}
\definecolor{zzttqq}{rgb}{0.6,0.2,0}
\definecolor{qqqqff}{rgb}{0,0,1}
\definecolor{ttttqq}{rgb}{0.2,0.2,0}

\definecolor{qqwwtt}{rgb}{0,0.4,0.2}
\definecolor{ubqqys}{rgb}{0.29,0,0.51}

\definecolor{wwttqq}{rgb}{0.4,0.2,0}
\definecolor{uuuuuu}{rgb}{0.27,0.27,0.27}
\definecolor{qqzzff}{rgb}{0,0.6,1}
\definecolor{xdxdff}{rgb}{0.49,0.49,1}
\definecolor{ccwwqq}{rgb}{0.8,0.4,0}

\definecolor{ttqqqq}{rgb}{0.2,0,0}

\definecolor{qqzzcc}{rgb}{0,0.6,0.8}

\newcommand{\R}{\mathbb{R}}

\newcommand{\x}{\chi}

\newcommand{\eref}[1]{$(\ref{#1})$}
\newcommand{\p}{\mathbb{\wp}}
\newcommand{\np}{p}

\newcommand{\nstd}{\vec{n}_{j}}
\newcommand{\etah}{ \hat{\bm{p}}}
\newcommand{\bphi}{ \bm{\phi}}
\newcommand{\bpsi}{ \bm{\psi}}

\newcommand{\nv}{\vec{n}}
\newcommand{\B}{\mathcal{B}}
\newcommand{\TF}{\mathbf{G}}

\newcommand{\st}{ {st}}

\newcommand{\Mpsi}{\bm{M}}

\newcommand{\idm}{\mathbf{I}}														
\newcommand{\vv}{{\mathbf{v}}}											
 
\newcommand{\rhov}{\rho {\mathbf{v}}}		
\newcommand{\rhoEE}{\rho E}      				
\newcommand{\rhoe}{\rho {e}}
\newcommand{\rhok}{\rho {k}}
\newcommand{\qh}{ \hat{\bm{Q}}}
\newcommand{\xx}{ \mathbf{x}}
\newcommand{\TT}{\mbox{\boldmath$T$}}
\newcommand{\QQ}{\mbox{\boldmath$R$}}

\newcommand{\Ni}{ N_i}
\newcommand{\Nj}{ N_j}
\newcommand{\phit}{ \tilde{\phi}}
\newcommand{\psit}{ \tilde{\psi}}
\newcommand{\dx}{ d\xx}
\newcommand{\dxt}{ d\xx \, dt \,}
\newcommand{\stens}{{\boldsymbol{\sigma}}}											

\newcommand{\D}{\bm{\mathcal{D}}}

\newcommand{\Dtilde}{{\bm{\mathcal{E}}}}
\newcommand{\LL}{ \bm{\mathcal{M}}}

\newcommand{\Q}{\bm{\mathcal{Q}}}
\newcommand{\RM}{\bm{\mathcal{R}}}
\newcommand{\LM}{\bm{\mathcal{L}}}
\newcommand{\M}{\bm{M}}
\newcommand{\MM}{\mathcal{\bm{M}}}

\newcommand{\ie}{ e}
\newcommand{\pph}{ \hat{\bm{p}}}

\newcommand{\rhoeh}{ \widehat{\bm{\rho \ie}}}
\newcommand{\rhokh}{ \widehat{\bm{\rho k}}}
\newcommand{\rhoEEh}{ \widehat{\bm{\rho E}}}
\newcommand{\rhovh}{ \widehat{\bm{\rho \mathbf{v}}}}

\newcommand{\hh}{ \hat{\mathbf{H}}}

\newcommand{\dtm}{ {d \rightarrow m} }
\newcommand{\mtd}{ {m \rightarrow d} }

\newcommand{\diff}[2]{\frac{\partial {#1} }{\partial {#2} } }

\usepackage{amssymb}

\journal{Journal of Computational Physics}

\begin{document}

\begin{frontmatter}



\title{A pressure-based semi-implicit space-time discontinuous Galerkin method on staggered unstructured meshes for the solution of the compressible Navier-Stokes equations at all Mach numbers }   

\author[1]{Maurizio Tavelli\fnref{label1}}
\author[1]{Michael Dumbser \corref{corr1} \fnref{label2}}
\address[1]{Department of Civil, Environmental and Mechanical Engineering, University of Trento, Via Mesiano 77, I-38123 Trento, Italy}

\fntext[label1]{\tt m.tavelli@unitn.it (M.~Tavelli)}
\fntext[label2]{\tt michael.dumbser@unitn.it (M.~Dumbser)}
\begin{abstract}
We propose a new arbitrary high order accurate semi-implicit space-time discontinuous Galerkin (DG) method for the solution of the two and three dimensional compressible Euler and Navier-Stokes 
equations on \textit{staggered} unstructured curved meshes. The method is pressure-based and semi-implicit and is able to deal with all Mach number flows.  

The new DG scheme extends the seminal ideas outlined in \cite{DumbserCasulli2016}, where a second order semi-implicit finite volume method for the solution of the compressible Navier-Stokes 
equations with a general equation of state was introduced on staggered Cartesian grids. Regarding the high order extension we follow \cite{3DSIINS}, where
a staggered space-time DG scheme for the incompressible Navier-Stokes equations was presented. 
In our scheme, the discrete pressure is defined on the primal grid, while the discrete velocity field and the density are defined on a face-based staggered dual grid. 
Then, the mass conservation equation, as well as the nonlinear convective terms in the momentum equation and the transport of kinetic energy in the energy equation are discretized
explicitly, while the pressure terms appearing in the momentum and energy equation are discretized implicitly. 
Formal substitution of the discrete momentum equation into the total energy conservation equation yields a linear system for only one unknown, namely the \textit{scalar} pressure. 
Here the equation of state is assumed linear with respect to the pressure. The enthalpy and the kinetic energy are taken explicitly and are then updated using a simple Picard procedure.
Thanks to the use of a staggered grid, the final pressure system is a very sparse block five-point system for three dimensional problems and it is a block four-point system in 
the two dimensional case. Furthermore, for high order in space and piecewise constant polynomials in time, the system is observed to be symmetric and positive definite. 
This allows to use fast linear solvers such as the conjugate gradient (CG) method. 
In addition, all the volume and surface integrals needed by the scheme depend only on the geometry and the polynomial degree of the basis and test functions and can therefore be precomputed and stored in a preprocessing stage. This leads to significant savings in terms of computational effort for the time evolution part. In this way also the extension to a fully curved isoparametric approach becomes natural and affects only the preprocessing step. The viscous terms and the heat flux are also discretized making use of the staggered grid by defining the 
viscous stress tensor and the heat flux vector on the dual grid, which corresponds to the use of a lifting operator, but on the dual grid. 
The time step of our new numerical method is limited by a CFL condition based only on the fluid velocity and not on the 
sound speed. This makes the method particularly interesting for low Mach number flows. 
Finally, a very simple combination of artificial viscosity and the \textit{a posteriori} MOOD technique allows to deal with shock waves and thus permits also to simulate high Mach number flows. 
We show computational results for a large set of two and three-dimensional benchmark problems, including both low and high Mach number flows and 
using polynomial approximation degrees up to $p=4$.  
\end{abstract}


\end{frontmatter}


\section{Introduction}
Computational fluid mechanics is a very important field for a wide set of applications that ranges from aerospace and mechanical engineering, energy production at the aid of gas, wind and water turbines over geophysical flows in oceans, rivers, lakes as well as atmospheric flows to blood flow in the human cardiovascular system. Although the field of application is extremely large, there exists one  universally accepted mathematical model of governing equations that can describe fluid flow in all the above mentioned circumstances. It can be derived from the conservation of mass, momentum and  total energy and is given by the well-known compressible Navier-Stokes equations. In their complete form they can describe a wide range of phenomena, including also the effects of momentum transport via molecular viscosity and heat conduction. 
The compressible Navier-Stokes equations also comprehend several simplified sub systems, like the compressible Euler equations in the inviscid case, or the incompressible Navier-Stokes equations as the zero  Mach number limit, where the Mach number is defined as usual by the ratio between the fluid velocity and the sound speed, see e.g. \cite{Klainermann81,Klainermann82,Munz2007}. Furthermore, from the incompressible Navier-Stokes equations the so-called shallow water equations can be derived by integration over the depth and by assuming a hydrostatic pressure. 
In this sense the applications of the compressible Navier-Stokes equations split into two main classes: low Mach number flows, typical for geophysical, environmental and biological applications, 
and high Mach number flows that are typical of industrial applications such as in aerospace and mechanical engineering. For the high Mach number case, the families of explicit 
\textit{density-based}  upwind finite difference and Godunov-type finite volume schemes are very popular, see for example 
\cite{Einfeldt1991,godunov,Harten:1983b,lax,LeVeque:2002a,munz94,osherandsolomon,roe,toro-book,Toro:1994}. 
Due to the elliptic behavior of the pressure in the incompressible limit, the use of purely explicit schemes introduces a very severe restriction on the maximum time step for low Mach number 
flows, since the CFL condition of explicit methods includes also the sound speed. 
This explains why semi-implicit \textit{pressure-based} schemes are more popular in this class of applications. In the past several 
semi-implicit numerical schemes have been developed for the incompressible case, see \cite{markerandcell,chorin1,chorin2,HirtNichols,Bell1989,CasulliCheng1992,patankar,vanKan,Casulli2014} 
and have been recently extended also to the new family of staggered semi-implicit discontinuous Galerkin schemes in \cite{DumbserCasulli,SINS,2STINS,3DSIINS}. Regarding semi-implicit schemes 
for the compressible case on staggered and collocated grids we refer the reader for example to the work presented in \cite{Casulli1984,Munz2003,klein,KleinMach,ParkMPV,Dolejsi1,Dolejsi2,Dolejsi3}. 

Very recently, a new weakly nonlinear semi-implicit finite volume scheme for the solution of the compressible Navier-Stokes and Euler equations with general equation of state (EOS) was presented by Dumbser and Casulli in \cite{DumbserCasulli2016}.

The aim of the present paper is to extend those ideas to higher order of accuracy in space and time and thus to develop a novel pressure-based semi-implicit staggered DG scheme for the solution of the  compressible Navier-Stokes and Euler equations in multiple space dimensions that is globally and locally conservative for mass momentum and total energy and that involves a CFL time step restriction that is only based on the local  flow velocity and not on the sound speed. Furthermore, we derive the method on a general unstructured curved staggered grid in order to fit also complex geometries. While the pressure is defined on a 
main tetrahedral (respectively triangular) grid, the density and the velocity fields are defined on a face-based staggered dual grid (respectively edge-based dual grid). The formal substitution of the  discrete momentum equation into the energy equation leads to a linear system for only one single unknown, namely the scalar pressure. 
The discrete form of the equations looks very similar to the high order staggered DG scheme proposed in \cite{3DSIINS} for the incompressible Navier-Stokes equations. 
Numerical evidence shows that the good properties of the resulting linear system for the pressure can be maintained also in the compressible framework, while in the incompressible case 
there is a rigorous analysis available, see \cite{SINS,2STINS,Capizzano2016}. The use of a staggered mesh as well as the semi-implicit solution strategy applied to the resulting discrete equations 
makes our new scheme totally different from the space-time DG schemes presented in \cite{spacetimedg1,spacetimedg2,KlaijVanDerVegt,Rhebergen2013}. 

High Mach number flows typically lead to shock waves where unlimited high order numerical schemes produce spurious oscillations - the well-known Gibbs phenomenon - which can also lead to 
unphysical quantities, i.e. negative values for pressure or density. It is then necessary to introduce a limiter in order to overcome the problem. Several types of limiting procedures have
been introduced in the past, such as WENO limiters \cite{QiuShu1,balsara2007,Zhu2008}, slope and moment limiting \cite{Krivodonova2007,cockburn1989, Biswas1994,Burbeau2001} and artificial 
viscosity (AV), see \cite{spacetimedg1,Hartman2002,Persson_06,Cesenek2013,Feistauer2007,Feistauer2010}. The concept of AV was already introduced in the 1950ies by Von Neumann and Richtmyer 
\cite{Neumann1950}, in order to deal with shock waves and high Mach number flows. 

A completely new way of limiting high order DG schemes was very recently developed by Dumbser et. al. in 
\cite{Dumbser2014,Zanotti2015,DGLimiter3}, where a novel \textit{a posteriori} sub-cell finite volume limiter was used to suppress spurious oscillations of the DG polynomials in the 
vicinity of shock waves and other discontinuities. Originally, the idea to use an \textit{a posteriori} approach was introduced by Clain, Diot and Loub\`ere 
\cite{Clain2011,Diot2012,Diot2013,Loubere2014} in the finite volume context with the so-called Multi-dimensional Optimal Order Detection (MOOD) method. For alternative 
\textit{a priori} subcell DG limiters, see the work of Sonntag and Munz \cite{Sonntag} and others \cite{CasoniHuerta1,CasoniHuerta2,MeisterOrtleb}. 
In this paper we propose for the first time to use the \textit{a posteriori} MOOD paradigm adapted to semi-implicit methods on staggered grids in combination with artificial viscosity. 

The rest of the paper is organized as follows: in Section \ref{sec_1} we present some basic definitions about the staggered meshes and the polynomial spaces used. 
In Section \ref{sec_semi_imp_dg} we derive the numerical method, discussing in particular the discretization of the nonlinear convective terms and of the viscous stress tensor 
in the momentum equation and of the heat flux in the energy equation. We also present the implementation of the new DG limiter based on the MOOD approach combined with artificial 
viscosity. 
In Section \ref{sec_NT} we show the numerical results obtained for a large set of two- and three-dimensional benchmark problems that cover a large range of Mach numbers. 
Finally, in Section \ref{sec.concl} we draw some conclusions and present an outlook to future research.

\section{Staggered space-time DG scheme for the compressible Navier-Stokes equations}
\label{sec_1}
\subsection{Governing equations}
The compressible Navier-Stokes equations can be derived by considering the conservation of mass, momentum and energy. In a compact vectorial form, the system reads 
\begin{eqnarray}
    \diff{\rho}{t}  +\nabla \cdot \left( \rhov \right) &=& 0 \label{eq:CNS_1}, \\ 
    \diff{\rhov}{t} +\nabla \cdot \left( \rho \vv \otimes \vv \right) + \nabla p   &=& 
		\nabla \cdot \stens 
		\label{eq:CNS_2}, \\
		\diff{\rhoEE}{t} + \nabla \cdot \left( \vv ( \rho E + p ) \right) & = & \nabla \cdot \left( \stens \vv + \lambda \nabla T \right) 
		\label{eq:CNS_3},
\end{eqnarray}
where the left hand side consists of the inviscid compressible Euler equations and the right hand side contains the viscous stress tensor and the heat flux. In the above system 
$\rho$ is the fluid density, $\vv$ is the velocity vector, $p$ is the pressure and  
$\rhoEE=\rhoe+\rhok$ is the total energy density; $\rhok = \frac{1}{2}\rho \vv^2$ is the  kinetic energy density and $e=e(p,\rho)$ represents the specific internal energy per unit mass and 
is given by the equation of state (EOS) as a function of the pressure $p$ and the density $\rho$, 
$H = e+\frac{p}{\rho}$ denotes the specific enthalpy, $\mu$ is the dynamic viscosity and $\lambda$ is the thermal conductivity coefficient. 
At the aid of the specific enthalpy $H$ and the specific kinetic energy $k$ the flux on the left hand side 
of the total energy equation can be rewritten as $\vv ( \rho E + p )  =  \rho \vv  ( k + H )$. With the usual Stokes hypothesis, the viscous stress tensor reads 
\begin{eqnarray}
	\stens = \mu (\nabla \vv +\nabla \vv^\top)-\frac{2}{3} \mu \left( \nabla \cdot \vv \right)\idm.
\label{eq:CNS_4}
\end{eqnarray}
The previous system can thus be rewritten as 
\begin{eqnarray}
    \diff{\rho}{t}   + \nabla \cdot \rhov &=& 0 \label{eq:CNS_1_2}, \\ 
    \diff{\rhov}{t}  + \nabla \cdot \mathbf{G} +\nabla p &=& 0 \label{eq:CNS_2_2}, \\
		\diff{\rhoEE}{t} + \nabla \cdot (\rho \vv k + \rho \vv H )&=&\nabla \cdot \left( \stens \vv + \lambda \nabla T \right) \label{eq:CNS_3_2},
\end{eqnarray}
where $\mathbf{G}=\mathbf{F}_{\rho \vv} - \stens$ contains the nonlinear convective terms and the viscous contribution in the momentum equation. 
Throughout this paper we will use the following notation for the convective terms: $\mathbf{F}_{\rho \vv} = \rho \vv \otimes \vv$ is an abbreviation for the nonlinear convective terms  
in the momentum equation, $\mathbf{F}_{\rho} = \rho \vv$ is the convective flux in the mass conservation equation and $\mathbf{F}_{\rho k}=\rho k \vv$ is the convective
flux of the kinetic energy. This type of splitting of the purely convective terms has been recently introduced by Toro and V\'azquez in \cite{ToroVazquez}. 
Furthermore, we will also use the abbreviation $\mathbf{w} = \boldsymbol{\sigma} \vv$ for the work of the viscous stress tensor in the energy equation  
and $\mathbf{q}=\lambda \nabla T$ for the heat flux vector. 

The system has to be closed with a suitable equation of state (EOS). 
Given a so called thermal equation of state $p = p(T,\rho)$ that links the pressure with the temperature and the density, as well as a caloric equation of state $e=e(T,\rho)$ that determines 
how the internal energy changes with temperature and density, we usually eliminate the temperature $T$ and thus obtain a single equation of state $e=e(p,\rho)$ that directly depends on the 
pressure $p$ and the density $\rho$. For an ideal gas, the thermal and the caloric EOS are given by 
\begin{eqnarray}
		p=RT\rho \qquad \mbox{and} \qquad e=c_v T,
\label{eq:CNS_5}
\end{eqnarray}
where $R=c_p-c_v$ is the specific gas constant; $c_v$ and $c_p$ are the heat capacities at constant volume and at constant pressure, respectively. 
From the previous relations one easily derives the equation of state as a function of only the pressure and the fluid density, given by 
\begin{eqnarray}
	e=\frac{p}{(\gamma-1)\rho},
\label{eq:CNS_6}
\end{eqnarray}
with $\gamma=\frac{c_p}{c_v}$ denoting the so-called ratio of specific heats.
For a generic cubic EOS the relation $e=e(p,\rho)$ can become highly nonlinear, see e.g.\cite{Vidal2001,DumbserCasulli2016}. For the methodology derived in the following, we 
allow a general equation of state $e(p,\rho)$, but with the constraint that it remains \textit{linear} with respect to $p$. This is the case, for example, for the ideal gas EOS 
or even for the well-known van der Waals EOS, which is only nonlinear in $\rho$ but not in $p$. 

\subsection{Staggered unstructured grid}
Throughout this paper we use the same unstructured spatially staggered mesh as the one used in \cite{SINS,2STINS,3DSIINS} for the two and three-dimensional case, respectively. 
In the following section we briefly summarize the grid construction and the main notation for the two dimensional triangular grid. After that, the primary and dual spatial elements 
are extended to the three dimensional case and also to the case of space-time control volumes.

\paragraph{Two space dimensions}
In the two-dimensional case the spatial computational domain $\Omega \subset \mathbb{R}^2$ is covered with a set of $\Ni$ non-overlapping triangular elements $\TT_i$ with $i=1 \ldots \Ni$. By denoting with $N_d$ the total number of edges, the $j-$th edge will be called $\Gamma_j$. $\B(\Omega)$ denotes the set of indices $j$ corresponding to boundary edges.
The three edges of each triangle $\TT_i$ constitute the set $S_i$ defined by $S_i=\{j \in [1,\Nj] \,\, | \,\, \Gamma_j \mbox{ is an edge of }\TT_i \}$. For every $j\in [1\ldots \Nj]-\B(\Omega)$ there exist two triangles $i_1$ and $i_2$ that share $\Gamma_j$. We assign arbitrarily a left and a right triangle called $\ell(j)$ and $r(j)$, respectively. The standard positive direction is assumed to be from left to right. Let $\nv_{j}$ denote the unit normal
vector defined on the edge $j$ and oriented with respect to the positive direction from left to right. For every triangular element $i$ and edge $j \in S_i$,
the neighbor triangle of element $\TT_i$ that share the edge $\Gamma_j$ is denoted by $\p(i,j)$.
\par
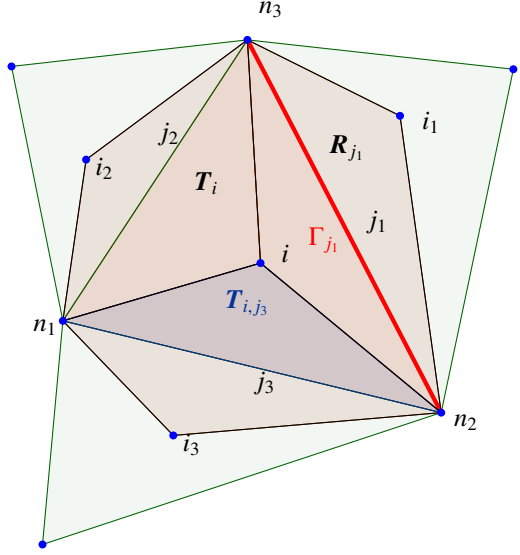
\begin{figure}[ht]
    \begin{center}
    \begin{tikzpicture}[line cap=round,line join=round,>=triangle 45,x=0.6373937677053826cm,y=0.6177884615384613cm]
\clip(2.11,-8.53) rectangle (16.23,3.95);
\fill[color=zzttqq,fill=zzttqq,fill opacity=0.1] (5.19,-3.03) -- (9,3) -- (13,-5) -- cycle;
\fill[color=qqwuqq,fill=qqwuqq,fill opacity=0.05] (9,3) -- (14.49,2.37) -- (13,-5) -- cycle;
\fill[color=qqwuqq,fill=qqwuqq,fill opacity=0.05] (9,3) -- (4.13,2.43) -- (5.19,-3.03) -- cycle;
\fill[color=qqwuqq,fill=qqwuqq,fill opacity=0.05] (5.19,-3.03) -- (4.77,-7.83) -- (13,-5) -- cycle;
\fill[color=zzttqq,fill=zzttqq,fill opacity=0.1] (9.27,-1.79) -- (9,3) -- (5.67,0.43) -- (5.19,-3.03) -- cycle;
\fill[color=zzttqq,fill=zzttqq,fill opacity=0.1] (9.27,-1.79) -- (9,3) -- (12.15,1.37) -- (13,-5) -- cycle;
\fill[color=zzttqq,fill=zzttqq,fill opacity=0.1] (5.19,-3.03) -- (7.47,-5.49) -- (13,-5) -- (9.27,-1.79) -- cycle;
\fill[color=qqttzz,fill=qqttzz,fill opacity=0.1] (13,-5) -- (5.19,-3.03) -- (9.27,-1.79) -- cycle;
\draw [color=zzttqq] (5.19,-3.03)-- (9,3);
\draw [color=zzttqq] (9,3)-- (13,-5);
\draw [color=zzttqq] (13,-5)-- (5.19,-3.03);
\draw [color=qqwuqq] (9,3)-- (14.49,2.37);
\draw [color=qqwuqq] (14.49,2.37)-- (13,-5);
\draw [color=qqwuqq] (13,-5)-- (9,3);
\draw [color=qqwuqq] (9,3)-- (4.13,2.43);
\draw [color=qqwuqq] (4.13,2.43)-- (5.19,-3.03);
\draw [color=qqwuqq] (5.19,-3.03)-- (9,3);
\draw [color=qqwuqq] (5.19,-3.03)-- (4.77,-7.83);
\draw [color=qqwuqq] (4.77,-7.83)-- (13,-5);
\draw [color=qqwuqq] (13,-5)-- (5.19,-3.03);
\draw (9.51,-1.19) node[anchor=north west] {$i$};
\draw (12.41,1.67) node[anchor=north west] {$i_1$};
\draw (5.67,0.67) node[anchor=north west] {$i_2$};
\draw (7.47,-5.25) node[anchor=north west] {$i_3$};
\draw (11.25,-0.47) node[anchor=north west] {$j_1$};
\draw (6.95,1.43) node[anchor=north west] {$j_2$};
\draw (8.93,-3.93) node[anchor=north west] {$j_3$};
\draw (4.39,-2.79) node[anchor=north west] {$n_1$};
\draw (13.07,-4.83) node[anchor=north west] {$n_2$};
\draw (9.05,4.03) node[anchor=north west] {$n_3$};
\draw (7.71,0.35) node[anchor=north west] {$\TT_i$};
\draw [color=zzttqq] (9.27,-1.79)-- (9,3);
\draw [color=zzttqq] (9,3)-- (5.67,0.43);
\draw [color=zzttqq] (5.67,0.43)-- (5.19,-3.03);
\draw [color=zzttqq] (5.19,-3.03)-- (9.27,-1.79);
\draw [color=zzttqq] (9.27,-1.79)-- (9,3);
\draw [color=zzttqq] (9,3)-- (12.15,1.37);
\draw [color=zzttqq] (12.15,1.37)-- (13,-5);
\draw [color=zzttqq] (13,-5)-- (9.27,-1.79);
\draw [color=zzttqq] (5.19,-3.03)-- (7.47,-5.49);
\draw [color=zzttqq] (7.47,-5.49)-- (13,-5);
\draw [color=zzttqq] (13,-5)-- (9.27,-1.79);
\draw [color=zzttqq] (9.27,-1.79)-- (5.19,-3.03);
\draw (10.49,1.15) node[anchor=north west] {$\QQ_{j_1}$};
\draw [color=ffqqqq](10.07,-0.83) node[anchor=north west] {$\Gamma_{j_1}$};
\draw [line width=1.6pt,color=ffqqqq] (9,3)-- (13,-5);
\draw [color=qqttzz] (13,-5)-- (5.19,-3.03);
\draw [color=qqttzz] (5.19,-3.03)-- (9.27,-1.79);
\draw [color=qqttzz] (9.27,-1.79)-- (13,-5);
\draw [color=qqttzz](8.35,-2.17) node[anchor=north west] {$\TT_{i,j_3}$};
\draw (5.19,-3.03)-- (5.67,0.43);
\draw (5.67,0.43)-- (9,3);
\draw (9,3)-- (9.27,-1.79);
\draw (9.27,-1.79)-- (5.19,-3.03);
\draw (9,3)-- (12.15,1.37);
\draw (12.15,1.37)-- (13,-5);
\draw (13,-5)-- (9.27,-1.79);
\draw (13,-5)-- (7.47,-5.49);
\draw (7.47,-5.49)-- (5.19,-3.03);
\begin{scriptsize}
\fill [color=qqqqff] (5.19,-3.03) circle (1.5pt);
\fill [color=qqqqff] (9,3) circle (1.5pt);
\fill [color=qqqqff] (13,-5) circle (1.5pt);
\fill [color=qqqqff] (9.27,-1.79) circle (1.5pt);
\fill [color=qqqqff] (14.49,2.37) circle (1.5pt);
\fill [color=qqqqff] (4.13,2.43) circle (1.5pt);
\fill [color=qqqqff] (4.77,-7.83) circle (1.5pt);
\fill [color=qqqqff] (12.15,1.37) circle (1.5pt);
\fill [color=qqqqff] (5.67,0.43) circle (1.5pt);
\fill [color=qqqqff] (7.47,-5.49) circle (1.5pt);
\end{scriptsize}
\end{tikzpicture}
    \caption{Example of a triangular mesh element with its three neighbors and the associated staggered edge-based dual control volumes, together with the notation
    used throughout the paper.}
    \label{fig.1}
		\end{center}
\end{figure}

 For every $j\in [1, \Nj]-\B(\Omega)$ the quadrilateral element associated to $j$ is called $\QQ_j$ and it is defined, in general, by the two centers of gravity of $\ell(j)$ and $r(j)$ and the two terminal nodes of $\Gamma_j$, see also \cite{Bermudez1998,USFORCE,2DSIUSW,StaggeredDG}. We denote by $\TT_{i,j}=\QQ_j \cap \TT_i$ the intersection element for every $i$ and $j \in S_i$. Figure $\ref{fig.1}$ summarizes the used notation, the primal triangular mesh and the dual quadrilateral grid. 
According to \cite{2STINS}, we will call the mesh of triangular elements $\{\TT_i \}_{i \in [1, \Ni]}$ the \textit{main grid} or \textit{primal grid} and the quadrilateral grid $\{\QQ_j \}_{j \in [1, N_d]}$ is termed the \textit{dual grid}. 


\paragraph{Three space dimensions}
The definitions given above are then readily extended to three space dimensions with the domain $\Omega \subset \mathbb{R}^3$. 
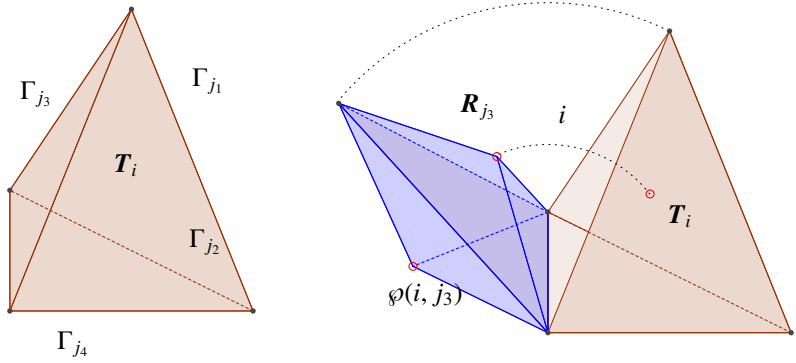
\begin{figure}[h!t]
    \begin{center}
    \begin{tikzpicture}[line cap=round,line join=round,>=triangle 45,x=0.8cm,y=0.8cm]
\clip(3.26,-5.21) rectangle (8.37,1.57);
\fill[color=zzttqq,fill=zzttqq,fill opacity=0.1] (4,-4) -- (8,-4) -- (4,-2) -- cycle;
\fill[color=zzttqq,fill=zzttqq,fill opacity=0.1] (4,-2) -- (6,1) -- (8,-4) -- cycle;
\fill[color=zzttqq,fill=zzttqq,fill opacity=0.1] (4,-2) -- (6,1) -- (4,-4) -- cycle;
\fill[color=zzttqq,fill=zzttqq,fill opacity=0.1] (6,1) -- (4,-4) -- (8,-4) -- cycle;
\fill[color=zzttqq,fill=zzttqq,fill opacity=0.1] (16,-4) -- (20,-4) -- (16,-2) -- cycle;
\fill[color=zzttqq,fill=zzttqq,fill opacity=0.1] (16,-2) -- (18,1) -- (20,-4) -- cycle;
\fill[color=zzttqq,fill=zzttqq,fill opacity=0.1] (18,1) -- (16,-4) -- (20,-4) -- cycle;
\fill[color=zzttqq,fill=zzttqq,fill opacity=0.1] (16,-2) -- (12.56,-0.2) -- (16,-4) -- cycle;
\fill[color=qqqqff,fill=qqqqff,fill opacity=0.1] (12.56,-0.2) -- (15.16,-1.08) -- (16,-4) -- cycle;
\fill[color=qqqqff,fill=qqqqff,fill opacity=0.1] (15.16,-1.08) -- (16,-2) -- (12.56,-0.2) -- cycle;
\fill[color=qqqqff,fill=qqqqff,fill opacity=0.1] (15.16,-1.08) -- (16,-4) -- (16,-2) -- cycle;
\fill[color=qqqqff,fill=qqqqff,fill opacity=0.1] (12.56,-0.2) -- (13.78,-2.9) -- (16,-4) -- cycle;
\fill[color=qqqqff,fill=qqqqff,fill opacity=0.1] (13.78,-2.9) -- (16,-2) -- (12.56,-0.2) -- cycle;
\fill[color=qqqqff,fill=qqqqff,fill opacity=0.1] (13.78,-2.9) -- (16,-4) -- (16,-2) -- cycle;
\fill[color=zzttqq,fill=zzttqq,fill opacity=0.15] (26,-6) -- (28,-3) -- (30,-8) -- cycle;
\draw [color=zzttqq] (4,-4)-- (8,-4);
\draw [color=zzttqq] (4,-2)-- (4,-4);
\draw [color=zzttqq] (4,-2)-- (6,1);
\draw [color=zzttqq] (6,1)-- (8,-4);
\draw [dash pattern=on 1pt off 1pt,color=zzttqq] (8,-4)-- (4,-2);
\draw [color=zzttqq] (4,-2)-- (6,1);
\draw [color=zzttqq] (6,1)-- (4,-4);
\draw [color=zzttqq] (4,-4)-- (4,-2);
\draw [color=zzttqq] (6,1)-- (4,-4);
\draw [color=zzttqq] (4,-4)-- (8,-4);
\draw [color=zzttqq] (8,-4)-- (6,1);
\draw [color=zzttqq] (16,-4)-- (20,-4);
\draw [color=zzttqq] (16,-2)-- (16,-4);
\draw [color=zzttqq] (16,-2)-- (18,1);
\draw [color=zzttqq] (18,1)-- (20,-4);
\draw [dash pattern=on 1pt off 1pt,color=zzttqq] (20,-4)-- (16,-2);
\draw [color=zzttqq] (18,1)-- (16,-4);
\draw [color=zzttqq] (16,-4)-- (20,-4);
\draw [color=zzttqq] (20,-4)-- (18,1);
\draw [color=zzttqq] (16,-2)-- (16.67,-2.33);
\draw [shift={(16,-3)},dotted]  plot[domain=1.11:2.46,variable=\t]({1*4.47*cos(\t r)+0*4.47*sin(\t r)},{0*4.47*cos(\t r)+1*4.47*sin(\t r)});
\draw [color=zzttqq] (12.56,-0.2)-- (16,-4);
\draw [color=zzttqq] (16,-4)-- (16,-2);
\draw [shift={(16,-3)},dotted]  plot[domain=0.66:1.98,variable=\t]({1*2.12*cos(\t r)+0*2.12*sin(\t r)},{0*2.12*cos(\t r)+1*2.12*sin(\t r)});
\draw [color=qqqqff] (12.56,-0.2)-- (15.16,-1.08);
\draw [color=qqqqff] (15.16,-1.08)-- (16,-4);
\draw [color=qqqqff] (16,-4)-- (12.56,-0.2);
\draw [color=qqqqff] (15.16,-1.08)-- (16,-2);
\draw [color=qqqqff] (12.56,-0.2)-- (15.16,-1.08);
\draw [color=qqqqff] (15.16,-1.08)-- (16,-4);
\draw [color=qqqqff] (16,-4)-- (16,-2);
\draw [color=qqqqff] (16,-2)-- (15.16,-1.08);
\draw [color=qqqqff] (12.56,-0.2)-- (13.78,-2.9);
\draw [color=qqqqff] (13.78,-2.9)-- (16,-4);
\draw [color=qqqqff] (16,-4)-- (12.56,-0.2);
\draw [dash pattern=on 1pt off 1pt,color=qqqqff] (16,-2)-- (12.56,-0.2);
\draw [color=qqqqff] (12.56,-0.2)-- (13.78,-2.9);
\draw [color=qqqqff] (13.78,-2.9)-- (16,-4);
\draw [color=qqqqff] (16,-4)-- (16,-2);
\draw [dash pattern=on 1pt off 1pt,color=qqqqff] (16,-2)-- (13.78,-2.9);
\draw [color=zzttqq] (26,-8)-- (30,-8);
\draw [color=zzttqq] (26,-6)-- (26,-8);
\draw [color=zzttqq] (26,-6)-- (28,-3);
\draw [color=zzttqq] (28,-3)-- (30,-8);
\draw [color=zzttqq] (30,-8)-- (26,-6);
\draw [color=zzttqq] (26,-6)-- (28,-3);
\draw [color=zzttqq] (28,-3)-- (26,-8);
\draw [color=zzttqq] (26,-8)-- (26,-6);
\draw [color=zzttqq] (28,-3)-- (26,-8);
\draw [color=zzttqq] (26,-8)-- (30,-8);
\draw [color=zzttqq] (30,-8)-- (28,-3);
\draw [dash pattern=on 1pt off 1pt,color=qqqqff] (30,-4)-- (26,-6);
\draw [color=qqqqff] (30,-8)-- (30,-4);
\draw [color=qqqqff] (28,-3)-- (30,-4);
\draw [color=qqqqff] (30,-4)-- (30,-8);
\draw [color=qqqqff] (30,-8)-- (28,-3);
\draw [color=qqqqff] (26,-6)-- (28,-3);
\draw [color=qqqqff] (28,-3)-- (30,-4);
\draw [dash pattern=on 1pt off 1pt,color=ffqqqq] (28.46,-5.11)-- (27.23,-5.8);
\draw [color=qqqqff] (28.67,-4.68)-- (30,-4);
\draw [color=ffqqqq] (27.23,-5.8)-- (27.86,-5.45);
\draw (5.56,-1.25) node[anchor=north west] {$\TT_i$};
\draw (6.85,0.18) node[anchor=north west] {$\Gamma_{j_1}$};
\draw (6.8,-2.49) node[anchor=north west] {$\Gamma_{j_2}$};
\draw (4.02,-0.04) node[anchor=north west] {$\Gamma_{j_3}$};
\draw (4.63,-4.18) node[anchor=north west] {$\Gamma_{j_4}$};
\draw (14.4,0.1) node[anchor=north west] {$\QQ_j$};
\draw (16.02,-0.08) node[anchor=north west] {$i$};
\draw (17.82,-1.74) node[anchor=north west] {$\TT_i$};
\draw (13.21,-3.03) node[anchor=north west] {$\mathbb{\wp}(i,j)$};
\draw [->,dash pattern=on 1pt off 1pt] (28.02,-5.62) -- (28.66,-5.29);
\draw (27.13,-5.78) node[anchor=north west] {$\ell(j)$};
\draw (28.99,-4.46) node[anchor=north west] {$r(j)$};
\draw (28.36,-5.41) node[anchor=north west] {$\vec{n_j}$};
\draw (30.3,-7.62) node[anchor=north west] {$\Gamma_j$};
\begin{scriptsize}
\fill [color=uuuuuu] (4,-4) circle (1.0pt);
\fill [color=uuuuuu] (8,-4) circle (1.0pt);
\fill [color=uuuuuu] (4,-2) circle (1.0pt);
\fill [color=uuuuuu] (6,1) circle (1.0pt);
\fill [color=uuuuuu] (16,-4) circle (1.0pt);
\fill [color=uuuuuu] (20,-4) circle (1.0pt);
\fill [color=uuuuuu] (16,-2) circle (1.0pt);
\fill [color=uuuuuu] (18,1) circle (1.0pt);
\fill [color=uuuuuu] (12.56,-0.2) circle (1.0pt);
\draw [color=ffqqqq] (17.68,-1.7) circle (1.5pt);
\draw [color=ffqqqq] (15.16,-1.08) circle (1.5pt);
\draw [color=ffqqqq] (13.78,-2.9) circle (1.5pt);
\fill [color=uuuuuu] (26,-8) circle (1.0pt);
\fill [color=uuuuuu] (30,-8) circle (1.0pt);
\fill [color=uuuuuu] (26,-6) circle (1.0pt);
\fill [color=uuuuuu] (28,-3) circle (1.0pt);
\fill [color=qqqqff] (30,-4) circle (1.5pt);
\draw[color=qqqqff] (0.19,5.94) node {$L$};
\draw [color=ffqqqq] (28.46,-5.11) circle (1.5pt);
\draw [color=ffqqqq] (27.23,-5.8) circle (1.5pt);
\end{scriptsize}
\end{tikzpicture}
		\begin{tikzpicture}[line cap=round,line join=round,>=triangle 45,x=0.8cm,y=0.8cm]
\clip(11.74,-4.85) rectangle (20.91,2.12);
\fill[color=zzttqq,fill=zzttqq,fill opacity=0.1] (4,-4) -- (8,-4) -- (4,-2) -- cycle;
\fill[color=zzttqq,fill=zzttqq,fill opacity=0.1] (4,-2) -- (6,1) -- (8,-4) -- cycle;
\fill[color=zzttqq,fill=zzttqq,fill opacity=0.1] (4,-2) -- (6,1) -- (4,-4) -- cycle;
\fill[color=zzttqq,fill=zzttqq,fill opacity=0.1] (6,1) -- (4,-4) -- (8,-4) -- cycle;
\fill[color=zzttqq,fill=zzttqq,fill opacity=0.1] (16,-4) -- (20,-4) -- (16,-2) -- cycle;
\fill[color=zzttqq,fill=zzttqq,fill opacity=0.1] (16,-2) -- (18,1) -- (20,-4) -- cycle;
\fill[color=zzttqq,fill=zzttqq,fill opacity=0.1] (18,1) -- (16,-4) -- (20,-4) -- cycle;
\fill[color=zzttqq,fill=zzttqq,fill opacity=0.1] (16,-2) -- (12.56,-0.2) -- (16,-4) -- cycle;
\fill[color=qqqqff,fill=qqqqff,fill opacity=0.1] (12.56,-0.2) -- (15.16,-1.08) -- (16,-4) -- cycle;
\fill[color=qqqqff,fill=qqqqff,fill opacity=0.1] (15.16,-1.08) -- (16,-2) -- (12.56,-0.2) -- cycle;
\fill[color=qqqqff,fill=qqqqff,fill opacity=0.1] (15.16,-1.08) -- (16,-4) -- (16,-2) -- cycle;
\fill[color=qqqqff,fill=qqqqff,fill opacity=0.1] (12.56,-0.2) -- (13.78,-2.9) -- (16,-4) -- cycle;
\fill[color=qqqqff,fill=qqqqff,fill opacity=0.1] (13.78,-2.9) -- (16,-2) -- (12.56,-0.2) -- cycle;
\fill[color=qqqqff,fill=qqqqff,fill opacity=0.1] (13.78,-2.9) -- (16,-4) -- (16,-2) -- cycle;
\fill[color=zzttqq,fill=zzttqq,fill opacity=0.15] (26,-6) -- (28,-3) -- (30,-8) -- cycle;
\draw [color=zzttqq] (4,-4)-- (8,-4);
\draw [color=zzttqq] (4,-2)-- (4,-4);
\draw [color=zzttqq] (4,-2)-- (6,1);
\draw [color=zzttqq] (6,1)-- (8,-4);
\draw [dash pattern=on 1pt off 1pt,color=zzttqq] (8,-4)-- (4,-2);
\draw [color=zzttqq] (4,-2)-- (6,1);
\draw [color=zzttqq] (6,1)-- (4,-4);
\draw [color=zzttqq] (4,-4)-- (4,-2);
\draw [color=zzttqq] (6,1)-- (4,-4);
\draw [color=zzttqq] (4,-4)-- (8,-4);
\draw [color=zzttqq] (8,-4)-- (6,1);
\draw [color=zzttqq] (16,-4)-- (20,-4);
\draw [color=zzttqq] (16,-2)-- (16,-4);
\draw [color=zzttqq] (16,-2)-- (18,1);
\draw [color=zzttqq] (18,1)-- (20,-4);
\draw [dash pattern=on 1pt off 1pt,color=zzttqq] (20,-4)-- (16,-2);
\draw [color=zzttqq] (18,1)-- (16,-4);
\draw [color=zzttqq] (16,-4)-- (20,-4);
\draw [color=zzttqq] (20,-4)-- (18,1);
\draw [color=zzttqq] (16,-2)-- (16.67,-2.33);
\draw [shift={(16,-3)},dotted]  plot[domain=1.11:2.46,variable=\t]({1*4.47*cos(\t r)+0*4.47*sin(\t r)},{0*4.47*cos(\t r)+1*4.47*sin(\t r)});
\draw [color=zzttqq] (12.56,-0.2)-- (16,-4);
\draw [color=zzttqq] (16,-4)-- (16,-2);
\draw [shift={(16,-3)},dotted]  plot[domain=0.66:1.98,variable=\t]({1*2.12*cos(\t r)+0*2.12*sin(\t r)},{0*2.12*cos(\t r)+1*2.12*sin(\t r)});
\draw [color=qqqqff] (12.56,-0.2)-- (15.16,-1.08);
\draw [color=qqqqff] (15.16,-1.08)-- (16,-4);
\draw [color=qqqqff] (16,-4)-- (12.56,-0.2);
\draw [color=qqqqff] (15.16,-1.08)-- (16,-2);
\draw [color=qqqqff] (12.56,-0.2)-- (15.16,-1.08);
\draw [color=qqqqff] (15.16,-1.08)-- (16,-4);
\draw [color=qqqqff] (16,-4)-- (16,-2);
\draw [color=qqqqff] (16,-2)-- (15.16,-1.08);
\draw [color=qqqqff] (12.56,-0.2)-- (13.78,-2.9);
\draw [color=qqqqff] (13.78,-2.9)-- (16,-4);
\draw [color=qqqqff] (16,-4)-- (12.56,-0.2);
\draw [dash pattern=on 1pt off 1pt,color=qqqqff] (16,-2)-- (12.56,-0.2);
\draw [color=qqqqff] (12.56,-0.2)-- (13.78,-2.9);
\draw [color=qqqqff] (13.78,-2.9)-- (16,-4);
\draw [color=qqqqff] (16,-4)-- (16,-2);
\draw [dash pattern=on 1pt off 1pt,color=qqqqff] (16,-2)-- (13.78,-2.9);
\draw [color=zzttqq] (26,-8)-- (30,-8);
\draw [color=zzttqq] (26,-6)-- (26,-8);
\draw [color=zzttqq] (26,-6)-- (28,-3);
\draw [color=zzttqq] (28,-3)-- (30,-8);
\draw [color=zzttqq] (30,-8)-- (26,-6);
\draw [color=zzttqq] (26,-6)-- (28,-3);
\draw [color=zzttqq] (28,-3)-- (26,-8);
\draw [color=zzttqq] (26,-8)-- (26,-6);
\draw [color=zzttqq] (28,-3)-- (26,-8);
\draw [color=zzttqq] (26,-8)-- (30,-8);
\draw [color=zzttqq] (30,-8)-- (28,-3);
\draw [dash pattern=on 1pt off 1pt,color=qqqqff] (30,-4)-- (26,-6);
\draw [color=qqqqff] (30,-8)-- (30,-4);
\draw [color=qqqqff] (28,-3)-- (30,-4);
\draw [color=qqqqff] (30,-4)-- (30,-8);
\draw [color=qqqqff] (30,-8)-- (28,-3);
\draw [color=qqqqff] (26,-6)-- (28,-3);
\draw [color=qqqqff] (28,-3)-- (30,-4);
\draw [dash pattern=on 1pt off 1pt,color=ffqqqq] (28.46,-5.11)-- (27.23,-5.8);
\draw [color=qqqqff] (28.67,-4.68)-- (30,-4);
\draw [color=ffqqqq] (27.23,-5.8)-- (27.86,-5.45);
\draw (5.56,-1.25) node[anchor=north west] {$\TT_i$};
\draw (6.85,0.18) node[anchor=north west] {$\Gamma_{j_1}$};
\draw (6.99,-2.49) node[anchor=north west] {$\Gamma_{j_2}$};
\draw (4.02,-0.04) node[anchor=north west] {$\Gamma_{j_3}$};
\draw (4.63,-4.18) node[anchor=north west] {$\Gamma_{j_4}$};
\draw (14.4,0.1) node[anchor=north west] {$\QQ_{j_3}$};
\draw (16.02,-0.08) node[anchor=north west] {$i$};
\draw (17.82,-1.74) node[anchor=north west] {$\TT_i$};
\draw (13.21,-3.03) node[anchor=north west] {$\mathbb{\wp}(i,j_3)$};
\draw [->,dash pattern=on 1pt off 1pt] (28.02,-5.62) -- (28.66,-5.29);
\draw (27.13,-5.78) node[anchor=north west] {$\ell(j)$};
\draw (28.99,-4.46) node[anchor=north west] {$r(j)$};
\draw (28.36,-5.41) node[anchor=north west] {$\vec{n_j}$};
\draw (30.3,-7.62) node[anchor=north west] {$\Gamma_j$};
\begin{scriptsize}
\fill [color=uuuuuu] (4,-4) circle (1.0pt);
\fill [color=uuuuuu] (8,-4) circle (1.0pt);
\fill [color=uuuuuu] (4,-2) circle (1.0pt);
\fill [color=uuuuuu] (6,1) circle (1.0pt);
\fill [color=uuuuuu] (16,-4) circle (1.0pt);
\fill [color=uuuuuu] (20,-4) circle (1.0pt);
\fill [color=uuuuuu] (16,-2) circle (1.0pt);
\fill [color=uuuuuu] (18,1) circle (1.0pt);
\fill [color=uuuuuu] (12.56,-0.2) circle (1.0pt);
\draw [color=ffqqqq] (17.68,-1.7) circle (1.5pt);
\draw [color=ffqqqq] (15.16,-1.08) circle (1.5pt);
\draw [color=ffqqqq] (13.78,-2.9) circle (1.5pt);
\fill [color=uuuuuu] (26,-8) circle (1.0pt);
\fill [color=uuuuuu] (30,-8) circle (1.0pt);
\fill [color=uuuuuu] (26,-6) circle (1.0pt);
\fill [color=uuuuuu] (28,-3) circle (1.0pt);
\fill [color=qqqqff] (30,-4) circle (1.5pt);
\draw[color=qqqqff] (0.19,5.94) node {$L$};
\draw [color=ffqqqq] (28.46,-5.11) circle (1.5pt);
\draw [color=ffqqqq] (27.23,-5.8) circle (1.5pt);
\end{scriptsize}
\end{tikzpicture}
    \caption{An example of a tetrahedral element of the primary mesh with $S_i=\{j_1, j_2, j_3, j_4\}$ (left) a non-standard dual face-based hexahedral element associated to the face $j_3$ (right).}
    \label{fig.MESH_1}
		\end{center}
\end{figure}
An example of the resulting main and dual grid in three space dimensions is reported in Figure \ref{fig.MESH_1}. The main grid consists of tetrahedral simplex elements, and the face-based dual
elements contain the three vertices of the common triangular face of two tetrahedra (a left and a right one), and the two barycenters of the two tetrahedra that share the same face. 
In three space dimensions the dual grid therefore consists of non-standard five-point hexahedral elements. The same face-based staggered dual mesh has also been used in 
\cite{Bermudez2014,USFORCE,USFORCE2}.  


\paragraph{Space-time extension}
In the time direction we cover the time interval $[0,T]$ with a sequence of times $0=t^0<t^1<t^2 \ldots <t^N<t^{N+1}=T$. We denote the time step by $\Delta t^{n+1} = t^{n+1}-t^{n} $ and 
the corresponding time interval by $T^{n+1}=[t^{n}, t^{n+1}]$ for $n=0 \ldots N$. In order to ease notation, sometimes we will use the abbreviation $\Delta t= \Delta t^{n+1}$. 
In this way the generic space-time element defined in the time interval $[t^n, t^{n+1}]$ is given by $\TT_i^\st = \TT_i \times T^{n+1}$ for the main grid, and $\QQ_j^\st=\QQ_j \times T^{n+1}$ for the dual grid.
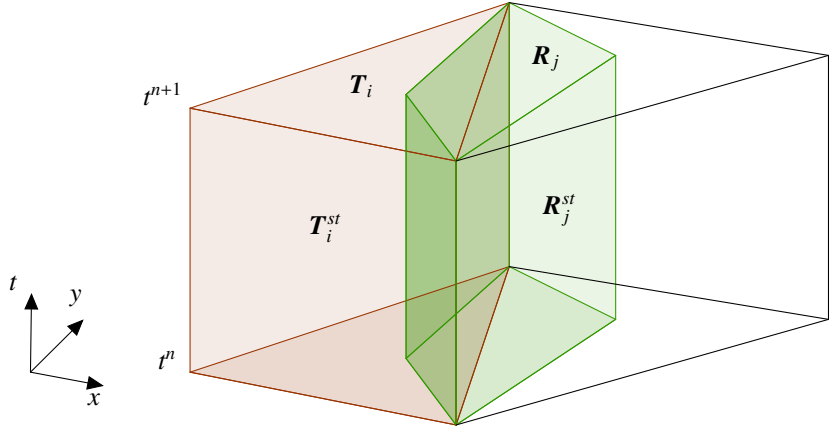
\begin{figure}[ht]
    \begin{center}
	\begin{tikzpicture}[line cap=round,line join=round,>=triangle 45,x=0.7cm,y=0.7cm]
\clip(-1.72,-8.88) rectangle (16.84,3.5);
\fill[color=zzttqq,fill=zzttqq,fill opacity=0.1] (9,-4) -- (3,-6) -- (8,-7) -- cycle;
\fill[color=ttzzqq,fill=ttzzqq,fill opacity=0.1] (9,-4) -- (7.06,-5.74) -- (8,-7) -- (11,-5) -- cycle;
\fill[color=zzttqq,fill=zzttqq,fill opacity=0.1] (9,1) -- (3,-1) -- (8,-2) -- cycle;
\fill[color=ttzzqq,fill=ttzzqq,fill opacity=0.1] (9,1) -- (7.06,-0.74) -- (8,-2) -- (11,0) -- cycle;
\fill[color=zzttqq,fill=zzttqq,fill opacity=0.1] (3,-6) -- (8,-7) -- (8,-2) -- (3,-1) -- cycle;
\fill[color=zzttqq,fill=zzttqq,fill opacity=0.1] (8,-7) -- (9,-4) -- (9,1) -- (8,-2) -- cycle;
\fill[color=ttzzqq,fill=ttzzqq,fill opacity=0.1] (8,-7) -- (11,-5) -- (11,0) -- (8,-2) -- cycle;
\fill[color=ttzzqq,fill=ttzzqq,fill opacity=0.25] (7.06,-0.74) -- (8,-2) -- (8,-7) -- (7.06,-5.74) -- cycle;
\fill[color=ttzzqq,fill=ttzzqq,fill opacity=0.2] (7.06,-5.74) -- (9,-4) -- (9,1) -- (7.06,-0.74) -- cycle;
\draw [color=zzttqq] (9,-4)-- (3,-6);
\draw [color=zzttqq] (3,-6)-- (8,-7);
\draw [color=zzttqq] (8,-7)-- (9,-4);
\draw [color=ttzzqq] (9,-4)-- (7.06,-5.74);
\draw [color=ttzzqq] (7.06,-5.74)-- (8,-7);
\draw [color=ttzzqq] (8,-7)-- (11,-5);
\draw [color=ttzzqq] (11,-5)-- (9,-4);
\draw [color=zzttqq] (9,1)-- (3,-1);
\draw [color=zzttqq] (3,-1)-- (8,-2);
\draw [color=zzttqq] (8,-2)-- (9,1);
\draw [color=ttzzqq] (9,1)-- (7.06,-0.74);
\draw [color=ttzzqq] (7.06,-0.74)-- (8,-2);
\draw [color=ttzzqq] (8,-2)-- (11,0);
\draw [color=ttzzqq] (11,0)-- (9,1);
\draw [color=zzttqq] (3,-6)-- (8,-7);
\draw [color=zzttqq] (8,-7)-- (8,-2);
\draw [color=zzttqq] (8,-2)-- (3,-1);
\draw [color=zzttqq] (3,-1)-- (3,-6);
\draw [color=zzttqq] (8,-7)-- (9,-4);
\draw [color=zzttqq] (9,-4)-- (9,1);
\draw [color=zzttqq] (9,1)-- (8,-2);
\draw [color=zzttqq] (8,-2)-- (8,-7);
\draw [color=ttzzqq] (8,-7)-- (11,-5);
\draw [color=ttzzqq] (11,-5)-- (11,0);
\draw [color=ttzzqq] (11,0)-- (8,-2);
\draw [color=ttzzqq] (8,-2)-- (8,-7);
\draw [color=ttzzqq] (7.06,-0.74)-- (8,-2);
\draw [color=ttzzqq] (8,-2)-- (8,-7);
\draw [color=ttzzqq] (8,-7)-- (7.06,-5.74);
\draw [color=ttzzqq] (7.06,-5.74)-- (7.06,-0.74);
\draw [color=ttzzqq] (7.06,-5.74)-- (9,-4);
\draw [color=ttzzqq] (9,-4)-- (9,1);
\draw [color=ttzzqq] (9,1)-- (7.06,-0.74);
\draw [color=ttzzqq] (7.06,-0.74)-- (7.06,-5.74);
\draw (5.82,-0.2) node[anchor=north west] {$\TT_i$};
\draw (9.24,0.36) node[anchor=north west] {$\QQ_j$};
\draw [->] (0,-6) -- (0.02,-4.5);
\draw [->] (0,-6) -- (1,-5);
\draw [->] (0,-6) -- (1.38,-6.26);
\draw (0.56,-4.26) node[anchor=north west] {$y$};
\draw (0.88,-6.22) node[anchor=north west] {$x$};
\draw (-0.58,-3.98) node[anchor=north west] {$t$};
\draw (1.94,-0.34) node[anchor=north west] {$t^{n+1}$};
\draw (2.22,-5.42) node[anchor=north west] {$t^{n}$};
\draw (5.06,-2.82) node[anchor=north west] {$\TT_i^{st}$};
\draw (9.44,-2.44) node[anchor=north west] {$\QQ_j^{st}$};
\draw (9,1)-- (15,0);
\draw (15,0)-- (8,-2);
\draw (15,-5)-- (15,0);
\draw (15,-5)-- (9,-4);
\draw (15,-5)-- (8,-7);
\end{tikzpicture}
    \caption{Example of space-time elements $\TT_i^\st$ (red) and $\QQ_j^\st$ (green) with $j \in S_i$.}
    \label{fig.st1}
		\end{center}
\end{figure}
Figure \ref{fig.st1} shows a graphical representation of the primary and dual space-time control volumes for the case of two space dimensions plus time. 

\subsection{Space-time basis functions}
\label{sec222}
According to \cite{2DSIUSW,2STINS,3DSIINS} we proceed as follows: in the two dimensional case, we first construct the polynomial basis up to a generic polynomial degree $p$ on some triangular and  quadrilateral reference elements with local coordinates $\xi$ and $\eta$. In particular, we take $T_{std}=\{(\xi,\eta) \in \R^{2} \,\, | \,\,  0 \leq \xi \leq 1, \,\,  0 \leq \eta \leq 1-\xi \}$ as the 
reference triangle and the unit square as the reference quadrilateral element $R_{std}=[0,1]^2$. Using the standard nodal approach of conforming continuous finite elements, we obtain 
$N_\phi=\frac{(p+1)(p+2)}{2}$ basis functions $\{\phi_k \}_{k \in [1,N_\phi]}$ on $T_{std}$ and $N_{\psi}=(p+1)^2$ basis functions on $R_{std}$. 
The connection between the reference coordinates $\boldsymbol{\xi}=(\xi,\eta)$ and the physical coordinates $\xx=(x,y)$ is performed by the maps $T_i:\TT_i \longrightarrow T_{std}$ for every $i =1 \ldots N_e$; $T_j:\QQ_j \longrightarrow R_{std}$ for every $j =1 \ldots N_d$ and its inverse, called $T_i^{-1}:\TT_i \longleftarrow T_{std}$ and $T_j^{-1}:\QQ_j \longleftarrow R_{std}$, respectively. The maps from the 
reference coordinates to the physical ones can be constructed following a classical sub-parametric or a complete iso-parametric approach. 

Regarding the basis functions in three space dimensions, we use the unit tetrahedron with vertices $(0,0,0)$, $(1,0,0)$, $(0,1,0)$ and $(0,0,1)$ 
to construct the basis polynomials for the main grid. We use again the standard nodal basis functions of conforming finite elements based on the 
reference element $T_{std}=\{(\xi,\eta,\zeta) \in \R^{3} \,\, | \,\,  0 \leq \xi \leq 1, \,\,  0 \leq \eta \leq 1-\xi, \,\, 0 \leq \zeta \leq 1-\xi-\eta \}$ 
and then using either a subparametric or an isoparametric map to connect the reference space $\boldsymbol{\xi}=(\xi,\eta,\zeta)$ to the physical space $\xx=(x,y,z)$ and vice-versa.  
For the non-standard five-point hexahedral elements of the dual mesh, we define the polynomial basis directly in the physical space using a simple modal basis function based on 
rescaled Taylor monomials, such as the ones defined in \cite{3DSIINS}. We thus obtain $N_\phi=N_\psi=\frac{(p+1)(p+2)(p+3)}{6}$ basis functions per element for both, 
the main grid and the dual mesh. 

Finally, we construct the time basis functions on a reference interval $I_{std}=[0,1]$ for polynomials of degree $p_\gamma$. 
In this case the resulting $N_\gamma=p_\gamma+1$ basis functions $\{\gamma_k\}_{k \in [1, N_\gamma]}$ are defined as the Lagrange interpolation polynomials passing through the Gauss-Legendre 
quadrature points for the unit interval. For every time interval $[t^n, t^{n+1}]$, the map between the reference interval and the physical one is simply given by $t=t^n+\tau \Delta t^{n+1}$ for 
every $\tau \in [0,1]$. 
Using the tensor product we can finally construct the basis functions on the space-time elements $\TT_i^\st$ and $\QQ_j^\st$ as $\tilde{\phi}(\boldsymbol{\xi},\tau)=\phi(\boldsymbol{\xi}) \cdot \gamma(\tau)$ and $\tilde{\psi}(\boldsymbol{\xi},\tau)=\psi\boldsymbol{\xi}) \cdot \gamma(\tau)$. The total number of basis functions becomes $N_\phi^\st=N_\phi \cdot N_\gamma$ and $N_\psi^\st=N_\psi \cdot N_\gamma$.

\section{Semi-implicit space-time DG scheme} 
\label{sec_semi_imp_dg} 
The discrete pressure $p_h$ is defined on the main grid, where we will use the notation $p_i(\xx,t)=p_h(\xx,t)|_{\TT_i^\st}$, while the discrete 
velocity vector field $\vv_h$ and the fluid density $\rho_h$ are defined on the dual grid, namely  $\vv_j(\xx,t)=\vv_h(\xx,t)|_{\QQ_j^\st}$ and $\rho_j(\xx,t)=\rho_h(\xx,t)|_{\QQ_j^\st}$. \par 
The discrete total energy density $(\rhoEE)_h$ is defined on the main grid while the discrete specific enthalpy $H_h$ and the discrete momentum density $(\rhov)_h$ are defined at the dual level.

The numerical solution of \eref{eq:CNS_1_2}-\eref{eq:CNS_3_2} is represented inside the space-time control volumes of the primal and the dual grid and for a time slice $T^{n+1}$ 
by piecewise space-time polynomials. The discrete pressure and the total energy on the main mesh as well as the momentum and the density on the dual mesh read as follows:  
\begin{eqnarray}
 	p_i(\xx,t)      & = & \sum\limits_{l=1}^{N_\phi^\st} \tilde{\phi}_l^{(i)}(\xx,t) \hat{p}_{l,i}^{n+1}=:\tilde{\bphi}^{(i)}(\xx,t)\hat{\mathbf{p}}_i^{n+1}, \nonumber \\ 
	\rho E_i(\xx,t) & = & \sum\limits_{l=1}^{N_\phi^\st} \tilde{\phi}_l^{(i)}(\xx,t) \widehat{\rho E}_{l,i}^{n+1}      =:\tilde{\bpsi}^{(j)}(\xx,t) \widehat{\boldsymbol{\rho} \mathbf{E}}_i^{n+1}, 
\label{eq:D_1}  
\end{eqnarray}
\begin{eqnarray}
	\rho \mathbf{v}_j(\xx,t)  & = & \sum\limits_{l=1}^{N_\psi^\st} \tilde{\psi}_l^{(j)}(\xx,t) \widehat{\rho \mathbf{v}}_{l,j}^{n+1}=:\tilde{\bpsi}^{(j)}(\xx,t)\widehat{\boldsymbol{\rho} \mathbf{v}}_j^{n+1}, \nonumber \\ 
	\rho_j(\xx,t) & = & \sum\limits_{l=1}^{N_\psi^\st} \tilde{\psi}_l^{(j)}(\xx,t) \hat{\rho}_{l,j}^{n+1}      =:\tilde{\bpsi}^{(j)}(\xx,t)\hat{\boldsymbol{\rho}}_j^{n+1}. 
\label{eq:D_2}
\end{eqnarray}
The same definitions hold also for all the other quantities on the main grid $\{Q_i\}_{i \in [1,\Ni]}$ and on the dual grid $\{Q_j\}_{j \in [1,\Nj]}$. 
Here and in the rest of the paper we use the convention that variables indexed by $j$ are defined on the dual grid, while the index $i$ is used for the quantities on the main grid. 
The vector of basis functions $\tilde{\bphi}^{(i)}(\xx,t)$ is generated via the map $T_i^{-1}$
from $\tilde{\bphi}(\boldsymbol{\xi}, \tau)$ on $T_{std}\times [0,1]$. The vector $\tilde{\bpsi}^{(j)}(\xx,t)$ is generated from $\tilde{\bpsi}(\boldsymbol{\xi}, \tau)$ on $R_{std} \times [0,1]$ in the two dimensional case and it is directly defined on the physical space for each element in the three dimensional case, see e.g. \cite{3DSIINS}.  

\subsection{Auxiliary variables}
\label{sec:auxvar}
First of all we define some average operators on staggered grids that will be used in the following. For all $i=1 \ldots \Ni$, $k=1\ldots N_\phi^\st$ we have the identity
\begin{eqnarray}
	\int\limits_{\TT_i^\st}{\phit_k^{(i)}\phit_l^{(i)}\dxt} \hat{Q}_{l,i}^{n+1} =\int\limits_{\TT_i^\st}{\phit_k^{(i)} Q_i(\xx,t) \dxt} = 
	\sum\limits_{j \in S_i} \int\limits_{\TT_{i,j}^\st}{\phit_k^{(i)} Q_j(\xx,t) \dxt}=\sum\limits_{j \in S_i}{\int\limits_{\TT_{i,j}^\st}{\phit_k^{(i)}\psit_l^{(j)}\dxt} \hat{Q}_{l,j}^{n+1} },
\label{eq:Aux1}
\end{eqnarray}
which is a high order average operator for a general quantity $\{Q_j\}_{j \in [1,\Nj]}$ from the dual grid to the main grid, where the same quantity will be denoted by 
$\{Q_i\}_{i \in [1,\Ni]}$. Recall that the index $i$ refers to the primary mesh, while the index $j$ refers to the staggered dual mesh. 
The previous equation $\eqref{eq:Aux1}$ is then written in a compact matrix form as  
\begin{eqnarray}
	\qh_i^{n+1} = \MM_i^{-1} \sum_{j \in S_i} \LL_{i,j} \qh_j^{n+1},
\label{eq:aux2}
\end{eqnarray}
where
\begin{eqnarray}
	\MM_i = \int\limits_{\TT_i^\st}{\phit_k^{(i)}\phit_l^{(i)}\dxt} \qquad \LL_{i,j}=\int\limits_{\TT_{i,j}^\st}{\phit_k^{(i)}\psit_l^{(j)}\dxt}.
\label{eq:aux3}
\end{eqnarray}
Equation \eqref{eq:aux2} defines the average operator $L_\dtm:\{\QQ_j\}_j \rightarrow \{\TT_i\}_i$. 
In a similar way one can derive the high order average operator $L_\mtd:\{\TT_i\}_i \rightarrow \{\QQ_j\}_j$ from the primary grid to the dual mesh as  
\begin{eqnarray}
	\int\limits_{\QQ_j^\st}{\psit_k^{(j)} Q_j(\xx,t) \dxt} = 
	\int\limits_{\TT_{\ell(j),j}^\st}{\psit_k^{(j)} Q_{\ell(j)}(\xx,t) \dxt} + 
	\int\limits_{\TT_{r(j),j}^\st}{\psit_k^{(j)} Q_{r(j)}(\xx,t) \dxt}   
\label{eq:Aux4a}
\end{eqnarray}
which in a more compact matrix form simply reads 
\begin{eqnarray}
	\qh_j^{n+1} = \bar{\MM}_j^{-1}\left( \LL^\top_{\ell(j),j} \qh_{\ell(j)}^{n+1} + \LL^\top_{r(j),j} \qh_{r(j)}^{n+1} \right),
\label{eq:aux4}
\end{eqnarray}
with the element mass matrix computed on the staggered dual mesh $\bar{\MM}_j$ given by 
\begin{eqnarray}
	\bar{\MM}_j = \int\limits_{\QQ_j^\st}{\psit_k^{(j)}\psit_l^{(j)}\dxt}. 
\label{eq:aux5}
\end{eqnarray}

Within the numerical scheme, we also need to compute derived quantities from the pressure and the conservative quantities. In two space dimensions we have a nodal basis on the dual mesh, 
hence for the degrees of freedom of the variables $\mathbf{v}_j$, $H_j$ and $\rhok_j$ we can use a standard point-wise evaluation, 
for instance   
\begin{eqnarray}
\hat{H}_{l,j}^{n+1} & = & \hat{e}_{l,j}^{n+1} + \frac{\hat{p}_{l,j}^{n+1}}{\hat{\rho}_{l,j}^{n+1}}, \qquad \,\, \textnormal{ with } \qquad \hat{e}_{l,j}^{n+1} = \ie(\hat{p}_{l,j}^{n+1},\hat{\rho}_{l,j}^{n+1}),
\label{eq:aux6}
\end{eqnarray}
\begin{eqnarray}
	\hat{\mathbf{v}}_{l,j}^{n+1} = \frac{\widehat{\rho \mathbf{v}}_{l,j}^{n+1}}{\hat{\rho}_{l,j}^{n+1}},  \qquad \widehat{\rho k}_{l,j}^{n+1}  = \frac{1}{2} \hat{\rho}_{l,j}^{n+1} \hat{\mathbf{v}}_{l,j}^{n+1} \cdot \hat{\mathbf{v}}_{l,j}^{n+1}.
\label{eq:aux7}
\end{eqnarray}
In three space dimensions we have a modal basis on the dual mesh, i.e. the degrees of freedom of the different variables have to be related to each other via standard $L_2$ projection. 

\subsection{Convective terms} 


In the framework of semi-implicit schemes \cite{Casulli1984,CasulliCheng1992,CasulliCattani,CasulliWalters2000,DumbserCasulli2016}, the nonlinear convective terms are typically discretized 
\textit{explicitly}, while the pressure terms are treated \textit{implicitly}. For more details on the relation to flux-vector splitting schemes see \cite{ToroVazquez}. 

It is therefore necessary to introduce a weak formulation of the convective part of the compressible Navier-Stokes equations on the main grid.  
The PDE of the purely convective part for a general quantity $Q$ is given by  
\begin{eqnarray}
	\frac{\partial Q}{\partial t} + \nabla \cdot \mathbf{F}_Q(Q) = 0,   
\label{eq:der1}
\end{eqnarray}
with the convective flux $\mathbf{F}_Q(Q) = Q \, \mathbf{v}$ for the quantity $Q$. In the compressible Navier-Stokes equations $Q$ represents either density $\rho$, momentum density $\rho \vv$
or kinetic energy density $\rho k$. Denoting the discrete solution for quantity $Q$ by $Q_h = \phi^{(i)}_l \hat Q_{l,i}^{n+1} = \boldsymbol{\phi}^{(i)} \hat{\mathbf{Q}}_{i}^{n+1}$, 
multiplying Equation \eqref{eq:der1} by $\phit_k^{(i)}$ and integrating over a space-time control volume $\TT_i^\st$, for every $k=1\ldots N_\phi^\st$ yields  
\begin{eqnarray}
	\int\limits_{\TT_i}{\phit_k^{(i)}(\xx,t^{n+1}_{-}) Q_h(\xx,t^{n+1}_{-}) \dx} - \int\limits_{\TT_i}{\phit_k^{(i)}(\xx,t^n_+) Q_h(\xx,t^{n}_{-}) \dx} - \int\limits_{\TT_i^\st}{\frac{\partial}{\partial t} \phit_k^{(i)} Q_h \dxt}  + && \nonumber \\ 
	\int \limits_{t^n}^{t^{n+1}} \int \limits_{\partial \TT_i}{\phit_k^{(i)} \mathcal{F}_Q(Q_h^-,Q_h^+) \cdot \vec{n} \, ds\,dt}-\int\limits_{\TT_i^\st}{\nabla \phit_k^{(i)} \cdot \mathbf{F}_Q(Q_h) \dxt}
	= 0, && 
\label{eq:der2}
\end{eqnarray}
where we have used the integration by parts formula in space and time and where a \textit{numerical flux} (Riemann solver) has been introduced to resolve the jumps at the element interfaces. 
The notation $t^n_\pm = \lim \limits_{\epsilon \to 0^\pm} t^n + \epsilon$ indicates whether the data are taken from within the current space-time slab, or from the previous one. 
This choice corresponds to a simple upwind flux in time direction due to the causality principle, see \cite{2STINS,3DSIINS}. For the spatial fluxes we use a simple Rusanov-type flux 
\begin{eqnarray}
		\mathcal{F}_Q(Q_h^-,Q_h^+) \cdot \vec{n} = \frac{1}{2}\left( \mathbf{F}_Q(Q_h^+) + \mathbf{F}_Q(Q_h^-) \right) \cdot \vec{n} - \frac{1}{2} s_{\max} (Q_h^+ - Q_h^-), 
\label{eq:der3}
\end{eqnarray}
with $s_{\max} = \left\| \vv \cdot \vec{n} \right\|$.
Equation \eqref{eq:der2} is written in a more compact matrix vector notation as 
\begin{equation}
	  \Mpsi_i  \hat{\mathbf{Q}}_{i}^{n+1} = \Mpsi_i^- \hat{\mathbf{Q}}_{i}^{n} - \widehat{\mathbf{FQ}}_i^{n+1}, 
\label{eq:der4}
\end{equation}
with 
\begin{eqnarray}
	  \Mpsi_i^+     = \int\limits_{\TT_i}{\tilde{\phi}_k^{(i)}(\xx,t^{n+1}_{-}) \tilde{\phi}_l^{(i)}(\xx,t^{n+1}_{-})  \dx},  \qquad \qquad   
    \Mpsi_i^-     = \int\limits_{\TT_i}{\tilde{\phi}_k^{(i)}(\xx,t^{n}_{+}) \tilde{\phi}_l^{(i)}(\xx,t^{n}_{-})    \dx},   && \nonumber \\  
    \Mpsi_i^\circ = \int\limits_{\TT_i^\st}{\diff{\tilde{\phi}_k^{(i)}}{t} \tilde{\phi}_l^{(i)} \dx dt},   \qquad \qquad                  
		\Mpsi_i       = \left(\Mpsi_i^+ - \Mpsi_i^\circ \right), && \nonumber \\
		\widehat{\mathbf{FQ}}_i^{n+1} = \int \limits_{t^n}^{t^{n+1}} \int \limits_{\partial \TT_i}{\phit_k^{(i)} \mathcal{F}_Q(Q_h^-,Q_h^+) \cdot \vec{n} \, ds\,dt}-\int\limits_{\TT_i^\st}{\nabla \phit_k^{(i)} \cdot \mathbf{F}_Q(Q_h) \dxt}.  && 		 
\end{eqnarray}

\subsection{Derivation of the semi-implicit space-time DG scheme}
We now derive a weak formulation of the system \eref{eq:CNS_1}-\eref{eq:CNS_3}. The presented algorithm is a rather natural extension of the scheme presented in 
\cite{3DSIINS} for the incompressible Navier-Stokes equations. It also follows the main ideas outlined in \cite{DumbserCasulli2016} for a simple 
and low order accurate staggered finite volume scheme. 

\paragraph{Discrete continuity equation}
The discretization of the continuity equation would be in principle immediate, since it contains only a 
convective flux. However, following the ideas put forward in \cite{DumbserCasulli2016}, the convective flux in the continuity equation should be discretized exactly in the same way as the 
convective flux in the momentum equation, in order to avoid unphysical oscillations for uniform pressure and velocity flows, see \cite{DumbserCasulli2016}. 
The momentum equation contains the divergence of the viscous stress tensor, which is naturally defined on the staggered \textit{dual grid} (see \cite{3DSIINS}) in order to yield a simple 
formulation for its divergence on the primary mesh. Therefore, we first have to average the discrete momentum from the dual mesh to the main grid, where then the convective and 
the viscous terms can be conveniently discretized. Hence, we must proceed in the same way for the continuity equation. Therefore, we first compute the density degrees of freedom at the old 
time level on the main grid by using the above-defined averaging operators  
\begin{equation}
  	\hat{\boldsymbol{\rho}}_i^{n} = \MM_i^{-1} \sum_{j \in S_i} \LL_{i,j} \hat{\boldsymbol{\rho}}_j^{n}.
		\label{eqn.avrho} 
\end{equation} 
Subsequently, we discretize the continuity equation on the main grid using the operator defined in \eqref{eq:der4} as 
\begin{equation}
  	\Mpsi_i \hat{\boldsymbol{\rho}}_i^{n+1} = \Mpsi_i^- \hat{\boldsymbol{\rho}}_i^{n} - \widehat{\mathbf{F}\boldsymbol{\rho}}_i^{n+1} 
	   \label{eqn.densitymain}
\end{equation}
and then we average back to the dual mesh using \eqref{eq:aux4} 
\begin{equation}
	\hat{\boldsymbol{\rho}}_j^{n+1} = \bar{\MM}_j^{-1}\left( \LL^\top_{\ell(j),j} \hat{\boldsymbol{\rho}}_{\ell(j)}^{n+1} + \LL^\top_{r(j),j}\hat{\boldsymbol{\rho}}_{r(j)}^{n+1} \right).
	\label{eqn.densitydual} 
\end{equation}
With \eqref{eqn.densitymain} and \eqref{eqn.densitydual} the new density degrees of freedom are known in all elements on the primary and the dual mesh. 
 
\paragraph{Discrete momentum equation}
For the discretization of the momentum equation, we first consider only the nonlinear convective terms and the viscous stress tensor. Suppose we have an initial guess
for $\hat{\vv}_i^{n+1}$. Following the ideas of \cite{3DSIINS}, the momentum is first averaged from the dual grid onto the main grid as 
 \begin{equation}
  	\widehat{\boldsymbol{\rho} \vv}_i^{n+1} = \MM_i^{-1} \sum_{j \in S_i} \LL_{i,j} \widehat{\boldsymbol{\rho} \vv} _j^{n+1}.
\end{equation} 
Since the density $\hat{\boldsymbol{\rho}}_i^{n+1}$ is already known from \eqref{eqn.avrho}, we also know the velocity field $\hat{\vv}_i^{n+1}$. The weak form of the discrete velocity 
gradient $\nabla_h \vv_j$ is then naturally defined on the \textit{dual mesh}, taking into account also the jump between 
$\vv_{\ell(j)}$ and $\vv_{r(j)}$ at the interface $\Gamma_j^\st$ in the sense of distributions \cite{3DSIINS}. 
Alternatively, this approach can also be interpreted as a Bassi-Rebay lifting operator \cite{BassiRebay}, but acting on the dual mesh. The discrete velocity gradient on the dual mesh 
is thus computed as  
\begin{eqnarray}
  \int\limits_{\QQ_j^\st}{\tilde{\psi}_k^{(j)} \nabla_h \vv_j \dxt} = 
	 \!\!\! \int\limits_{\TT_{\ell(j),j}^\st} \!\! {\tilde{\psi}_k^{(j)} \nabla \vv_{\ell(j)} \dxt} + 
	 \!\!\! \int\limits_{\TT_{r(j),j}^\st}    \!\! {\tilde{\psi}_k^{(j)} \nabla \vv_{r(j)} \, \dxt} + 
	 \int\limits_{\Gamma_j^\st}{\tilde{\psi}_k^{(j)} \left(\vv_{r(j)}-\vv_{\ell(j)}\right) \otimes \nstd \, ds dt}.
	\label{eqn.viscstress} 
\end{eqnarray} 
Once the discrete velocity gradient is known on the dual mesh, the viscous stress tensor can then be immediately computed directly from its definition \eqref{eq:CNS_4} as 
$\boldsymbol{\sigma}_j = \boldsymbol{\sigma}(\nabla_h \vv_j)$. A discrete version of the divergence of the nonlinear convective and viscous flux terms $\mathbf{G}=\mathbf{F}_{\rho \vv} - \stens$ 
on the main grid is then given by 
\begin{eqnarray}
	\widehat{\mathbf{G} \vv}_i^{n+1} = \int\limits_{\TT_i^\st}{\phit_k^{(i)} \nabla \cdot \mathbf{\TF}_h \dxt}  &=&
	\int \limits_{t^n}^{t^{n+1}} \int \limits_{\partial \TT_i}{\phit_k^{(i)} \mathcal{F}_{\rho \vv}({\rho \vv}_h^-,{\rho \vv}_h^+) \cdot \vec{n} \, ds\,dt} -
	\int\limits_{\TT_i^\st}{\nabla \phit_k^{(i)} \cdot \mathbf{F}_{\rho \vv}((\rho \vv)_h) \dxt} +  \nonumber \\ 
	&& \sum_{j \in S_i}\left( \int\limits_{\Gamma_j^\st}   \phit_k^{(i)} \boldsymbol{\sigma}_j \cdot \vec{n}_{i,j} \, ds dt 
		-\int\limits_{\TT_{ij}^\st}{\nabla \phit_k^{(i)} \cdot \boldsymbol{\sigma}_j \, \dxt}  \right). 
\label{eq:visconv}
\end{eqnarray}
Note that in \eqref{eq:visconv} a numerical flux function is only used for the hyperbolic convective terms and not for the parabolic viscous terms, unlike in DG schemes on 
collocated meshes, where a numerical flux is also needed for parabolic operators \cite{BassiRebay,CBS-convection-diffusion,Gassner2007}. 
With the average operator \eqref{eq:aux4}, we can also map this discrete operator back onto the staggered dual mesh as
\begin{equation}
 \widehat{\mathbf{G} \vv}_j^{n+1} = \bar{\MM}_j^{-1}\left( \LL^\top_{\ell(j),j} \widehat{\mathbf{G} \vv}_{\ell(j)}^{n+1} + \LL^\top_{r(j),j} \widehat{\mathbf{G} \vv}_{r(j)}^{n+1} \right).
\label{eq:visconv2} 
\end{equation} 

In what follows, we will also consider the pressure terms. Multiplication of the momentum equation \eref{eq:CNS_2} by $\tilde{\bpsi}$ and integrating over a space time control volume $\QQ_j^\st$ yields
\begin{equation}
\int\limits_{\QQ_j^\st}{\tilde{\psi}_k^{(j)}\left( \diff{(\rhov)_h}{t} + \nabla \cdot \mathbf{\TF}_h \right)  \dxt}+\int\limits_{\QQ_j^\st}{\tilde{\psi}_k^{(j)} \nabla p_h \, \dxt}=0.
\label{eq:CS_5}
\end{equation}
The term $\nabla \cdot \mathbf{\TF}_h$ is the contribution of the discrete nonlinear convective and viscous terms on the dual grid as given by \eqref{eq:visconv} above, hence we assume it as known. 
By taking into account the jumps of $\nabla p_h$ inside the dual elements one can write (see e.g. \cite{3DSIINS}), 
\begin{eqnarray}
\int\limits_{\QQ_j^\st}{\tilde{\psi}_k^{(j)}\left( \diff{(\rhov)_j}{t}+\nabla \cdot \mathbf{\TF}_h \right)  \dxt}
+\hspace{-3mm} \int\limits_{\TT_{\ell(j),j}^\st} \!\! {\tilde{\psi}_k^{(j)} \nabla \np_{\ell(j)} \dxt} 
+\hspace{-3mm} \int\limits_{\TT_{r(j),j}^\st}    \!\! {\tilde{\psi}_k^{(j)} \nabla \np_{r(j)} \, \dxt} + 
\int\limits_{\Gamma_j^\st}{\tilde{\psi}_k^{(j)} \left(\np_{r(j)}-\np_{\ell(j)}\right) \nstd \, ds dt}=0.
\label{eq:CS_9}
\end{eqnarray}
Integrating the first term of \eref{eq:CS_9} by parts in time one obtains  
\begin{eqnarray}
\int\limits_{\QQ_j^\st}{\tilde{\psi}_k^{(j)} \diff{(\rhov)_j}{t} \dxt} &=& 
 \int\limits_{\QQ_j}{\tilde{\psi}_k^{(j)}(\xx,t^{n+1}_{-}) (\rhov)_j(\xx,t^{n+1}_{-}) \dx }  - \int\limits_{\QQ_j}{\tilde{\psi}_k^{(j)}(\xx,t^{n}_{+}) (\rhov)_j(\xx,t^{n}_{-}) \dx }   
- \int\limits_{\QQ_j^\st}{\diff{\tilde{\psi}_k^{(j)}}{t} (\rhov)_j \dxt}.  
\label{eq:CS_11_2}
\end{eqnarray}
In Equation \eqref{eq:CS_11_2} one can again recognize the fluxes between the current space-time element $\QQ_j \times T^{n+1}$, the future space-time slab and the past space-time elements, as well as an internal space-time volume contribution that connects the layers inside the space-time element $\QQ_j^\st$ in an asymmetric way. This formulation includes directly the initial condition of the momentum in a weak sense. It can also be interpreted as using an upwind flux in time direction. Note that for $p_\gamma=0$ the basis functions are constant in time and so the last integral in 
\eqref{eq:CS_11_2} vanishes. 
By substituting Equation \eref{eq:CS_11_2} into \eref{eq:CS_9} and using \eqref{eq:D_2}, we obtain the following weak formulation of the momentum equation:
\begin{eqnarray}
 \left( \int\limits_{\QQ_j}{\tilde{\psi}_k^{(j)}(\xx,t^{n+1}_{-}) \tilde{\psi}_l^{(j)}(\xx,t^{n+1}_{-}) \dx }   
- \int\limits_{\QQ_j^\st}{\diff{\tilde{\psi}_k^{(j)}}{t} \tilde{\psi}_l^{(j)} \dxt} \right) (\rhov)_{l,j}^{n+1}   
-  \int\limits_{\QQ_j}{\tilde{\psi}_k^{(j)}(\xx,t^{n}_{+}) \tilde{\psi}_l^{(j)}(\xx,t^{n}_{-}) \dx }   (\rhov)_{l,j}^{n}  & & \nonumber \\ 
+ \int\limits_{\QQ_j^\st}{\tilde{\psi}_k^{(j)} \nabla \cdot \mathbf{\TF}_h \dxt} 
+\int\limits_{\TT_{\ell(j),j}^\st}{\tilde{\psi}_k^{(j)} \nabla \tilde{\phi}_{l}^{(\ell(j))}  \dxt} \,  \, \hat \np_{l,\ell(j)}^{n+1}
+\int\limits_{\TT_{r(j),j}^\st}{   \tilde{\psi}_k^{(j)} \nabla \tilde{\phi}_{l}^{(r(j))}     \dxt} \,  \, \hat \np_{l,r(j)}^{n+1}    & &  \nonumber \\
 +\int\limits_{\Gamma_j^\st}{\tilde{\psi}_k^{(j)} \tilde{\phi}_{l}^{(r(j))}    \nstd ds dt} \,  \hat \np_{l,r(j)}^{n+1}
 -\int\limits_{\Gamma_j^\st}{\tilde{\psi}_k^{(j)} \tilde{\phi}_{l}^{(\ell(j))} \nstd ds dt} \,  \hat \np_{l,\ell(j)}^{n+1} = 0.  & & 
\label{eq:CS_11}
\end{eqnarray}

\noindent For every $j=1 \ldots \Nj$, Equation \eref{eq:CS_11} is then rewritten in a compact matrix form as
\begin{eqnarray}
    \Mpsi_j (\rhovh)_j^{n+1} - \Mpsi_j^-(\rhovh)_j^{n} + \widehat{\mathbf{G v}}_j^{n+1} + \RM_j \etah_{r(j)}^{n+1} - \LM_j \etah_{\ell(j)}^{n+1} =0, \label{eq:CS_12_1} 
\end{eqnarray}
where we have introduced several matrices:
\begin{eqnarray}
	  \Mpsi_j^+ = \int\limits_{\QQ_j}{\tilde{\psi}_k^{(j)}(\xx,t^{n+1}_{-}) \tilde{\psi}_l^{(j)}(\xx,t^{n+1}_{-})  \dx},  &&  
    \Mpsi_j^- = \int\limits_{\QQ_j}{\tilde{\psi}_k^{(j)}(\xx,t^{n}_{+}) \tilde{\psi}_l^{(j)}(\xx,t^{n}_{-})    \dx},   \\ 
    \Mpsi_j^\circ = \int\limits_{\QQ_j^\st}{\diff{\tilde{\psi}_k^{(j)}}{t} \tilde{\psi}_l^{(j)} \dx dt},  &&  
		\Mpsi_j  = \left(\Mpsi_j^+ - \Mpsi_j^\circ \right),  
\end{eqnarray}
%

\begin{equation}
	\RM_{j}=\int\limits_{\Gamma_j^\st}{\tilde{\psi}_k^{(j)} \tilde{\phi}_{l}^{(r(j))}\nstd ds dt}+\int\limits_{\TT_{r(j),j}^\st}{\tilde{\psi}_k^{(j)} \nabla \tilde{\phi}_{l}^{(r(j))}  \dx dt},
\label{eq:MD_4}
\end{equation}

\begin{equation}
	\LM_{j}=\int\limits_{\Gamma_j^\st}{\tilde{\psi}_k^{(j)} \tilde{\phi}_{l}^{(\ell(j))}\nstd ds dt}-\int\limits_{\TT_{\ell(j),j}^\st}{\tilde{\psi}_k^{(j)} \nabla \tilde{\phi}_{l}^{(\ell(j))}  \dx dt}.
\label{eq:MD_5}
\end{equation}
%
 
\noindent Following \cite{3DSIINS} the matrices $\LM$ and $\RM$ are then generalized by introducing a new matrix $\Q_{i,j}$, defined as
\begin{equation}
	\Q_{i,j}=\int\limits_{\TT_{i,j}^\st}{\tilde{\psi}_k^{(j)} \nabla \tilde{\phi}_{l}^{(i)}  \dxt}-\int\limits_{\Gamma_j^\st}{\tilde{\psi}_k^{(j)} \tilde{\phi}_{l}^{(i)} \sigma_{i,j} \nstd ds dt},
\label{eq:MD_6}
\end{equation}
where $\sigma_{i,j}$ is a sign function defined by
\begin{equation}
	\sigma_{i,j}=\frac{r(j)-2i+\ell(j)}{r(j)-\ell(j)}.
\label{eq:SD_1}
\end{equation}
In this way $\Q_{\ell(j),j}=-\LM_j$ and $\Q_{r(j),j}=\RM_j$, and then equation \eref{eq:CS_12_1} becomes in terms of $\Q$
\begin{equation}
	\Mpsi_j (\rhovh)_j^{n+1} - \Mpsi_j^-(\rhovh)_j^{n} + \widehat{\mathbf{Gv}}_j^{n+1}  + \Q_{r(j),j}  \etah_{r(j)}^{n+1} + \Q_{\ell(j),j} \etah_{\ell(j)}^{n+1} =0.
\label{eq:CS_12_2}
\end{equation}
\paragraph{Discrete energy equation} A weak form of the energy equation \eqref{eq:CNS_3_2} is obtained again by multiplication of \eqref{eq:CNS_3_2} by a space-time test function $\phit_k^{(i)}$ and integrating over a space-time control volume $\TT_i^\st$. For every $k=1\ldots N_\phi^\st$ we get 
\begin{eqnarray}
		\int\limits_{\TT_i^\st}{\phit_k^{(i)}\diff{(\rhoEE)_h}{t} \dxt} + \int\limits_{\TT_i^\st}{\phit_k^{(i)}\nabla \cdot (k \rho \vv)_h  \dxt} + \int\limits_{\TT_i^\st}{\phit_k^{(i)} \nabla \cdot  
		(H \rho \vv)_h \dxt}&=&\int\limits_{\TT_i^\st}{\phit_k^{(i)}\nabla \cdot \left[ \sigma  \vv + \lambda \nabla T \right]_h \dxt} \label{eq:CNS_15}.
\end{eqnarray}
The discrete enthalpy $H_h$, the momentum $(\rho \vv)_h$ and the stress tensor $\boldsymbol{\sigma}_h$ are already naturally given on the dual mesh. 
The discrete work of the stress tensor can thus be obtained by direct evaluation of the product as $\mathbf{w}_h = \boldsymbol{\sigma}_h \vv_h$. For the computation of the 
discrete heat flux vector $\mathbf{q}_h$ we follow the same strategy as for the discrete velocity gradient in the viscous stress tensor, thus obtaining 
\begin{eqnarray}
  \int\limits_{\QQ_j^\st}{\tilde{\psi}_k^{(j)} \mathbf{q}_j \dxt} = 
	 \!\! \int\limits_{\TT_{\ell(j),j}^\st} \!\! {\tilde{\psi}_k^{(j)} \lambda \nabla T_{\ell(j)} \dxt} + 
	 \!\! \int\limits_{\TT_{r(j),j}^\st}    \!\! {\tilde{\psi}_k^{(j)} \lambda \nabla T_{r(j)} \, \dxt} + 
	 \int\limits_{\Gamma_j^\st}{\tilde{\psi}_k^{(j)} \lambda \left(T_{r(j)}-T_{\ell(j)}\right) \nstd \, ds dt}.
\end{eqnarray} 
Note that the temperature $T$ is computed from the pressure and the density via the equation of state, which we assume to be linear in $p$. The weak form of the energy equation thus becomes 
\begin{eqnarray}
		\int\limits_{\TT_i^\st}{\phit_k^{(i)}\diff{(\rhoEE)_h}{t} \dxt} + \int\limits_{\TT_i^\st}{\phit_k^{(i)} \nabla \cdot (k \rho \vv)_h  \dxt} 
		+ \sum_{j \in S_i}\left( \int\limits_{\Gamma_j^\st} \phit_k^{(i)} (H \rho \vv)_j \cdot \vec{n}_{i,j} \, ds dt 
		-\int\limits_{\TT_{ij}^\st}{\nabla \phit_k^{(i)} \cdot (H \rho \vv)_j \dxt}  \right) && \\ 
		= \sum_{j \in S_i}\left( \int\limits_{\Gamma_j^\st} \phit_k^{(i)} \mathbf{w}_j \cdot \vec{n}_{i,j} \, ds dt 
		-\int\limits_{\TT_{ij}^\st}{\nabla \phit_k^{(i)} \cdot \mathbf{w}_j \, \dxt}  \right) 
		+ \sum_{j \in S_i}\left( \int\limits_{\Gamma_j^\st}   \phit_k^{(i)} \mathbf{q}_j \cdot \vec{n}_{i,j} \, ds dt 
		-\int\limits_{\TT_{ij}^\st}{\nabla \phit_k^{(i)} \cdot \mathbf{q}_j \, \dxt}  \right). 
			\label{eq:energy2}
\end{eqnarray}
In order to compute the convective term for the kinetic energy, we simply average the quantity $k$ from the dual to the main grid and apply the discrete convection operator \eqref{eq:der4}, 
which we denote by $\widehat{\mathbf{F} \mathbf{k}}^{n+1}_i$. 
By using the same reasoning as for the momentum equation and writing all the quantities in \eqref{eq:CNS_15} in terms of the basis functions we obtain the following compact expression for 
the discrete energy equation: 
\begin{eqnarray}
	\M_i(\rhoEEh)_i^{n+1} - \M_i^-(\rhoEEh)_i^n + \widehat{\mathbf{F} \mathbf{k}}^{n+1}_i + \sum_{j \in S_i}{\Dtilde_{ij} \hh_j^{n+1} (\rhovh)_j^{n+1}}
		&=& \sum_{j \in S_i}{\D_{ij} \left( \hat{\mathbf{w}}_j^{n+1}  + \hat{\mathbf{q}}_j^{n+1} \right) } 
\label{eq:CNS_17}
\end{eqnarray}
with
\begin{eqnarray}
	\Dtilde_{ij}&=&  \int\limits_{\Gamma_j^\st}{\phi_k^{(i)}\psi_l^{(j)}\psi_r^{(j)}\vec{n}_{i,j}ds dt}-\int\limits_{\TT_{ij}^\st}{\nabla \phi_k^{(i)}\psi_l^{(j)}\psi_r^{(j)}\dxt},
\label{eq:Dtilde}
\end{eqnarray}
%
\begin{eqnarray}
	\D_{ij}&=&  \int\limits_{\Gamma_j^\st}{\phi_k^{(i)}\psi_l^{(j)}\vec{n}_{i,j}ds \, dt}-\int\limits_{\TT_{ij}^\st}{\nabla \phi_k^{(i)}\psi_l^{(j)}\dxt}.
\label{eq:D}
\end{eqnarray}
Note that the matrix $\D_{ij}$ defined here is exactly the same as the one defined in the incompressible case for the divergence of the velocity field, see e.g. \cite{3DSIINS}.
For the multiplication of the three-dimensional tensor $\Dtilde_{ij}$ with the two vectors of degrees of freedom $\hh_j^{n+1}$ and $(\rhovh)_j^{n+1}$ we use the convention 
(see \cite{DumbserCasulli} and \cite{2DSIUSW}) 
\begin{equation}
   \Dtilde_{ij} \hh_j^{n+1} (\rhovh)_j^{n+1} = 
	  \left( \int\limits_{\Gamma_j^\st}{\phi_k^{(i)}\psi_l^{(j)}\psi_r^{(j)}\vec{n}_{i,j}ds dt}-\int\limits_{\TT_{ij}^\st}{\nabla \phi_k^{(i)}\psi_l^{(j)}\psi_r^{(j)}\dxt} \right)
		\hat{H}_{l,j}^{n+1} \widehat{\rho \mathbf{v}}_{r,j}^{n+1}, 
\end{equation} 
with summation over the repeated indices $l$ and $r$. 

\paragraph{Pressure system} 
The resulting discrete momentum and energy equations can be summarized as follows: 
\begin{eqnarray}
			\Mpsi_j (\rhovh)_j^{n+1} + \Q_{r(j),j}  \etah_{r(j)}^{n+1} + \Q_{\ell(j),j} \etah_{\ell(j)}^{n+1} 
			&=& \Mpsi_j^-(\rhovh)_j^{n} - \widehat{\mathbf{Gv}}_j^{n+1},
  \label{eq:CNS_18_1} \\
			\M_i(\rhoEEh)_i^{n+1} + \sum_{j \in S_i}{\Dtilde_{ij} \hh_j^{n+1} (\rhovh)_j^{n+1}}
		  &=& \M_i^-(\rhoEEh)_i^n - \widehat{\mathbf{F} \mathbf{k}}^{n+1}_i  + \sum_{j \in S_i}{\D_{ij} \left( \hat{\mathbf{w}}_j^{n+1}  + \hat{\mathbf{q}}_j^{n+1} \right) }. 
 		\label{eq:CNS_18_2} 
\end{eqnarray}
Here the pressure in the momentum equation is discretized implicitly as well as the momentum in the energy equation. 
In the case of $p_\gamma=0$ one can take $\etah_{\cdot}^{n+\theta}$ and $(\rhovh)_j^{n+\theta}$ in the momentum and energy equation, respectively, 
in order to recover second order of accuracy in time with a Crank-Nicolson time discretization by setting $\theta=0.5$. 
The density from the continuity equation is already available from \eqref{eqn.densitymain} and \eqref{eqn.densitydual} on the main grid and on the dual mesh, respectively.
Furthermore, the total energy $(\rhoEEh)_i^{n+1}=(\rhoeh)_i^{n+1}+(\rhokh)_i^{n+1}$ can be expressed in terms of the internal energy and the kinetic energy. 
Formal substitution of the discrete momentum equation \eqref{eq:CNS_18_1} into the discrete energy equation \eqref{eq:CNS_18_1} leads to the following expression:
\begin{eqnarray}
\M_i \left[ (\rhoeh)_i^{n+1}+(\rhokh)_i^{n+1} \right] - \sum_{j \in S_i}{\Dtilde_{ij} \hh_j^{n+1} \left[ \M_j^{-1}\left( \Q_{r(j)j} \pph_{r(j)}^{n+1}  + \Q_{\ell(j)j} \pph_{\ell(j)}^{n+1} \right) \right]}= \M_i^-(\rhoEEh)_i^n - \widehat{\mathbf{F} \mathbf{k}}^{n+1}_i  + \sum_{j \in S_i}{\D_{ij} \left( \hat{\mathbf{w}}_j^{n+1}  + \hat{\mathbf{q}}_j^{n+1} \right) } && \nonumber \\
 - \sum_{j \in S_i}{\Dtilde_{ij} \hh_j^{n+1} \left[ \M_j^{-1}\left( \Mpsi_j^-(\rhovh)_j^{n} - \widehat{\mathbf{Gv}}_j^{n+1} \right) \right]},
\label{eq:CNS_19_1}
\end{eqnarray} 
with $(\rhoeh)_i^{n+1} = \boldsymbol{\rho} \mathbf{e}\left( \hat{\mathbf{p}}_i^{n+1}, \hat{\boldsymbol{\rho}}_i^{n+1} \right)$, to be understood 
as a componentwise evaluation of the internal energy density $\rho e$ for each component of the input vectors. 
Due to the dependency of the enthalpy on the pressure, and due to the presence of the kinetic energy in the total energy, 
the system \eqref{eq:CNS_19_1} is highly nonlinear. In order to obtain a \textit{linear} system for the pressure, a simple but very efficient 
Picard iteration is used, as suggested in \cite{DumbserCasulli2016} and \cite{2STINS,3DSIINS}. 

\paragraph{Picard iteration and summary of the semi-implicit algorithm}
In \cite{CasulliZanolli2010} and \cite{DumbserCasulli2016}, but also in \cite{Dumbser2008,DumbserZanotti} the nonlinearities appearing in 
semi-implicit schemes and in locally implicit schemes for nonlinear PDE systems were successfully \textit{avoided} by the introduction 
of a simple Picard iteration. This means that certain nonlinear terms are discretized at the previous Picard iteration, and thus 
become essentially explicit. In the nonlinear pressure system \eqref{eq:CNS_19_1} above, we therefore take all convective and diffusive terms 
on the right hand side, as well as the enthalpy and the kinetic energy appearing on the left hand side at the previous (old) Picard iteration. 
The entire staggered semi-implicit space-time DG algorithm for the compressible Navier-Stokes equations \eqref{eq:CNS_1}-\eqref{eq:CNS_3}  
can therefore be summarized for each Picard iteration $m=1,\ldots, N_{Pic}$ as follows:  
\begin{eqnarray}
  	\Mpsi_i \hat{\boldsymbol{\rho}}_i^{n+1,m+1} = \Mpsi_i^- \hat{\boldsymbol{\rho}}_i^{n} - \widehat{\mathbf{F}\boldsymbol{\rho}}_i^{n+1,m}, 
		& 
		\qquad \hat{\boldsymbol{\rho}}_j^{n+1,m+1} = \bar{\MM}_j^{-1}\left( \LL^\top_{\ell(j),j} \hat{\boldsymbol{\rho}}_{\ell(j)}^{n+1,m+1} + \LL^\top_{r(j),j}\hat{\boldsymbol{\rho}}_{r(j)}^{n+1,m+1} \right), 
	   \label{eqn.conti.pic}
\end{eqnarray}
\begin{eqnarray}
\M_i \, \left[ {\boldsymbol{\rho}} \mathbf{e}\left( \hat{\mathbf{p}}_i^{n+1,m+1}, \hat{\boldsymbol{\rho}}_i^{n+1,m+1} \right)  \right] - \sum_{j \in S_i}{\Dtilde_{ij} \hh_j^{n+1,m} \left[ \M_j^{-1}\left( \Q_{r(j)j} \pph_{r(j)}^{n+1,m+1}  + \Q_{\ell(j)j} \pph_{\ell(j)}^{n+1,m+1} \right) \right]} = \M_i^-(\rhoEEh)_i^n && \nonumber \\ 
 - \M_i  (\rhokh)_i^{n+1,m} - \widehat{\mathbf{F} \mathbf{k}}^{n+1,m}_i  + \sum_{j \in S_i}{\D_{ij} \left( \hat{\mathbf{w}}_j^{n+1,m}  + \hat{\mathbf{q}}_j^{n+1,m} \right) } 
 - \sum_{j \in S_i}{\Dtilde_{ij} \hh_j^{n+1,m} \left[ \M_j^{-1}\left( \Mpsi_j^-(\rhovh)_j^{n} - \widehat{\mathbf{Gv}}_j^{n+1,m} \right) \right]},
\label{eqn.pressure.pic}
\end{eqnarray} 
\begin{eqnarray}
			\Mpsi_j (\rhovh)_j^{n+1,m+1} 
			&=& \Mpsi_j^-(\rhovh)_j^{n} - \widehat{\mathbf{Gv}}_j^{n+1} - \Q_{r(j),j}  \etah_{r(j)}^{n+1,m+1} - \Q_{\ell(j),j} \etah_{\ell(j)}^{n+1,m+1},
  \label{eqn.mom.pic} \\
			\M_i(\rhoEEh)_i^{n+1} 
		  &=& \M_i^-(\rhoEEh)_i^n - \widehat{\mathbf{F} \mathbf{k}}^{n+1}_i  - \sum_{j \in S_i}{\Dtilde_{ij} \hh_j^{n+1} (\rhovh)_j^{n+1,m+1}} + \sum_{j \in S_i}{\D_{ij} \left( \hat{\mathbf{w}}_j^{n+1,m}  + \hat{\mathbf{q}}_j^{n+1,m} \right) }. 
 		\label{eqn.energy.pic} 
\end{eqnarray}
As initial guess for all quantities at iteration $m=0$ we simply take the known data from the previous time $t^n$. Typically, we use $N_{Pic} \in \left\{ p_{\gamma}+1, p_{\gamma}+2 \right\}$, 
since for ODE the Picard iteration allows to gain one order in time per iteration. Numerical evidence allows us to conjecture that this behavior also holds for the present algorithm, 
but further theoretical analysis on this topic is still necessary in the future. The algorithm can thus be summarized in the following main steps: 
\begin{enumerate}
\item Solve the discrete continuity equation \eqref{eqn.conti.pic} for the new density on the primary and the dual mesh at the new Picard iteration $m+1$. 
\item Compute the nonlinear convective terms as well as the contribution of the viscous stress tensor to the momentum equation $\widehat{\mathbf{Gv}}_j^{n+1,m}$. 
      This is done by first averaging the momentum to the main grid and by subsequently defining the discrete viscous stress tensor on the dual mesh by \eqref{eqn.viscstress}. 
			This step corresponds to the use of a lifting operator, but acting on the \textit{dual mesh} rather than on the primary mesh. It also avoids the use of a 
			Riemann solver in the discretization of parabolic terms. Then, compute the discrete divergence operator on the main grid \eqref{eq:visconv} and subsequently 
			average back the result to the dual mesh \eqref{eq:visconv2}. 
\item Solve the linear system for the pressure given by \eqref{eqn.pressure.pic}. The scalar pressure field at the new Picard iteration $m+1$ is the only 
unknown therein, while all the other terms on the right hand side of \eqref{eqn.pressure.pic} are known from the previous Picard level $m$. 
\item Once the new pressure is known, the discrete momentum density is readily updated by Equation \eqref{eqn.mom.pic}. 
\item With the new pressure and the new momentum density, the total energy density can be updated by Equation \eqref{eqn.energy.pic}. 
\end{enumerate} 
The algorithm above consists in an \textit{explicit} discretization of all convective terms at the previous Picard iteration $m$,  
while the pressure terms are discretized \textit{implicitly} at the new iteration level $m+1$. Also the enthalpy $H$ and the kinetic energy in the 
total energy equation \eqref{eqn.pressure.pic} are taken at the old Picard iteration $m$. This is why we call the method a \textit{semi-implicit} 
scheme. All equations are written in a flux form at the discrete level, hence the scheme is locally and globally conservative for all conserved
quantities: mass, momentum and total energy. If necessary, the heat flux can be discretized implicitly by evaluating $\hat{\mathbf{q}}_j^{n+1,m+1}$
at the new iteration $m+1$. Also in the discretization of the nonlinear convective and viscous terms, the stress tensor can be taken implicitly
if needed, rewriting \eqref{eq:der4} and \eqref{eq:visconv} as 
\begin{equation} 
\Mpsi_i  \widehat{\boldsymbol{\rho} \vv}_{i}^{n+1,m+\frac{1}{2}} = \Mpsi_i^- \widehat{\boldsymbol{\rho} \vv}_{i}^{n} - \widehat{\mathbf{F} \boldsymbol{\rho} \vv}_i^{n+1,m} + 
\sum_{j \in S_i} \D_{ij} \left( \hat{\boldsymbol{\sigma}}_j^{n+1,m+\frac{1}{2}} \right), 
\end{equation} 
Being still embraced by the outer Picard loop, this approach corresponds to a simple but high order splitting, see \cite{3DSIINS} for details.  

Experimentally we also observe that the pressure system \eqref{eqn.pressure.pic} is symmetric and positive definite for the special case $p_\gamma=0$. 
Hence we can use the conjugate gradient (CG) method \cite{cgmethod} to solve the system. In contrast to the incompressible case \cite{Fambri2016}, the system as it 
is written here cannot be made symmetric in the general case $p_\gamma > 0$, due to a non symmetric contribution of the time derivatives contained
in $\M_i$ and $\M_j$. Hence, in this case we have to use a more general linear solver, such as the GMRES algorithm \cite{GMRES}.

\subsection{Time step restriction} 
Here we give a brief discussion about the CFL time restriction. Since we take the nonlinear convective term explicitly in order avoid the introduction of a non-linearity into the main system 
for the pressure, the maximum time step is restricted by a CFL-type restriction that can become rather severe for DG schemes. However, since the pressure is discretized
implicitly, the CFL condition is based only on the local flow velocity $\vv$ and not on the sound speed $c$: 
\begin{equation}
	\Delta t_{max}=\frac{\textnormal{CFL}}{2p+1}\cdot \frac{h_{\min}}{2|v_{\max}|},
\label{eq:CNS_38}
\end{equation}
with CFL$<1/d$ and where $d$ denotes the number of space dimensions. Furthermore, $h_{\min}$ is the smallest insphere diameter (in 3D) or incircle radius (in 2D) and $v_{max}$ is the maximum 
convective speed. 
In the case of an explicit discretization of the viscous stress tensor and the heat flux, we also must consider the eigenvalues of the viscous operator (see \cite{ADERNSE}) which are 
given by 
\begin{equation}
	\lambda_\nu=\max{\left(\frac{4}{3}\frac{\mu}{\rho},\frac{\gamma \mu}{Pr \rho}\right)},
\label{eq:CNS_38_2}
\end{equation} 
For large viscosities, the time step must obey a quadratic restriction proportional to $h_{\min}^2 / (2p+1)^2 \cdot 1 / \lambda_\nu$, and this condition on the maximum time step can become 
rather severe for large viscosities $\mu$. In particular, if we employ a limiter based on artificial viscosity, then this condition becomes too restrictive for our purpose. In this case,
an implicit treatment of the viscous stresses and the heat flux is necessary. 
However, the explicit treatment of the viscous terms remains a simple and computationally efficient solution for small viscosities and in the absence of the AV limiter.

\subsection{A posteriori DG limiter with artificial viscosity}
\label{sec_DGLIMITER}
Since we want to apply the staggered DG scheme also to high Mach number flows with shock waves, we need to introduce a limiter. 
In this section we describe how to apply the novel \textit{a posteriori} MOOD concept proposed in \cite{Clain2011,Diot2012,Diot2013}, 
which has already been successfully applied to DG schemes in \cite{Dumbser2014,Zanotti2015,DGLimiter3}, to the staggered DG method described in this paper. 
Since the sub-cell limiting procedure as used in \cite{Dumbser2014,Zanotti2015,DGLimiter3} is not easily generalizable to \textit{staggered} unstructured meshes, 
we opt for a simpler artificial viscosity method, combined with the \textit{a posteriori} detection criteria of the MOOD approach. The limiter will therefore be composed of two steps: 
the first one where we identify \textit{a posteriori} those cells where a certain set of detection criteria is not satisfied (troubled cells), and a second step where we apply the 
artificial viscosity limiter on those cells. In the following, we will denote the vector of conservative variables by $\mathbf{Q}=\left(\rho, \rho \mathbf{v}, \rho E \right)$. 

\subsubsection{Detection criteria} 
For each time step the algorithm described above produces a candidate solution $\mathbf{Q}_h^*(\xx,t^{n+1})$ at the new time $t^{n+1}$ that is based on the solution $\mathbf{Q}_h(\xx,t^n)$
at the old time $t^{n}$. 
Once we have the candidate solution, we can easily check \textit{a posteriori} the physical and numerical admissibility of this candidate solution. In particular, we require that the discrete 
solution satisfies some physical criteria such as the positivity of density and pressure, i.e. $\rho_h^*(\xx,t^{n+1})>0$ and $p_h^*(\xx,t^{n+1})>0$ $\forall \xx \in \Omega$. In addition, we require 
that it also satisfies a numerical admissibility criterion, namely a relaxed discrete version of the maximum principle. According to \cite{Dumbser2014,Zanotti2015,DGLimiter3} we use the following 
relaxed discrete maximum principle (DMP) component-wise for all the conserved variables $\mathbf{Q}$: 
\begin{eqnarray}
\min_{\mathbf{y} \in \mathcal{V}_i}\left(\mathbf{Q}_h(\mathbf{y},t^n) \right)-\delta \leq \mathbf{Q}_h^*(\xx,t^{n+1}) \leq \max_{\mathbf{y} \in \mathcal{V}_i}\left(\mathbf{Q}_h(\mathbf{y},t^n) \right)+\delta \qquad 
\forall \xx \in \TT_i ,
\label{eq:CNS_39}
\end{eqnarray}
where $\delta$ is a small relaxation parameter and the set $\mathcal{V}_i$ contains the current element $\TT_i$ and its Voronoi neighbors, i.e. those elements which share a common node with $\TT_i$.  
%
In our numerical experiments we use 
\begin{equation}
 \delta = \max \left( \delta_0, \epsilon \left( \max_{\mathbf{y} \in \mathcal{V}_i}\left(\mathbf{Q}_h(\mathbf{y},t^n) \right) - \min_{\mathbf{y} \in \mathcal{V}_i}\left(\mathbf{Q}_h(\mathbf{y},t^n) \right) \right) \right), 
\end{equation} 
with $\delta_0 = 2 \cdot 10^{-3}$ and $\epsilon= 10^{-3}$. 
Furthermore, since the evaluation of the global extrema of the solution inside each element $\TT_i$ can be rather expensive, we limit the evaluation of the maximum and minimum 
only to the known degrees of freedom stored inside each element. 

\subsubsection{Limiter}
If the detector described above is not activated anywhere in the entire domain $\Omega$, i.e. if the set containing the list of troubled cells  
$\mathcal{T}^{n+1} = \{\TT_i \in \Omega \,\, | \,\, \TT_i \mbox{ is troubled}\}$ is empty, then we accept the candidate  solution $\mathbf{Q}_h^*(\xx,t^{n+1})$ as is 
and set $\mathbf{Q}_h(\xx,t^{n+1})=\mathbf{Q}_h^*(\xx,t^{n+1})$. Otherwise, if the detector found some troubled cells and thus $\mathcal{T}^{n+1}$ is not empty, 
we discard the candidate solution $\mathbf{Q}_h^*(x,t^{n+1})$ and we add some artificial viscosity $\mu_{add,i}$ to the physical viscosity $\mu$ in the troubled cells
so that for all $\TT_i \in \mathcal{T}^{n+1}$ the effective viscosity used in the Navier-Stokes equations is $\mu_i=\mu+\mu_{add,i}$. The additional viscosity 
$\mu_{add,i}$ is computed in such a way that the mesh Reynolds number in the troubled cells becomes   
\begin{eqnarray}
	Re_i=\frac{\rho s_{max} h_s}{\mu_i}=1,
\label{eq:CNS_43}
\end{eqnarray}
where $h_s=\frac{h_i}{2p+1}$ is a rescaled mesh spacing based on the cell size $h_i$ (circumcircle or circumsphere diameter) and the polynomial degree $p$. 
Here $s_{max}=\left\| \mathbf{v} \right\| + c$ is the maximum eigenvalue of the hyperbolic part of the system, which contains the norm of the local flow velocity 
and the sound speed $c = \sqrt{\gamma p/\rho}$. The corresponding artificial heat conduction coefficient $\lambda_i$ in the heat flux is obtained by setting the 
Prandtl number to $Pr_i = \mu_i \gamma c_v / \lambda_i = 1$. Finally, we extend the new viscosity to all the dual elements as follows 
\begin{equation}
 \mu_j = \max\limits_{i \in [\ell(j), r(j)]} \mu_i.
\end{equation}
Note that the value of this artificial viscosity can become rather high, since now $s_{max}$ contains also the sound speed. Consequently, the time step restriction 
induced by \eqref{eq:CNS_38_2} for an explicit discretization of the heat flux and the viscous stresses becomes in general too severe and thus the implicit treatment 
of the viscous terms and of the heat flux becomes more advantageous. Once the artificial viscosity coefficients have been computed for the entire mesh, the time step 
is \textit{recomputed}, but with the artificial viscosity activated in all troubled cells. 

\section{Numerical tests}
\label{sec_NT}
\subsection{Two dimensional numerical tests}
In this section we report several classical test cases and academic benchmarks for the two dimensional compressible Navier-Stokes equations using the two dimensional version of 
the proposed staggered DG algorithm. The obtained numerical solutions are compared against exact or numerical reference solutions available in the literature. We will also use a 
complete isoparametric description of the geometry when it is necessary. For more details see e.g. \cite{3DSIINS}. If not stated otherwise, we will always assume that the ratio 
of specific heats is given by $\gamma=1.4$ for all test cases. 

\subsubsection{Isentropic vortex}
Here we consider a smooth two dimensional isentropic vortex, whose initial condition is given by a perturbation of a constant state as 
$\rho(\xx,0) = 1 + \delta \rho$, $\mathbf{v}(\xx,0) = \mathbf{v}_\infty + \delta \mathbf{v}$ and $p(\xx,0) = 1 + \delta p$, 
where $\mathbf{v}_\infty = (u_\infty,v_\infty)$ determines the velocity of the background field; the perturbations $\delta$ are given by 
\begin{eqnarray}
	\delta u= \frac{\epsilon}{2 \pi}e^{\frac{1-r^2}{2}}(-(y-y_0)), & & 
	\delta v= \frac{\epsilon}{2 \pi}e^{\frac{1-r^2}{2}}(x-x_0), \nonumber \\
	\delta \rho= (1+\delta T)^{\frac{1}{\gamma-1}}-1, & & 
	\delta p= (1+\delta T)^{\frac{1}{\gamma-1}}-1, \qquad \delta T=-\frac{(\gamma-1)\epsilon^2}{8\gamma \pi^2}e^{1-r^2},  
\label{eq:NT_1_2}
\end{eqnarray}
where $\gamma=1.4$ is the ratio of specific heats, $r=\sqrt{(x-x_0)^2+(y-y_0)^2}$ is the radial coordinate and $\epsilon=5$ is the vortex strength. 
The exact solution is simply given by a rigid body translation of the initial state with velocity $\mathbf{v}_\infty$. For this test we take $\Omega=[-10,10]^2$ and apply 
periodic boundary conditions everywhere. Furthermore we set $\mathbf{v}_\infty=(1,1)$, $(x_0,y_0)=(0,0)$, $\mu=0$ and $t_{\textnormal{end}}=1.0$. We consider a set of 
successively refined grids and polynomial degrees of $p=p_\gamma=1 \ldots 4$. 
\begin{table} 
\caption{$L_2$ error norms $\epsilon_p$, $\epsilon_u$  and $\epsilon_\rho$ and related numerical convergence rates 
$\mathcal{O}_p$, $\mathcal{O}_u$ and $\mathcal{O}_\rho$ obtained for the pressure, the velocity component $u$ and the density $\rho$ 
for the isentropic vortex problem using polynomial degrees in space and time from one to four.}  
\begin{center}%
\begin{tabular}{c ccc p{3mm}p{3mm}p{3mm} | c ccc p{3mm}p{3mm}p{3mm}}
\hline 
				\multicolumn{7}{c|}{$p_\gamma=p=1$}  &  \multicolumn{7}{c}{$p_\gamma=p=2$} \\  
\hline 
$N_i$ &	$\epsilon_p$ & $\epsilon_u$ & $\epsilon_\rho$  & $\mathcal{O}_p$ & $\mathcal{O}_u$ & $\mathcal{O}_\rho$ & 
$N_i$ &	$\epsilon_p$ & $\epsilon_u$ & $\epsilon_\rho$  & $\mathcal{O}_p$ & $\mathcal{O}_u$ & $\mathcal{O}_\rho$ \\ 
\hline 
5684  &  3.5E-02 & 3.3E-02 &  1.7E-02  &     &     &       &    5684   & 5.5E-03 & 6.4E-03 & 3.5E-03  &     &     &      \\
7684  &  2.7E-02 & 2.5E-02 &  1.3E-02  & 1.9 & 1.9 & 1.8   &    7424   & 4.5E-03 & 5.0E-03 & 2.8E-03  & 1.5 & 1.8 & 1.7  \\
9396  &  2.1E-02 & 2.0E-02 &  1.0E-02  & 2.0 & 2.0 & 1.8   &    9396   & 3.7E-03 & 4.0E-03 & 2.3E-03  & 1.7 & 1.8 & 1.7  \\
11600 &  1.7E-02 & 1.6E-02 &  8.8E-03  & 1.9 & 1.9 & 1.8   &    11600  & 3.0E-03 & 3.2E-03 & 1.9E-03  & 1.9 & 2.0 & 1.7   \\
\hline 
				\multicolumn{7}{c|}{$p_\gamma=p=3$}  &  \multicolumn{7}{c}{$p_\gamma=p=4$} \\  
\hline 
$N_i$ &	$\epsilon(p)$ & $\epsilon(u)$ & $\epsilon(\rho)$  & $\mathcal{O}_p$ & $\mathcal{O}_u$ & $\mathcal{O}_\rho$ & 
$N_i$ &	$\epsilon(p)$ & $\epsilon(u)$ & $\epsilon(\rho)$  & $\mathcal{O}_p$ & $\mathcal{O}_u$ & $\mathcal{O}_\rho$ \\ 
\hline 
1856  & 	 2.4E-03 & 3.0E-03 &  1.7E-03  &      &     &     &     116   & 6.7E-02  &  9.1E-02  &  4.5E-02	&       &     &      \\  
2900  & 	 1.1E-03 & 1.3E-03 &  7.7E-04  & 3.6  & 3.6 & 3.5 &     464   & 5.5E-03  &  7.6E-03  &  4.0E-03  &  3.6  & 3.6 & 3.5 \\ 
4176  & 	 5.4E-04 & 6.5E-04 &  3.7E-04  & 4.0  & 3.9 & 3.9 &    1044   & 1.1E-03  &  1.4E-03  &  8.1E-04  &  3.9  & 4.1 & 3.9 \\
5684  & 	 3.0E-04 & 3.6E-04 &  2.1E-04  & 3.8  & 3.8 & 3.8 &    1856   & 3.4E-04  &  4.2E-04  &  2.4E-04  &  4.2  & 4.2 & 4.2 \\
\hline 
	\label{tab:NT_1}
\end{tabular}
\end{center}
\end{table} 
In Table \ref{tab:NT_1} the error in standard $L_2$ norm is reported for the pressure, the velocity field and the density, using polynomial approximation degrees in space and time $p=p_\gamma$ from $1$ to $4$. In this case we observe an optimal convergence rate for odd polynomial degrees and a suboptimal one for even polynomial degrees. One can also check that by taking $p_\gamma=0$ the  resulting order tends to $1$, as expected. In other words, for this smooth unsteady problem it is very important to use high order polynomials also in time, since the time error significantly affects 
the total error. On the contrary, the resulting matrix for the pressure in the case of $p_\gamma>0$ becomes non-symmetric and this issue cannot be fixed using the reasoning of \cite{Fambri2016}. 
Hence, one has to use less efficient iterative linear solvers, such as the GMRES algorithm \cite{GMRES}. In that case also the theoretical analysis carried out in  
\cite{Capizzano2016} for the incompressible case does not hold any more. 
%
\subsubsection{Smooth acoustic wave propagation in 2D}
\label{sec_2Dsmooth}
Here we consider another simple 2D test problem with smooth solution, characterized by an acoustic wave traveling in radial direction. The aim of this test case is to verify the correct 
propagation speed of sound waves, since the unstructured grid used here is \textit{not} orthogonal, unlike the unstructured meshes used in classical low order semi-implicit finite volume or finite 
difference schemes \cite{CasulliWalters2000,BoscheriDumbser}. A further difficulty in this test case is that the effective Mach number is very low and the time step size is chosen 
according to a CFL condition based on the flow velocity and not according to a CFL condition based on the sound speed. 
The initial condition is given by $\rho(\xx,0)=1$, $\mathbf{v}(\xx,0)=0$ and $p(\xx,0) = 1 + e^{-\alpha r^2}$, with $r^2=x^2+y^2$. Due to the angular symmetry of the problem, we can obtain a 
reference solution by simply solving an 
equivalent one-dimensional PDE in the radial direction with a geometrical source term, see e.g. \cite{toro-book}. In particular, the computation of the reference solution was performed 
using a second order TVD finite volume scheme with the Osher-type flux \cite{OsherUniversal} on a radial grid consisting of $10^4$ cells. The computational domain is given by $\Omega=[-2,2]^2$, 
with periodic boundary conditions everywhere. Furthermore, we set $\alpha=40$; $(p, p_\gamma)=(3,0)$; $\mu=0$; $N_i=5616$ and $t_{end}=1.0$. 
\begin{figure}[ht!]
    \begin{center}
   \includegraphics[width=0.48\textwidth]{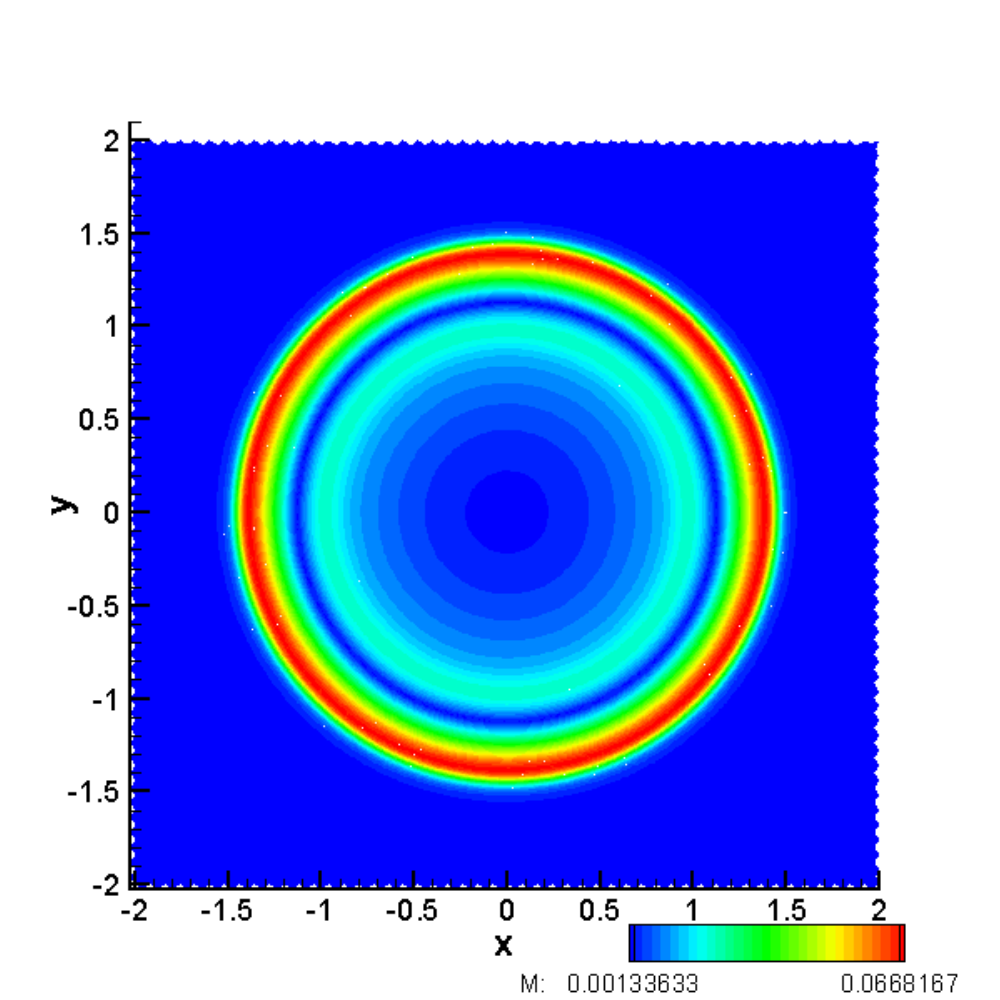}
   \includegraphics[width=0.48\textwidth]{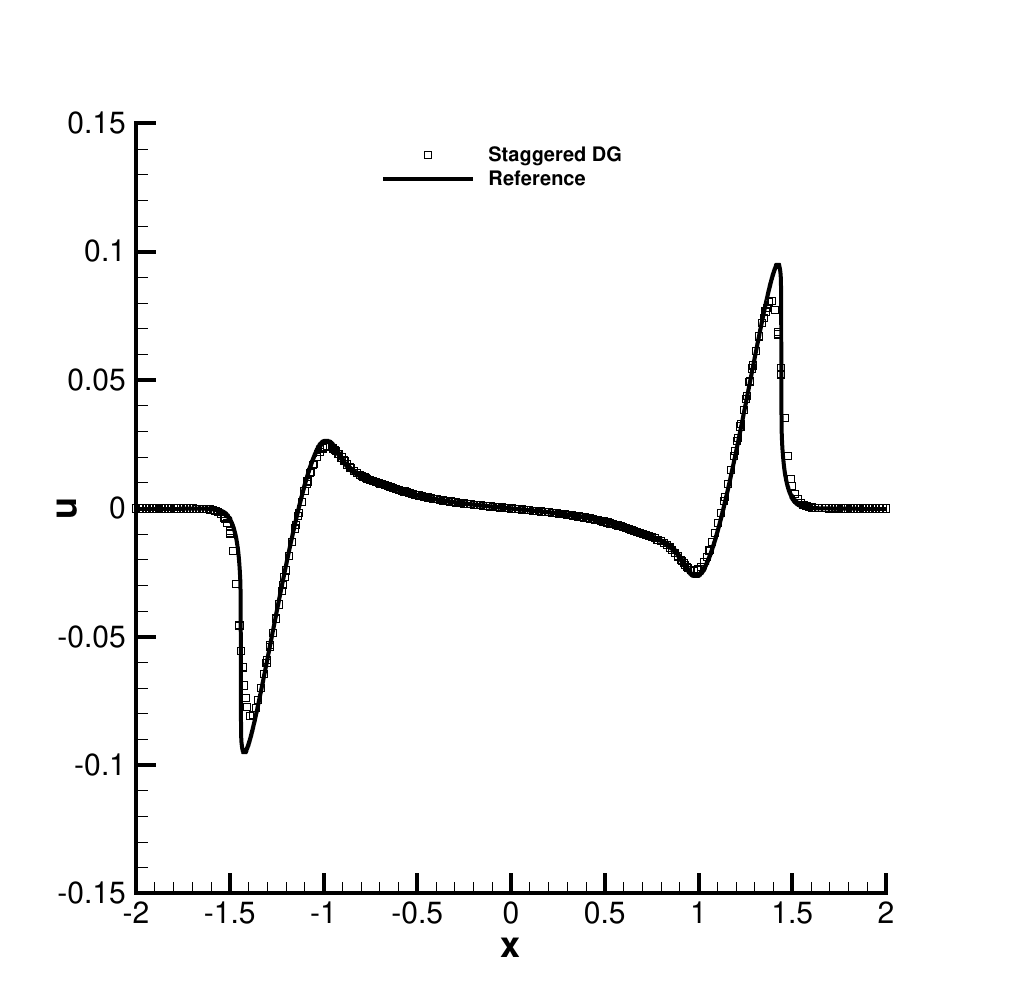}
   \includegraphics[width=0.48\textwidth]{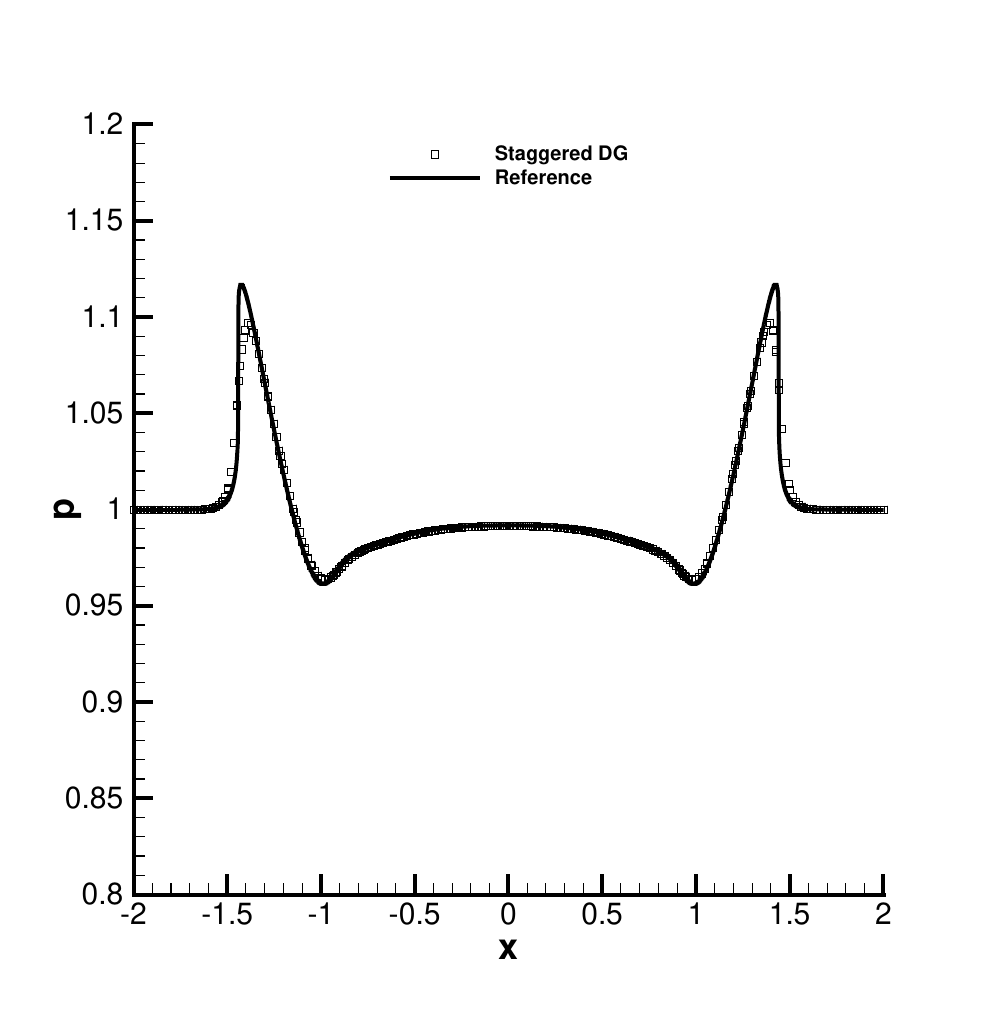}
   \includegraphics[width=0.48\textwidth]{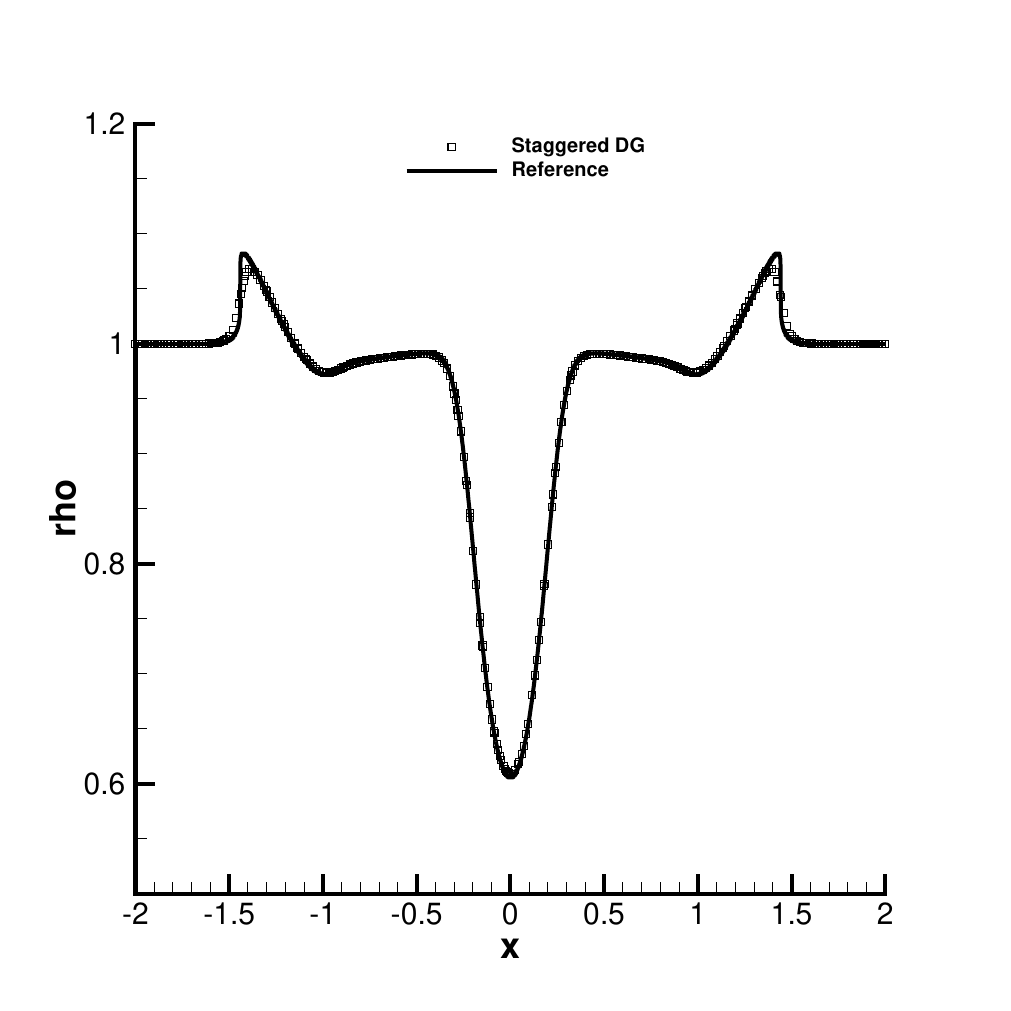}
    \caption{Numerical solution for the smooth acoustic wave propagation problem in 2D at time $t=1.0$ and comparison against the reference solution. 
		From top left to bottom right: Mach number contours, velocity $u$, pressure $p$ and density $\rho$.} 
    \label{fig.NT_2_1}
	\end{center}
\end{figure}
Figure $\ref{fig.NT_2_1}$ shows the obtained numerical results for the Mach number contours in the domain $\Omega$ and a comparison of the numerical solution along the line $y=0$ with the reference solution for all the main variables of the system. We can also observe in the plot of the Mach number that the resulting solution is symmetric with respect to the angular direction, which is not 
trivial, since the grid is fully unstructured and not symmetric. A good agreement between the reference solution and our numerical results can be observed, in particular we can see how the position 
of the acoustic wave is well reproduced by our numerical scheme. 

\subsubsection{Heat conduction}
\begin{figure}[ht!]
    \begin{center}
   \includegraphics[width=0.34\textwidth]{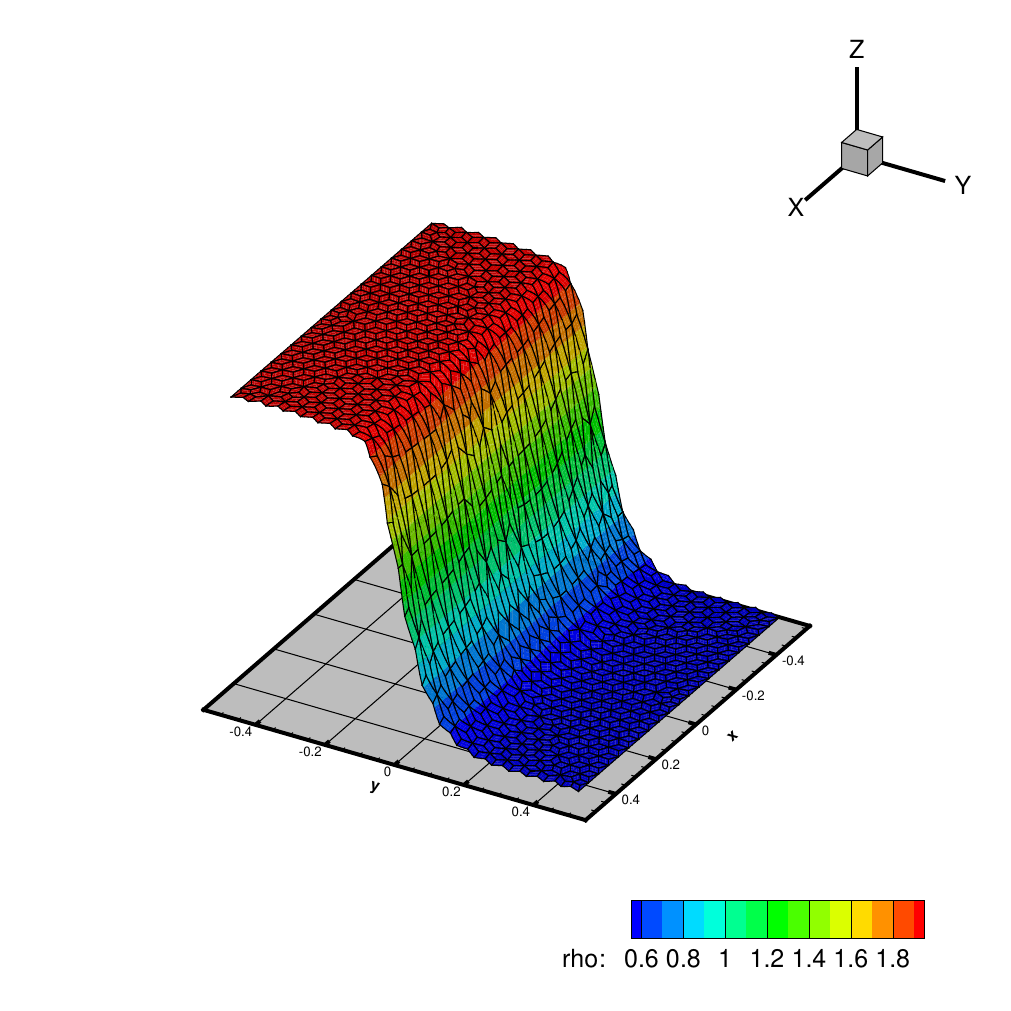}
   \includegraphics[width=0.3\textwidth]{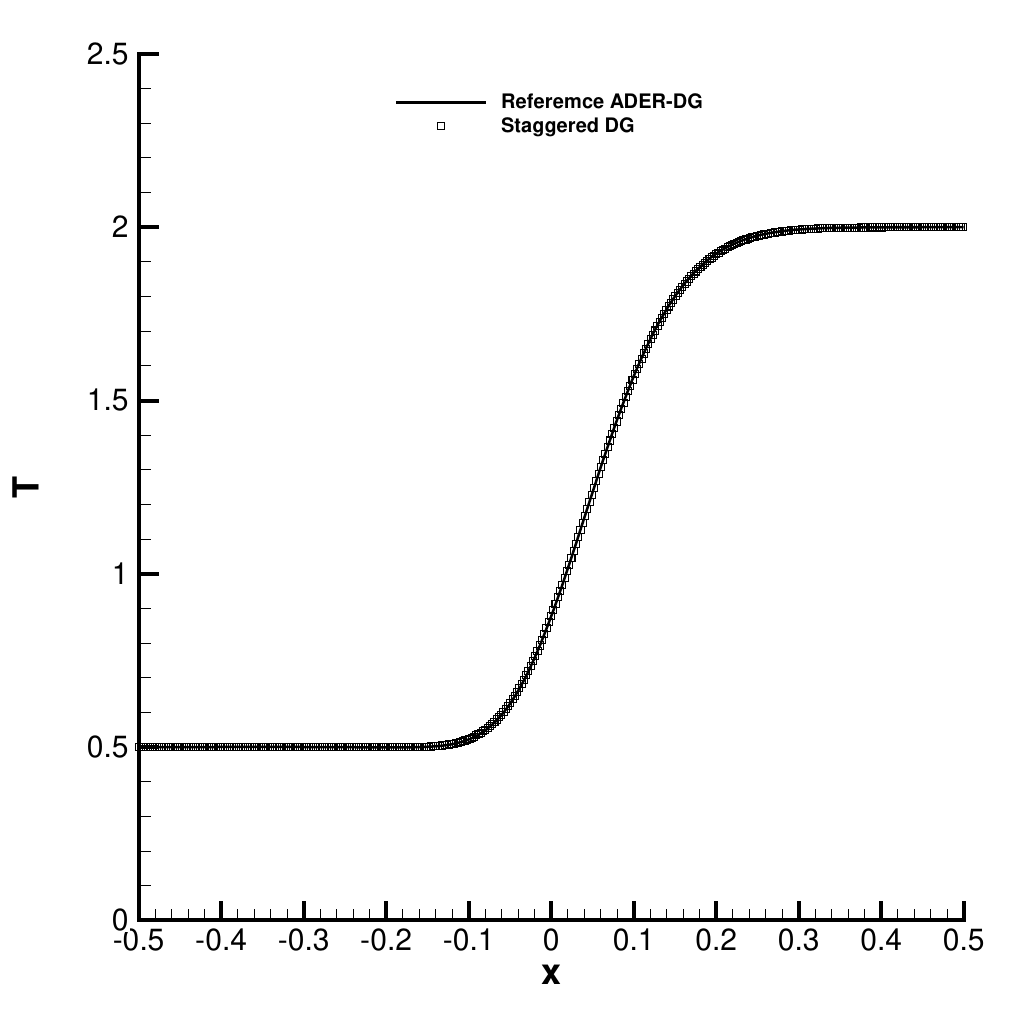}
   \includegraphics[width=0.3\textwidth]{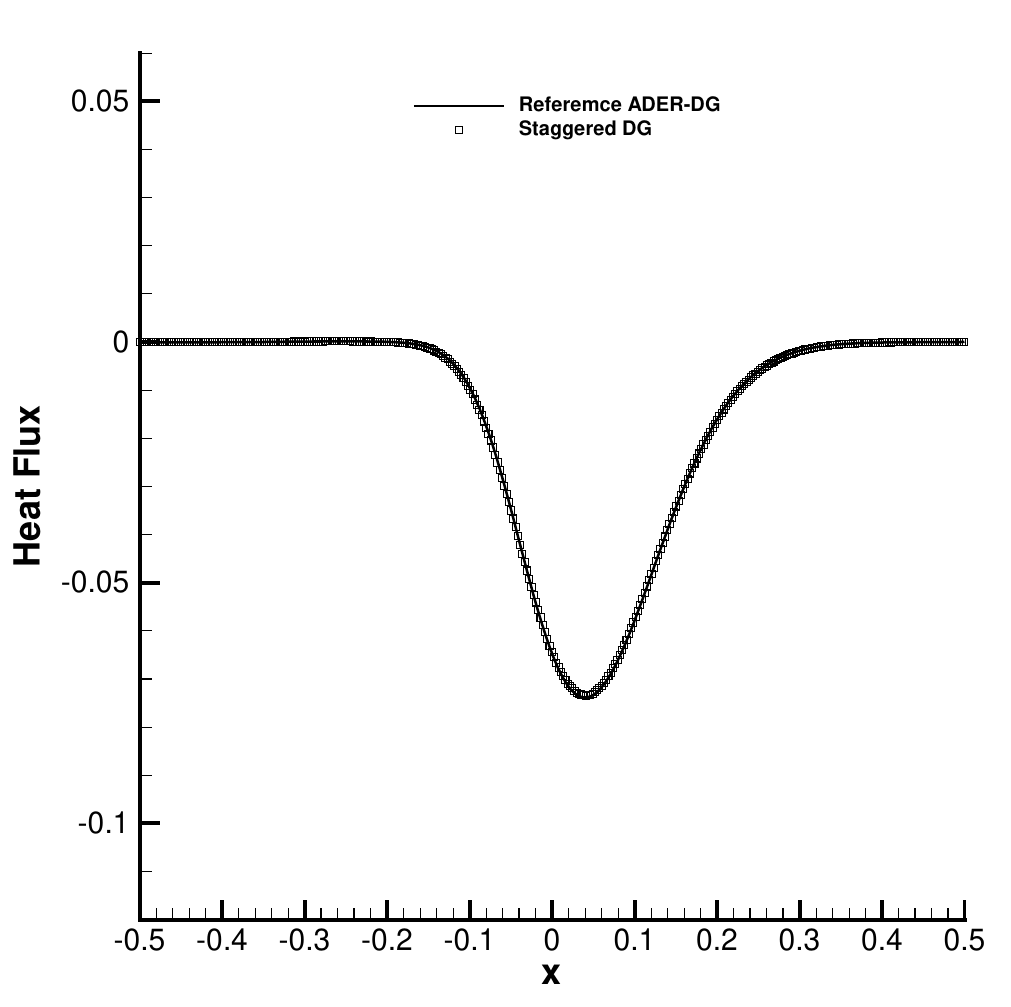}
    \caption{Heat conduction test case at time $t=1$. From left to right: 3D plot of the density and the dual mesh; temperature profile and heat flux component $q_1=-\lambda T_x$
		along the line $y=0$ compared with the 1D reference solution.}  
    \label{fig.NT_3_1}
	\end{center}
\end{figure}

Up to now we have checked the algorithm only in the absence of viscosity and heat conduction. Here we want to study a simple test problem dominated by heat transfer via heat conduction. 
The chosen initial condition is taken from \cite{DumbserPeshkov2016} and is briefly summarized below: 
\begin{equation}
\rho(\xx,0) = \left\{
\begin{array}{ll}
	2   & x <0 \\
	0.5 &  x \geq 0
\end{array}\right.
\label{eq:NT_3_1}
\end{equation}
$\mathbf{v}(\xx,0) = 0$ and $p(\xx,0)=1$. The parameters of the gas are $\gamma=1.4$, $c_v=2.5$, $\mu=10^{-2}$ and $\lambda=10^{-2}$. The computational domain is given by $\Omega=[-0.5,0.5]^2$ and the reference solution is the one used in \cite{DumbserPeshkov2016} computed using the ADER-DG scheme \cite{ADERNSE}. Finally we use $(p, p_\gamma)=(3,0)$, $t_{end}=1.0$ and the domain is discretized with 
$\Ni=902$ triangular elements. Figure $\ref{fig.NT_3_1}$ shows the numerical results obtained with the proposed staggered DG scheme on the computational domain $\Omega$ with a plot of the dual 
grid as well as a comparison with the reference solution on the cut along the line $y=0$ for the temperature, the heat flux component $q_1=-\lambda T_x$, and the density. Overall, a good agreement 
is obtained also in this simple problem dominated by heat conduction. 

\subsubsection{Viscous shock}
\label{sec_NT_VS}
\begin{figure}[ht!]
    \begin{center}
   \includegraphics[width=0.96\textwidth]{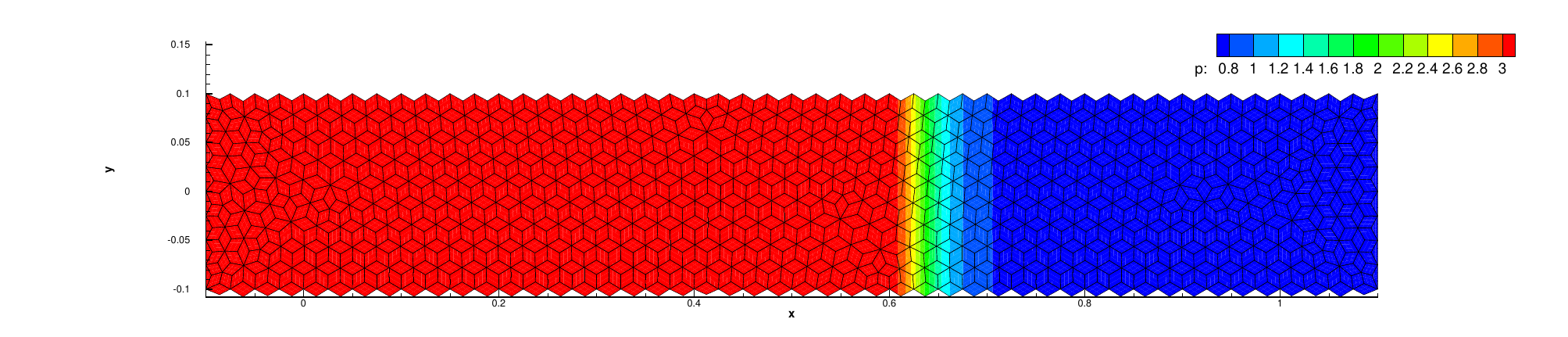}
    \caption{Viscous shock test problem. Sketch of the dual mesh and final pressure contours $p$ at $t=t_{end}$.} 
    \label{fig.NT_4_1}
	\end{center}
\end{figure}

\begin{figure}[ht!]
    \begin{center}
   \includegraphics[width=0.32\textwidth]{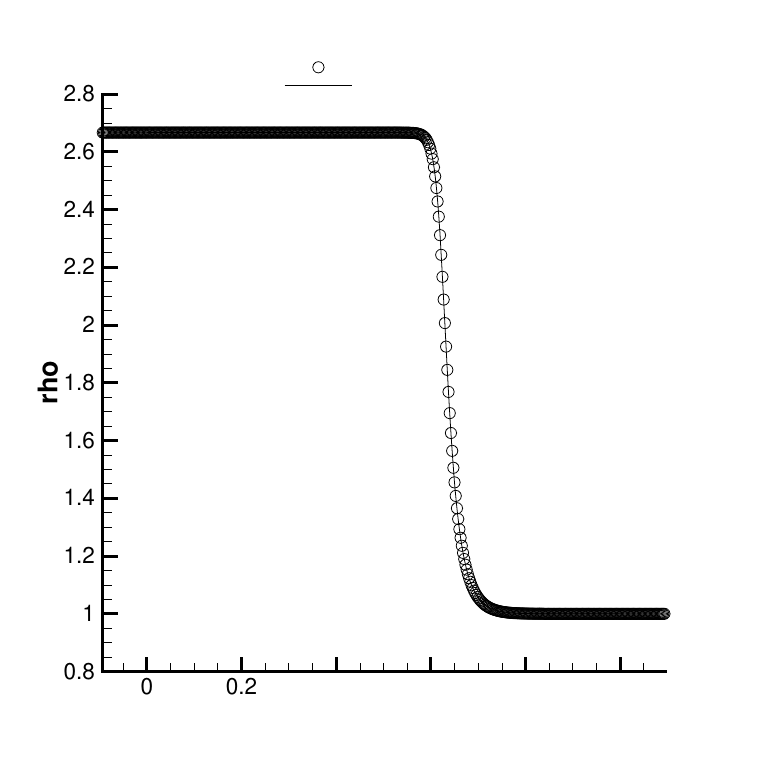}
	\includegraphics[width=0.32\textwidth]{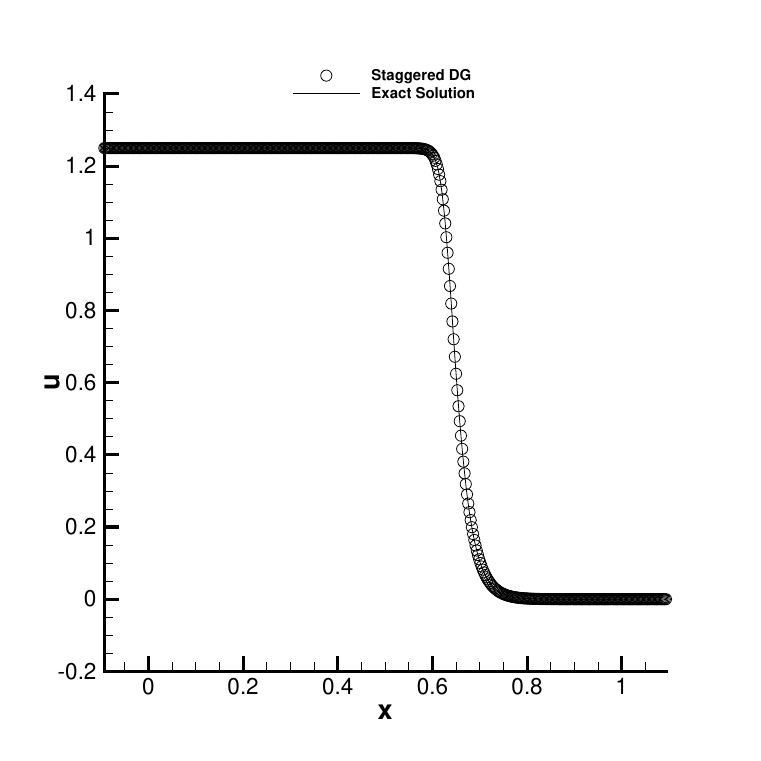}
	\includegraphics[width=0.32\textwidth]{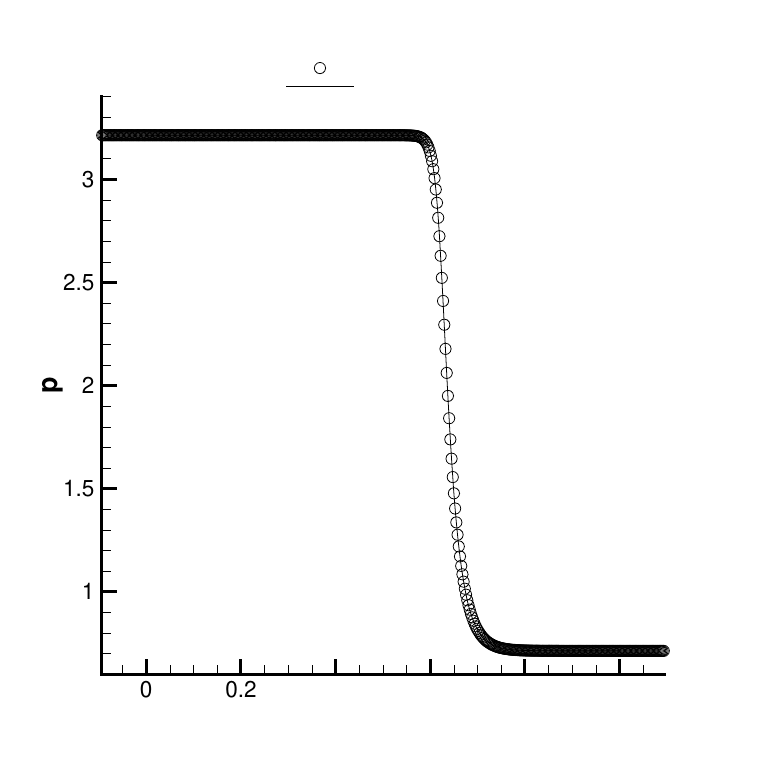}
    \caption{Viscous shock test problem  at time $t=1$. Comparison of the main flow variables $\rho$, $u$, and $p$ with the exact solution at $t=t_{end}$.}  
    \label{fig.NT_4_2}
	\end{center}
\end{figure}

\begin{figure}[ht!]
    \begin{center}
   \includegraphics[width=0.32\textwidth]{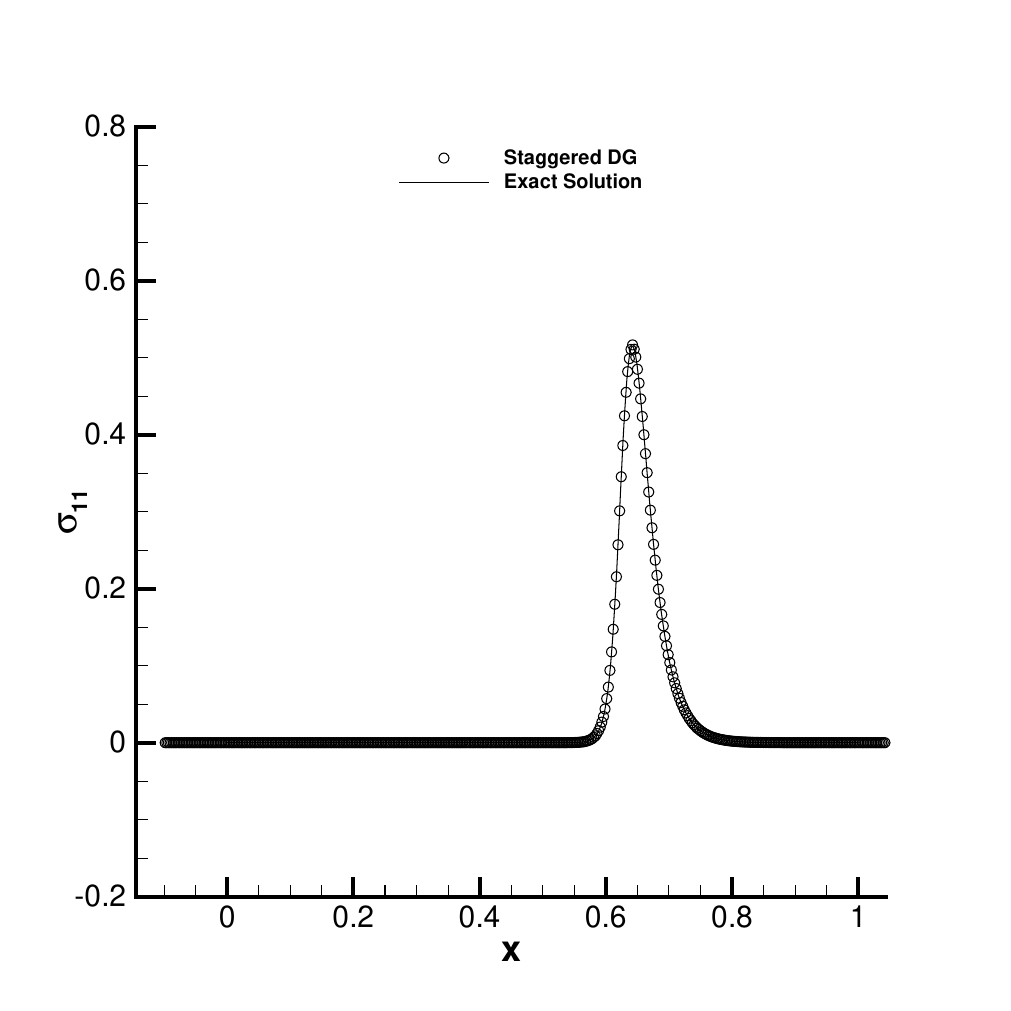}
	\includegraphics[width=0.32\textwidth]{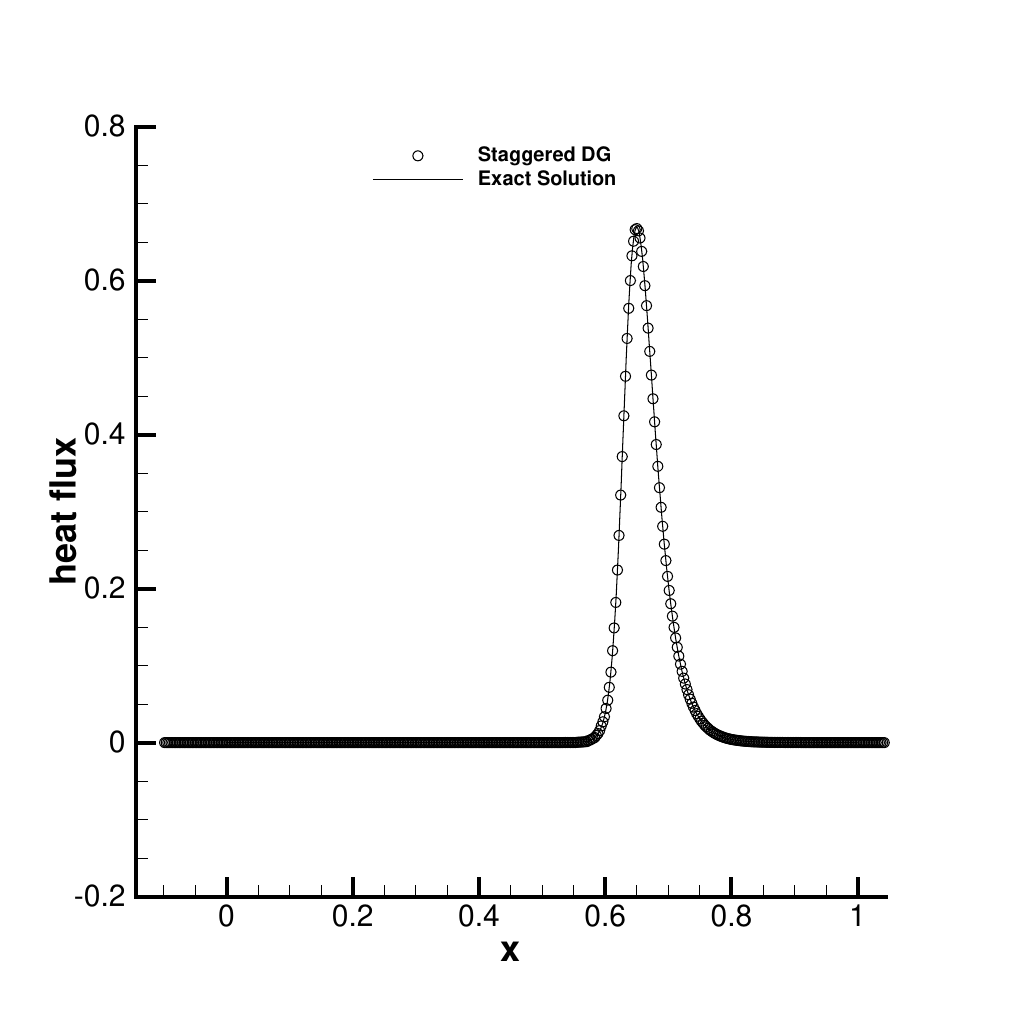}
    \caption{Viscous shock test problem at $t=1$. Comparison of the viscous stress and the heat flux with the exact solution.} 
    \label{fig.NT_4_3}
	\end{center}
\end{figure}

All previous examples involved only low Mach number flows. Here we consider still a smooth test problem, but in the supersonic regime, i.e. $M>1$. 
In the case of a Prandtl number of Pr$=0.75$, there exists an exact traveling wave solution of the compressible Navier-Stokes equation that was first 
found by Becker \cite{Becker1923} and that consists in a viscous shock profile. 
For the special case of a stationary viscous shock wave at Prandtl number $Pr=0.75$ and constant viscosity,  
the compressible Navier-Stokes equations can be reduced to one single ordinary differential 
equation that can be integrated analytically. The exact solution for the 
dimensionless velocity $\bar u = \frac{u}{M_s \, c_0}$ of this stationary shock wave with shock Mach number $M_s$ is then given by the root of 
the following equation, see \cite{Becker1923,BonnetLuneau,ADERNSE}:
\begin{equation} 
 \label{eqn.alg.u} 
  \frac{|\bar u - 1|}{|\bar u - \lambda^2|^{\lambda^2}} = \left| \frac{1-\lambda^2}{2} \right|^{(1-\lambda^2)} 
  \exp{\left( \frac{3}{4} \textnormal{Re}_s \frac{M_s^2 - 1}{\gamma M_s^2} x \right)},
\end{equation}
with
\begin{equation}
  \lambda^2 = \frac{1+ \frac{\gamma-1}{2}M_s^2}{\frac{\gamma+1}{2}M_s^2}.
\end{equation}
From eqn. \eqref{eqn.alg.u} one obtains the dimensionless velocity $\bar u$ as a function of $x$. 
The form of the viscous profile of the dimensionless pressure $\bar p = \frac{p-p_0}{\rho_0 c_0^2 M_s^2}$ is given by
the relation 
\begin{equation}
 \label{eqn.alg.p} 
  \bar p = 1 - \bar u +  \frac{1}{2 \gamma}
                         \frac{\gamma+1}{\gamma-1} \frac{(\bar u - 1 )}{\bar u} (\bar u - \lambda^2).  
\end{equation}
Finally, the profile of the dimensionless density $\bar \rho = \frac{\rho}{\rho_0}$ is found from the 
integrated continuity equation: $\bar \rho \bar u = 1$. In order to obtain an unsteady shock wave traveling 
into a medium at rest, it is sufficient to superimpose a constant velocity field $u = M_s c_0$ to the solution of the 
stationary shock wave found in the previous steps.
Our test problem consists in a shock wave initially centered in $x=0.25$ and traveling with a shock Mach number of $M_s=2.0$ in a domain $\Omega=[-0.1,1.1]\times[-0.1,0.1]$. 
The fluid parameters are $\gamma=1.4$, $\mu=2\cdot 10^{-2}$ and $Pr=0.75$. For this test we use $(p,p_\gamma)=(3,2)$, $t_{end}=0.2$ and a grid consisting of only $\Ni=870$ elements. 
Periodic boundary conditions are applied in the $y$-direction. Furthermore, we use an implicit discretization of the viscous stress tensor and of the heat flux. A sketch of the dual 
grid and resulting pressure profile at $t_{end}$ is shown in Figure $\ref{fig.NT_4_1}$. Figure $\ref{fig.NT_4_2}$ shows the obtained numerical results compared with the exact solution 
for the main variables $\rho$, $u$ and $p$. Finally, in Figure $\ref{fig.NT_4_3}$ we report the value of the first component of the viscous stress tensor $\sigma_{11}$ and the classical 
Fourier heat flux compared with the exact solution. A very good agreement can be observed also for these quantities which involve derivatives of the main flow variables.   
Note that this benchmark problem is very useful to test new numerical methods since it involves all terms of the compressible Navier-Stokes system \eqref{eq:CNS_1_2}-\eqref{eq:CNS_3_2}, 
in particular viscosity, heat conduction, nonlinear convection and the pressure forces. 

\subsubsection{Compressible mixing layer}
\begin{figure}[!htbp]
\begin{center}
\begin{tabular}{l} 
\includegraphics[width=0.7\textwidth]{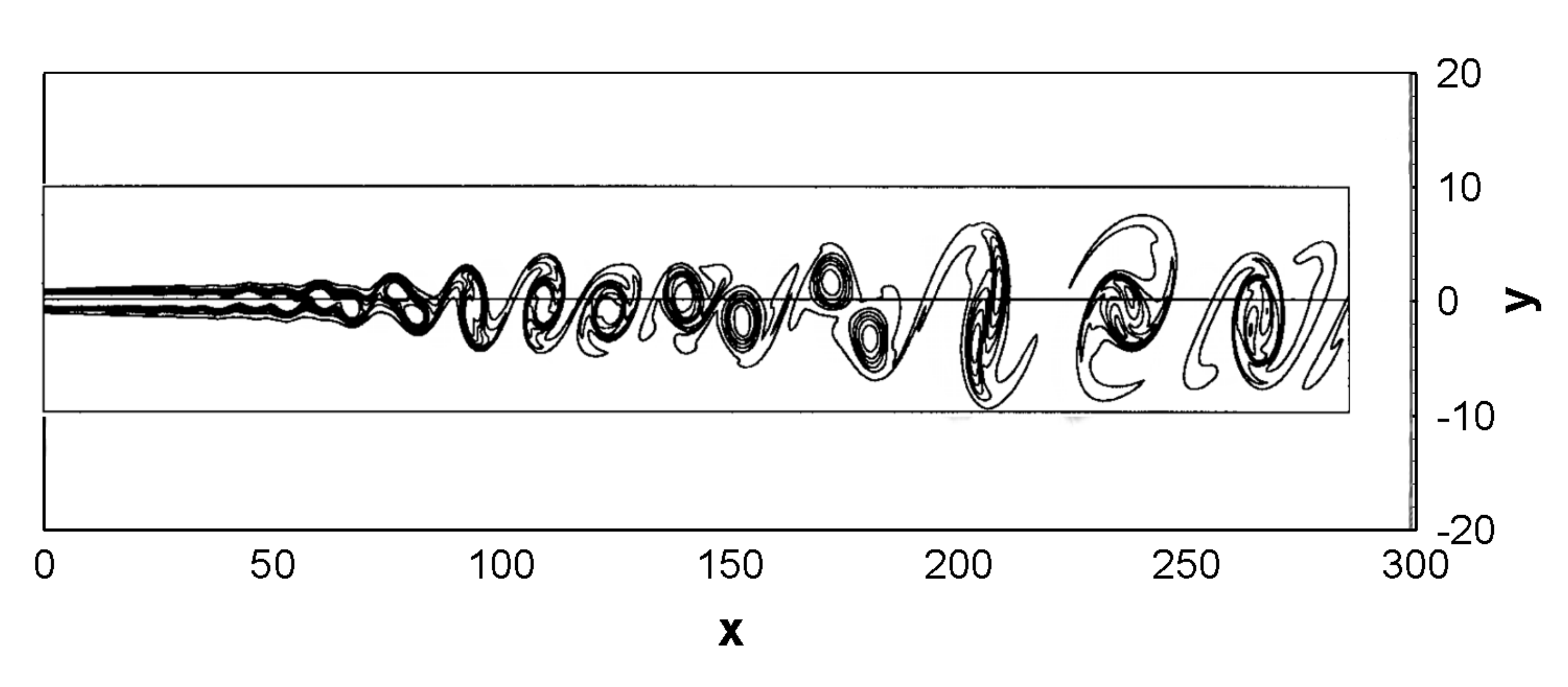}       \\ 
\includegraphics[width=0.7\textwidth]{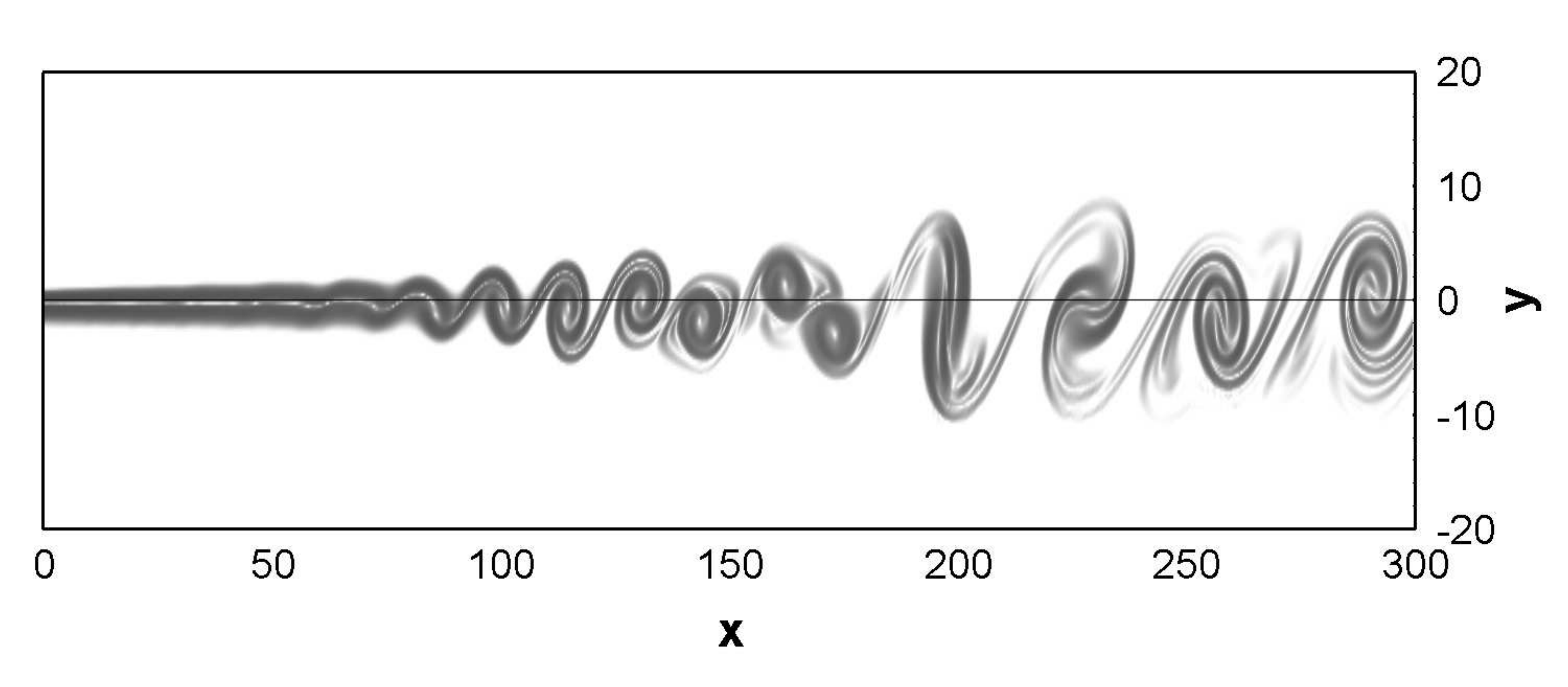} \\  
\includegraphics[width=0.68\textwidth]{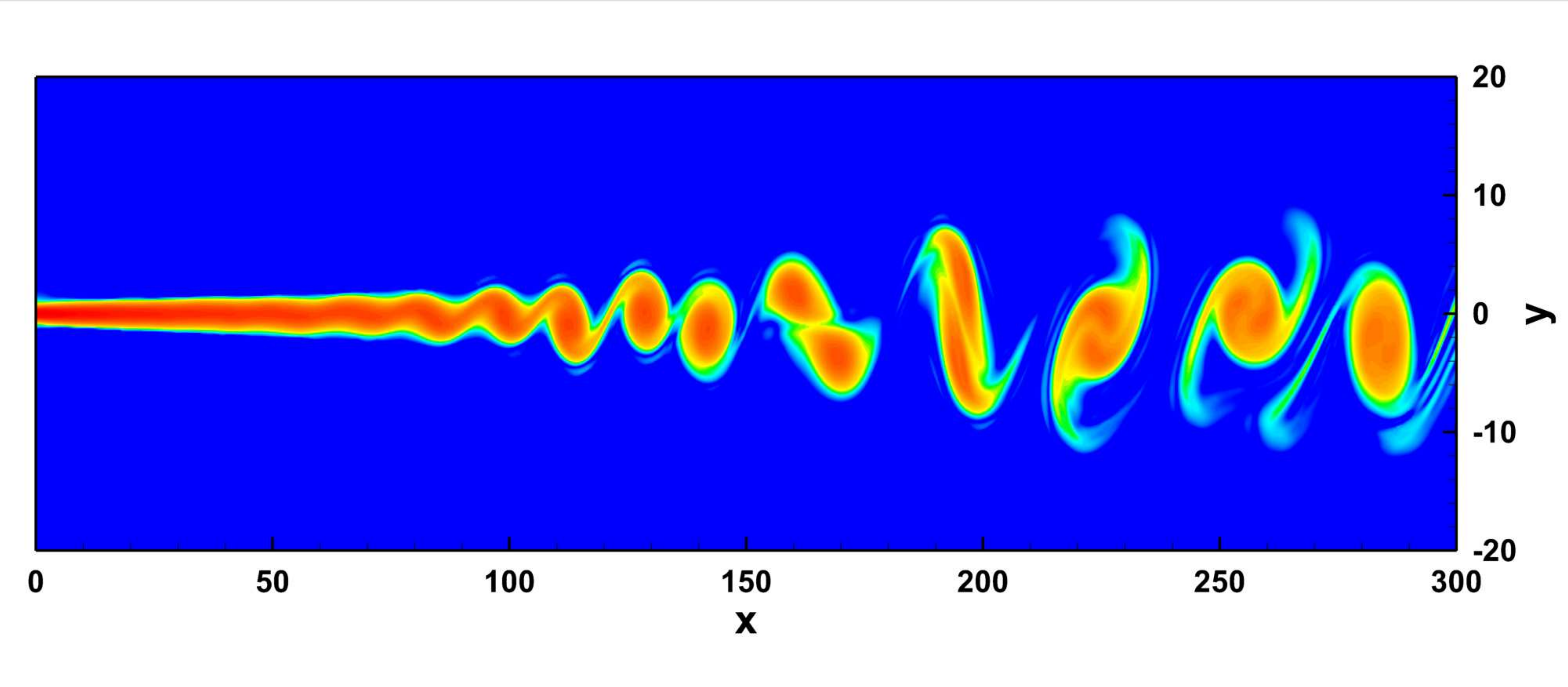}         \\ 
\includegraphics[width=0.665\textwidth]{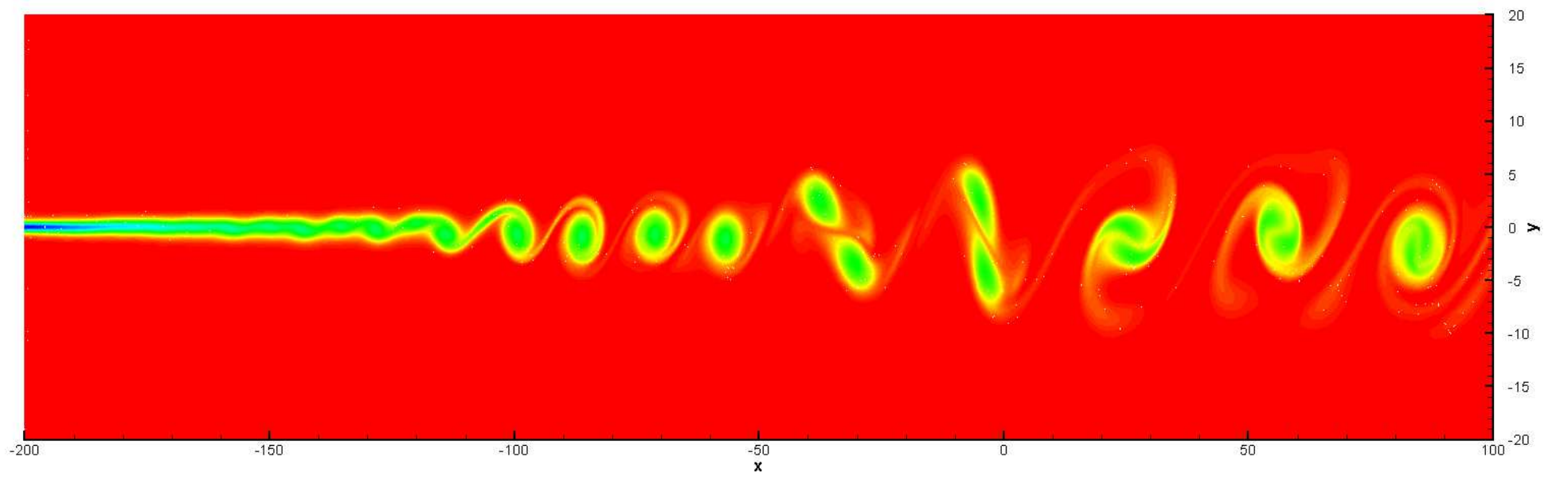}  
\end{tabular} 
\caption{Numerical results for the vorticity distribution in the compressible mixing layer at $t=t_{end}$. From top to bottom: Solution obtained for the compressible Navier-Stokes equations with a compact finite difference scheme  \cite{colonius}; solution obtained with a sixth order $P_NP_M$ scheme ($P_3P_5$), see \cite{ADERNSE}; third order ADER-WENO finite volume scheme applied to the unified Godunov-Peshkov-Romenski (GPR) 
model of continuum mechanics, see \cite{DumbserPeshkov2016}; the proposed semi-implicit staggered DG scheme. 
} 
\label{fig.NT_5_1}
\end{center}
\end{figure}

Here we consider an unsteady compressible mixing layer, originally studied by Colonius et al in \cite{colonius} and recently also investigated in \cite{ADERNSE, DumbserPeshkov2016}. 
The initial setup consists in two fluid layers moving with different velocities. For $y \rightarrow +\infty$ we set $u_{+\infty}=0.5$ while for $y \rightarrow -\infty$ we set $u_{-\infty}=0.25$. 
In the domain $\Omega$ we impose a smooth transition between the two velocities. The initial condition is given by 
\begin{equation*}
	u(\xx,0) = u_0 = \frac{1}{8}\tanh (2y) +\frac{3}{8}, \qquad v(\xx,0) = v_0 = 0, \qquad \rho(\xx,0) = \rho_0 = 1, \qquad p(\xx,0) = p_0 = \frac{1}{\gamma}. 
\label{eq:NT_5_1}
\end{equation*}
The vorticity thickness and its associated Reynolds number are given such as in \cite{DumbserPeshkov2016} as
\begin{eqnarray}
\theta = \frac{u_\infty-u_{-\infty}}{ \max{ \left( \left. \frac{\partial u}{\partial y} \right|_{x=0} \right)} } := 1 \qquad
\textnormal{Re}_{\theta} = \frac{\rho_0  u_\infty \theta}{\mu} = 500.
\label{eq:NT_5_2}
\end{eqnarray}
According to \cite{DumbserPeshkov2016} we introduce a perturbation at the inflow boundary at $x=0$ given by 
\begin{equation}
	\rho(0,y,t) = \rho_0 + 0.05 \delta(y,t), \qquad 
	\mathbf{v}(0,y,t) = \mathbf{v}_0 + \left(
	\begin{array}{c}
		1.0 \\ 0.6
	\end{array} \right) \delta(y,t), \qquad 
p(0,y,t) = p_0 + 0.2\delta(y,t), \nonumber
\label{eq:NT_5_3}
\end{equation}
where $\delta=\delta(y,t)$ is a periodic function given by
\begin{equation}
  \delta(y,t) = -10^{-3} \exp( -0.25 y^2)   \left(    \cos \left( \omega t \right) + \cos \left( \frac{1}{2} \omega t - 0.028 \right) + \cos \left( \frac{1}{4} \omega t + 0.141 \right) + \cos \left( \frac{1}{8} \omega t + 0.391 \right) \right),
\end{equation}
and $\omega=0.3147876$. The fluid parameters are set to $\gamma=1.4$, $\mu=10^{-3}$ and $Pr=+\infty$ in order to neglect the heat conduction. The computational domain is given by $\Omega=[-200,200]\times [-50,50]$ covered with $\Ni=19742$ triangles. We set $(p,p_\gamma)=(3,0)$ and $t_{end}=1596.8$. In this test problem, the viscous terms are discretized explicitly. The vorticity pattern obtained with 
our scheme at the final time is compared with some  reference solutions \cite{colonius,ADERNSE,DumbserPeshkov2016} in Figure \ref{fig.NT_5_1}. A good qualitative agreement between our numerical solution 
and the other reference solutions can be observed.  
\subsubsection{Two dimensional lid-driven cavity flow}
\label{sec.2DCavity}
\begin{figure}[ht!]
    \begin{center}
   \includegraphics[width=0.32\textwidth]{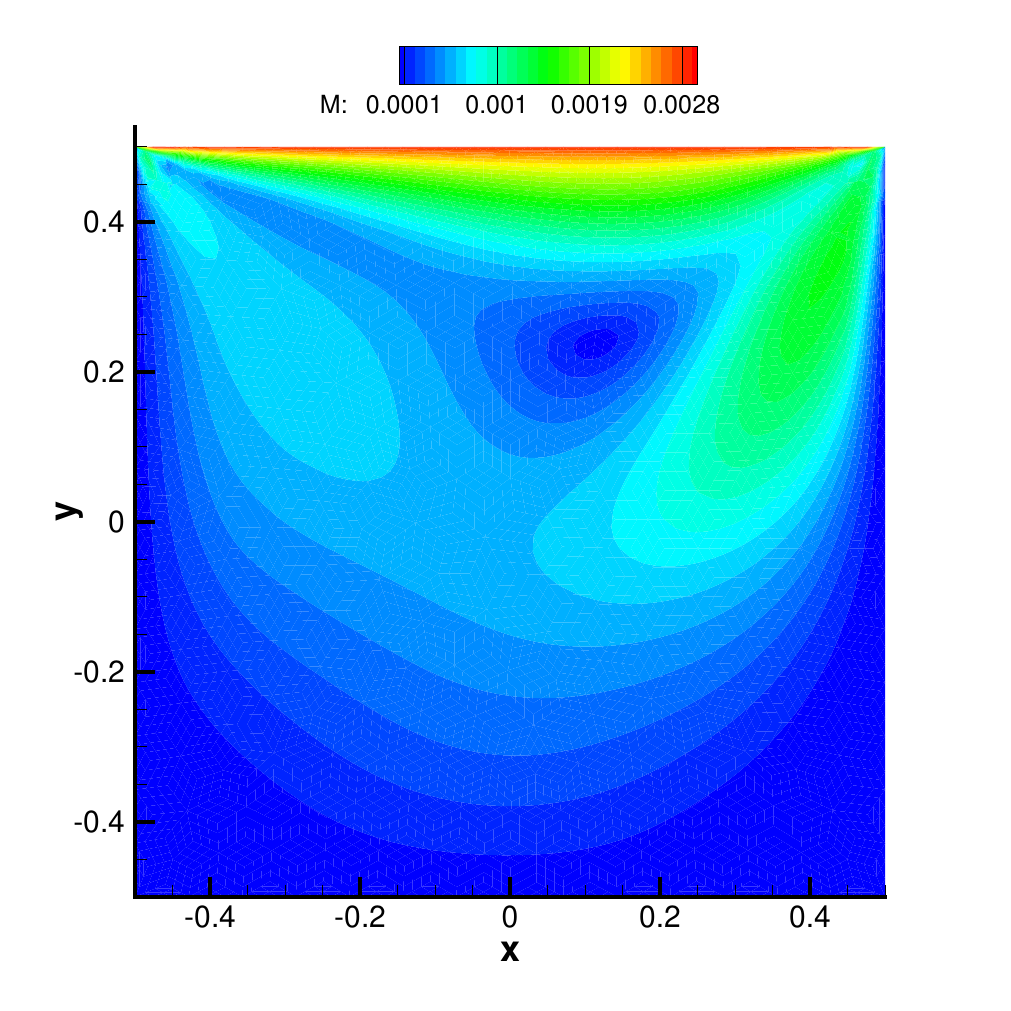}
	\includegraphics[width=0.32\textwidth]{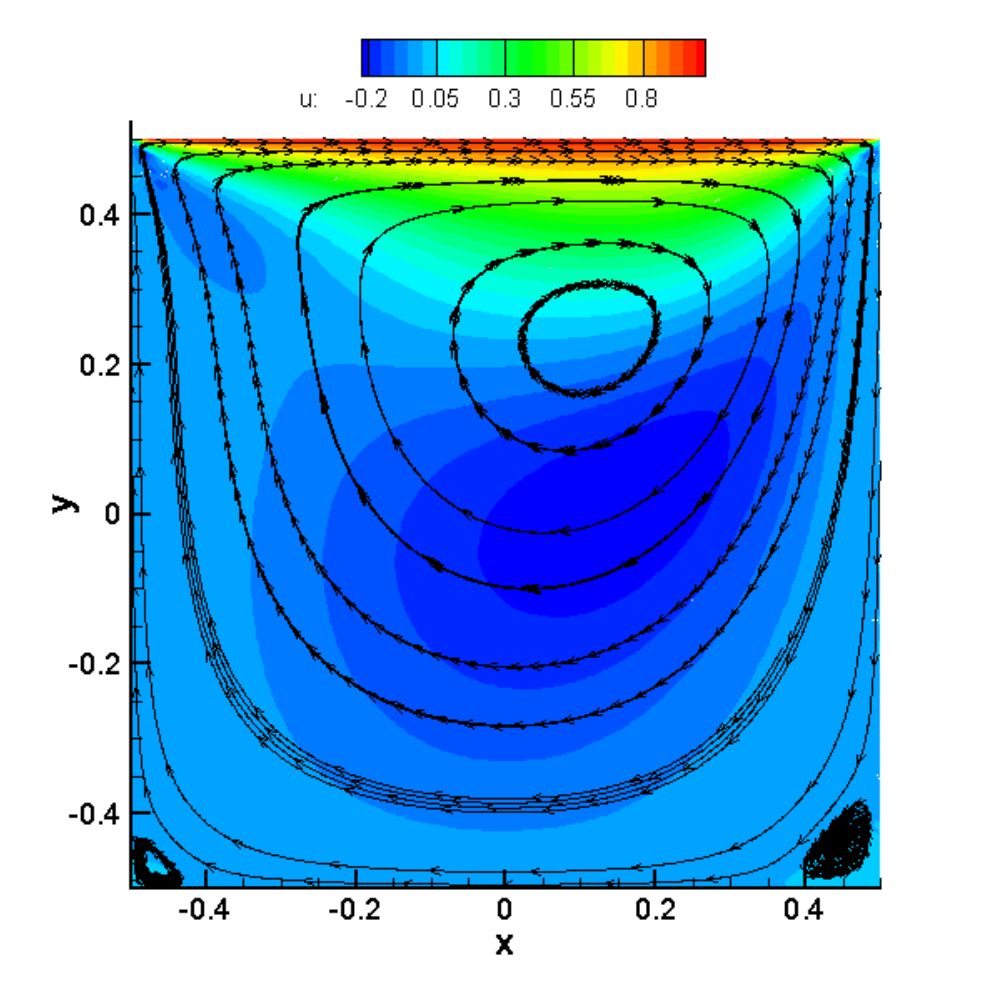}
	\includegraphics[width=0.32\textwidth]{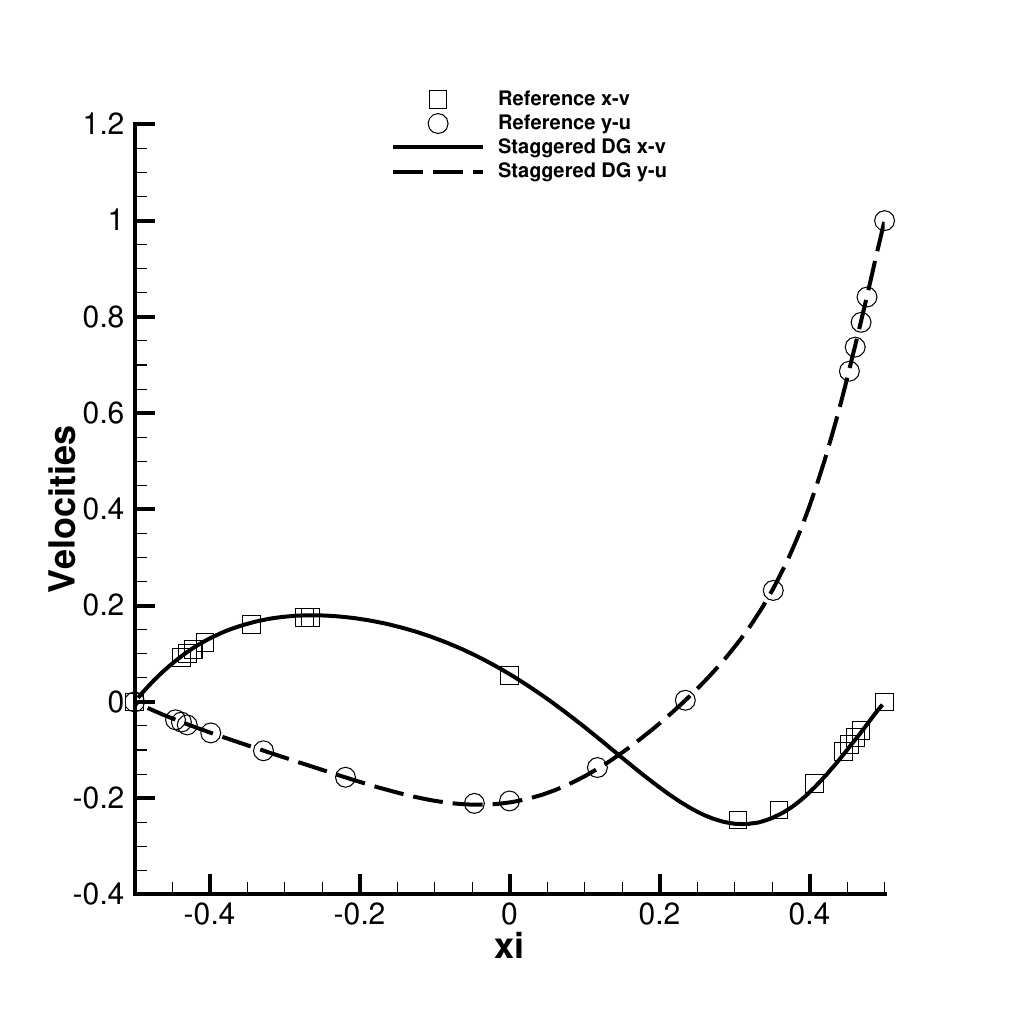}
    \caption{Lid driven cavity flow. From left to right: Mach number contours; horizontal velocity with streamlines; velocity profiles compared with the data given by Ghia et al in \cite{Ghia1982} at $Re=100$.}
    \label{fig.NT_6_1}
	\end{center}
\end{figure}

\begin{figure}[ht!]
    \begin{center}
   \includegraphics[width=0.32\textwidth]{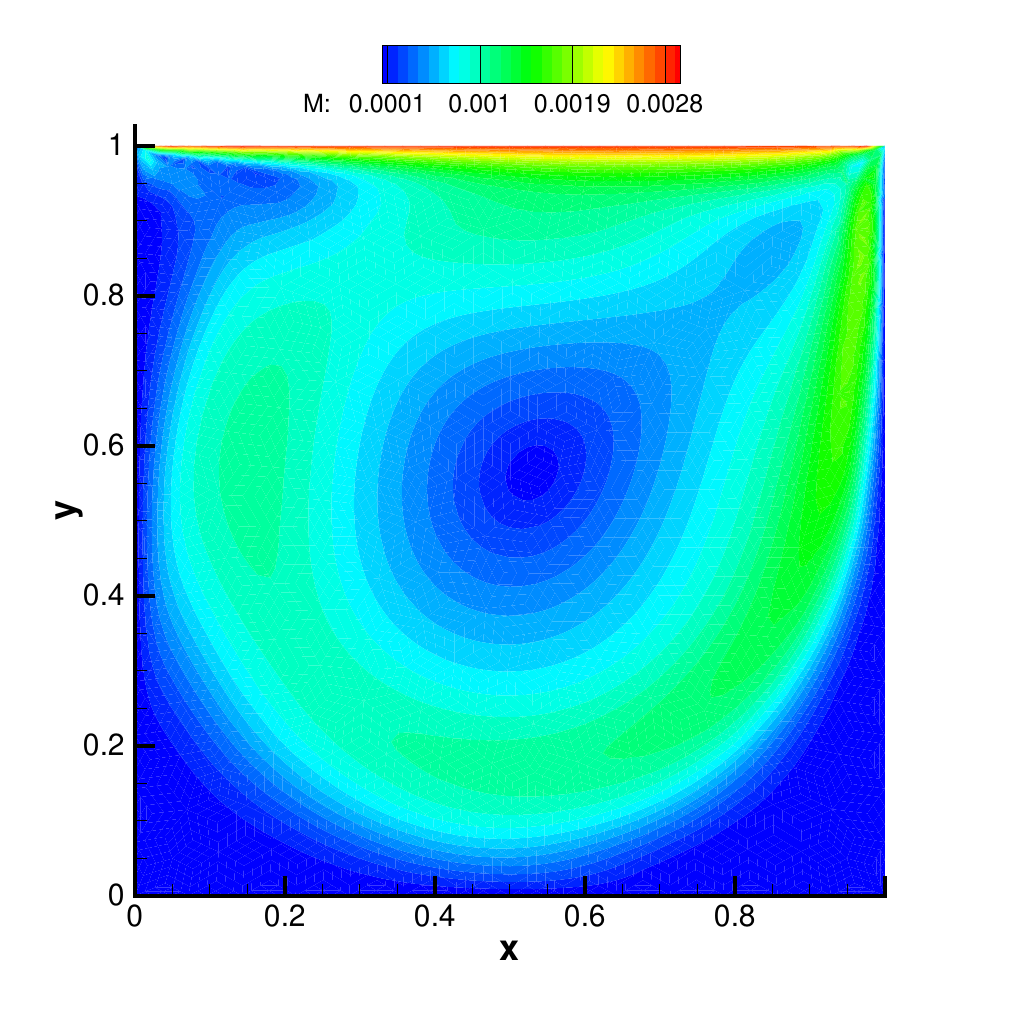}
	\includegraphics[width=0.32\textwidth]{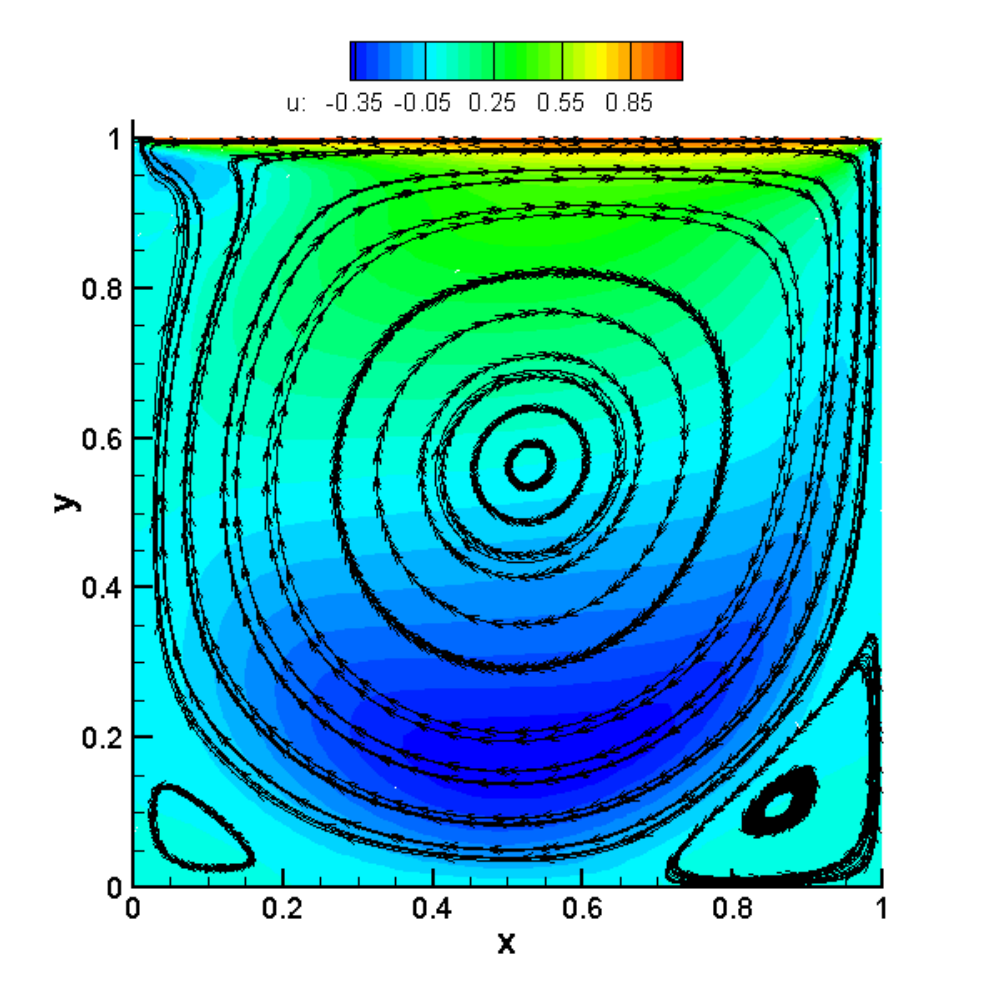}
	\includegraphics[width=0.32\textwidth]{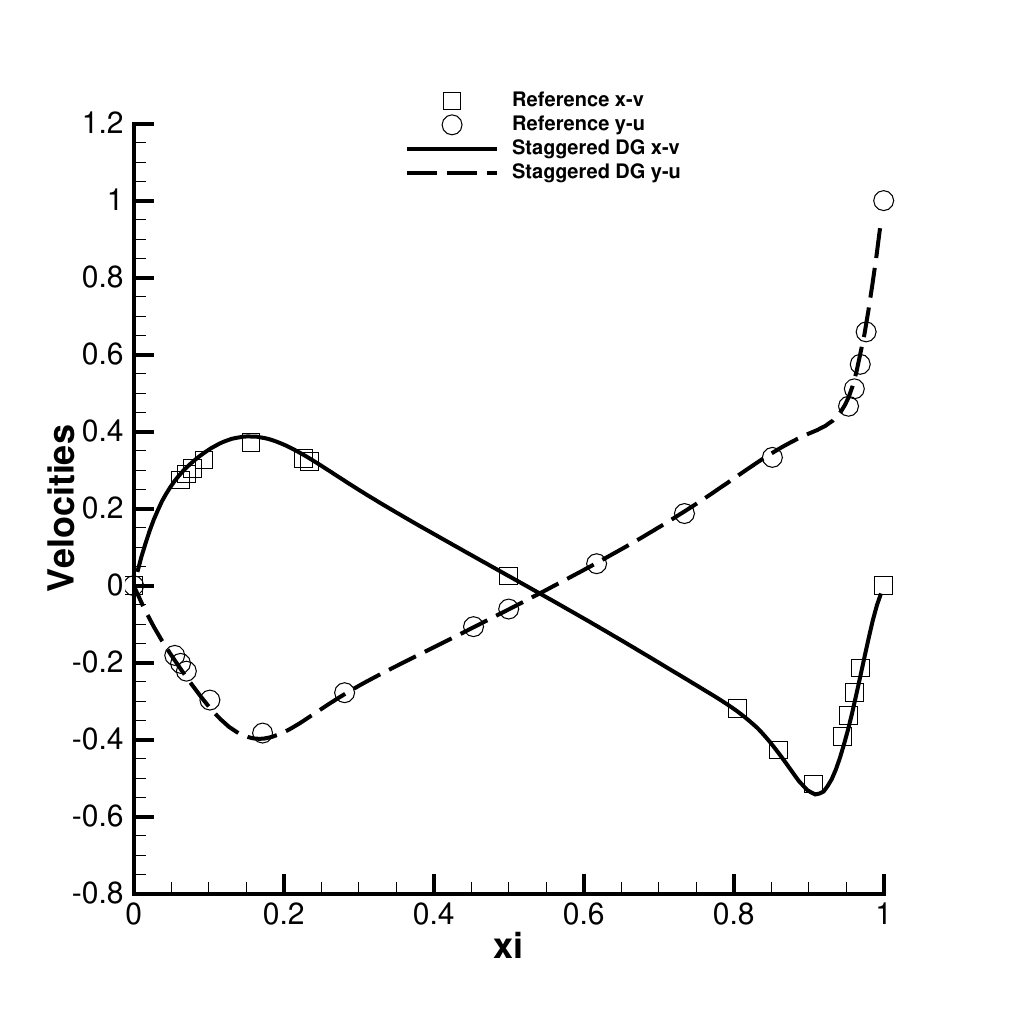}
    \caption{From left to right: Mach number; horizontal velocity with streamlines; and velocity cut compared with the data given by Ghia et al in \cite{Ghia1982} at $Re=1000$.}
    \label{fig.NT_6_2}
	\end{center}
\end{figure}

Another classical two dimensional test problem is the Lid-driven cavity flow. For the incompressible case some well-known standard reference solutions are available at several Reynolds numbers, 
see \cite{Ghia1982}.
The computational domain is $\Omega=[-0.5,0.5]\times [-0.5, 0.5]$ and the flow inside the cavity is driven by a velocity $\mathbf{v}=(1,0)$ imposed at the upper boundary, i.e. at $y=0.5$. 
No-slip wall boundary conditions are then imposed on the rest of the boundaries. As initial condition we take $\mathbf{v}(\xx,0)=0$, $\rho(\xx,0)=1$ and $p(\xx,0)=10^5$ in order to approach 
the incompressible limit. 
We use a grid with only $\Ni=440$ triangles, setting $(p,p_\gamma)=(3,0)$ and we run the simulation until the steady state is reached. The resulting Mach number and velocity contours as
well as the velocity profiles obtained by a cut along the lines $x=0$ and $y=0$ are compared with the reference solution given in \cite{Ghia1982} and are reported in Figure \ref{fig.NT_6_1} for $Re=100$ and in Figure \ref{fig.NT_6_2} for $Re=1000$.
Our numerical solution fits the incompressible reference solution of Ghia \cite{Ghia1982} very well for both Reynolds numbers under consideration. The resulting recirculation zones 
are visualized by the streamline plots in Figures \ref{fig.NT_6_1} and \ref{fig.NT_6_2} and show the formation of the main recirculation but also some secondary corner vortices, 
as already observed in several papers, see \cite{Erturk2005,2STINS,Fambri2016}. 

\subsubsection{Low Reynolds number flow past a circular cylinder}
\begin{figure}[ht!]
    \begin{center}
		\includegraphics[width=0.7\textwidth]{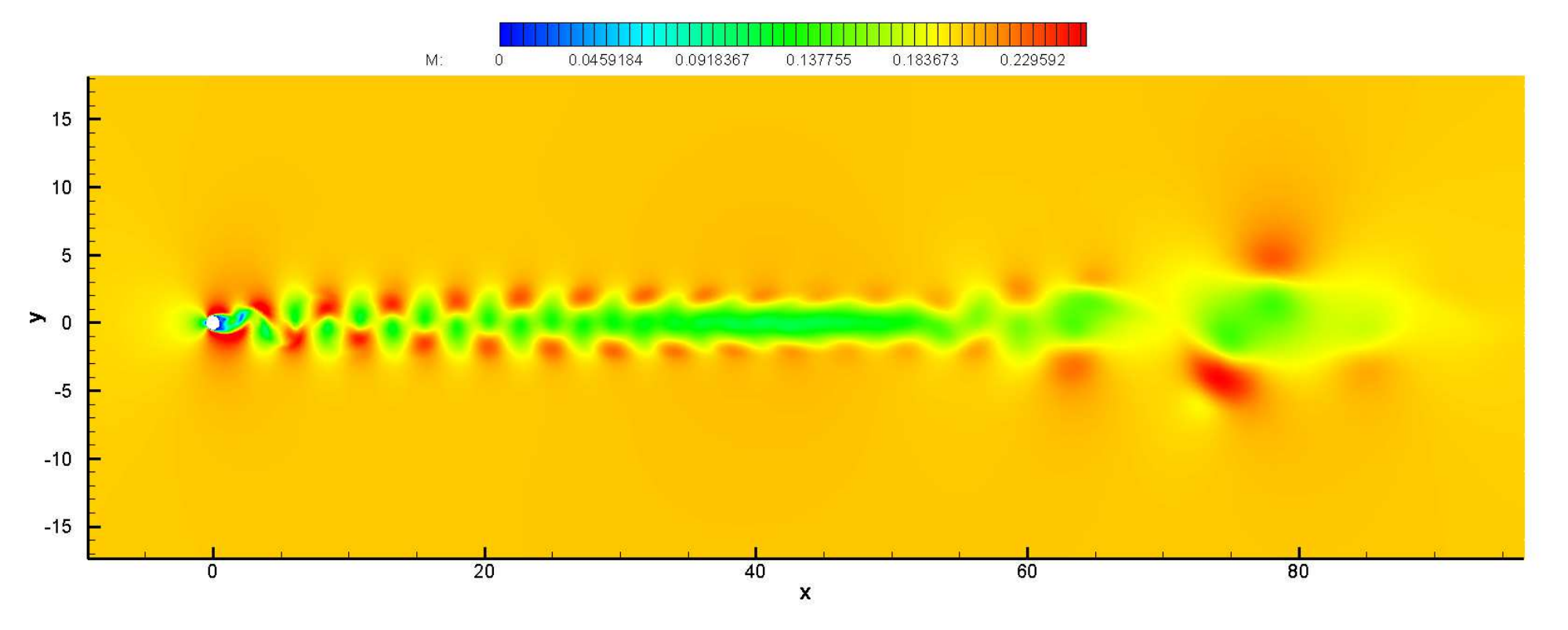}
		\includegraphics[width=0.7\textwidth]{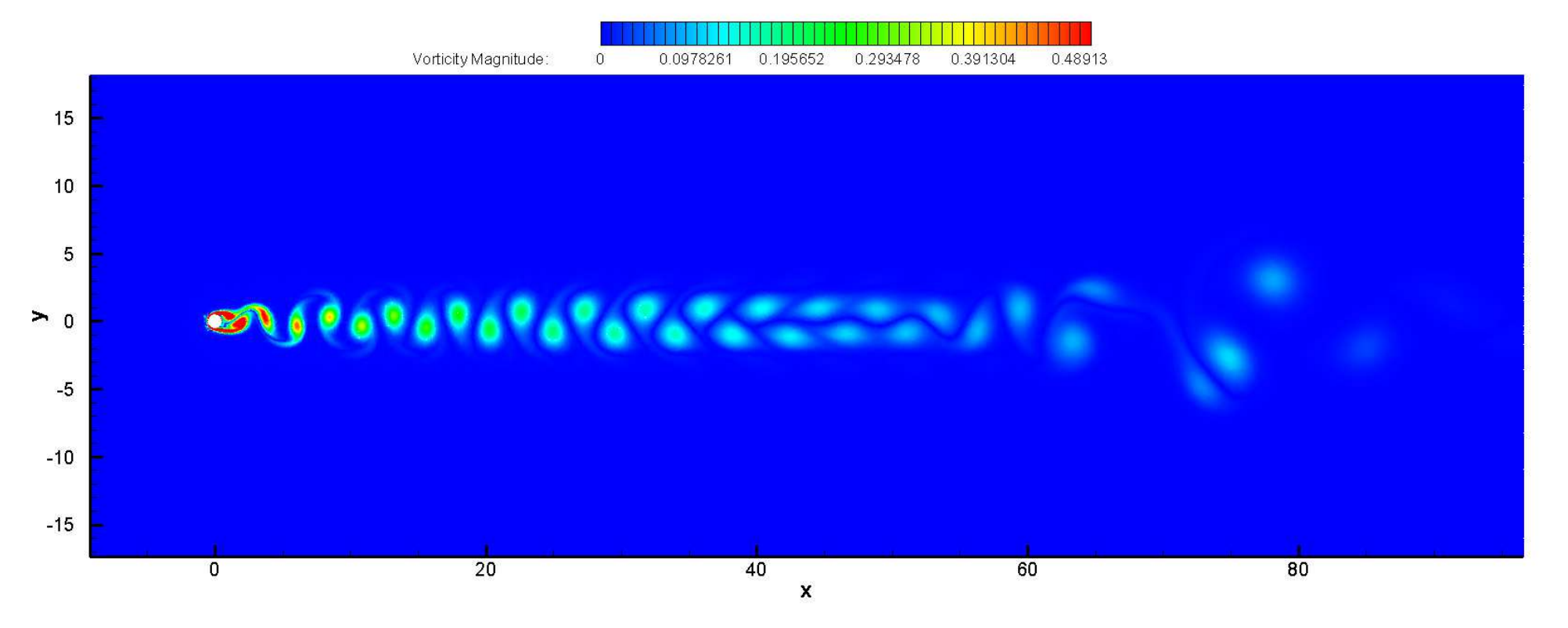}
    \caption{Mach number profile of the von Karman street at $t=500$ (top) and vorticity magnitude (bottom).}
    \label{fig.NT_7_1}
	\end{center}
\end{figure}

\begin{figure}[ht!]
    \begin{center}
		\includegraphics[width=0.6\textwidth]{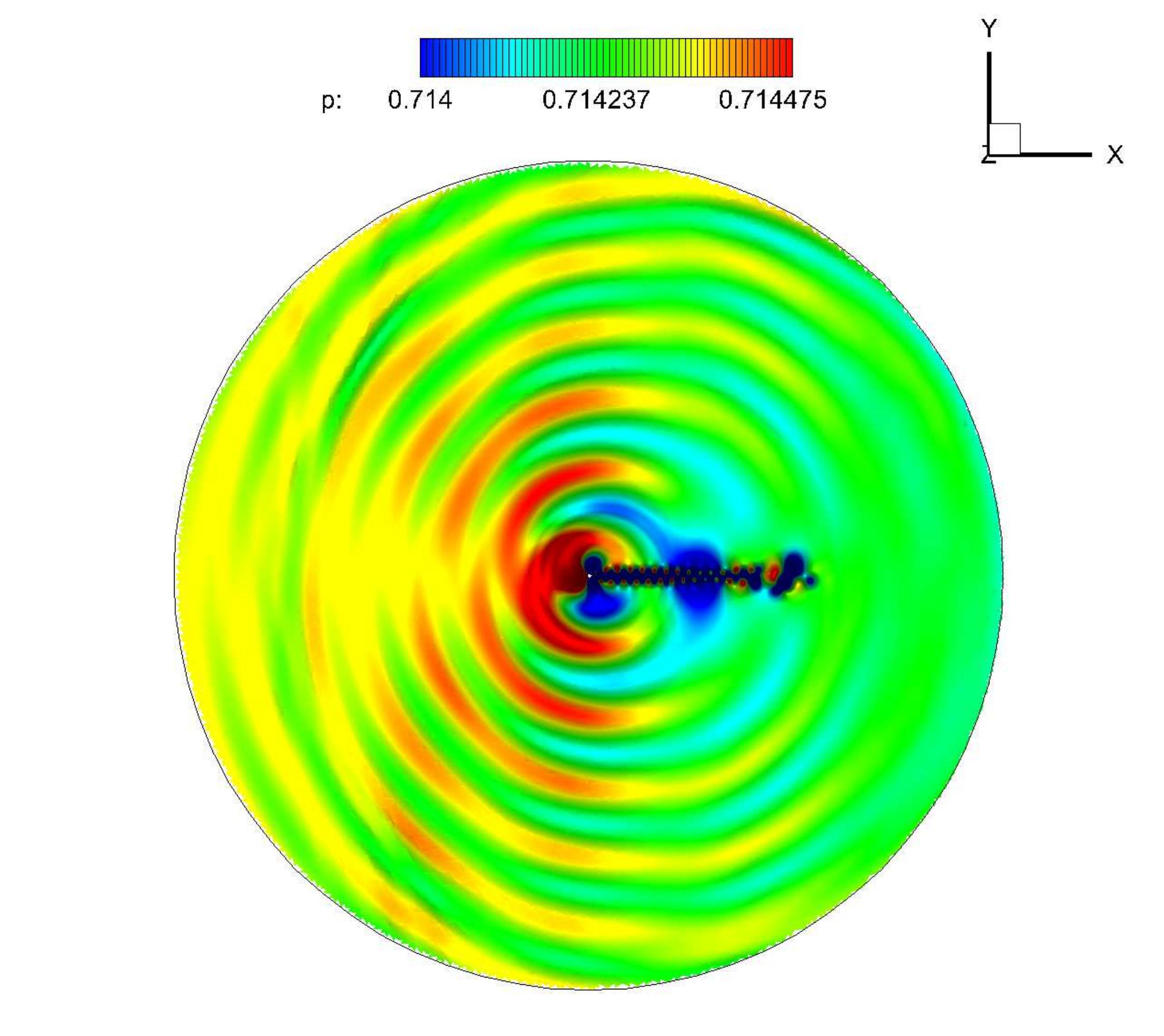}
    \caption{Sound pressure field generated by the von Karman street.}
    \label{fig.NT_7_2}
	\end{center}
\end{figure}

In this section we consider the flow past a circular cylinder. In this case, the use of an isoparametric finite element approach is mandatory to represent the curved geometry of the cylinder wall, as already discussed in \cite{BassiRebay,2DSIUSW}. 
We consider a circular cylinder of diameter $d=1$, the Mach number of the incoming flow is $M=0.2$, the Reynolds number is $Re=150$, the Prandtl number is set to $Pr=1$ and $\gamma=1.4$. As initial  condition we simply take the constant far field $(\rho_\infty,p_\infty,u_\infty,v_\infty)=(1,1/\gamma,0.2,0)$ everywhere. The simulation is performed on the domain $\Omega=C_{200} \backslash C_{0.5}$,  where $C_r=\left\{ \xx \, | \, \left\| \xx \right\| \leq r \right\}$ indicates the area delimited by the circle of radius $r$. We cover $\Omega$ with $\Ni=12494$ triangles and use a sponge layer of thickness $20$ in order to allow acoustic waves to leave the domain without generating spurious reflections at the outer boundary. We furthermore use $(p,p_\gamma)=(3,0)$ and $t_{end}=500$. The resulting Mach number profile and vorticity magnitude at the final time is reported in Figure \ref{fig.NT_7_1} where we can clearly see the von Karman vortex street appearing behind the circular cylinder. Figure \ref{fig.NT_7_2} shows the sound pressure field generated by the von Karman street using the same contour colors as the ones already employed in \cite{ADERNSE}.  The resulting frequency for the sound waves obtained in our numerical experiment is $f=0.036688$ that corresponds to a Strouhal number of $St=\frac{fd}{u_{\infty}}=0.1834$, which is exactly the same obtained by M\"uller in \cite{Muller2008} and in good agreement with the one obtained in \cite{ADERNSE}, where $St=0.182$ was reported.

\subsubsection{Viscous single Mach reflection problem}
\begin{figure}[ht!]
    \begin{center}
		\includegraphics[width=0.8\textwidth]{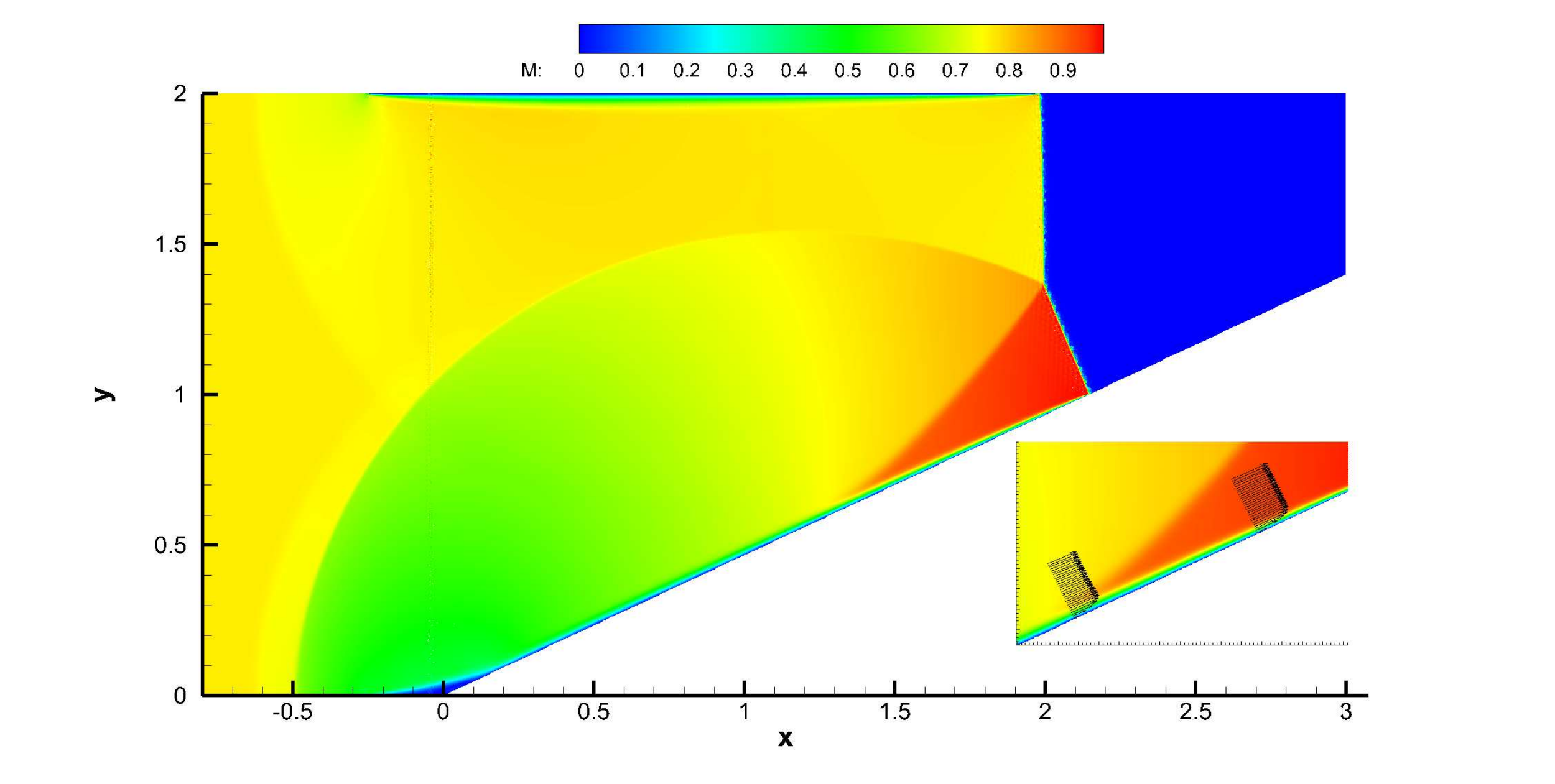}
    \caption{Mach number profile for the viscous single Mach reflection problem.}
    \label{fig.NT_8_1}
	\end{center}
\end{figure}

We now consider a well known test case in the supersonic regime. We consider a viscous shock wave traveling at a shock Mach number $M_s=1.7$ that hits a wedge with an angle of inclination of 
$\alpha=25^{\circ}$. The considered Reynolds number based on the shock speed is $Re_0=\frac{\rho_0 c_0 M_s L}{\mu_0}=3400$, where $\rho_0=1$, $p_0=1/\gamma$  and $\gamma=1.4$. By taking $Pr=0.75$ we have again the smooth analytical representation of the viscous shock discovered by Becker in \cite{Becker1923} and discussed in the viscous shock test case, see Section \ref{sec_NT_VS}. For this test we consider $\Omega=[-0.8,3]\times[0,2]-W_{25}$ where $W_{\alpha}$ is wedge of inclination $\alpha$. We use again $(p,p_\gamma)=(3,0)$ and $\Ni=94560$ triangles. As initial condition we use the exact viscous shock profile initially placed at $x=-0.04$ and defined in Section \ref{sec_NT_VS}. Observe that the exact solution in this case is smooth for sufficiently fine grid resolution and hence no limiter is necessary. We run the simulation up to $t_{end}=1.2$ so that the exact shock location at the final time must be $x=2$. As boundary condition we impose no-slip boundary condition for the upper and lower boundaries corresponding to $x\geq -0.25$; slip wall boundary for $x<0.25$; and Dirichlet boundary conditions with the initial state at the inflow and outflow. The resulting Mach number profile at the final time is reported in Figure \ref{fig.NT_8_1} with a zoom of the boundary layer generated by the presence of the viscous no-slip wall boundary conditions. A good agreement with experimental results presented in \cite{toro-book} and other numerical experiments presented in \cite{ADERNSE} can be observed also in this case. 

\subsubsection{Viscous double Mach reflection problem}
\begin{figure}[ht!]
    \begin{center}
		\includegraphics[width=0.8\textwidth]{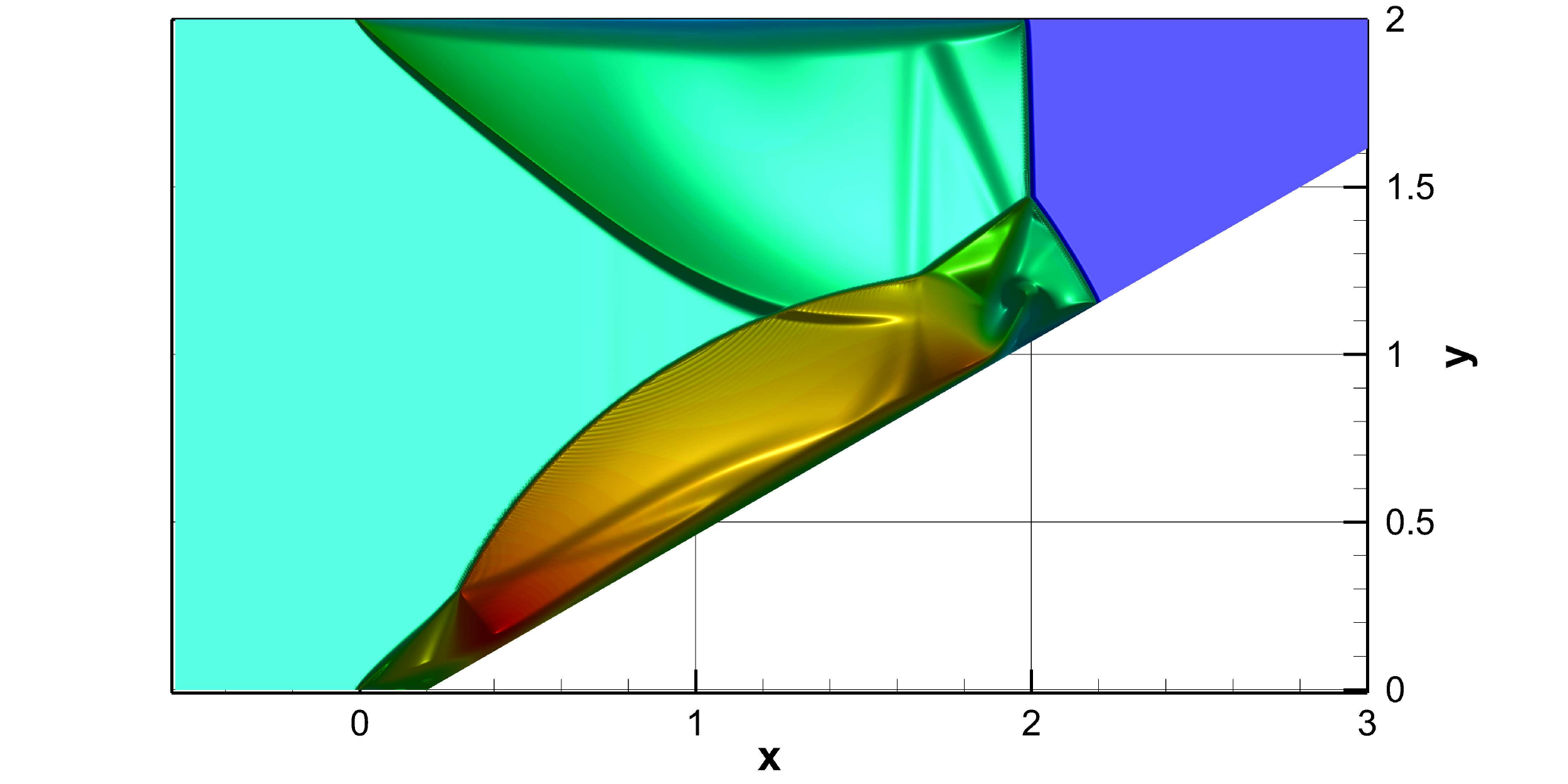}
    \caption{Density profile for the viscous double Mach reflection problem.}
    \label{fig.NT_9_1}
	\end{center}
\end{figure}

\begin{figure}[ht!]
    \begin{center}
		\includegraphics[width=0.49\textwidth]{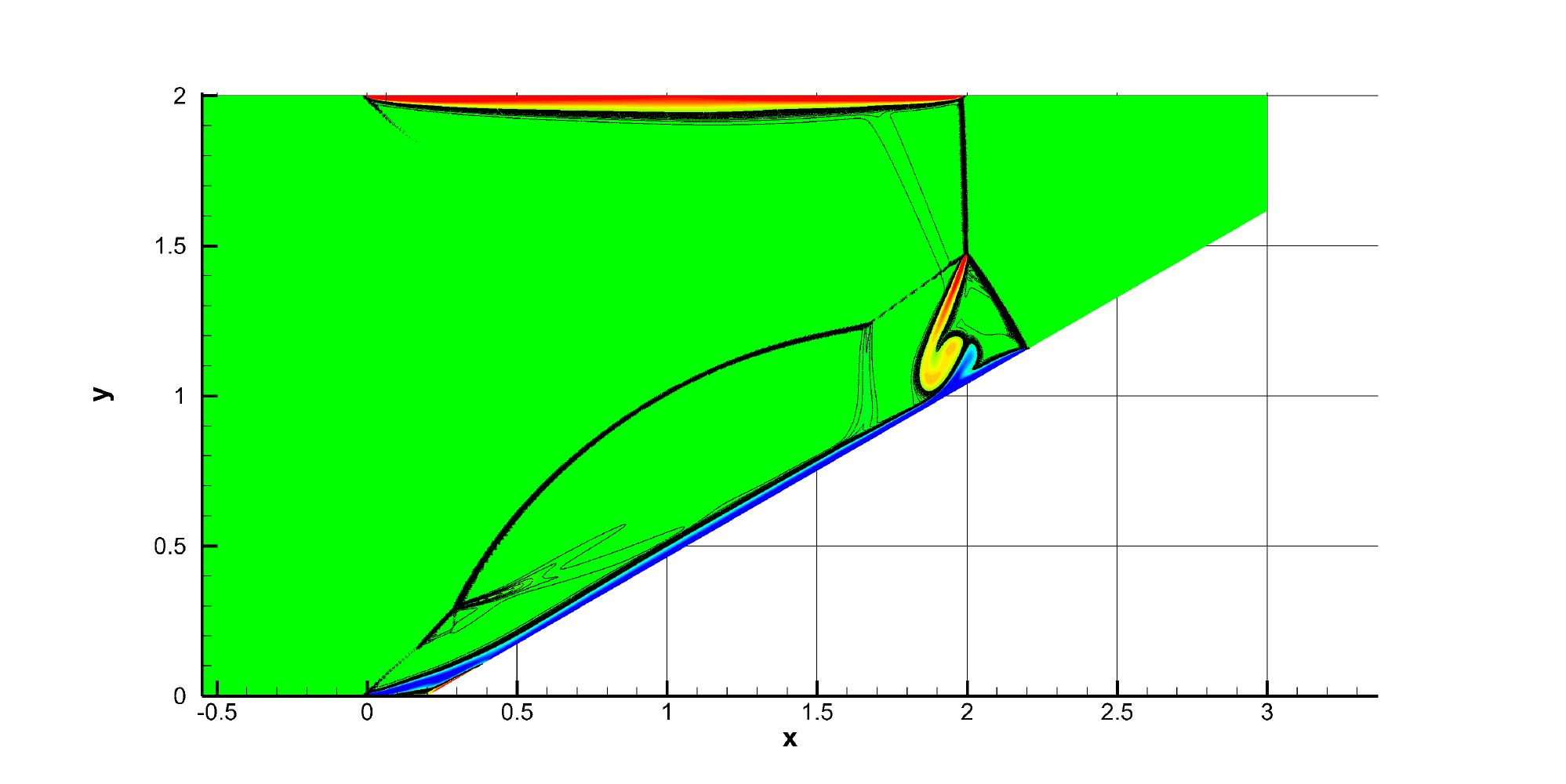}
		\includegraphics[width=0.49\textwidth]{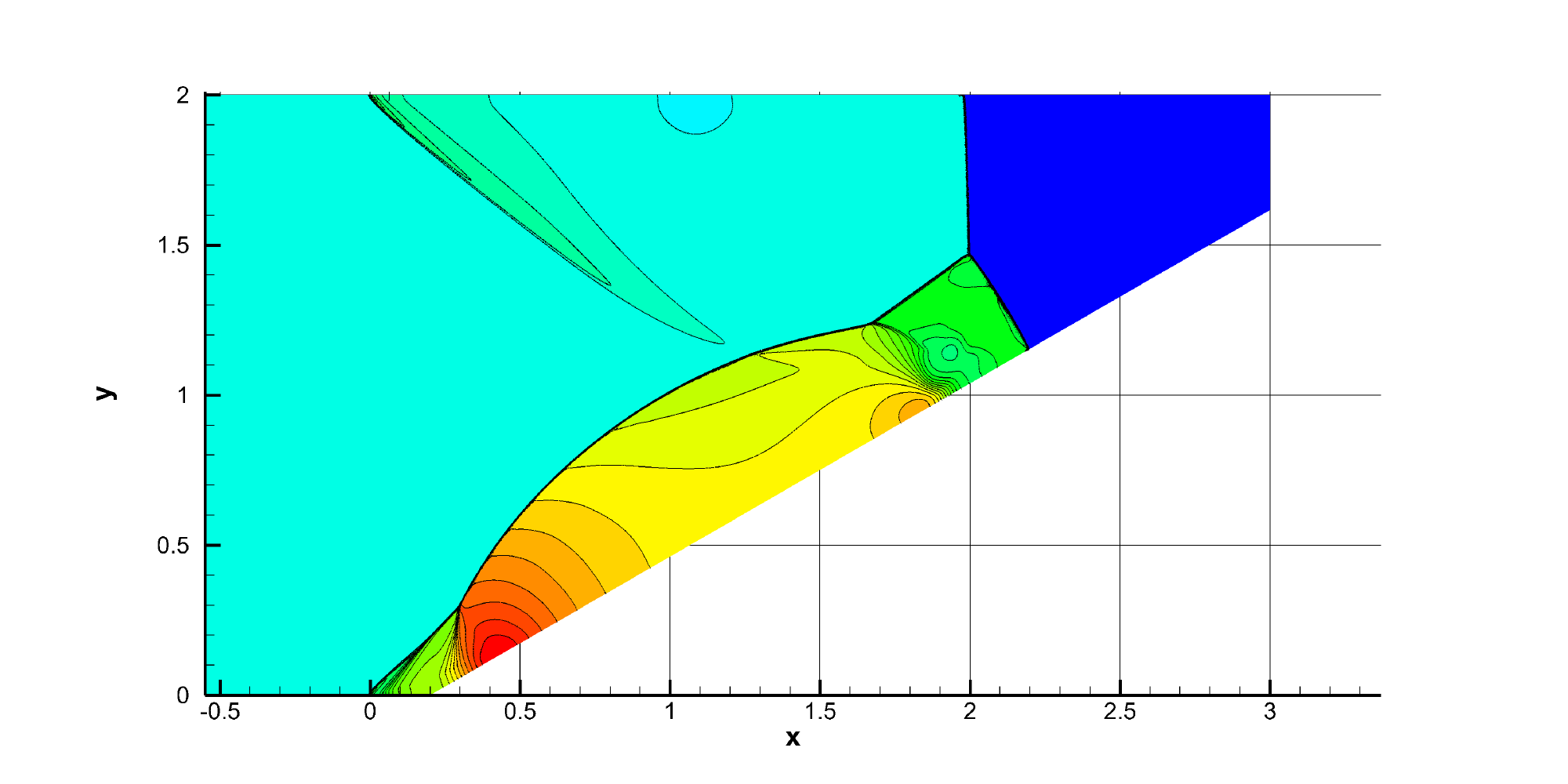} \\
		\includegraphics[width=0.49\textwidth]{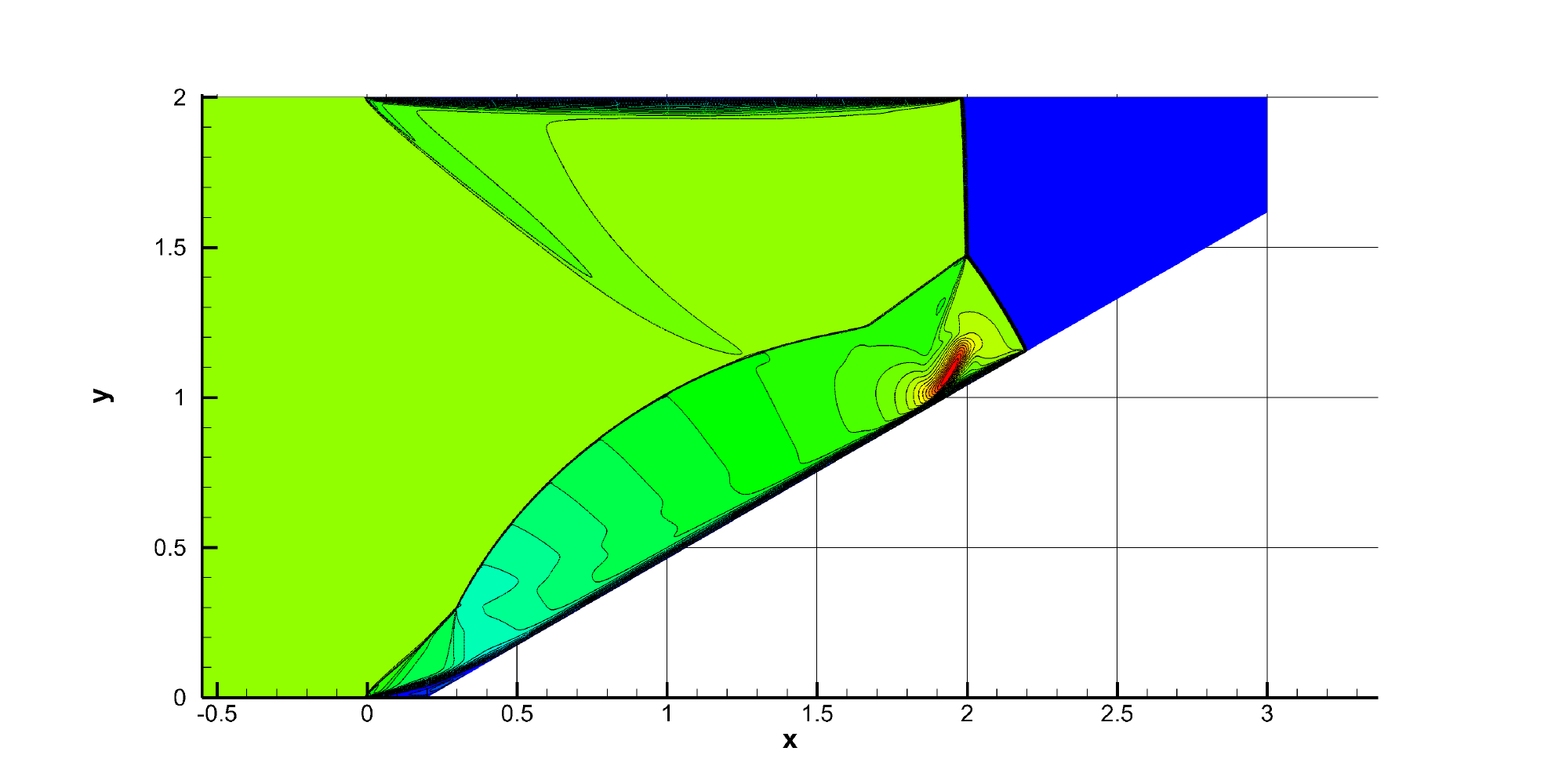}
		\includegraphics[width=0.49\textwidth]{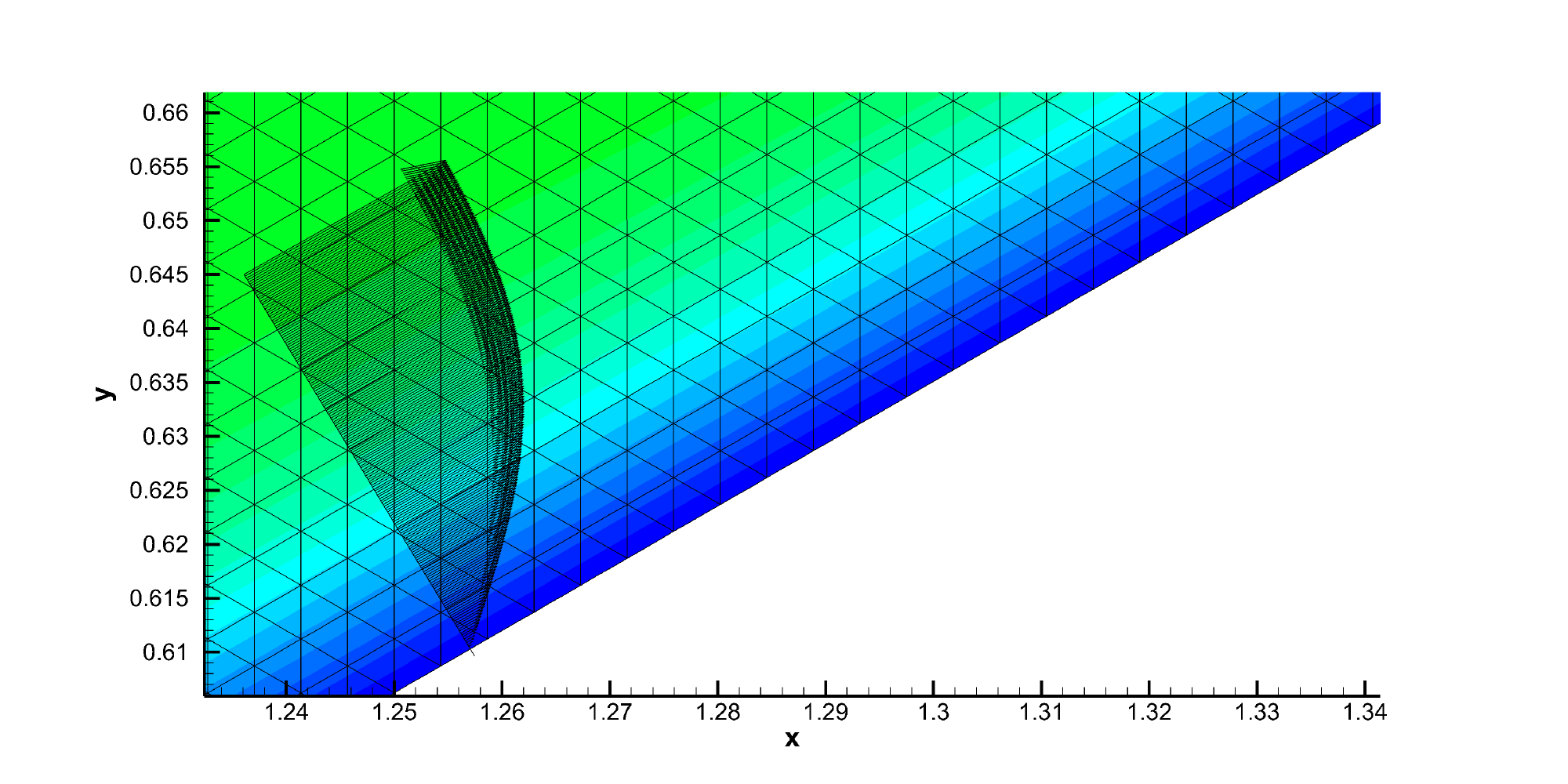}
    \caption{From top left to bottom right: vorticity, pressure, Mach number and boundary layer zoom for the viscous double Mach reflection problem.}
    \label{fig.NT_9_2}
	\end{center}
\end{figure}

Here we investigate a similar test case as before, but with a much stronger incident shock wave. The problem under consideration is a viscous version of the double Mach reflection problem that was studied 
in \cite{ADERNSE} and which has been originally proposed in the inviscid case by Woodward and Colella in \cite{woodwardcol84}. The shock Mach number is now $M_s=10.0$ and the wedge surface has an inclination of $\alpha=30^{\circ}$. The gas parameters are taken as in the single Mach reflection problem, but now $Re_0=800$. The domain $\Omega=[-0.55,3]\times[0,2]-W_{30}$ is covered with $\Ni=438344$ triangles; $(p,p_\gamma)=(3,0)$ and $t_{end}=0.2$. The initial shock position is $x=0$ so that the exact shock location at the final time is $x=2$. Again we use no-slip wall boundary conditions for $x\geq0$ and slip boundary for $x<0$. The resulting density contour is reported in Figure \ref{fig.NT_9_1} and looks very similar to the one computed in \cite{ADERNSE}. Figure \ref{fig.NT_9_2} shows the resulting vorticity, Mach number and pressure contours. On the upper boundary we can clearly see the boundary layer introduced from $x=0$ on, similar to the one shown in the previous example. A zoom on the wedge boundary is then reported with the main grid at the bottom right of Figure \ref{fig.NT_9_2}. Also in this case, since we are considering a smooth viscous shock, no limiter is necessary.

\subsubsection{Riemann problems}
\begin{figure}[ht!]
		\begin{tabular}{ccc} 
		\includegraphics[width=0.32\textwidth]{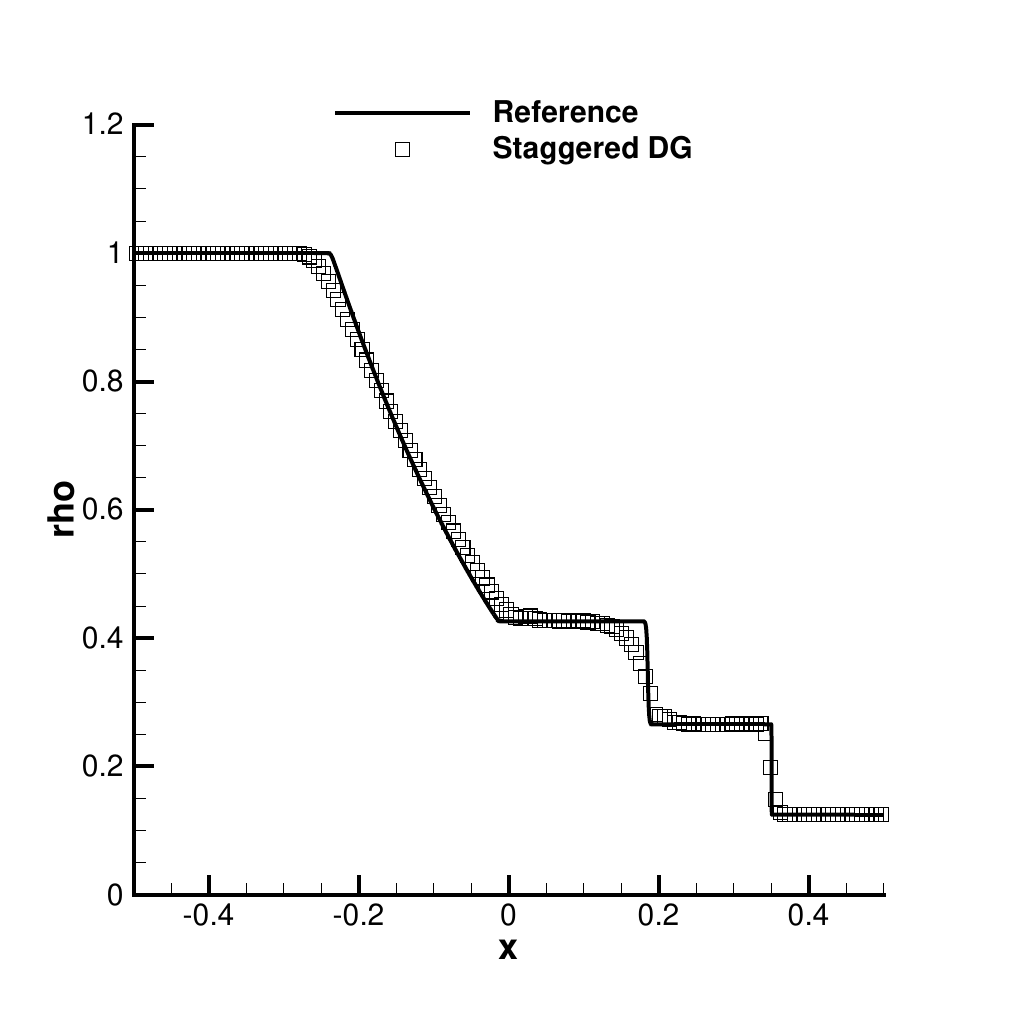}   & 
		\includegraphics[width=0.32\textwidth]{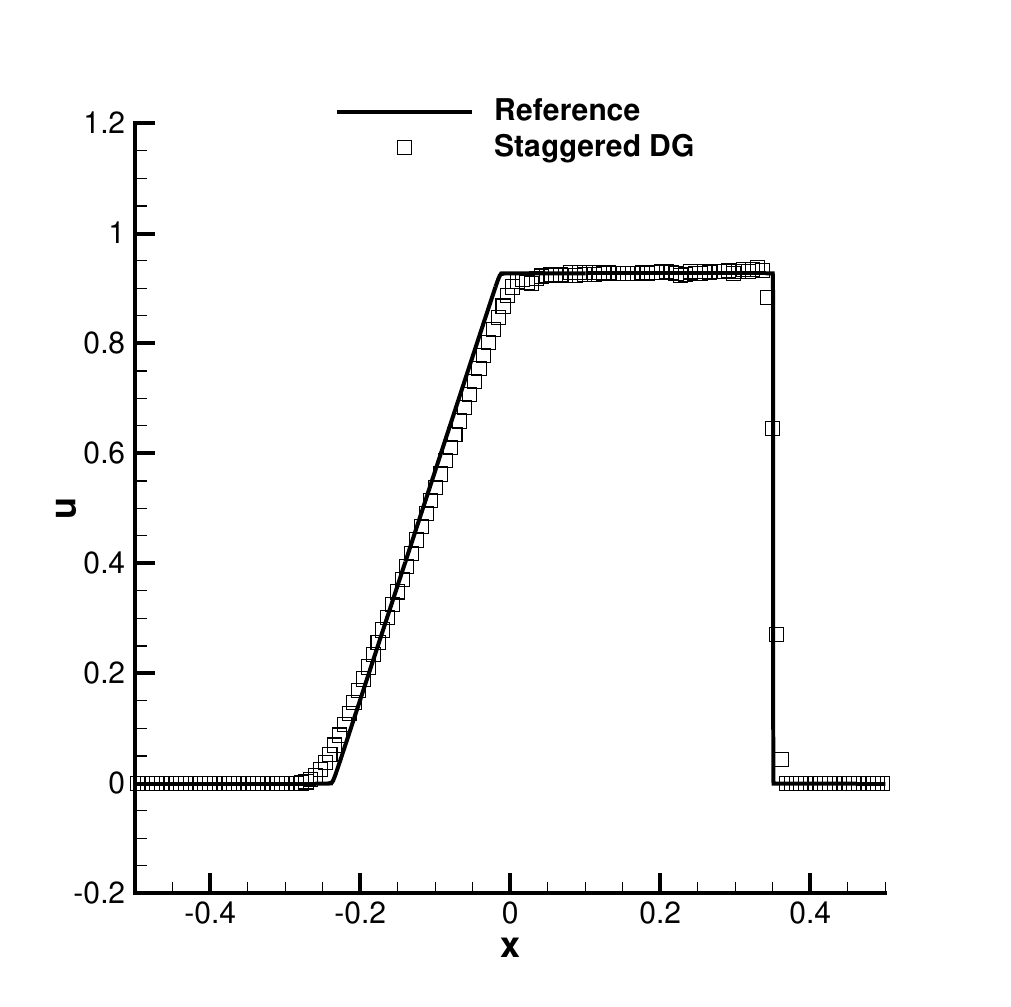}     & 
		\includegraphics[width=0.32\textwidth]{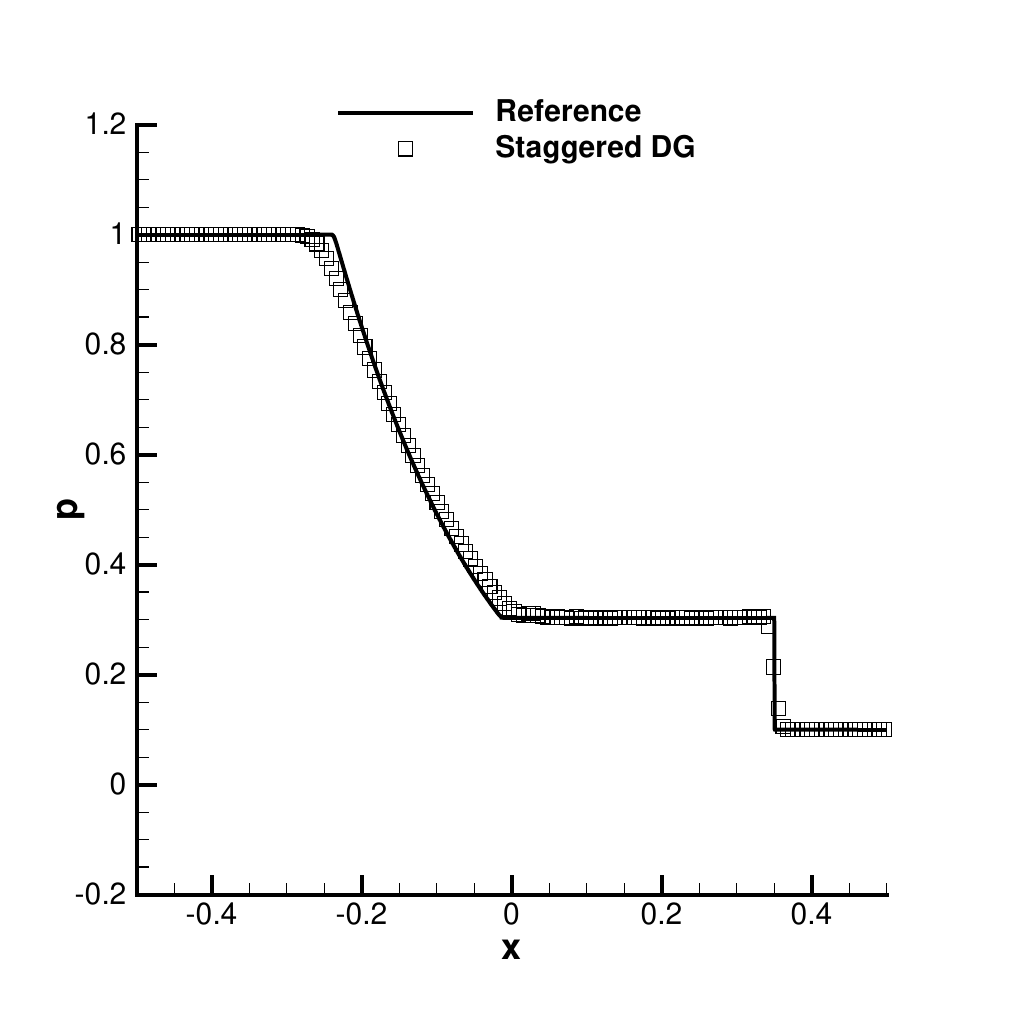}     \\  
		\multicolumn{3}{c}{
		\includegraphics[width=0.55\textwidth]{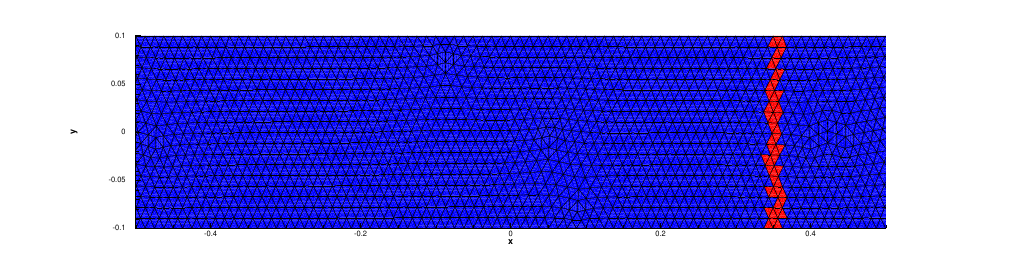}
		}
		\end{tabular} 
    \caption{Riemann problem RP1 at time $t=0.2$ using $p=4$. Top row from left to right: 1D cut through the density, velocity and pressure profile compared against the exact solution. 
		         Bottom row: computational grid, limited cells highlighted in red and unlimited cells in blue.}
    \label{fig:NT_10_1}
\end{figure}

\begin{figure}[ht!]
		\begin{tabular}{ccc} 
		\includegraphics[width=0.32\textwidth]{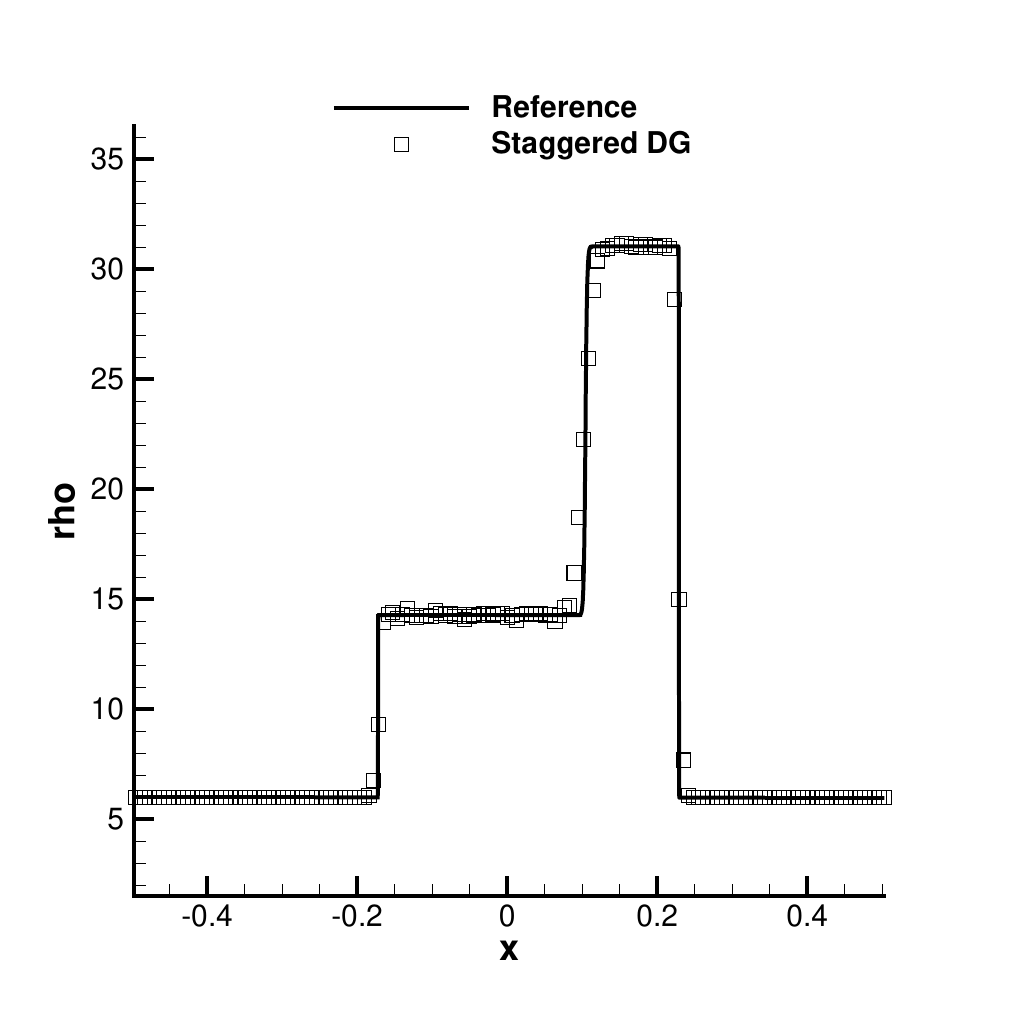}   & 
		\includegraphics[width=0.32\textwidth]{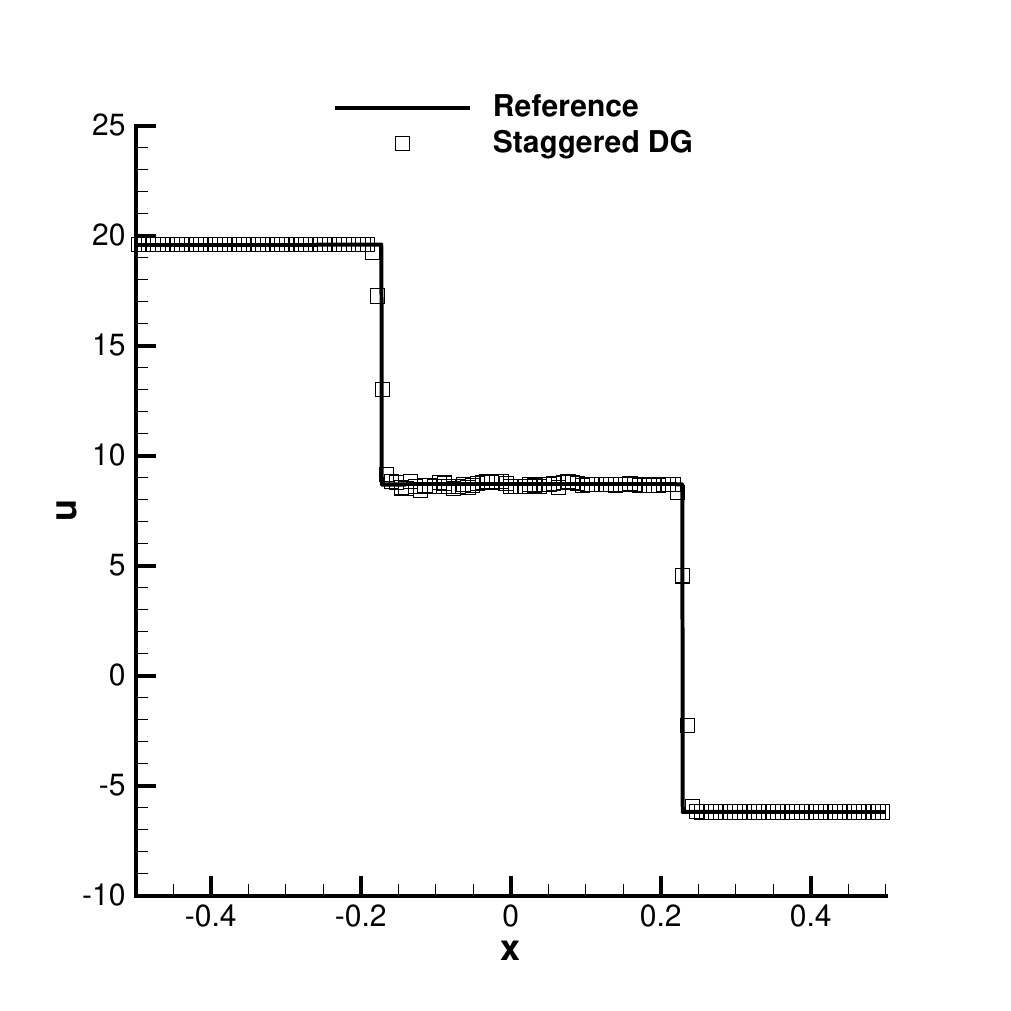}     & 
		\includegraphics[width=0.32\textwidth]{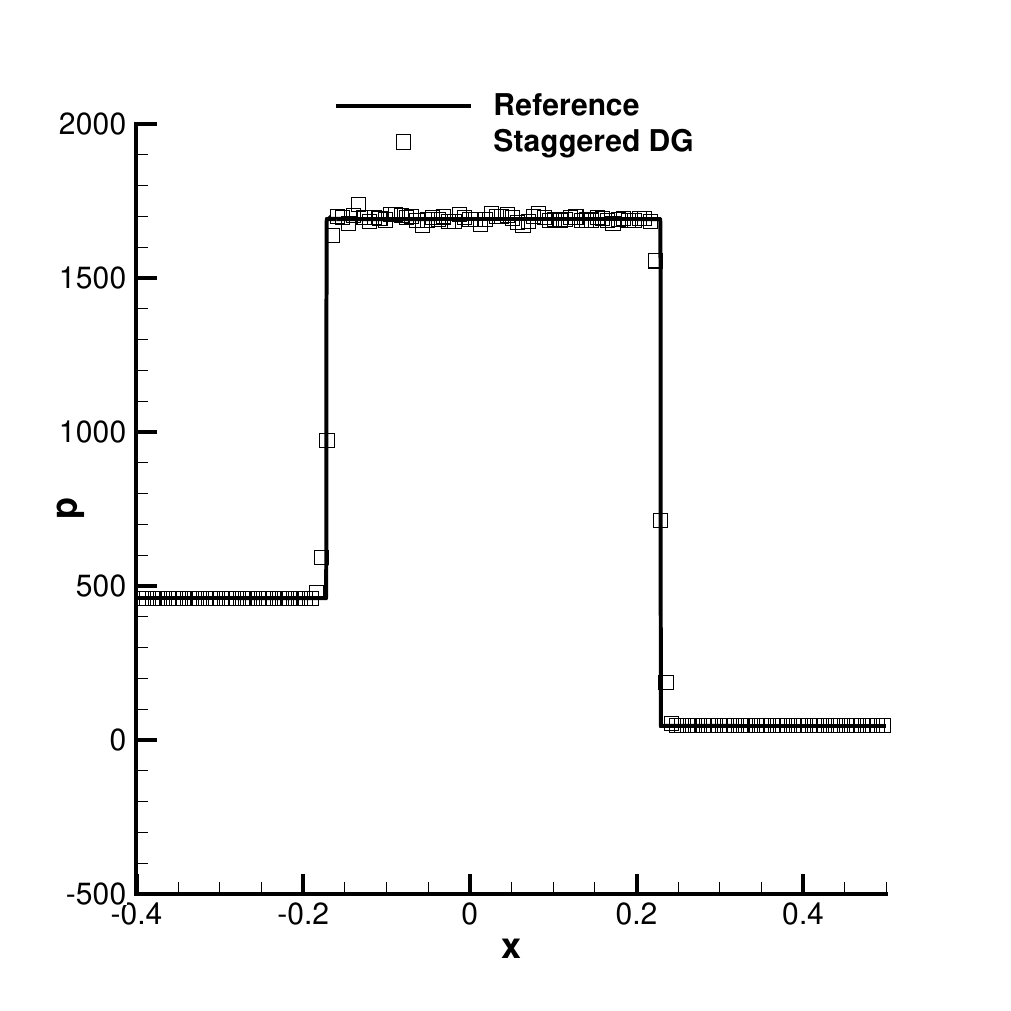}     \\  
		\multicolumn{3}{c}{
		\includegraphics[width=0.55\textwidth]{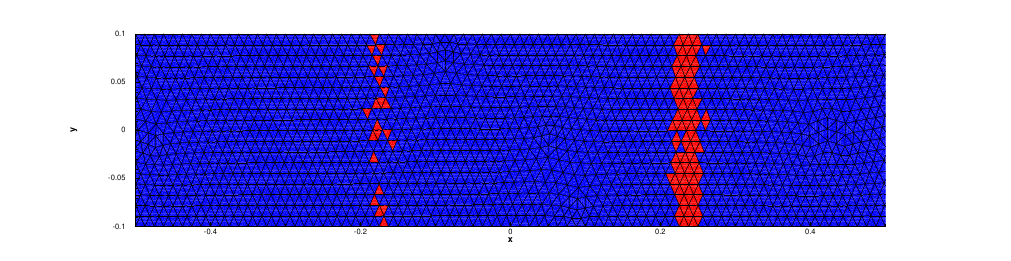}
		}
		\end{tabular} 
    \caption{Riemann problem RP2 at time $t=0.035$ using $p=4$. Top row from left to right: 1D cut through the density, velocity and pressure profile compared against the exact solution. 
		         Bottom row: computational grid, limited cells highlighted in red and unlimited cells in blue.}
    \label{fig:NT_10_2}
\end{figure}

Up to now we have treated only smooth problems where the introduction of the limiter discussed in Section \ref{sec_DGLIMITER} was not necessary. Therefore, we now want to study some classical 
Riemann problems where the use of the limiter is mandatory in order to achieve non-oscillatory and physically admissible solutions. 
\begin{table}%
\caption{Initial states and solution times for the considered Riemann problems.}
 \begin{center}
\begin{tabular}{cccccccc}
\hline RP &	$\rho_L$  & $u_L$	& $p_L$ &	$\rho_R$  & $u_R$	& $p_R$ & $t_{end}$\\ \hline
1 	& 1.0 	& 0.0	&1.0	&0.125	&0.0	&0.1	&0.2\\  
2 	& 5.99924 	& 19.5975	&460.894	&5.99242	&-6.19633	&46.095	&0.035\\ 
3 	& 0.445 	& 0.698	&3.528	&0.5	&0.0	&0.571	&0.15\\ 
4 	& 1.0 	& 2.0	&0.1	&1.0	&-2.0	&0.1	&0.8 \\
5 	& 1.0 	& 0.75 &1.0	&0.125	&0.0	&0.1	&0.2 \\ \hline
\end{tabular}
\label{tab.NT_10_1}
 \end{center}
\end{table}
The initial conditions for the Riemann problems under consideration are taken from \cite{toro-book} and are reported in Table \ref{tab.NT_10_1}. For all the test cases we 
take $\Omega=[-0.5, 0.5]\times[-0.1,0.1]$, covered with $\Ni=3480$ triangles and setting $(p,p_\gamma)=(4,0)$. As previously discussed in Section \ref{sec_DGLIMITER}, 
in this case it is mandatory to use an implicit treatment of the viscous terms and of the heat flux, since in the AV approach the artificial viscosity 
and heat conduction coefficients used in the vicinity of shock waves and discontinuities can become very high. 
In order to avoid unphysical oscillations at the initial time, where the originally discontinuous initial data are projected onto the piecewise polynomial approximation space via 
classical $L_2$ projection, we \textit{smooth} the initial condition as follows:
\begin{eqnarray}
	Q=\frac{1}{2}(Q_R+Q_L)+\frac{1}{2}(Q_R-Q_L) \, \textnormal{erf}\left(\frac{x-x_0}{\epsilon_0}\right),
\label{eq:NT_10_th}
\end{eqnarray}
where $Q_L$ and $Q_R$ are the left and the right state of the Riemann problem, respectively; $\textnormal{erf}$ is the classical error function; $x_0$ is the initial position of the discontinuity 
and $\epsilon$ is a small relaxation parameter. 
The first considered Riemann problem is the classical {\itshape Sod problem} initially proposed by Sod in \cite{Sod1978}. The used numerical parameters are the ones corresponding to (RP1) in 
Table $\ref{tab.NT_10_1}$; $x_0=0$ and $\epsilon_0 = 2\cdot 10^{-2}$. In Figure \ref{fig:NT_10_1}  we report our numerical results for the main variables $(\rho, u, p)$ at the final time $t_{end}$ 
and compare them with the exact solution provided in \cite{toro-book}. We also show a plot of the limited cells at $t=t_{end}=0.2$. We can observe that the position of the shock is well reproduced 
by our numerical scheme and no significant oscillations appear, thanks to the use of the MOOD+AV limiter. Furthermore, in this simple case, the limiter acts only at the shock wave, as expected. 
Next, we consider a very difficult Riemann problem taken from \cite{toro-book} and reported as (RP2) in Table \ref{tab.NT_10_1}. Here, we set $x_0=0$ and $\epsilon_0=1\cdot 10^{-2}$. 
The resulting numerical solution for all main flow variables as well as the limited cells are reported in Figure \ref{fig:NT_10_2}. 
Also in this case the limiter acts only on a very small number of elements in the computational domain. It is important to emphasize that this test case does not run without the use 
of a limiter, where negative values of the pressure appear after only a few time steps. Hence, in this case the use of the limiter is mandatory to obtain a solution. 

From Figure \ref{fig:NT_10_2} one can note that the exact solution is well reproduced by our numerical scheme. However, in this case one can also observe the production 
of some small spurious oscillations that require some further investigations. At this point we would like to 
stress that in all computational results shown in this paper, we did never play with the AV coefficients, as is usually done in the literature. We only apply 
the MOOD detector, together with an artificial viscosity that leads to a mesh Reynolds number of unity. 

\begin{figure}[ht!]
		\begin{tabular}{ccc} 
		\includegraphics[width=0.32\textwidth]{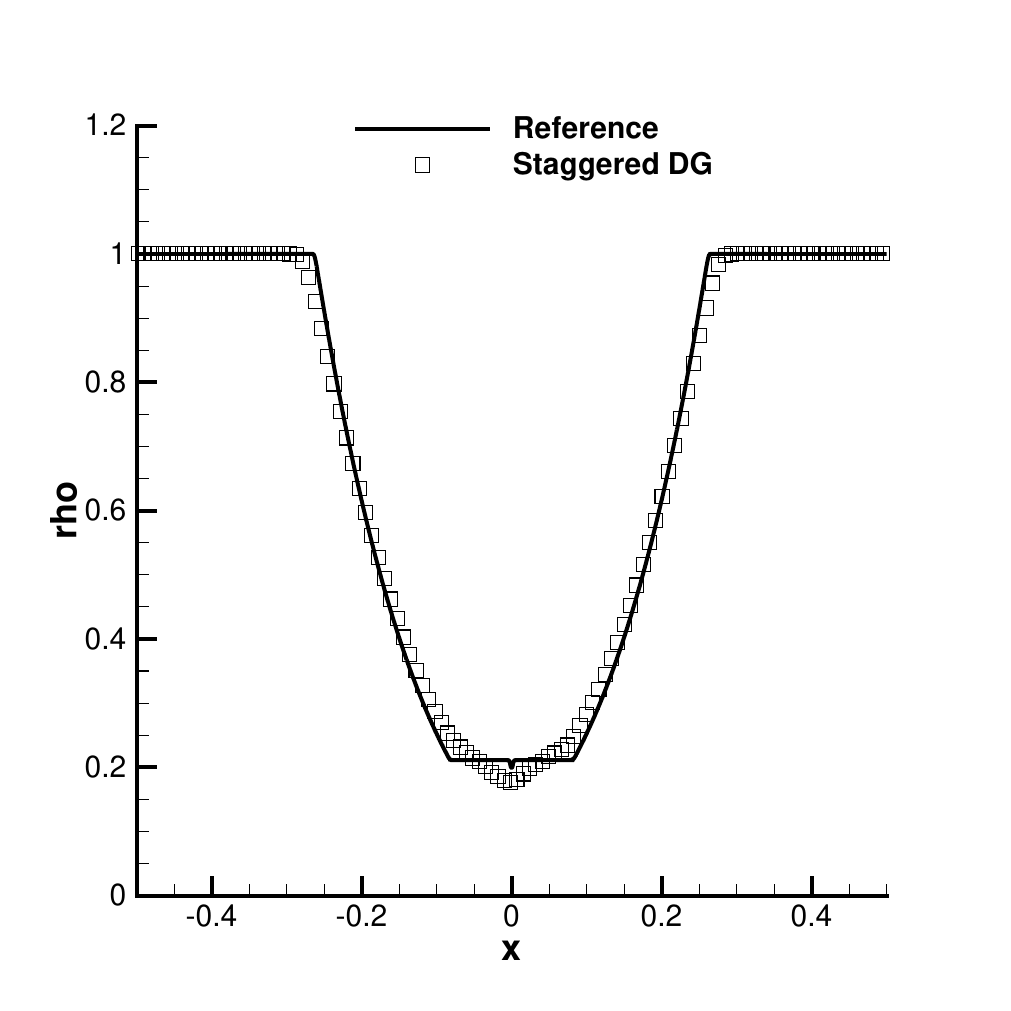}   & 
		\includegraphics[width=0.32\textwidth]{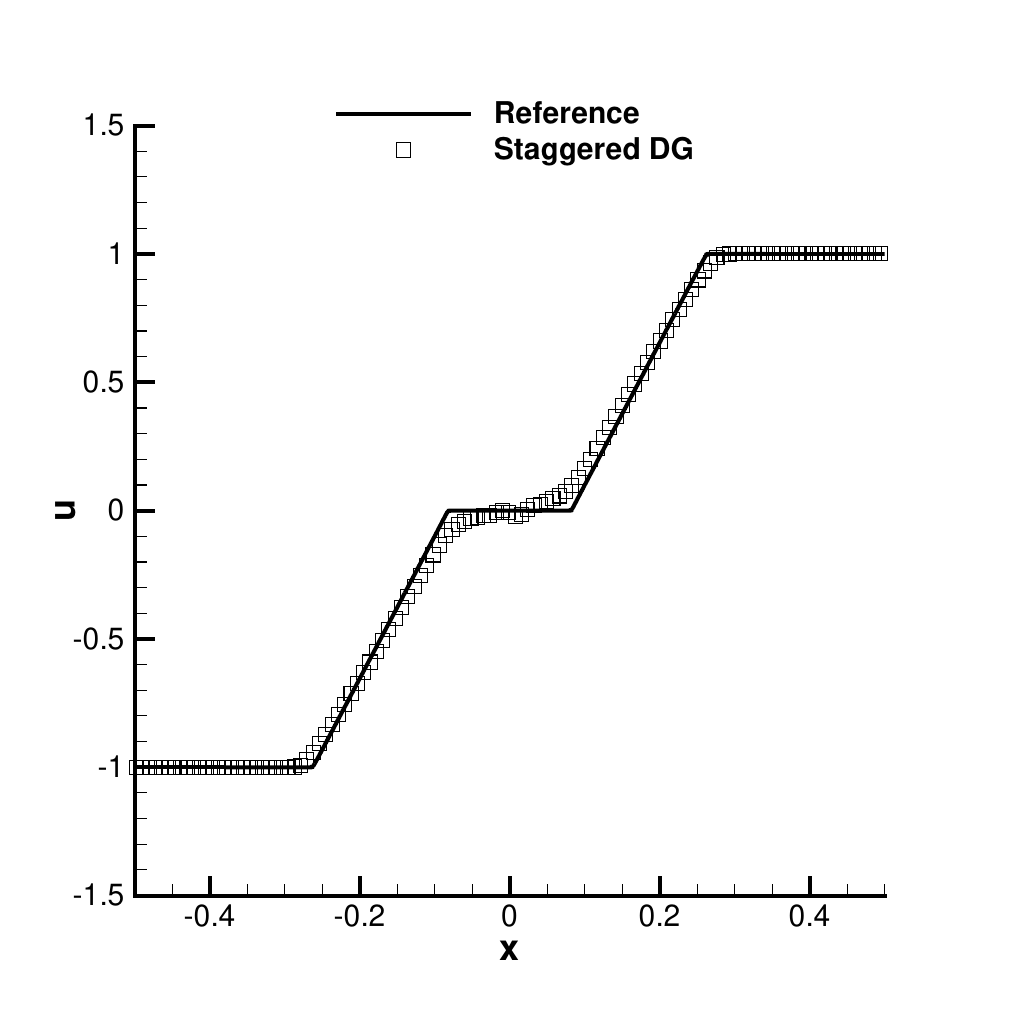}     & 
		\includegraphics[width=0.32\textwidth]{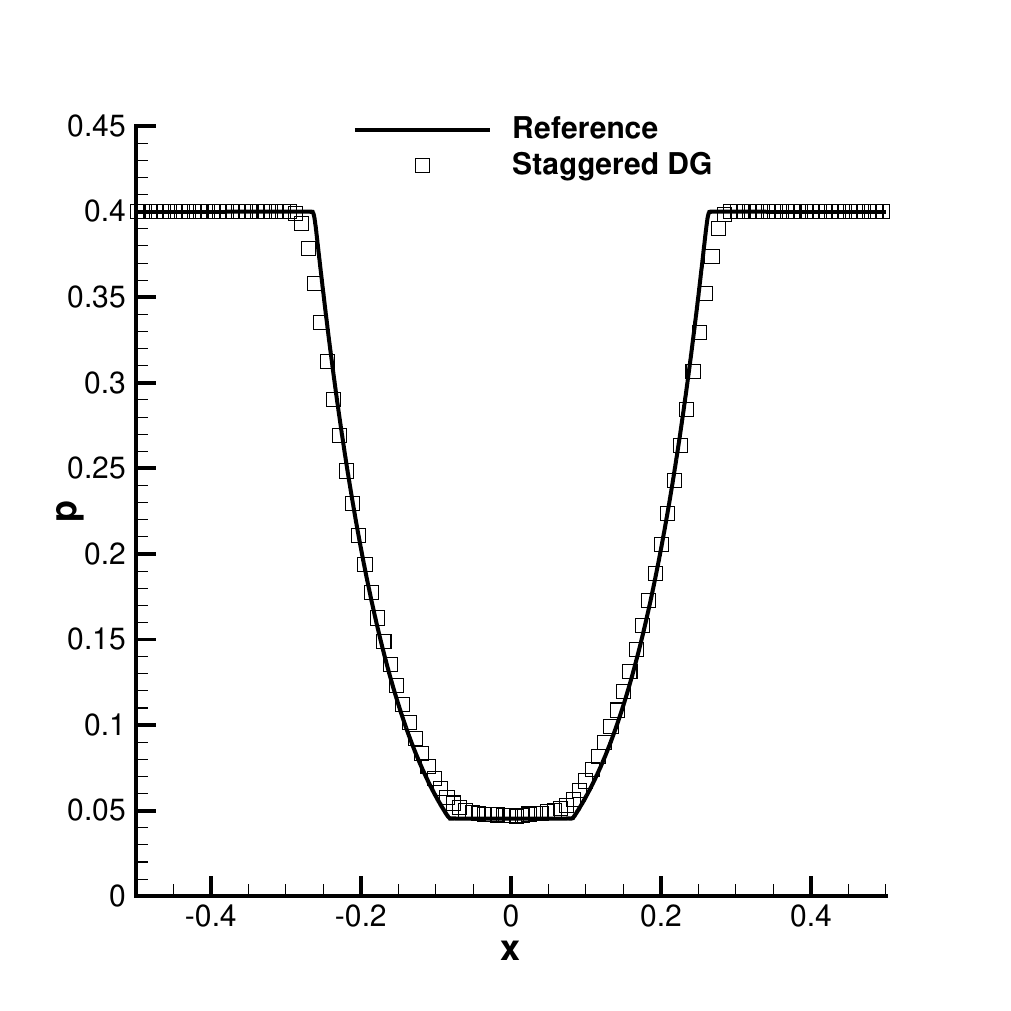}     \\  
		\multicolumn{3}{c}{
		\includegraphics[width=0.55\textwidth]{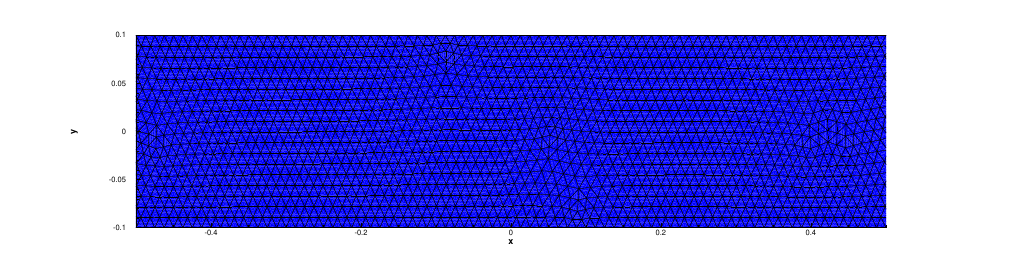}
		}
		\end{tabular} 
    \caption{Riemann problem RP3 at time $t=0.15$ using $p=4$. Top row from left to right: 1D cut through the density, velocity and pressure profile compared against the exact solution. 
		         Bottom row: computational grid, limited cells highlighted in red and unlimited cells in blue.}
    \label{fig:NT_10_3}
\end{figure}

\begin{figure}[ht!]
    \begin{center}
		\includegraphics[width=0.32\textwidth]{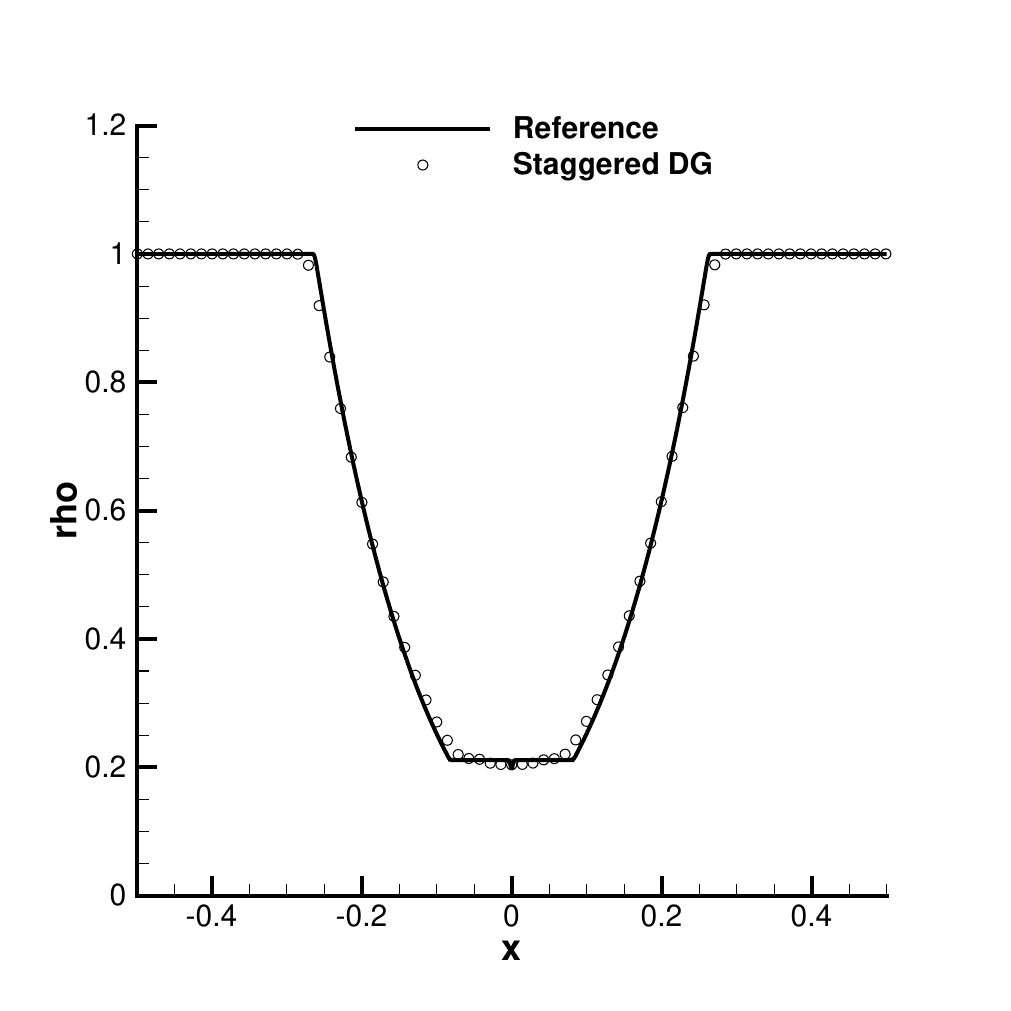}
		\includegraphics[width=0.32\textwidth]{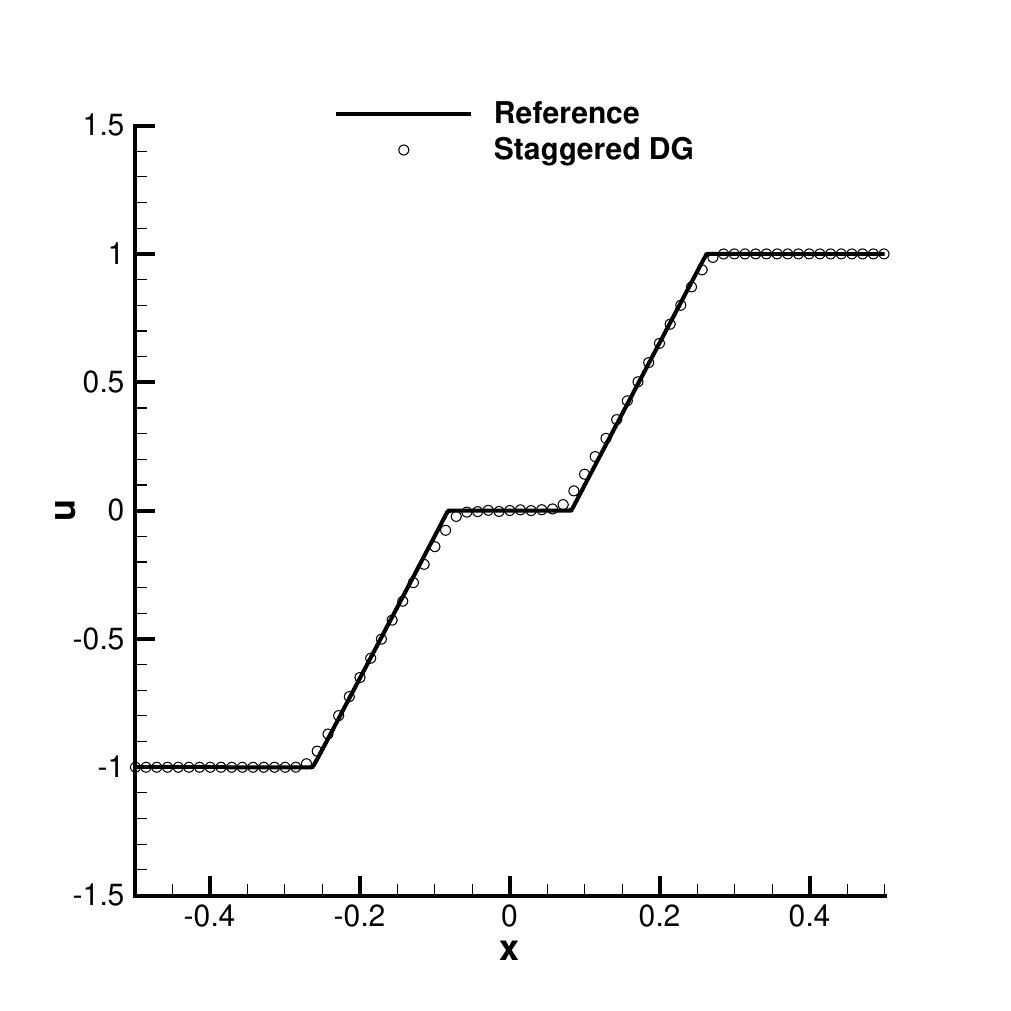}
		\includegraphics[width=0.32\textwidth]{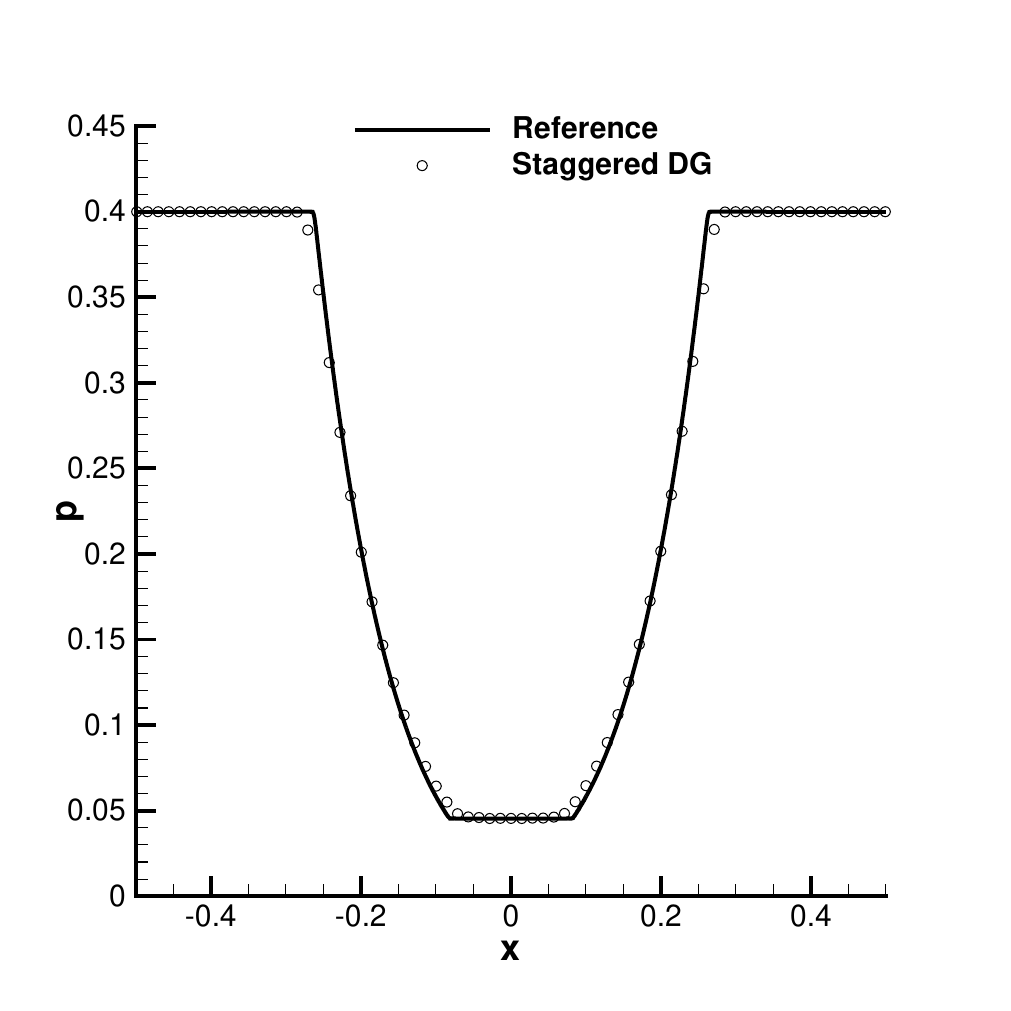}
    \caption{Same as in Fig. \ref{fig:NT_10_3}, but without the limiter.}
    \label{fig:NT_10_3_2}
	\end{center}
\end{figure}

The next test consists in a double rarefaction problem (RP3), where we have used $x_0=0$ and $\epsilon_0=5\cdot 10^{-3}$. Since no shocks appear, this test problem can be considered as smooth, 
except from the initial condition. Figure \ref{fig:NT_10_3} shows the numerical results obtained for this test problem. Also in this case the exact solution is well reproduced, except from an 
unphysical value of the density in the origin. At the final time we can observe that no cells are limited, as expected. However, in the very beginning of the simulation the limiter is activated 
due to the steep gradients present in the solution. The same simulation but with deactivated limiter reported in Figure \ref{fig:NT_10_3_2} shows instead a very good agreement with the exact 
solution also at the origin, since no artificial viscosity is added to the solution. 

\begin{figure}[ht!]
		\begin{tabular}{ccc} 
		\includegraphics[width=0.33\textwidth]{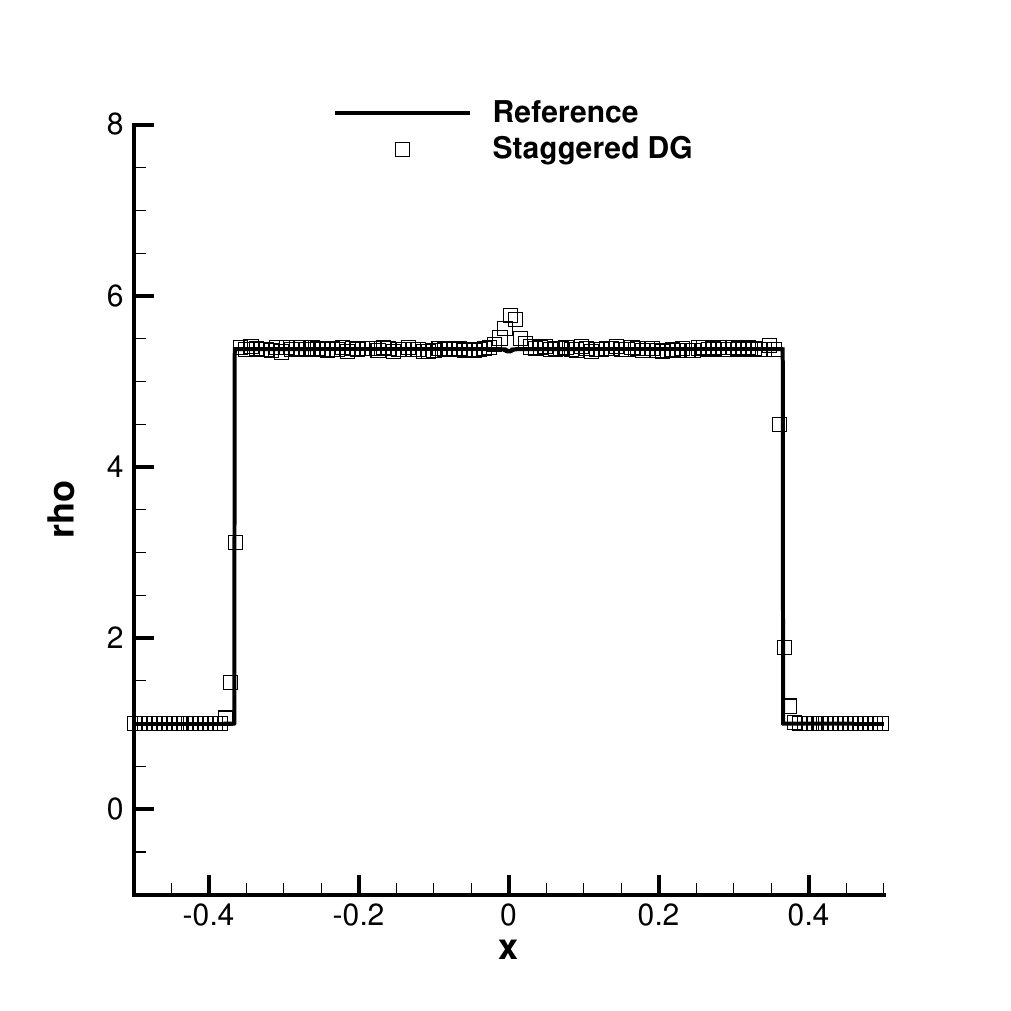}   & 
		\includegraphics[width=0.33\textwidth]{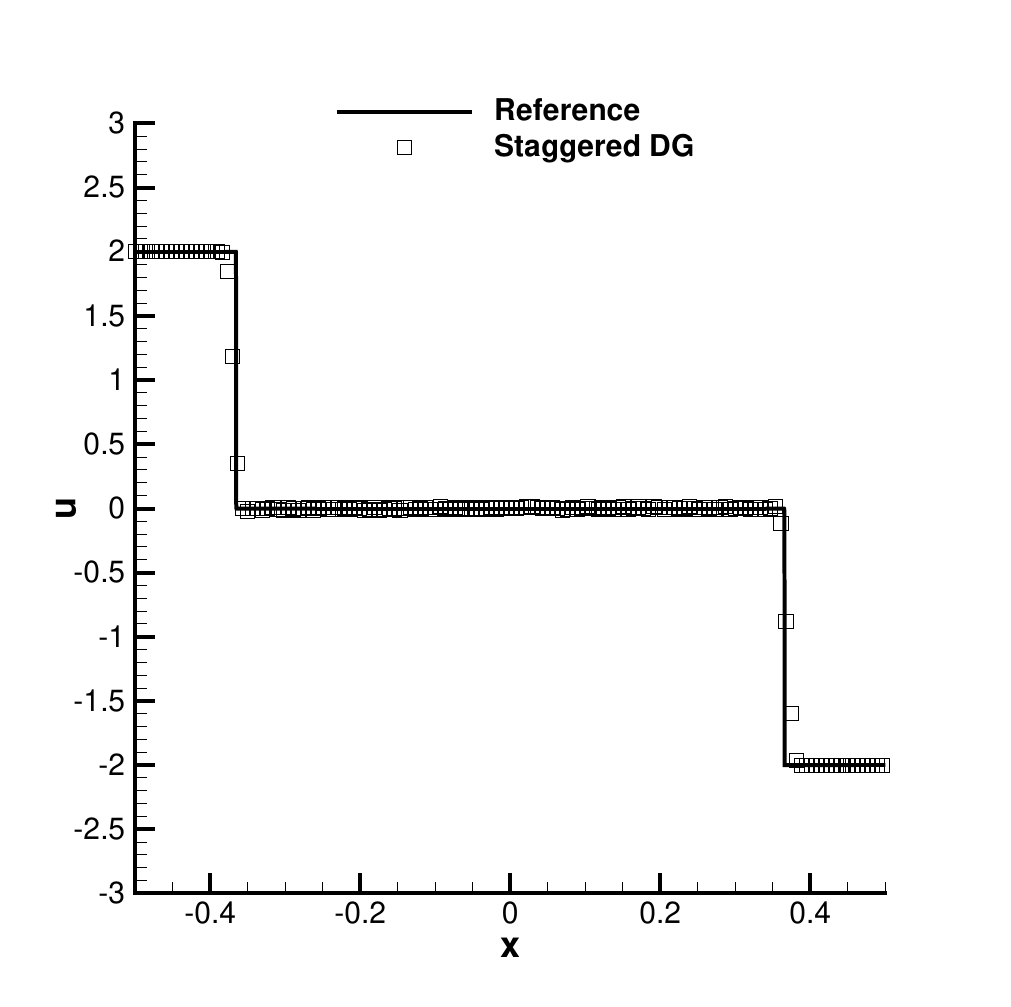}     & 
		\includegraphics[width=0.33\textwidth]{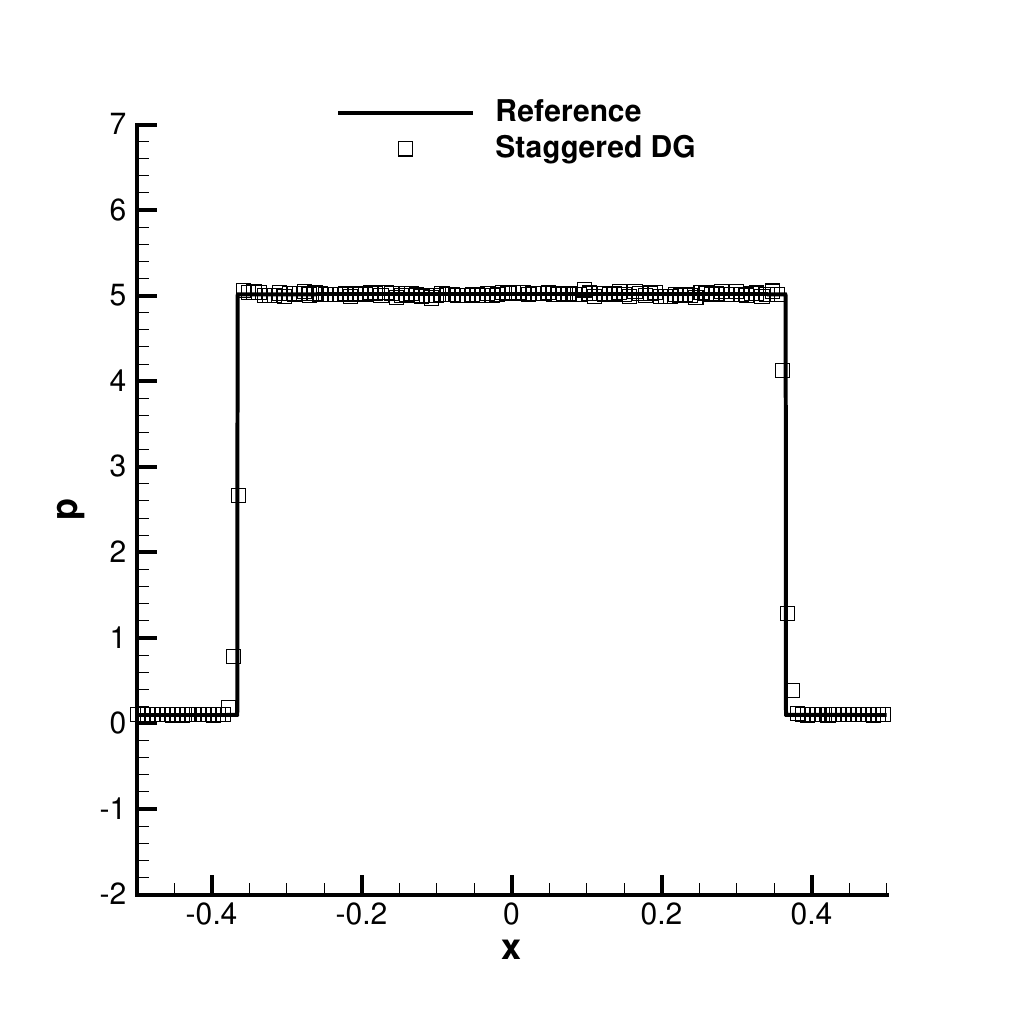}     \\  
		\multicolumn{3}{c}{
		\includegraphics[width=0.55\textwidth]{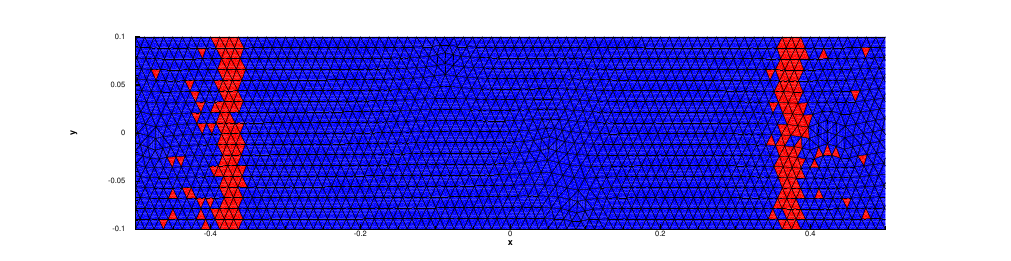}
		}
		\end{tabular} 
    \caption{Riemann problem RP4 at time $t=0.8$ using $p=4$. Top row from left to right: 1D cut through the density, velocity and pressure profile compared against the exact solution. 
		         Bottom row: computational grid, limited cells highlighted in red and unlimited cells in blue.}
    \label{fig:NT_10_4}
\end{figure}

\begin{figure}[ht!]
		\begin{tabular}{ccc} 
		\includegraphics[width=0.33\textwidth]{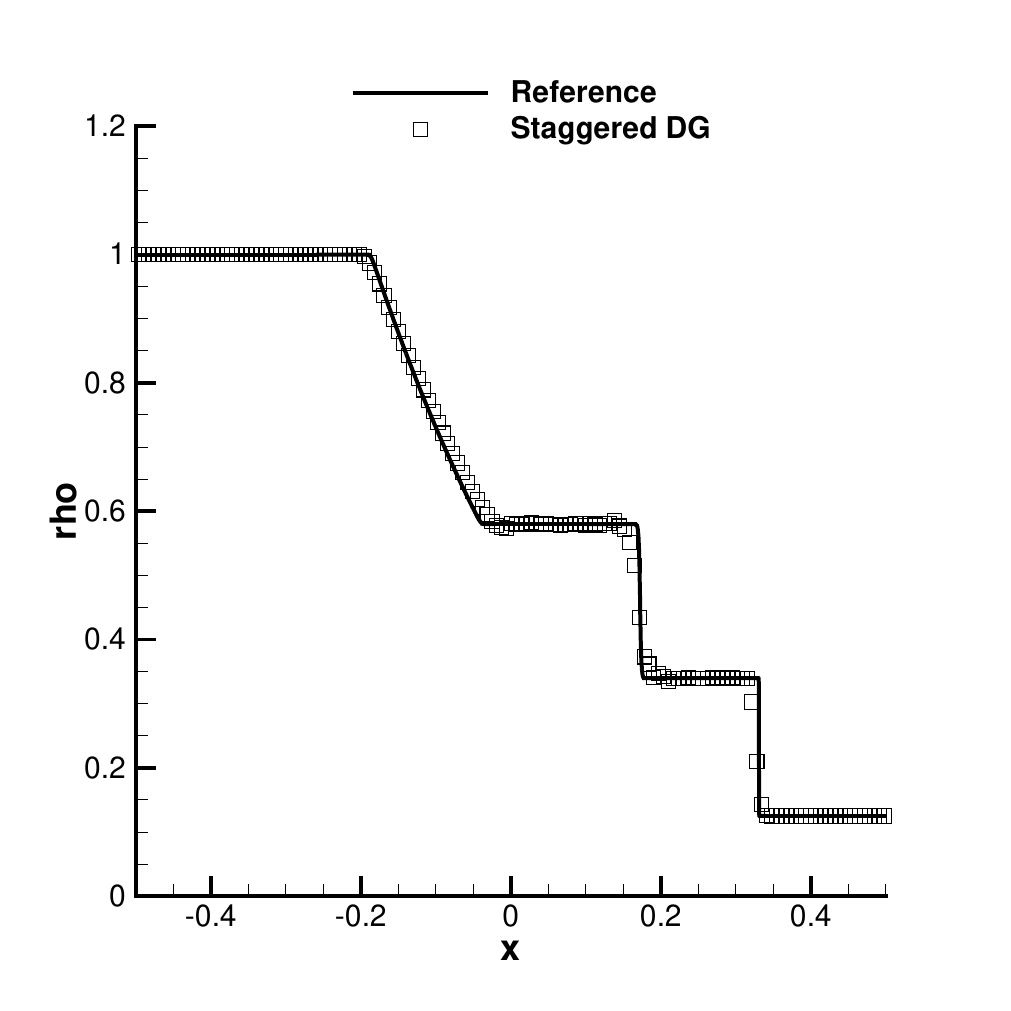}   & 
		\includegraphics[width=0.33\textwidth]{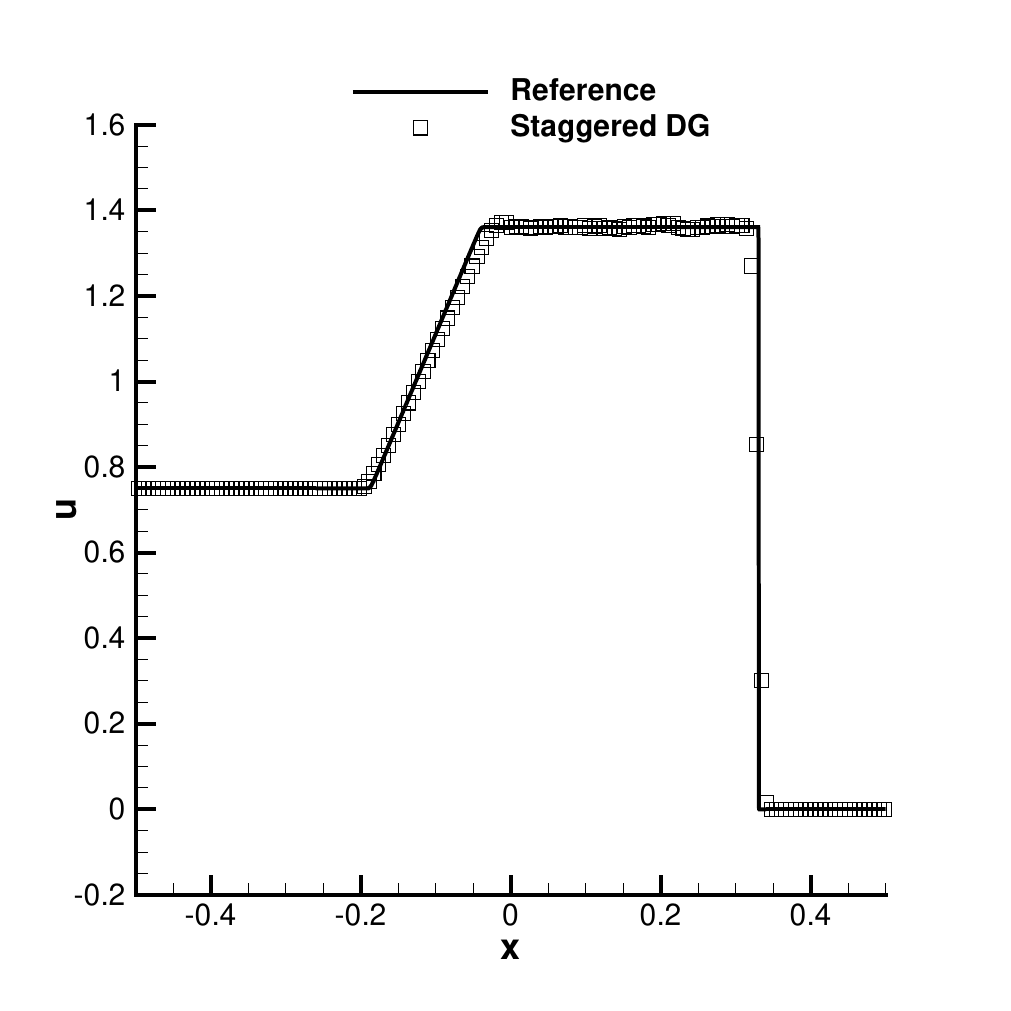}     & 
		\includegraphics[width=0.33\textwidth]{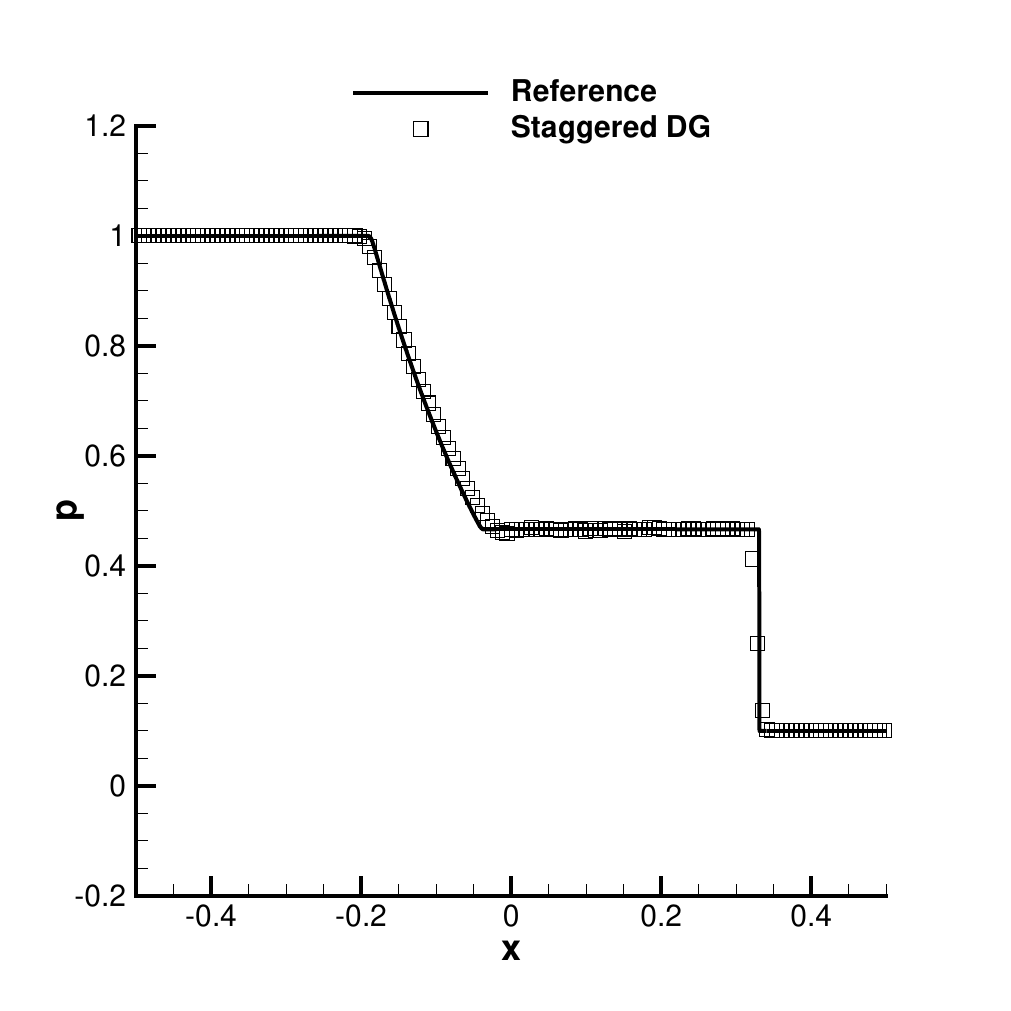}     \\  
		\multicolumn{3}{c}{
		\includegraphics[width=0.55\textwidth]{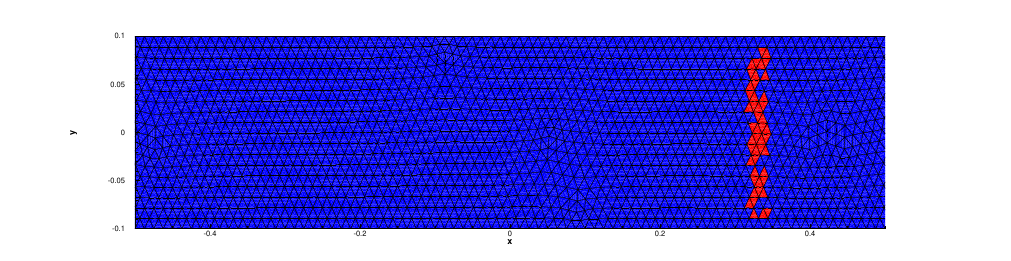}
		}
		\end{tabular} 
    \caption{Riemann problem RP5 at time $t=0.2$ using $p=4$. Top row from left to right: 1D cut through the density, velocity and pressure profile compared against the exact solution. 
		         Bottom row: computational grid, limited cells highlighted in red and unlimited cells in blue.}
    \label{fig:NT_10_5}
\end{figure}

Now we consider two strong colliding shocks (RP4) setting $x_0=0$ and $\epsilon_0=1\cdot 10^{-2}$. Figure \ref{fig:NT_10_4} reports the numerical results obtained with our staggered DG scheme. 
In this case we obtain a similar unexpected behavior for the density in the origin as already observed in the second order case, see \cite{DumbserCasulli2016}. So it seems that also the high 
order extension of the pressure-based scheme \cite{DumbserCasulli2016} is affected by a similar local wall-cooling error. Apart from the origin, the solution shows some small oscillations 
but a relatively good shape if compared with the exact reference solution, see Figure \ref{fig:NT_10_4}.  
The last Riemann problem that we consider here is a modified version of the Sod problem originally proposed by Toro in \cite{toro-book}. The initial state and the final time for this test are 
reported in Table \ref{tab.NT_10_1} as (RP5); we furthermore use $x_0=-0.1$ and $\epsilon_0 = 1\cdot 10^{-2}$. The purpose of this test is to show whether the proposed numerical method is 
affected by the well-known sonic glitch problem that typically appears at the sonic point inside transonic rarefaction waves. Figure \ref{fig:NT_10_5} shows the results obtained with the 
staggered DG method compared with the exact solution. A good agreement can be observed also in this case with only some small oscillations in the region between the shock and the rarefaction. 
Note that no sonic glitch is present inside the rarefaction fan, as already observed for the low order scheme \cite{DumbserCasulli2016}. 

\subsubsection{Circular explosion problem}
\begin{figure}[ht!]
    \begin{center}
		\includegraphics[width=0.45\textwidth]{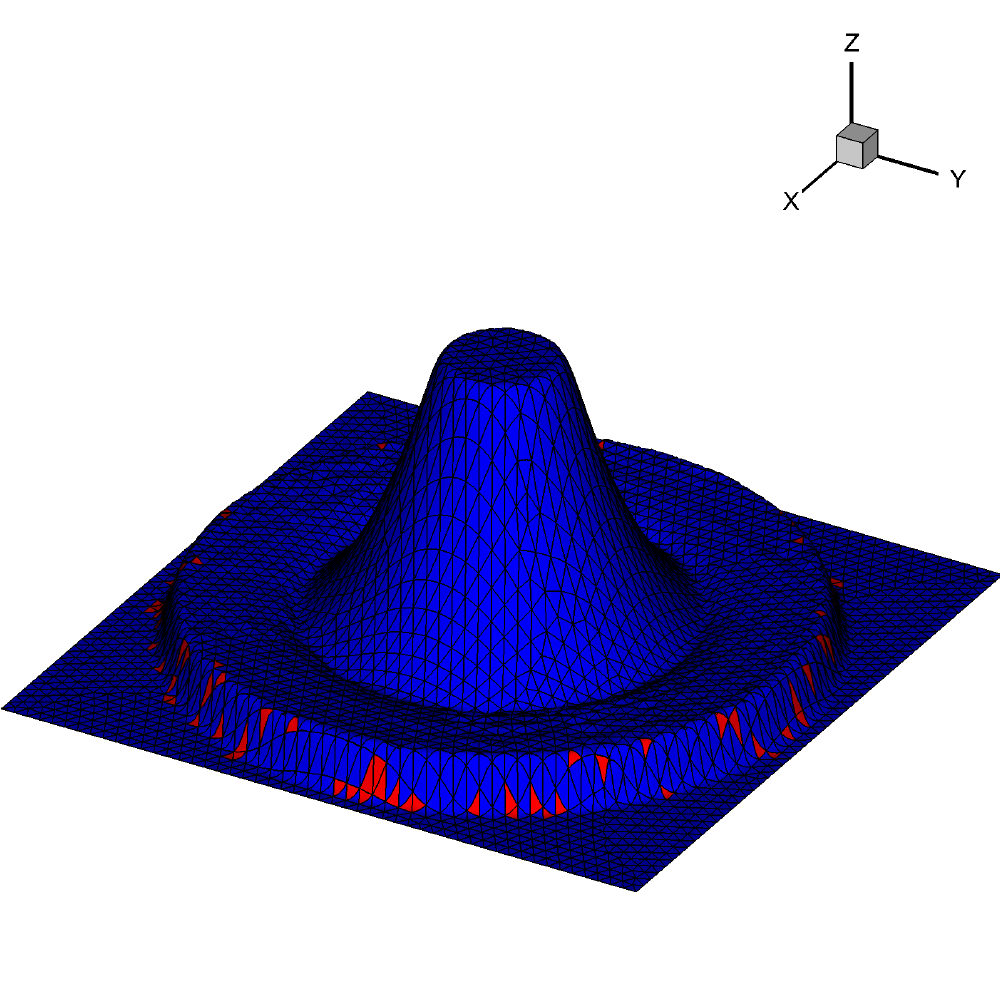}	
    \includegraphics[width=0.32\textwidth]{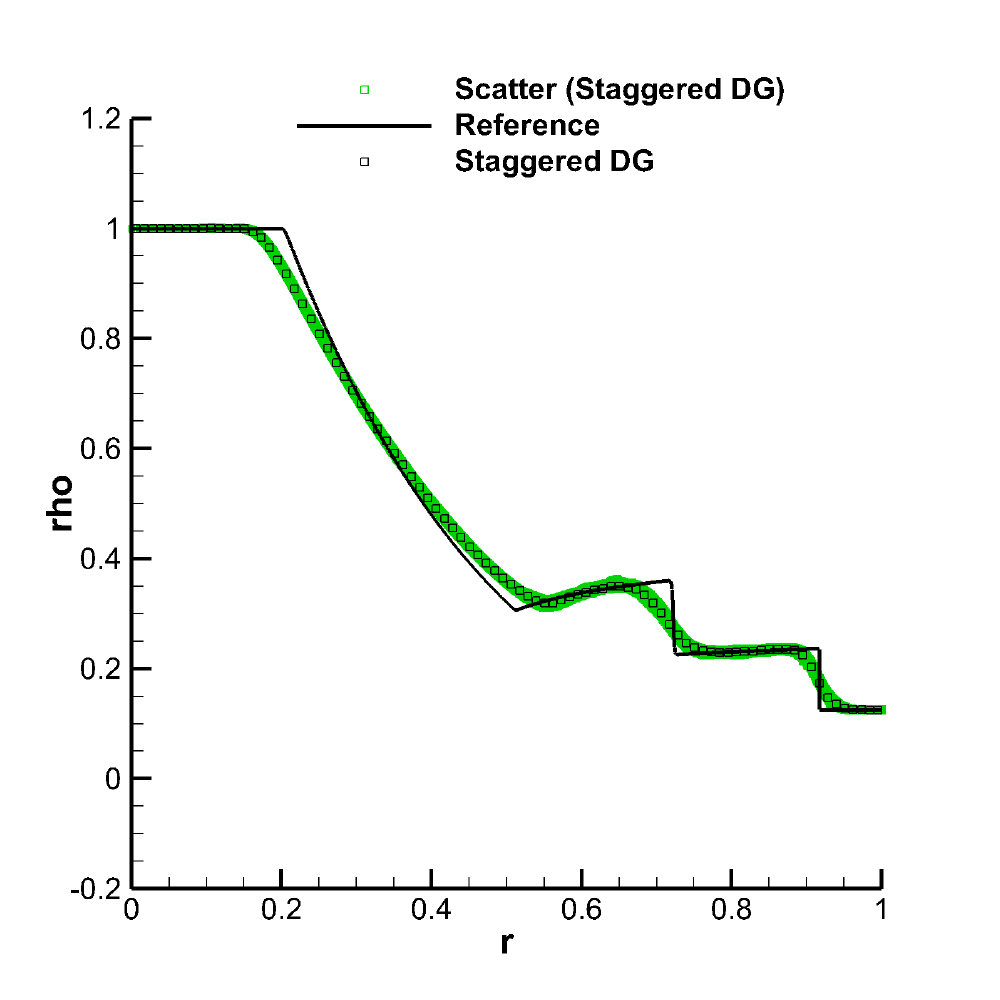} \\ 
		\includegraphics[width=0.32\textwidth]{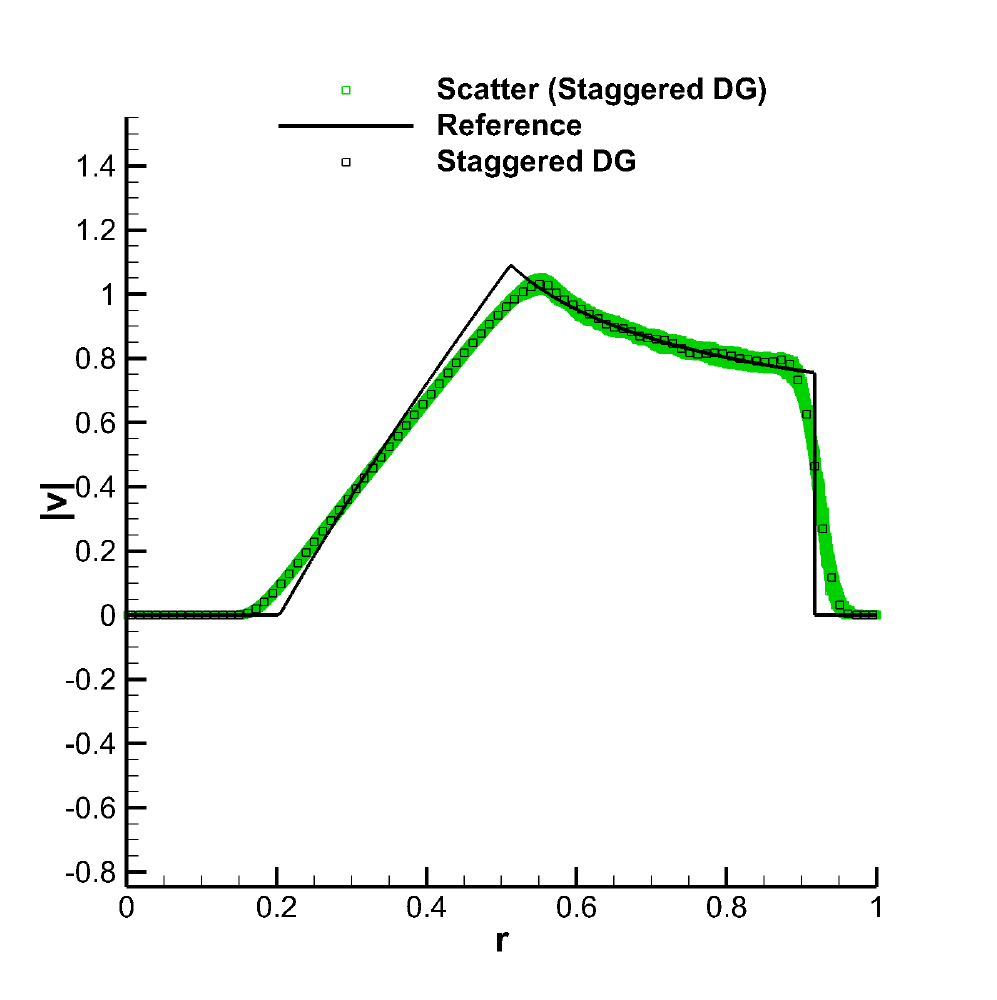}
		\includegraphics[width=0.32\textwidth]{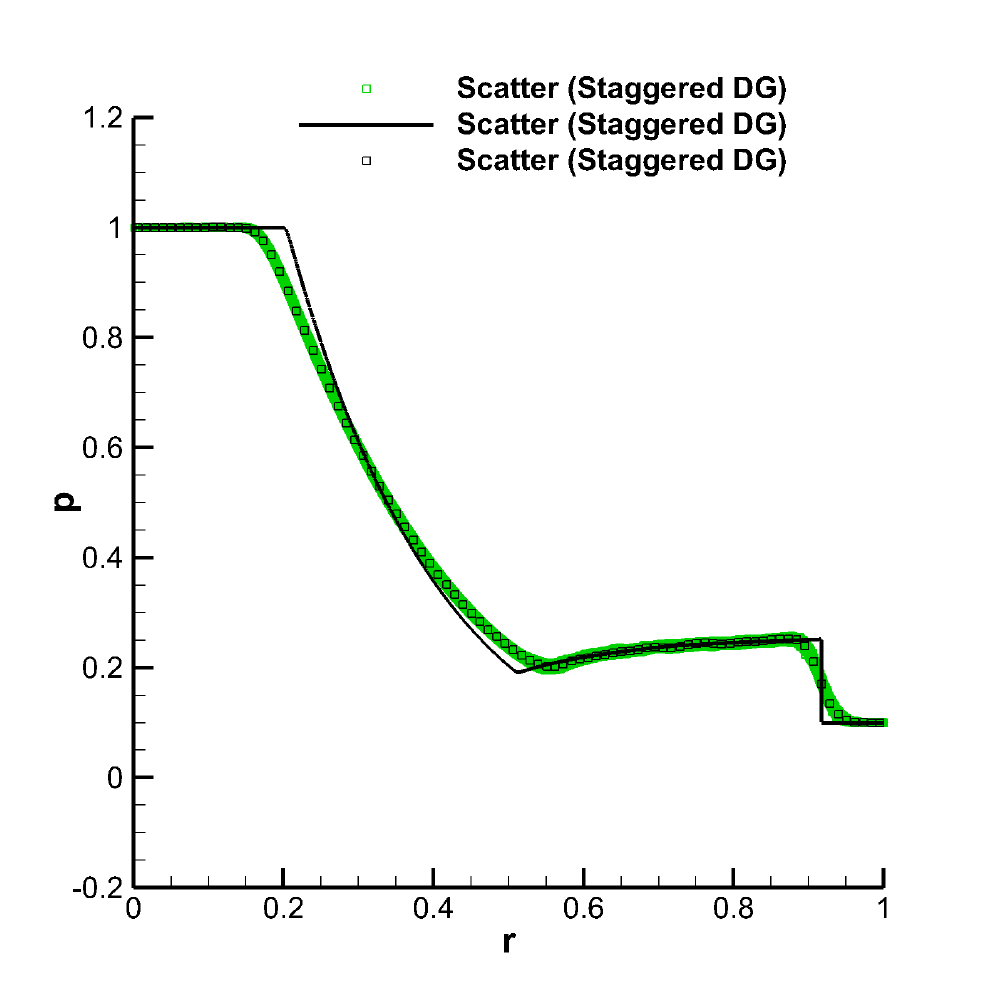}
    \caption{Two-dimensional explosion problem at time $t=0.25$. 3D plot of the density with the computational mesh and the limited cells highlighted in red, while unlimited cells are plotted blue (top left). Scatter plots and 1D profiles of density (top right), velocity magnitude (bottom left) and pressure (bottom right).}
    \label{fig:ep2d}
	\end{center}
\end{figure}

While the previous examples have shown the application of the limiter to essentially one dimensional problems, we consider here a two dimensional explosion problem. 
In particular we take as initial condition: 
\begin{eqnarray}
	(\rho_0,u_0,v_0,p_0)=\left\{\begin{array}{ll}
		(1,0,0,1)       & r \leq 0.5 \\
		(0.125,0,0,0.1) & r > 0.5
	\end{array}\right.
\label{eq:NT_11_1}
\end{eqnarray}
The computational domain $\Omega=[-1,1]^2$ is covered with $\Ni=5616$ triangles; we set $(p,p_\gamma)=(3,0)$ and transmissive boundary conditions are applied everywhere. We perform the simulation 
up to $t_{end}=0.25$. A reference solution can be obtained by considering the angular symmetry of the problem, hence the 2D Euler equations in polar coordinates reduce to a 1D system in radial
direction with geometric source terms, see \cite{toro-book}. For the solution of the 1D system in radial direction, which will serve as reference solution, we use a second order TVD scheme 
with Osher-type flux \cite{OsherUniversal} on a very fine grid with $10^4$ cells. Figure $\ref{fig:ep2d}$ shows a 1D cut as well as a scatter plot of the numerical results for the main flow 
variables compared with the the reference solution. Overall we can observe a good agreement and the limiter seems to act only at the discontinuity, as also shown in Figure $\ref{fig:ep2d}$. 


\subsection{Three dimensional numerical tests}
In this section we report some numerical results for the three-dimensional version of the proposed algorithm. In particular we will apply the method to some low to moderate Mach number flows 
such as the three-dimensional propagation of a smooth acoustic wave, the 3D lid-driven cavity flow, the classical Taylor Green vortex and the flow past a sphere at $M=0.5$. In this last case 
we use a fully isoparametric description of the geometry with curved elements. 

\subsubsection{Smooth acoustic wave propagation in 2D}
\begin{figure}[ht!]
    \begin{center}
		\includegraphics[width=0.4\textwidth]{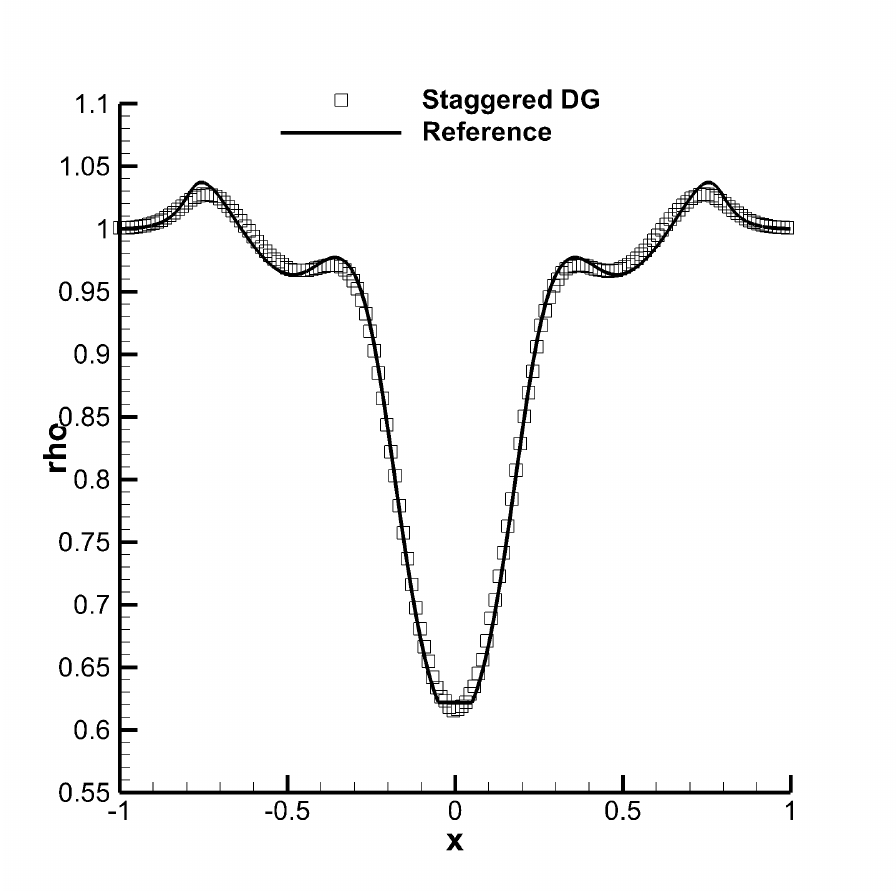}	
		\includegraphics[width=0.4\textwidth]{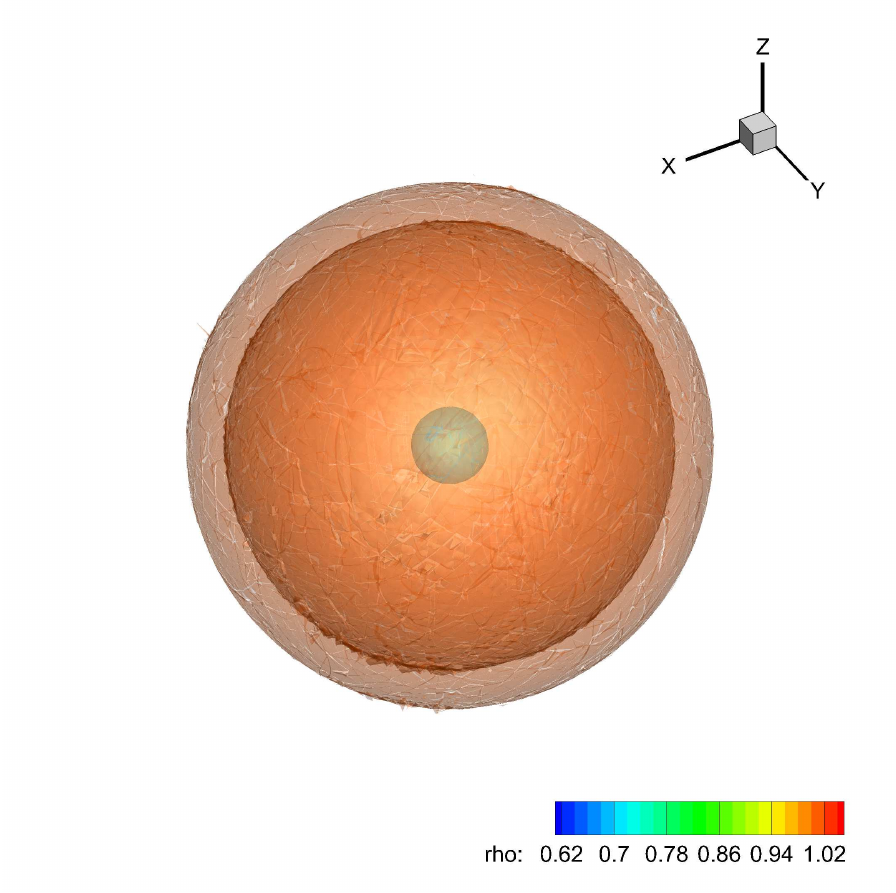}	\\
		\includegraphics[width=0.4\textwidth]{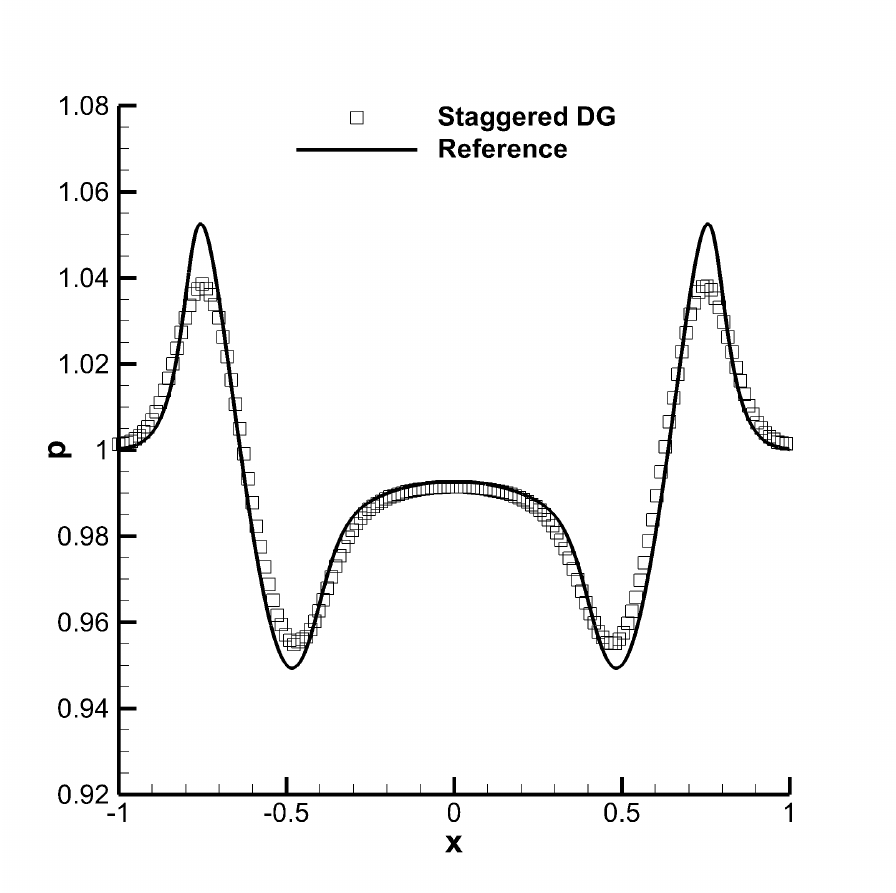}	
		\includegraphics[width=0.4\textwidth]{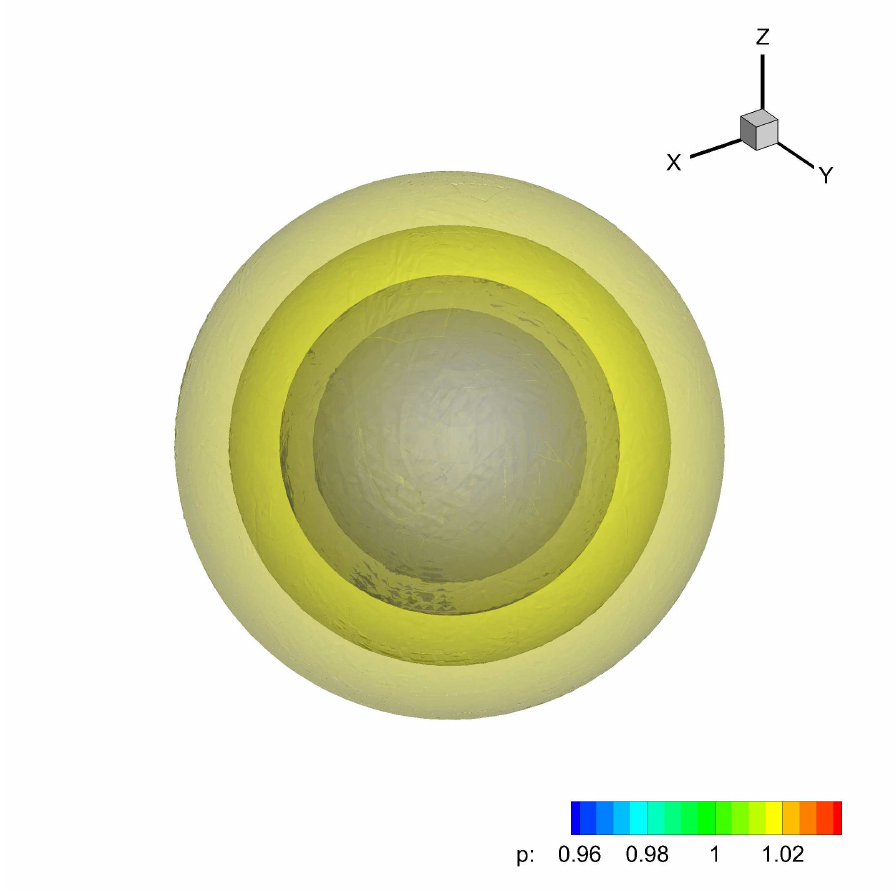}	\\
		\includegraphics[width=0.4\textwidth]{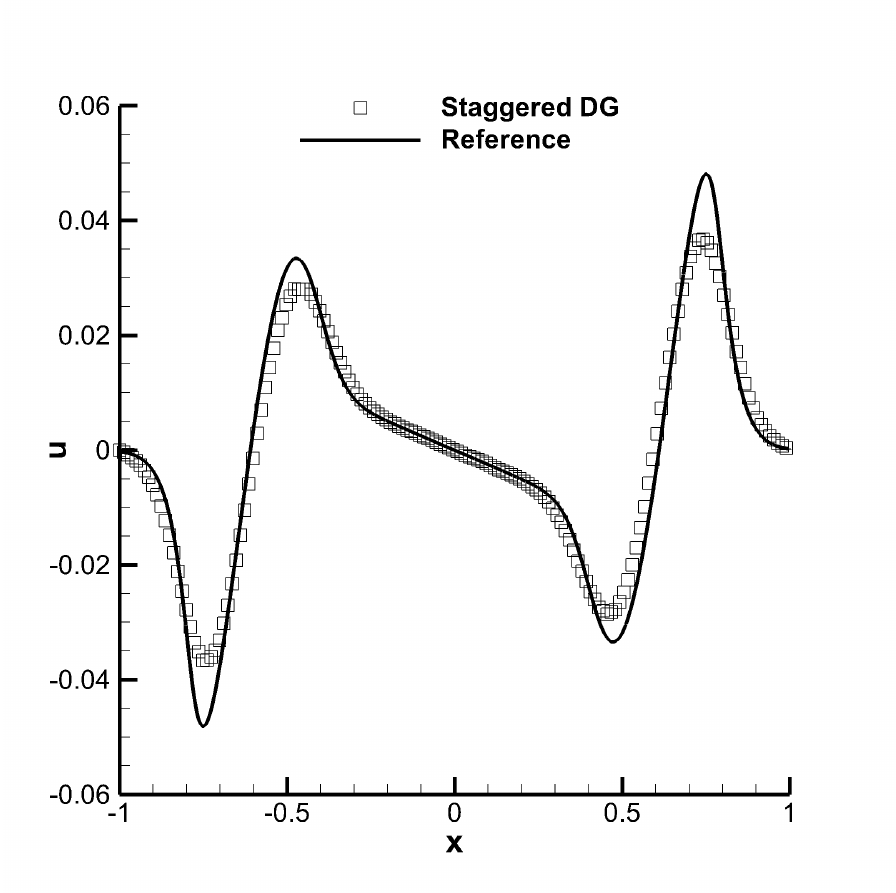}	
		\includegraphics[width=0.4\textwidth]{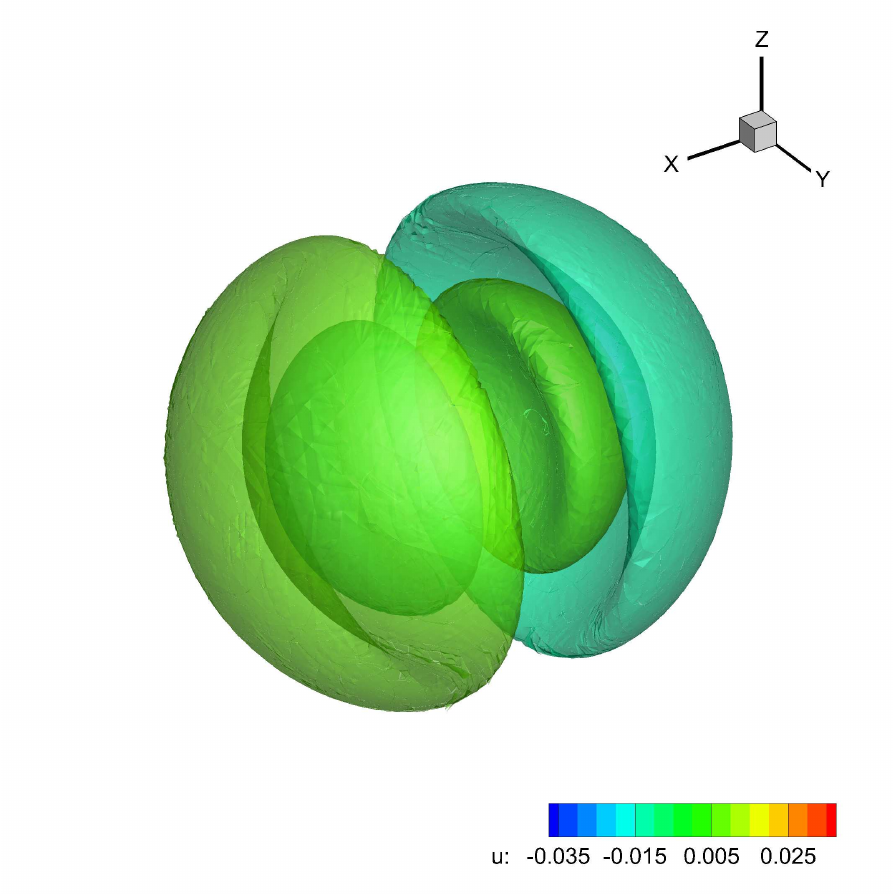}	\\
    \caption{Density, pressure and velocity component $u$; Isosurface (left column) and a cut of the solution compared with the reference one (right one).}
    \label{fig:NT_112_1}
	\end{center}
\end{figure}

The first test that we consider is a simple three dimensional propagation of a smooth acoustic wave, equivalent to the test presented in Section \ref{sec_2Dsmooth} for the two 
dimensional case. The initial condition is the same as the one given in Section \ref{sec_2Dsmooth} and the problem is inviscid ($\mu=0$). 
Due to the symmetry of the problem, we can use again a second order TVD scheme  on $10^4$ cells applied to a reduced 1D system with geometric source terms \cite{toro-book} in order to 
compute the reference solution for the three-dimensional case. 
The computational domain is $\Omega=[-1,1]^3$ with transmissive boundaries everywhere. The domain $\Omega$ is covered with $\Ni=26082$ tetrahedra; $(p,p_\gamma)=(3,0)$ and $t_{end}=0.5$. 
A sketch of the iso-surfaces for the main variables as well as a one dimensional cut along the line $y=0$ and $z=0$ is compared with the reference solution in Figure \ref{fig:NT_12_1}. 
A reasonably good agreement of the numerical solution with the reference has been obtained also in this simple test, where the three dimensional effects change the shape of the solution
compared to the 2D simulation. 
Furthermore, we can observe from Figure \ref{fig:NT_112_1} that the resulting numerical solution is symmetric, despite the use of an unstructured mesh and the wave propagation speed is 
correct. We only note some numerical dissipation close to the velocity and pressure peaks. 

\subsubsection{3D Lid-driven cavity}
\begin{figure}[ht!]
    \begin{center}
		\includegraphics[width=0.45\textwidth]{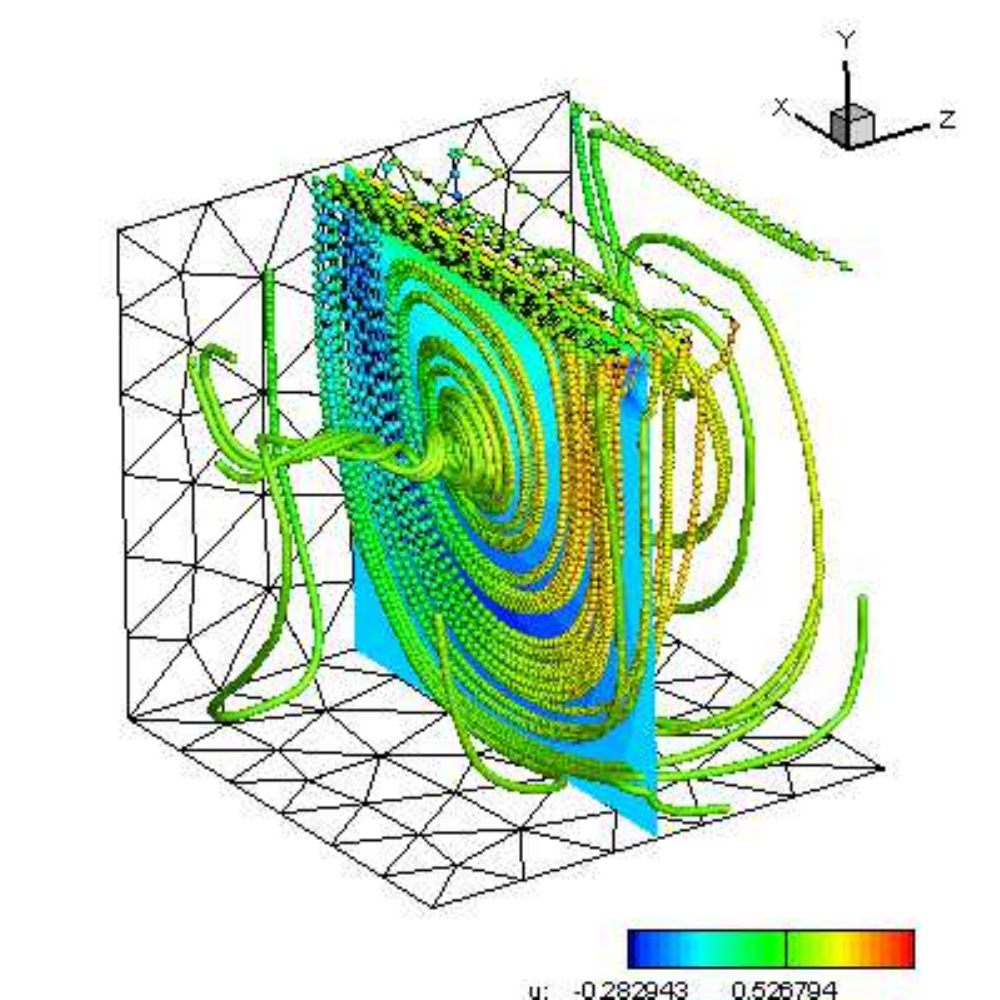}	
		\includegraphics[width=0.45\textwidth]{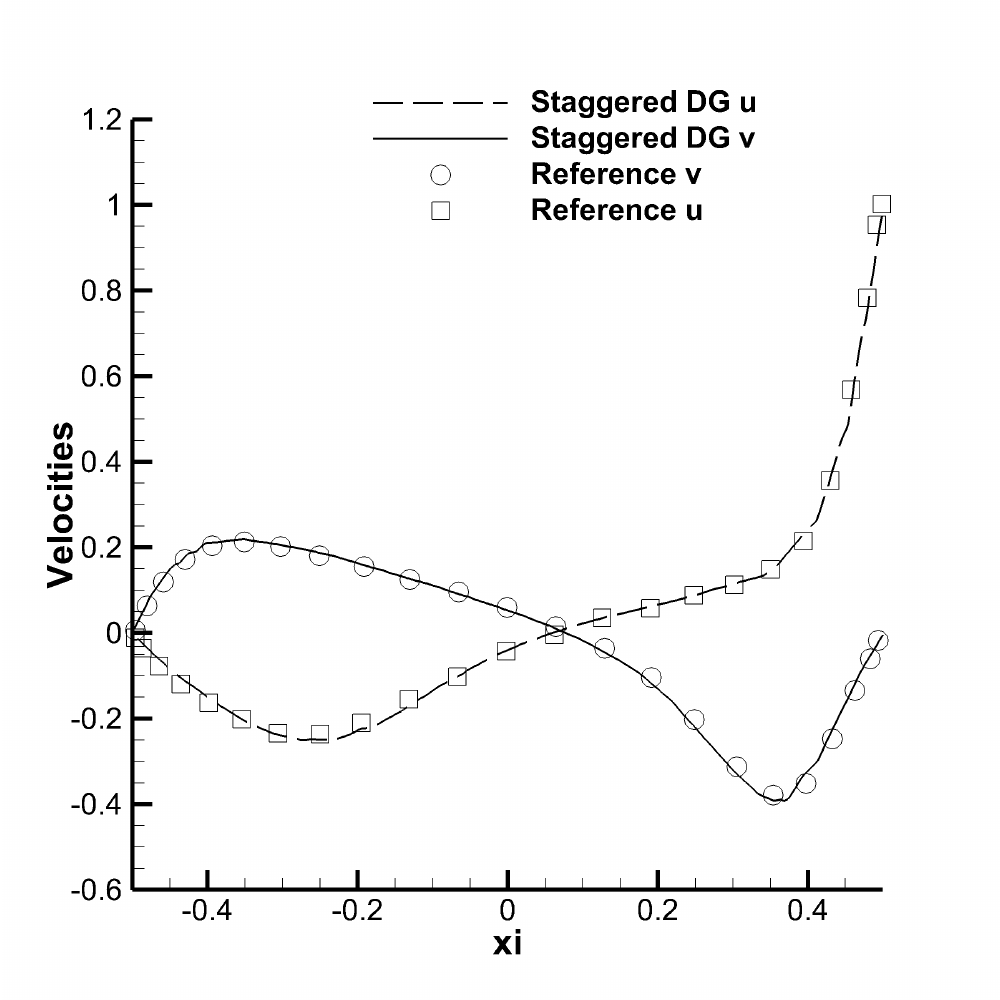}	\\
		\includegraphics[width=0.32\textwidth]{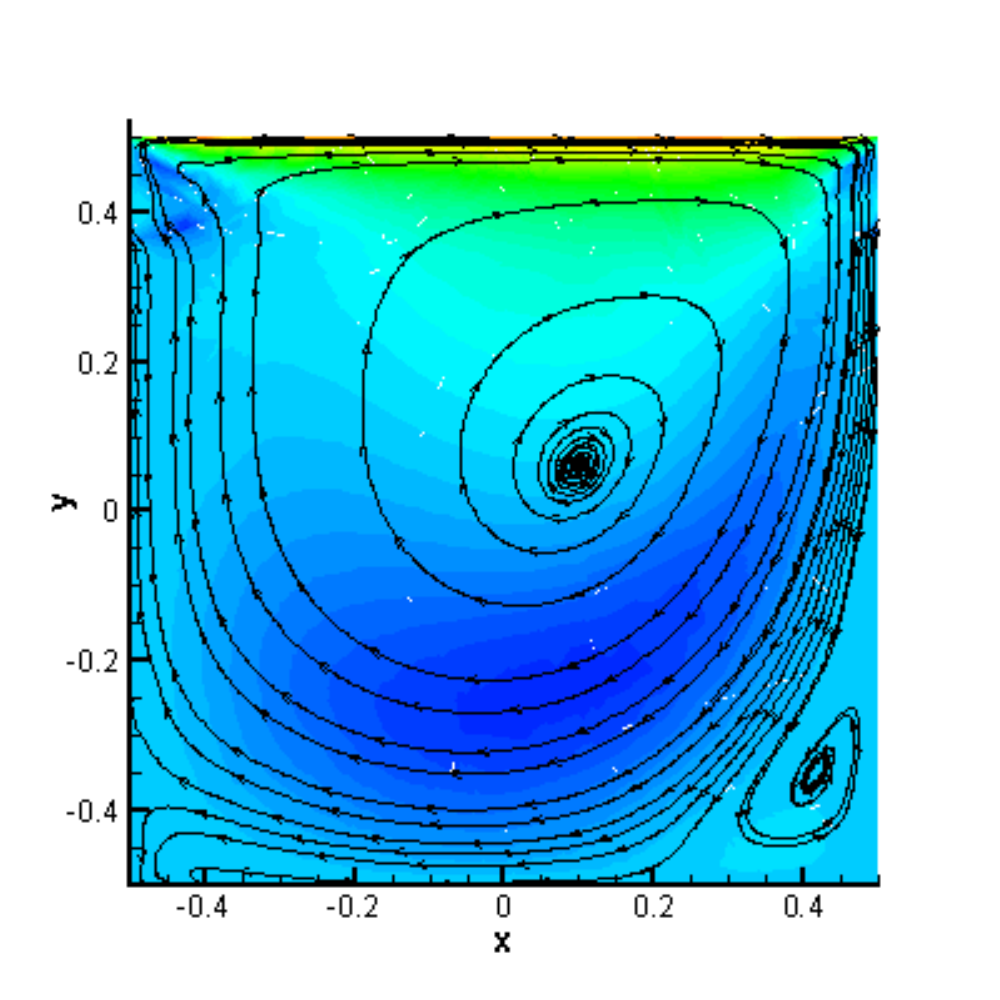}	
		\includegraphics[width=0.32\textwidth]{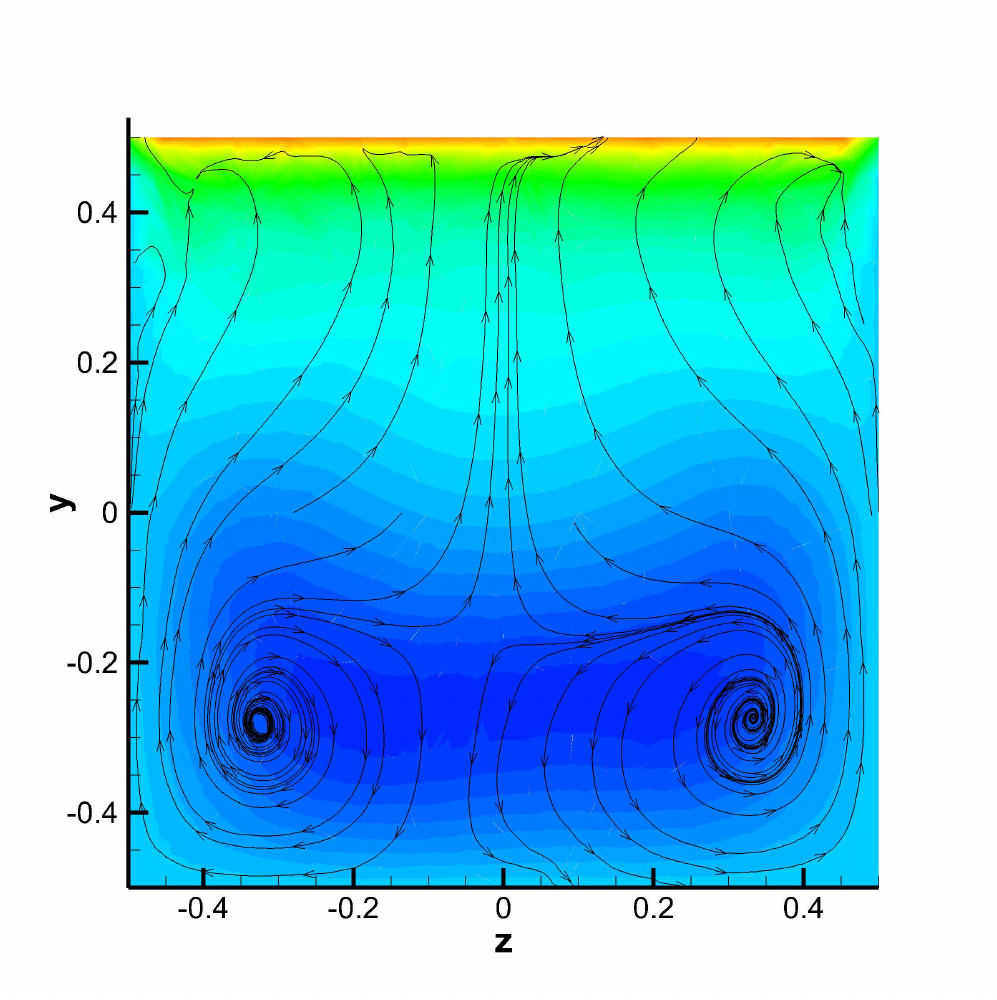}
		\includegraphics[width=0.32\textwidth]{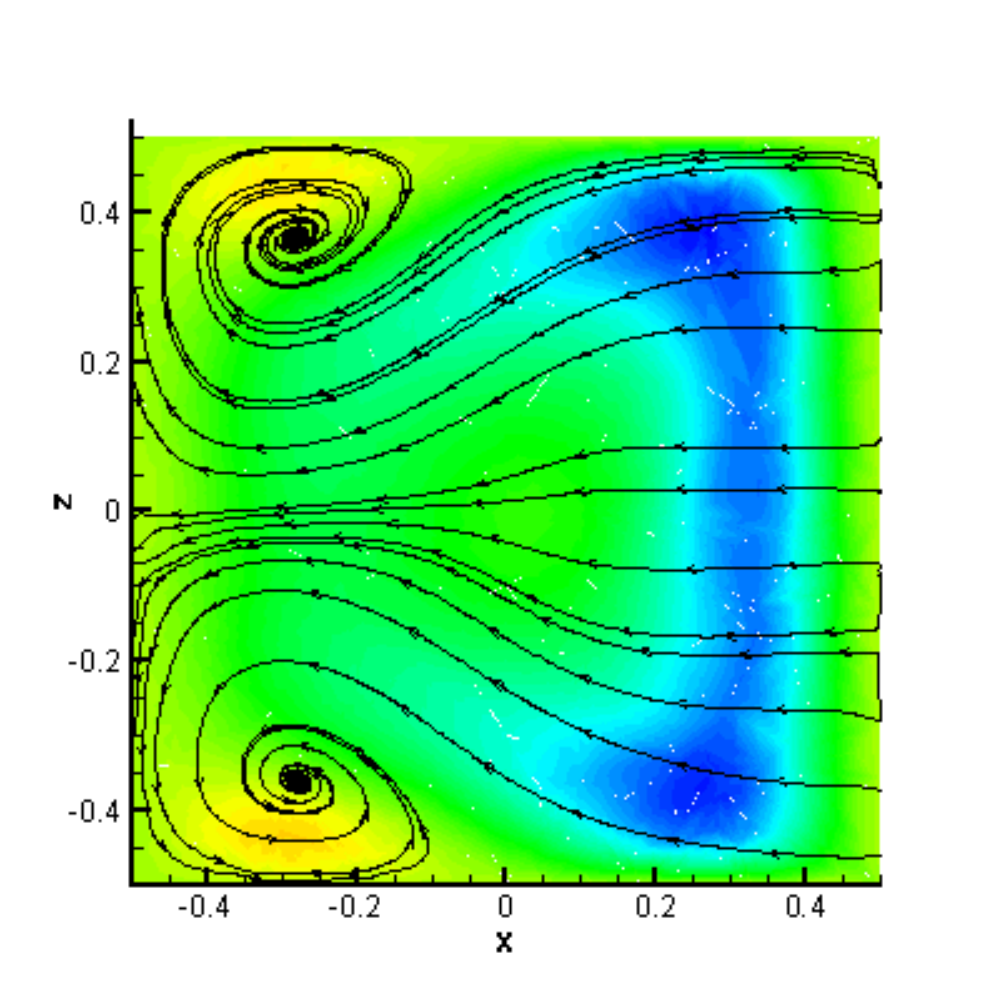}
    \caption{Numerical results for the lid-driven cavity at $Re=400$, from top left to bottom right: Three-dimensional plot with stream traces, comparison with the reference data obtained by Albensoeder et. al. in \cite{Albensoeder2005}, sketch of the main slices $\{z=0\}$, $\{z=0\}$ and $\{y=0\}$ at $t=t_{end}$.}
    \label{fig:NT_12_1}
	\end{center}
\end{figure}

\begin{figure}[ht!]
    \begin{center}
		\includegraphics[width=0.45\textwidth]{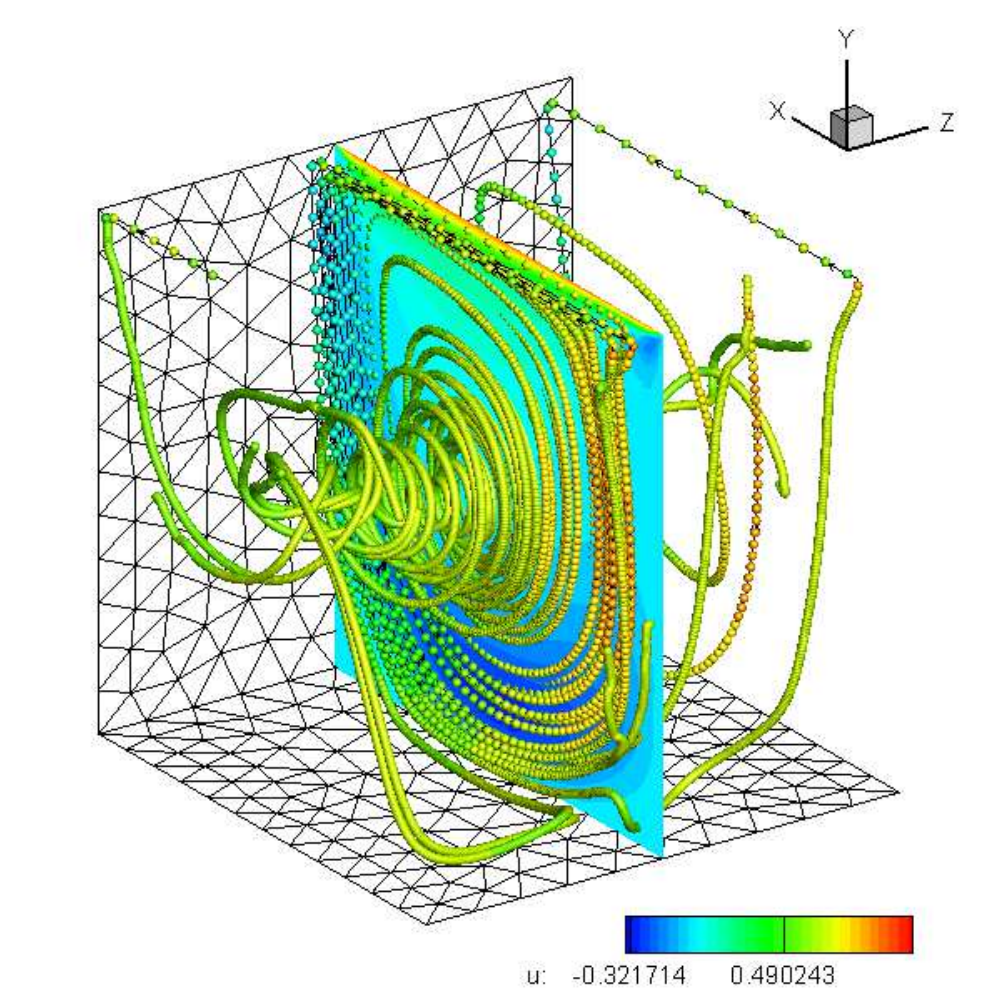}	
		\includegraphics[width=0.45\textwidth]{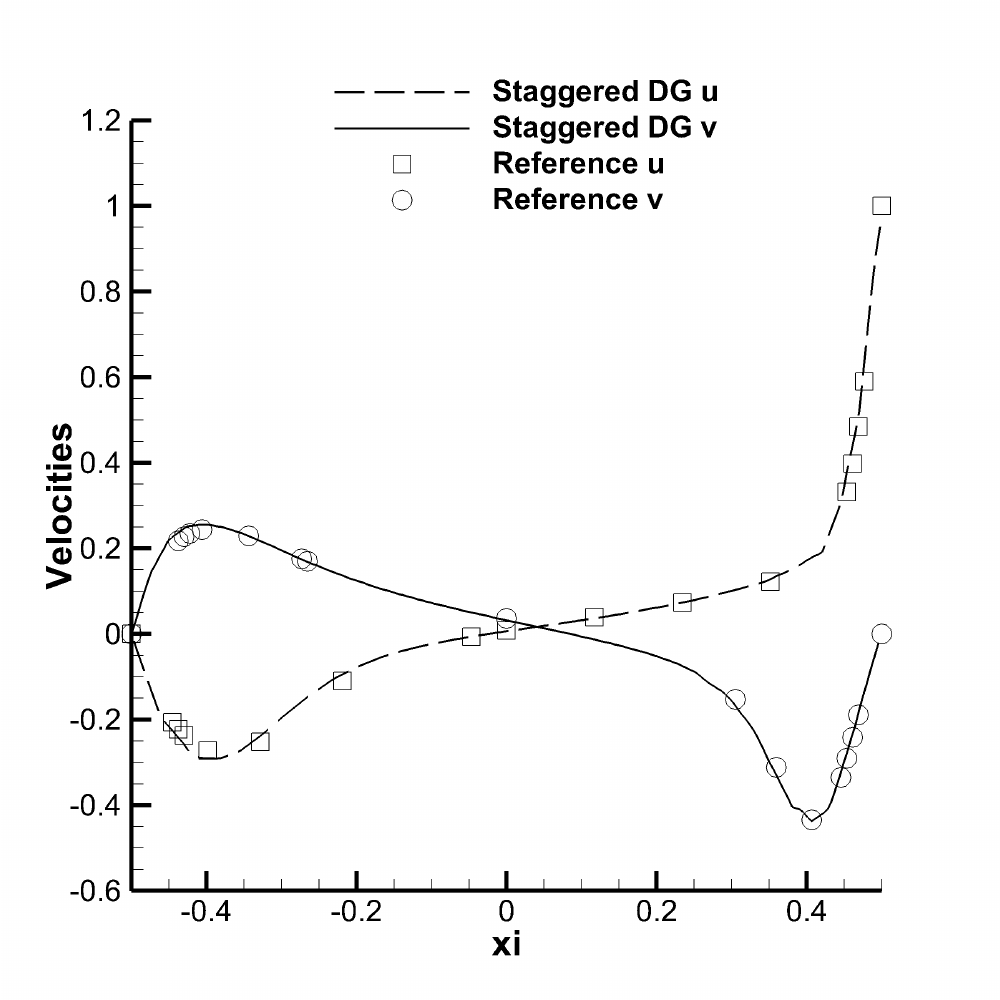}	\\
		\includegraphics[width=0.32\textwidth]{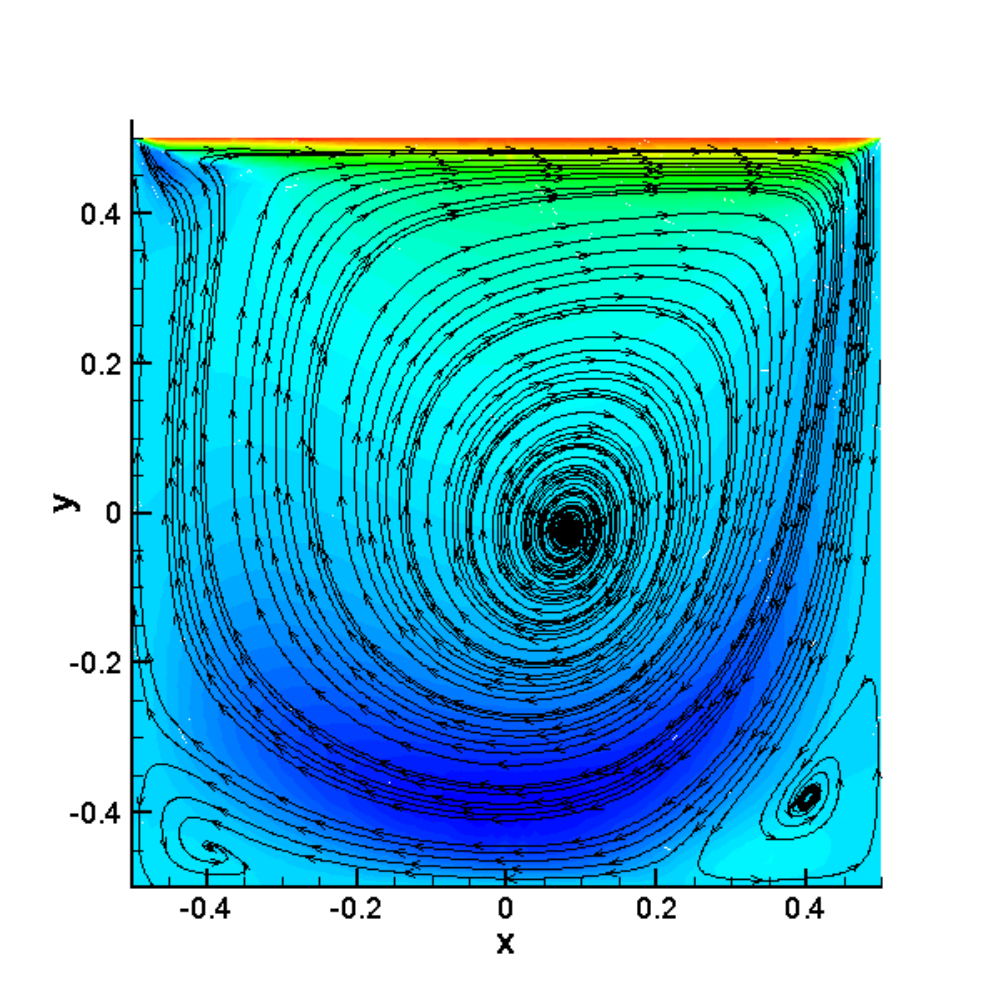}	
		\includegraphics[width=0.32\textwidth]{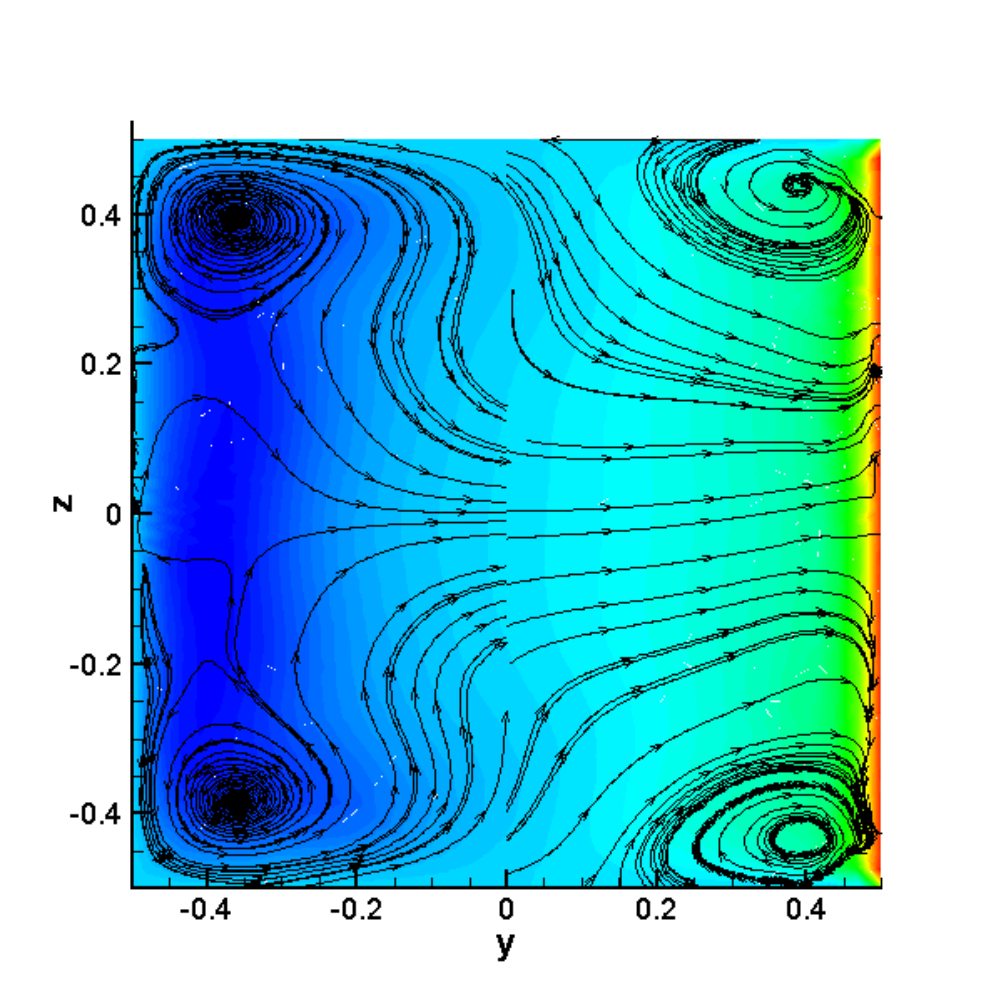}
		\includegraphics[width=0.32\textwidth]{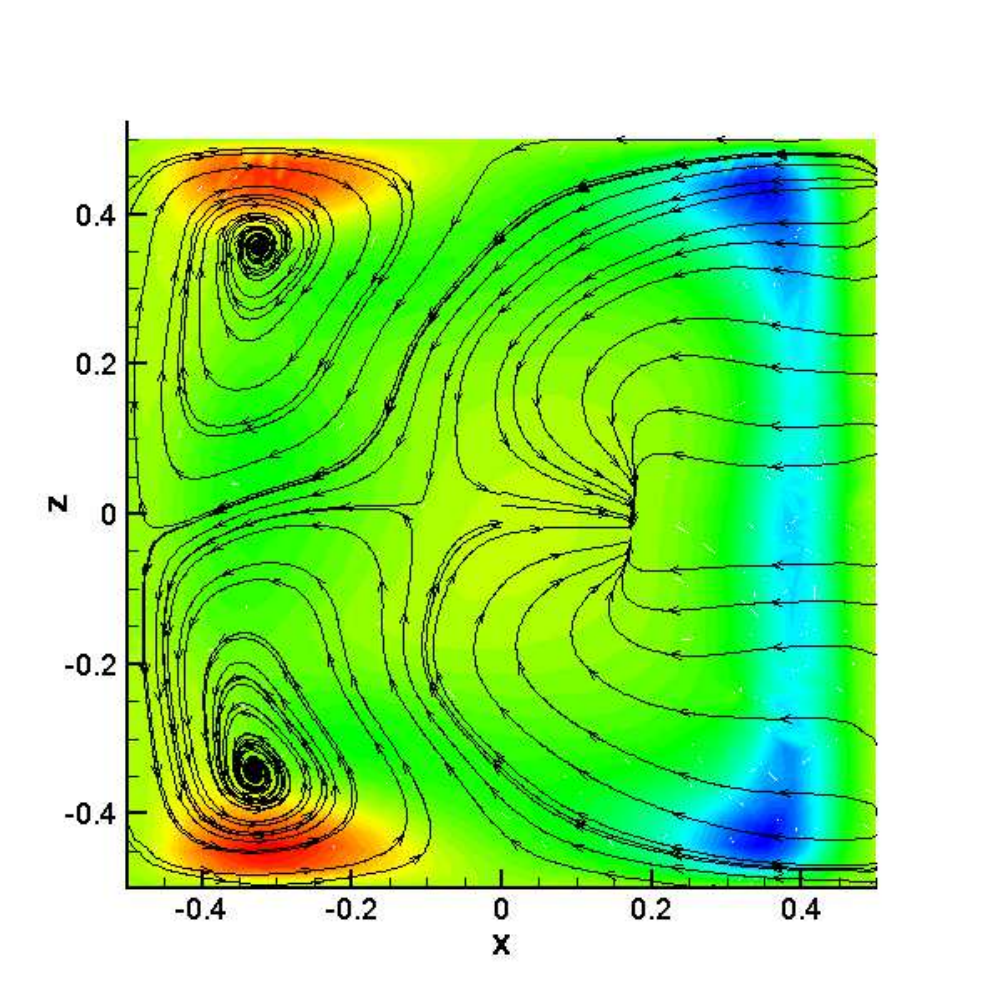}
    \caption{Numerical results for the lid-driven cavity at $Re=1000$, from top left to bottom right: Three-dimensional plot with stream traces, comparison with the reference data obtained by Albensoeder et. al. in \cite{Albensoeder2005}, sketch of the main slices $\{z=0\}$, $\{z=0\}$ and $\{y=0\}$ at $t=t_{end}$.}
    \label{fig:NT_12_2}
	\end{center}
\end{figure}

Here we consider a low Mach number lid-driven cavity flow in three space dimensions. Two dimensional results have already been presented in Section \ref{sec.2DCavity}. For the three dimensional case several reference solutions are available, see e.g. \cite{Albensoeder2005,Hwar1987,Aidun1991}. As computational domain we take a simple cubic cavity $\Omega=[-0.5, 0.5]^3$. The initial condition is  
given by $\mathbf{v}=0$; $\rho=1$; and $p=10^5$. At the upper boundary ($y=0.5$) we impose a constant velocity of $\mathbf{v}=(1,0,0)$, while no-slip wall boundary conditions are imposed for the 
remaining boundaries. 
We consider two different Reynolds numbers, namely $Re=400$ and $Re=1000$. The number of tetrahedra is $\Ni=1380$ for $Re=400$ and $\Ni=11040$ for $Re=1000$, respectively. 
Furthermore, we take $p=3, p_\gamma=0$ and $t_{\textnormal{end}}=36$ for $Re=400$ and $t_{\textnormal{end}}=50$ for $Re=1000$, respectively.  
In Figure \ref{fig:NT_12_1} the results for $Re=400$ are shown, while in Figure \ref{fig:NT_12_2} the results for $Re=1000$ are reported. 
A very good agreement with the numerical reference data available in the literature can be observed in both cases, as well as several secondary Taylor-G{\"o}rtler like vortex structures 
that can be recognized in the secondary planes. The observed structures are in good agreement with the ones reported in the literature, see e.g. \cite{Albensoeder2005} 

\subsubsection{3D Taylor-Green vortex}
In this section we investigate another typical benchmark, namely the 3D Taylor-Green vortex. In this test case a very simple initial solution degenerates quickly to a turbulent flow 
with very complex small scale structures. We take the initial condition as given in \cite{Oriol2015}   
\begin{eqnarray}
	&& \mathbf{v}(\xx,0) = \left( \, \, \sin(x)\cos(y)\cos(z),  \, \, -\cos(x)\sin(y)\cos(z), \, \, 0 \right),  \nonumber \\ 
	&& \rho(\xx,0) = \rho_0, \qquad 
	p(\xx,0) = p_0 + \frac{\rho_0}{16}\left(\cos(2x)+\cos(2y) \right)\left(\cos(2z)+2) \right), 
\label{eq:NT_13_1}
\end{eqnarray}
with $\Omega=[-\pi,\pi]^3$ and periodic boundary conditions everywhere. For this test case a very accurate reference solution is available and is given by a well-resolved DNS carried out by 
Brachet et al. \cite{Brachet1983} for an incompressible fluid. In order to approximate the incompressible case we take $\rho_0=1$ and $p_0=10^3$. We set $(p,p_\gamma)=(4,0)$, $t_{end}=10$  
and two different values of $\mu$ so that the Reynolds numbers under consideration are $Re=400$ and $Re=800$. The computational domain is covered with $\Ni=61824$ tetrahedra for $Re=400$ 
and $\Ni=120750$ elements for $Re=800$. 
\begin{figure}[ht!]
    \begin{center}
		\includegraphics[width=0.6\textwidth]{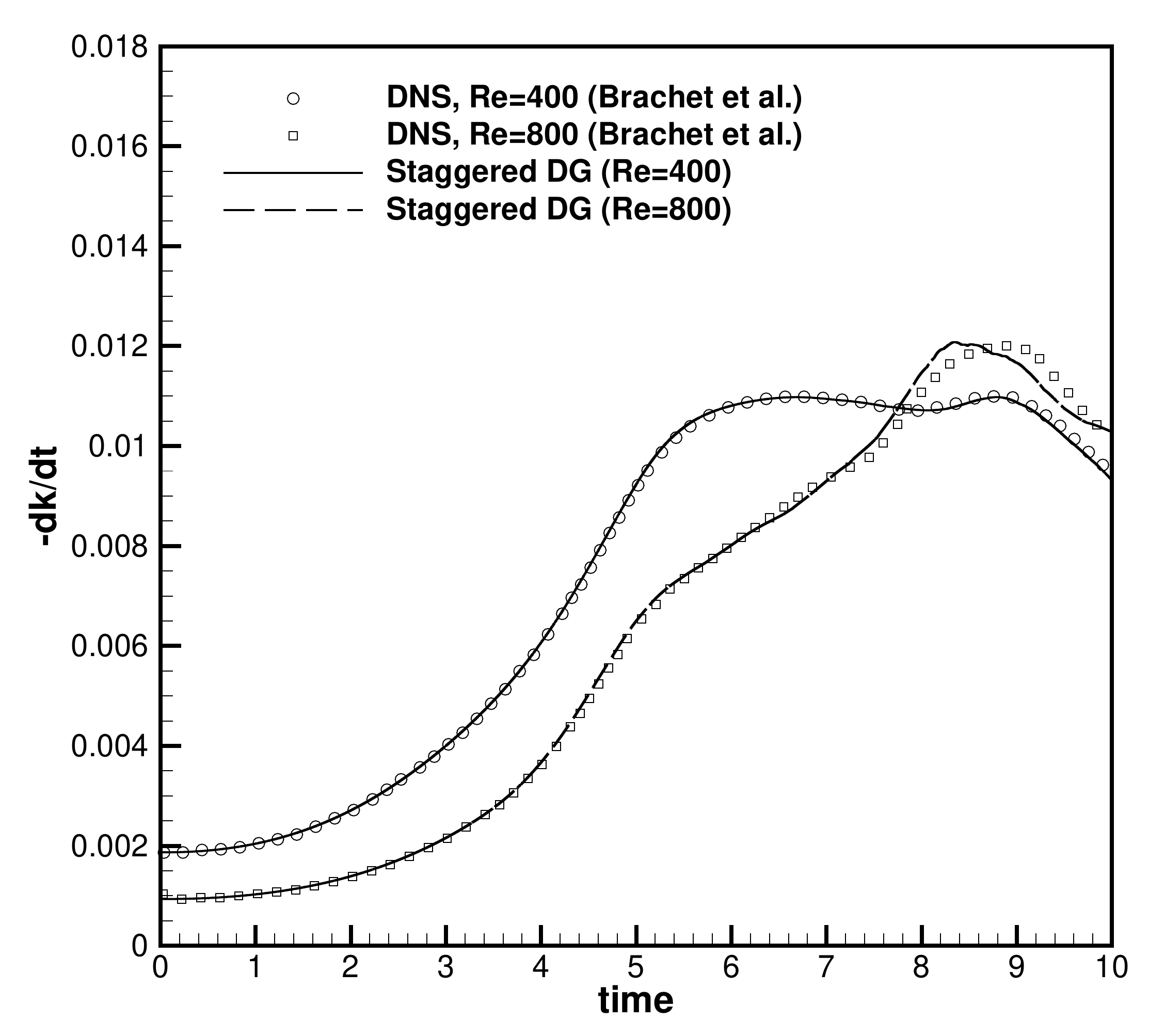}	
    \caption{Time evolution of the kinetic energy dissipation rate $-dk/dt$ for the 3D Taylor-Green vortex, compared with available DNS data of Brachet et al \cite{Brachet1983} 
		for $Re=400$ and $800$.}
    \label{fig:NT_13_1}
	\end{center}
\end{figure}

\begin{figure}[ht]
    \begin{center}
		\includegraphics[width=0.3\textwidth]{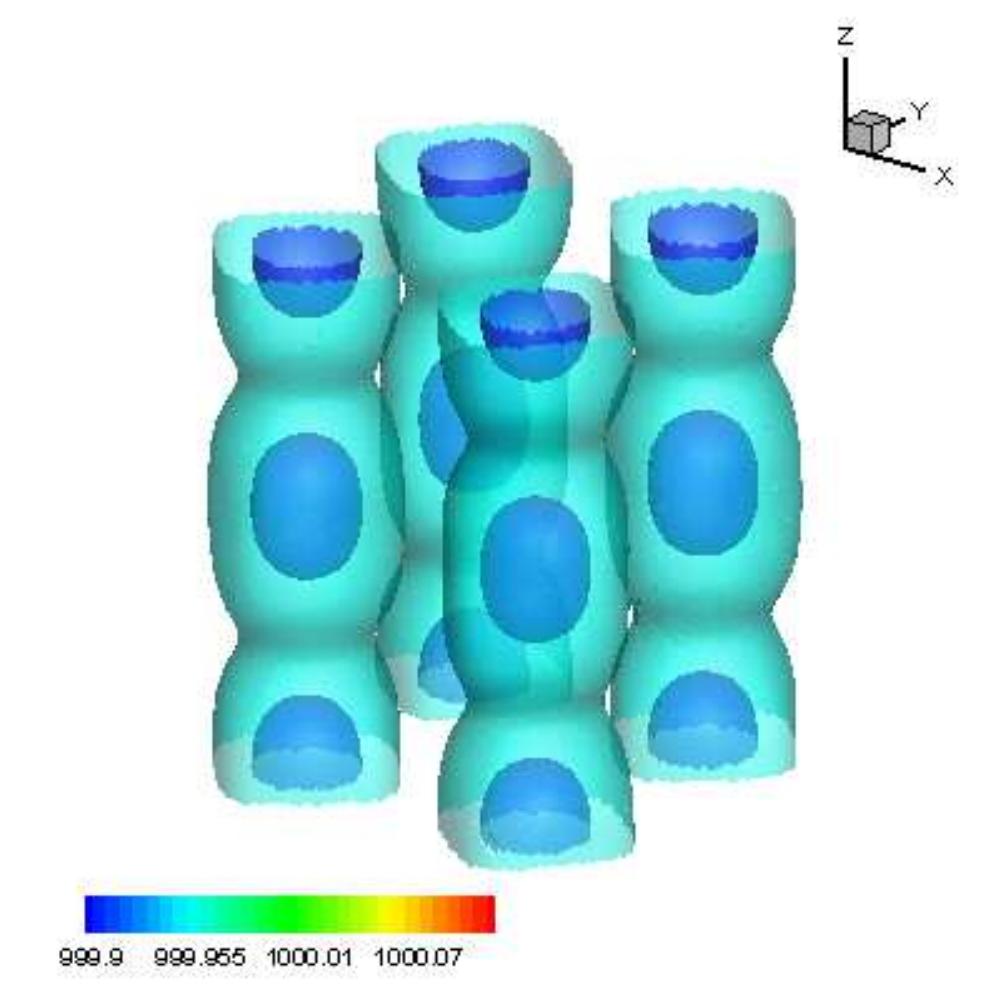}
		\includegraphics[width=0.3\textwidth]{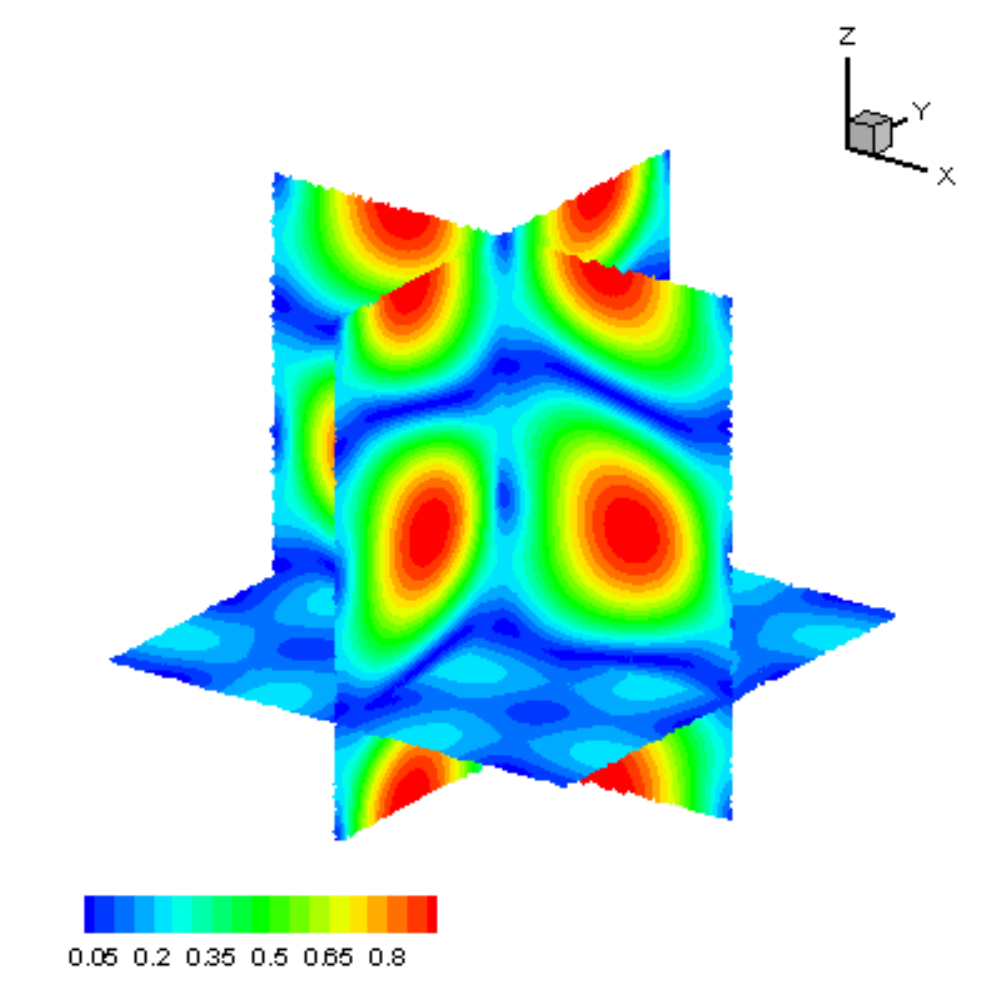}
		\includegraphics[width=0.3\textwidth]{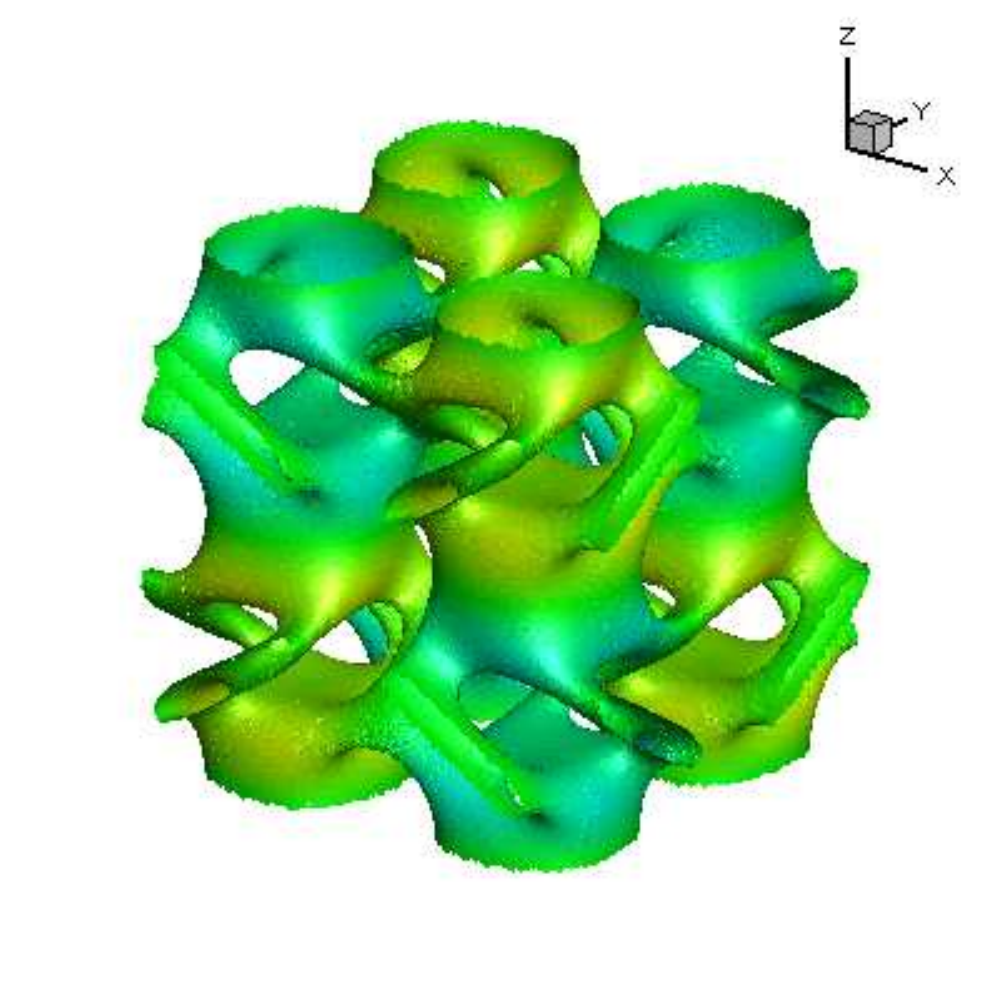} \\
		\includegraphics[width=0.3\textwidth]{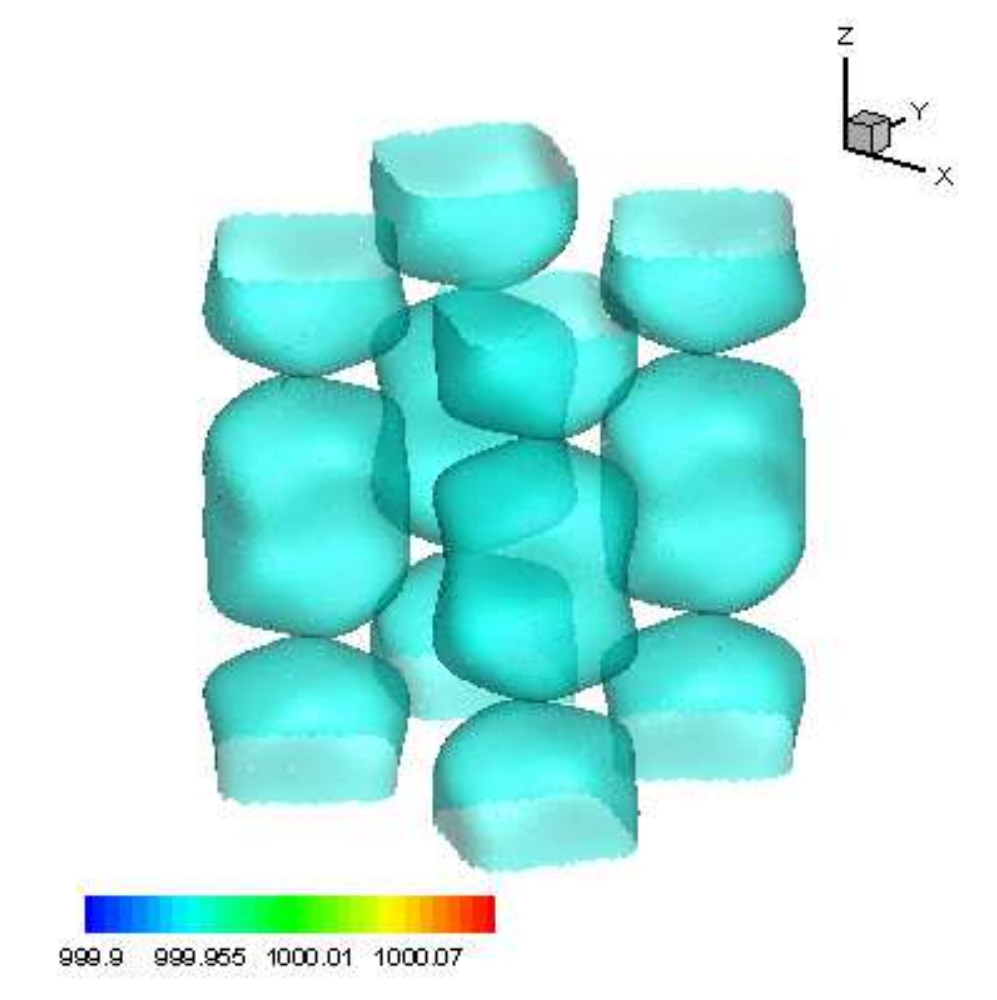}
		\includegraphics[width=0.3\textwidth]{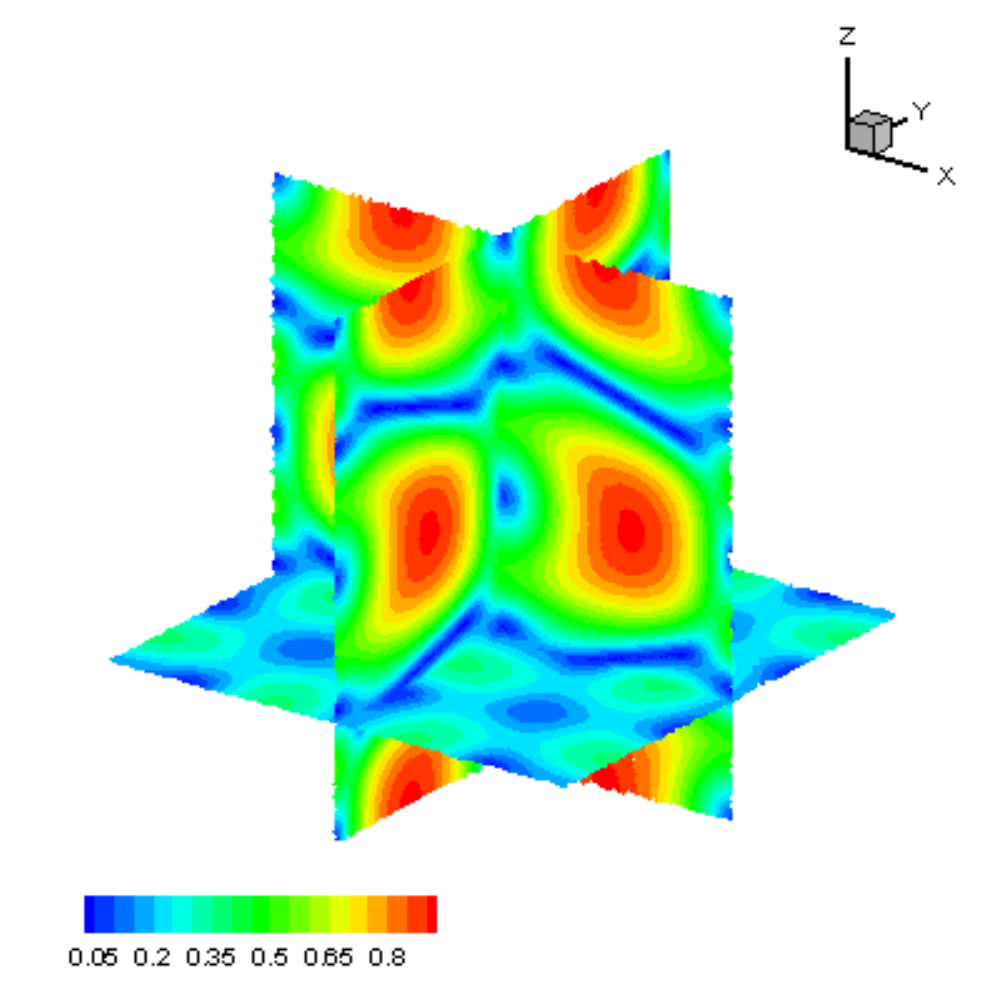}
		\includegraphics[width=0.3\textwidth]{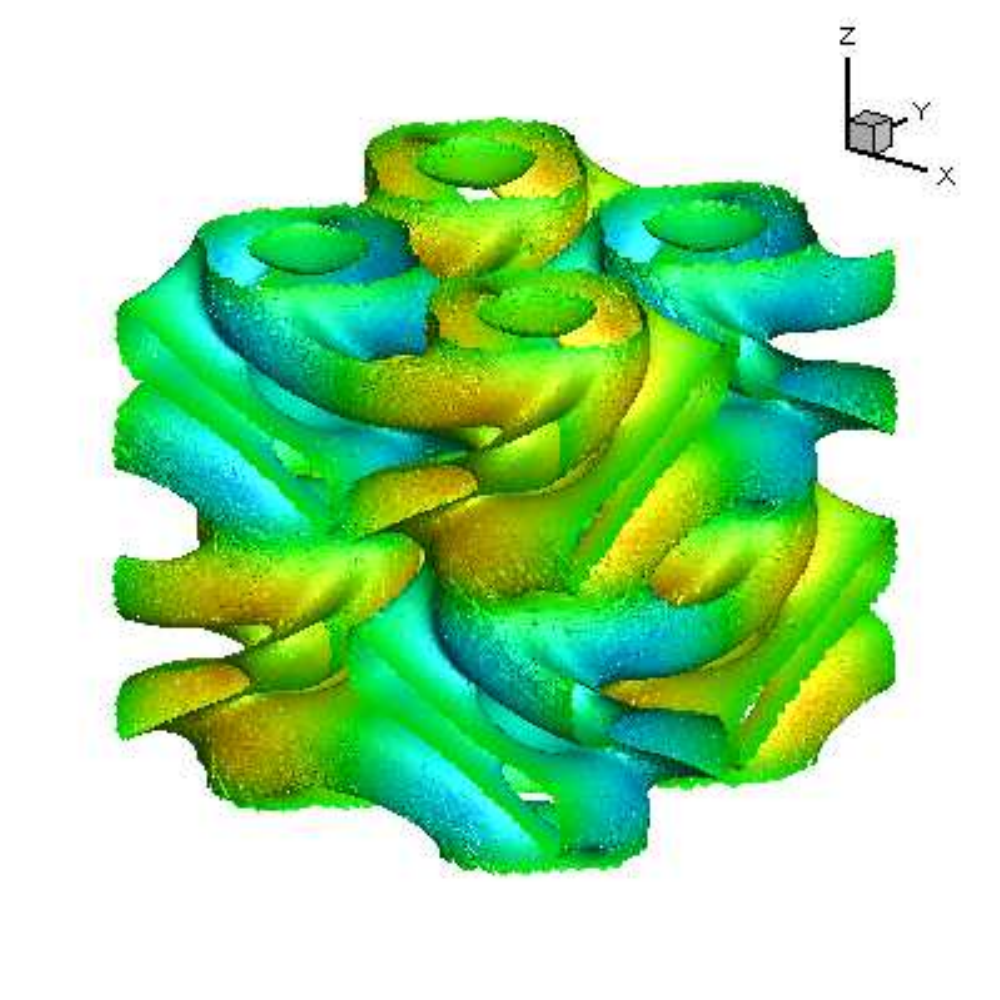} \\
			\includegraphics[width=0.3\textwidth]{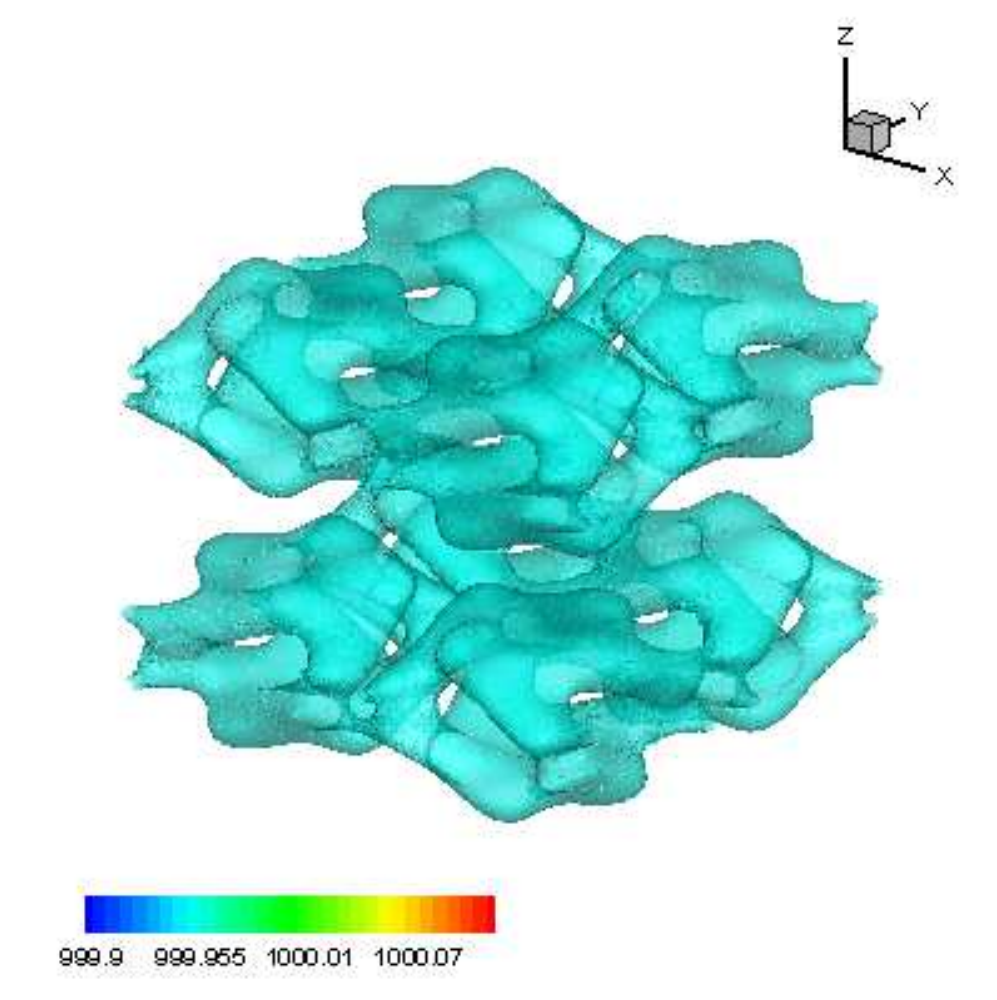}
		\includegraphics[width=0.3\textwidth]{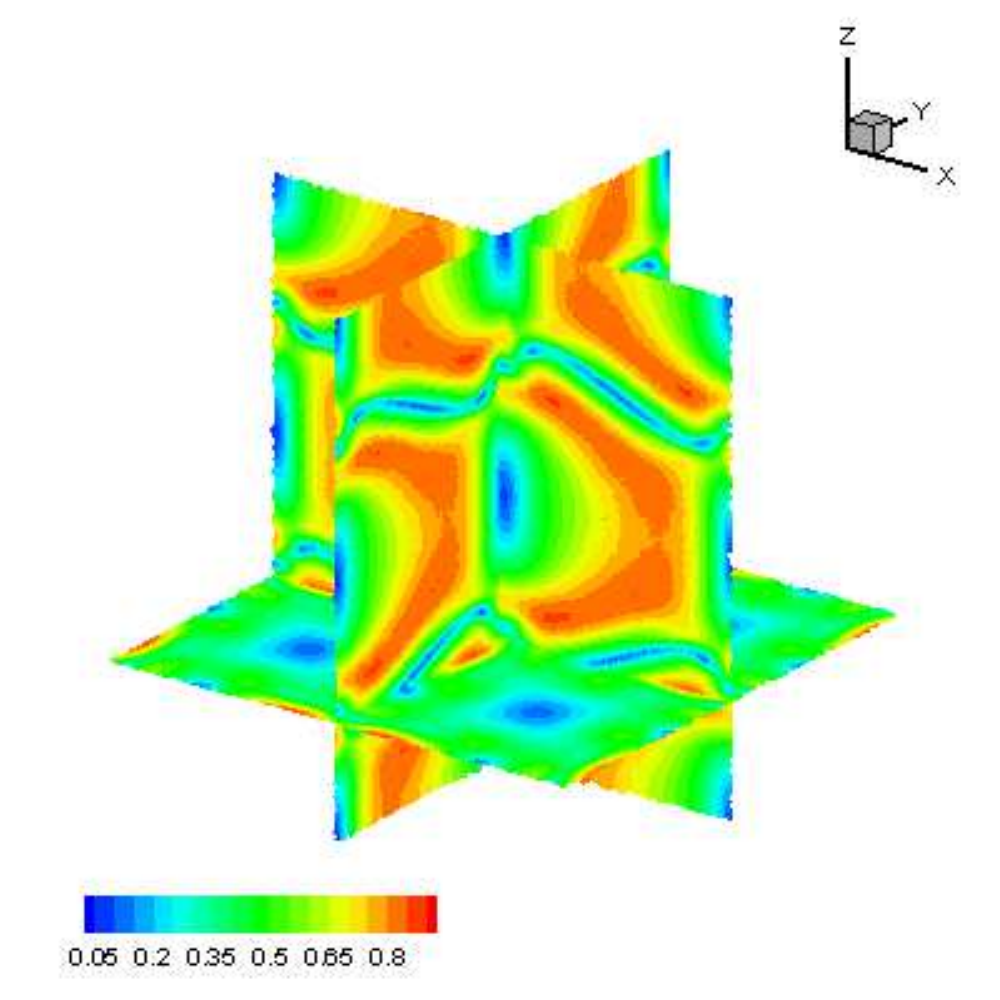}
		\includegraphics[width=0.3\textwidth]{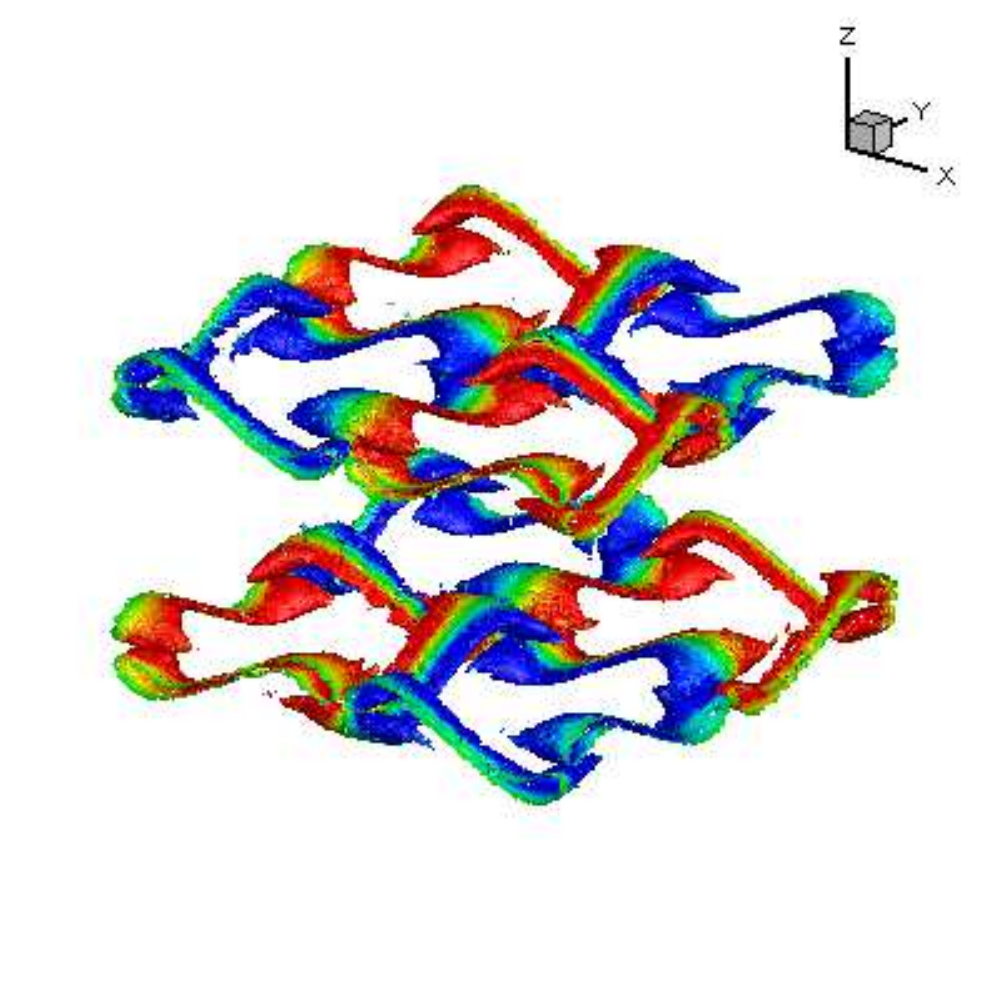} \\
			\includegraphics[width=0.3\textwidth]{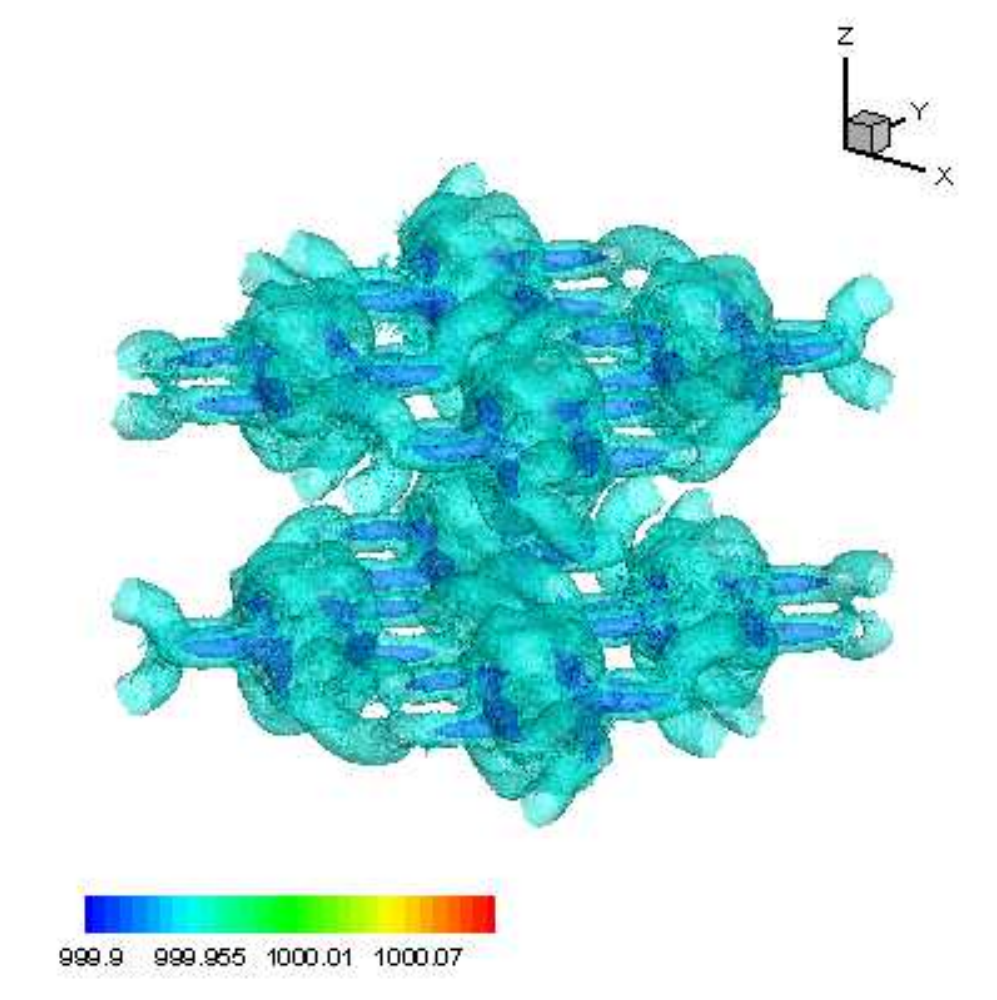}
		\includegraphics[width=0.3\textwidth]{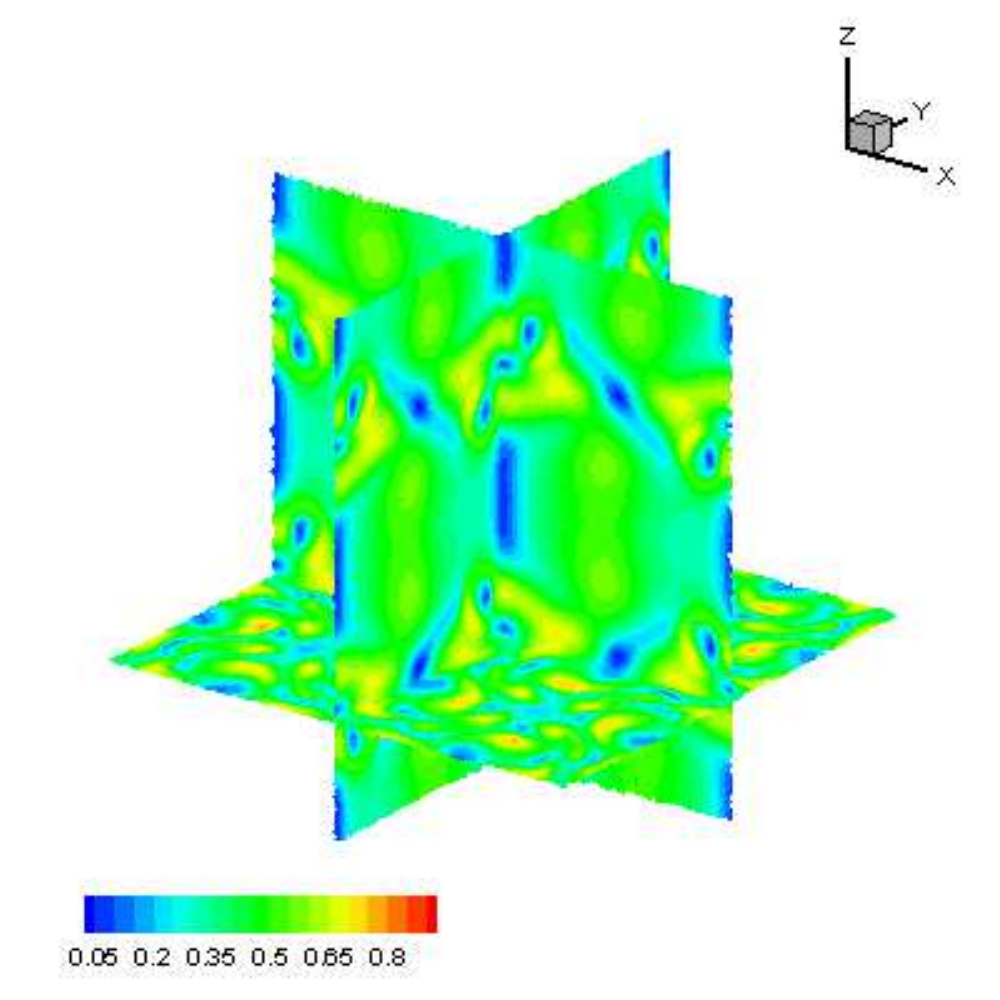}
		\includegraphics[width=0.3\textwidth]{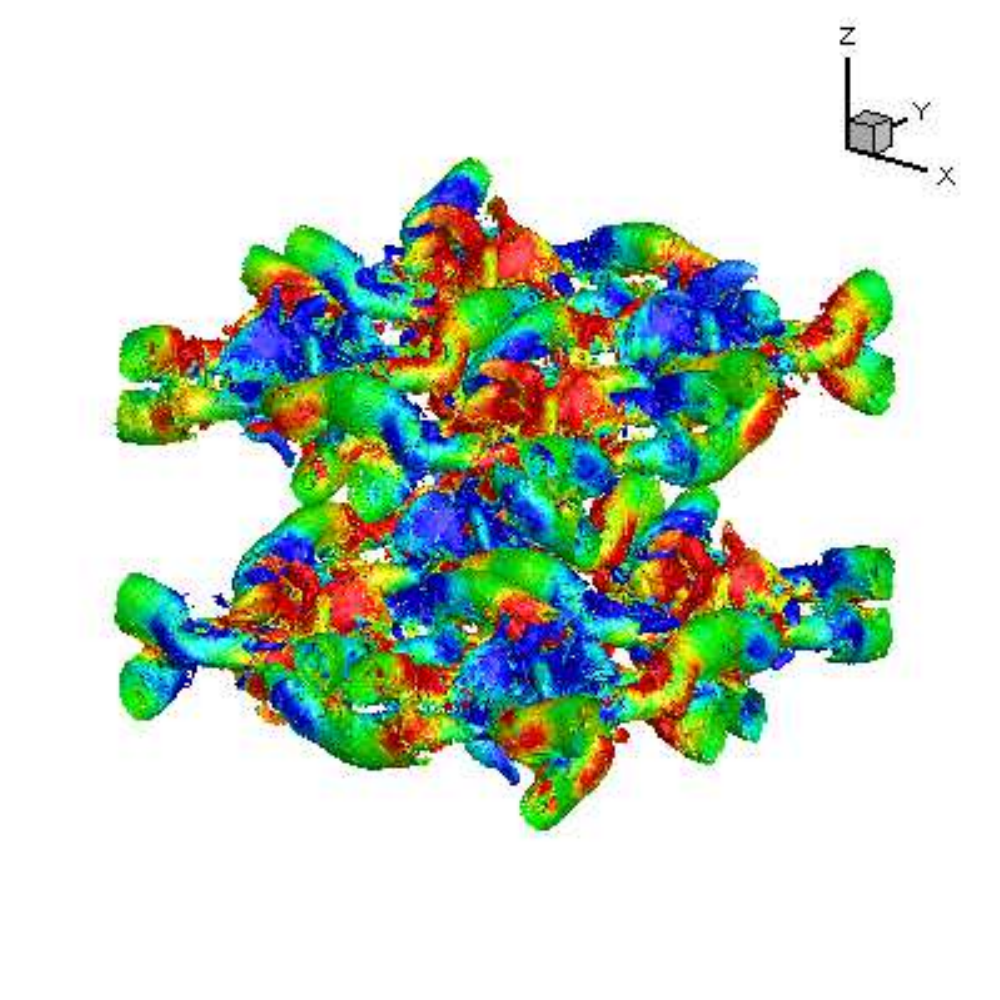} \\
    \caption{3D Taylor-Green vortex at $Re=400$. From left to right: Pressure isosurfaces, velocity magnitude and vorticity isosurfaces at times $t=0.9$ (top), $t=1.9$, $t=3.7$ and $t=8.2$ (bottom).}
    \label{fig:NT_13_2}
		\end{center}
\end{figure}
The resulting numerical solution for the pressure, velocity magnitude and vorticity colored by the helicity is reported  in Figure \ref{fig:NT_13_2} for $Re=400$ at several intermediate times, 
showing the development of the small-scale structures. In Figure \ref{fig:NT_13_1} the time series of the calculated total kinetic energy dissipation rates is compared with the DNS data of 
Brachet et al. \cite{Brachet1983}. Also for this rather complex test case, the resulting numerical results fit well with the reference solution and show only a slightly higher dissipation 
rate at $Re=800$. 
At $Re=400$ our numerical solution seems to represent the energy decay quite very well also at later times ($t>7$), where the smaller scales dissipate the major part of the energy, 
see e.g. \cite{Brachet1983} for a detailed discussion. 

\subsubsection{Flow around a sphere}
\begin{figure}[ht]
    \begin{center}
    \includegraphics[width=0.6\textwidth]{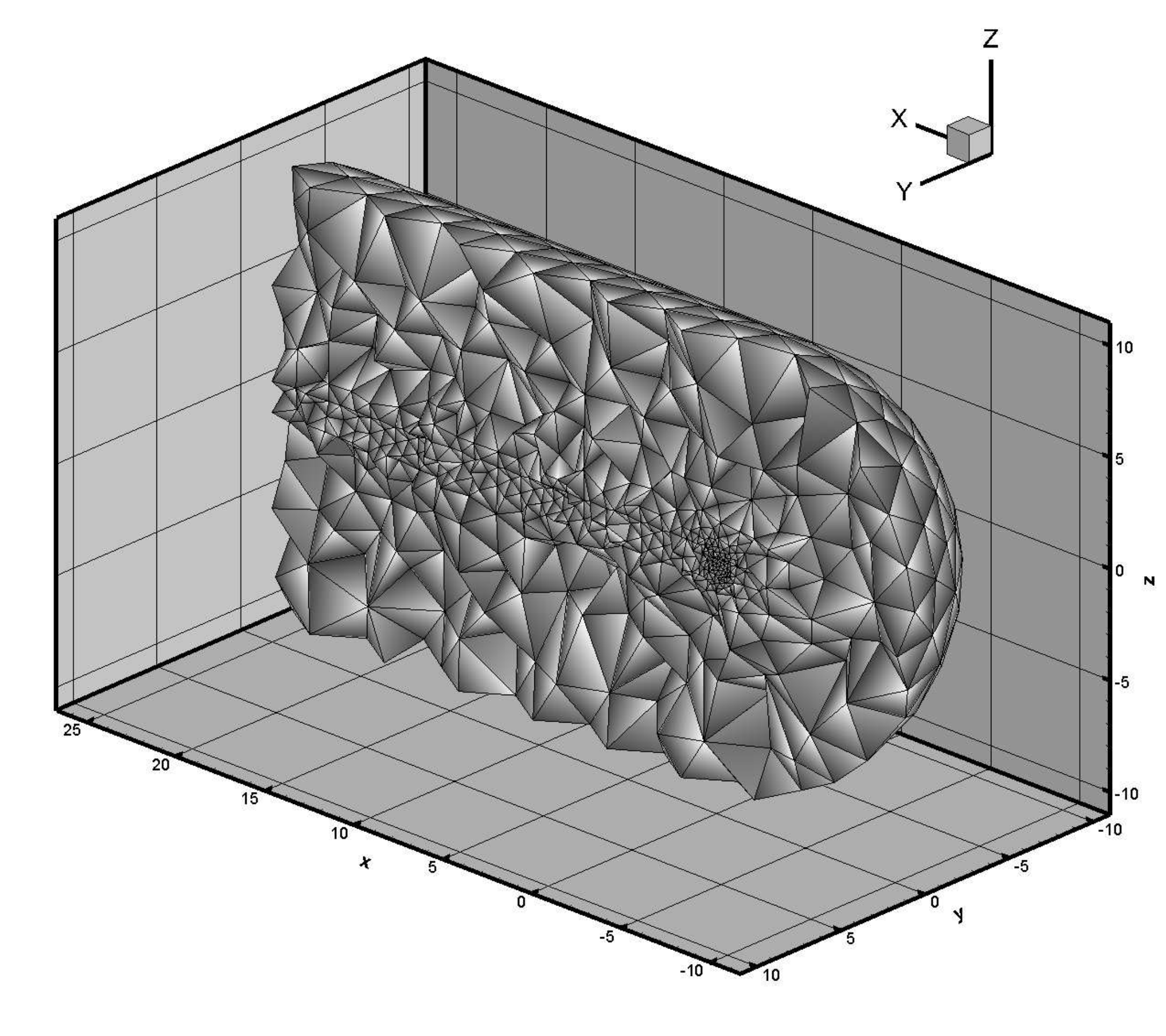}
    \caption{Flow around a sphere. Cut view of the computational domain with $\Ni=14403$.}
    \label{fig.Sph.1}
		\end{center}
\end{figure}

In this last numerical test we consider the flow around a sphere at a moderate Mach number. In particular we take as computational domain 
$\Omega=\mathcal{S}_{10} \cup \mathcal{C}_{10,25} \backslash \mathcal{S}_{0.5}$, where $\mathcal{S}_r$ is a generic sphere with center in the origin and radius $r$; 
$\mathcal{C}_{r,H}$ is a cylinder with circular basis on the $yz$-plane, radius $r$ and height $H$. We use a very coarse grid that is composed of a total 
number of only $\Ni=14403$ tetrahedra. A sketch of the grid is shown in Figure $\ref{fig.Sph.1}$. 

We start from an initially uniform flow with velocity $\mathbf{v}=(u_\infty,0,0)$, $u_\infty=0.5$, density $\rho = \rho_\infty=1$ and a pressure of $p=p_\infty=\frac{1}{\gamma}$, 
so that we obtain a moderate Mach number at infinity given by $M_{\infty}=0.5$. We impose $u_\infty$ on $\mathcal{S}_{10}\cap \{x\leq 0\}$ as boundary condition; transmissive 
boundary conditions on $\mathcal{C}_{10,25}$ and no-slip wall boundary conditions on $\mathcal{S}_{0.5}$. We use polynomial approximation degrees $(p,p_\gamma)=(3,0)$. The 
Reynolds number of the flow is chosen as $Re=300$, while the Prandtl number is set to $Pr=0.75$. The final simulation time is set to $t_{end}=200$. 

\begin{figure}[ht]
    \begin{center}
		\includegraphics[width=0.8\textwidth]{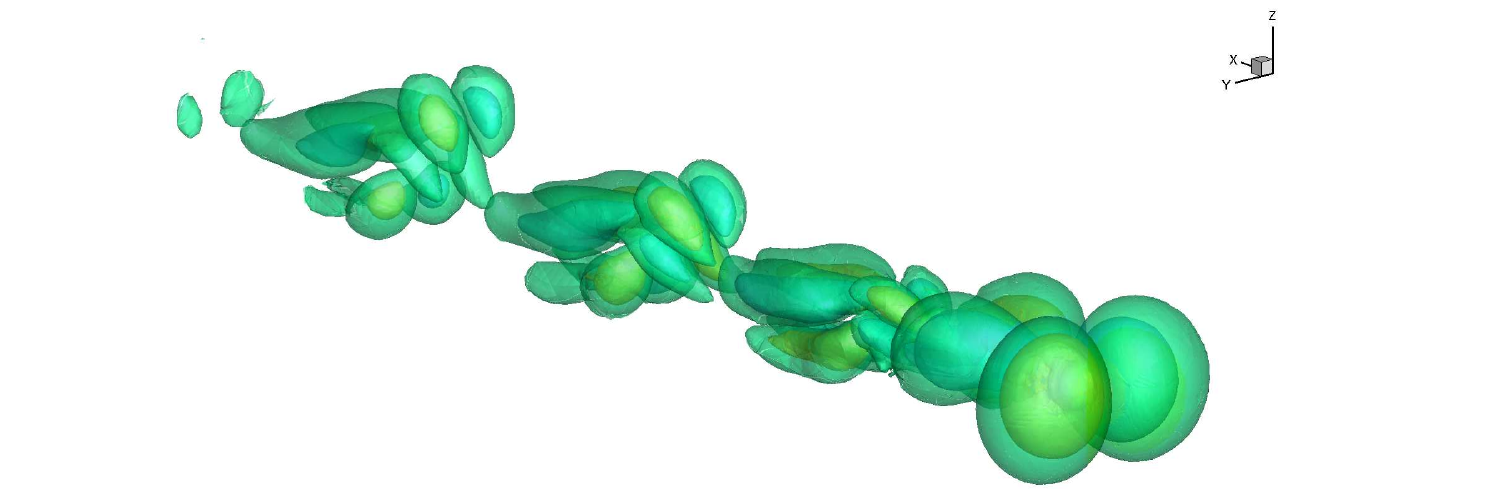}
    \includegraphics[width=0.8\textwidth]{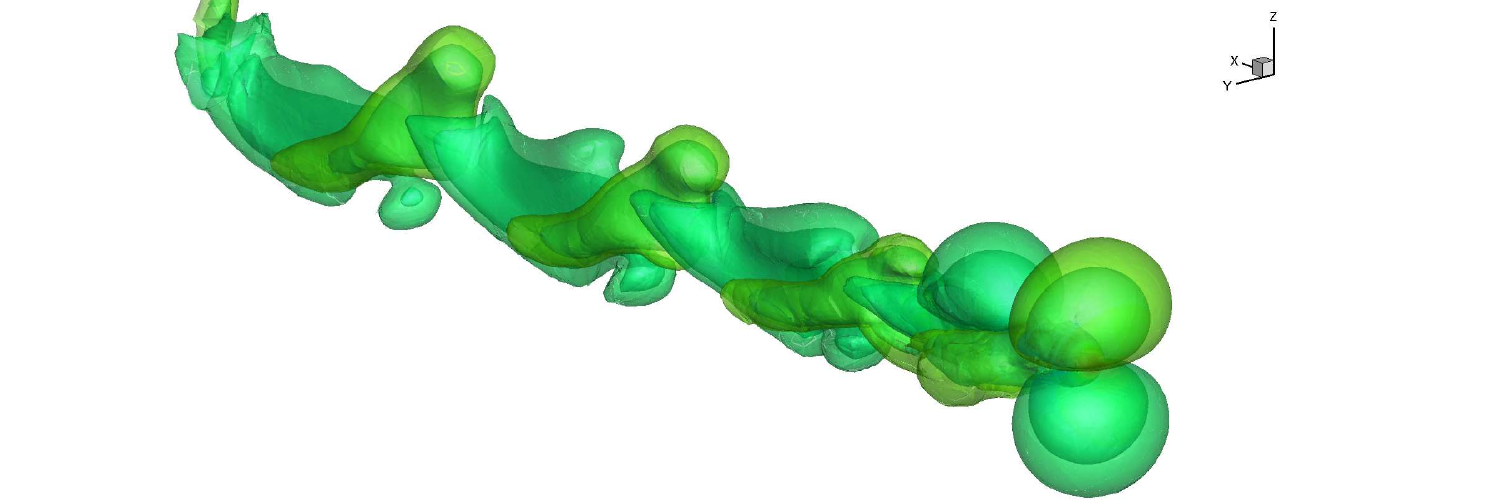}
		\includegraphics[width=0.8\textwidth]{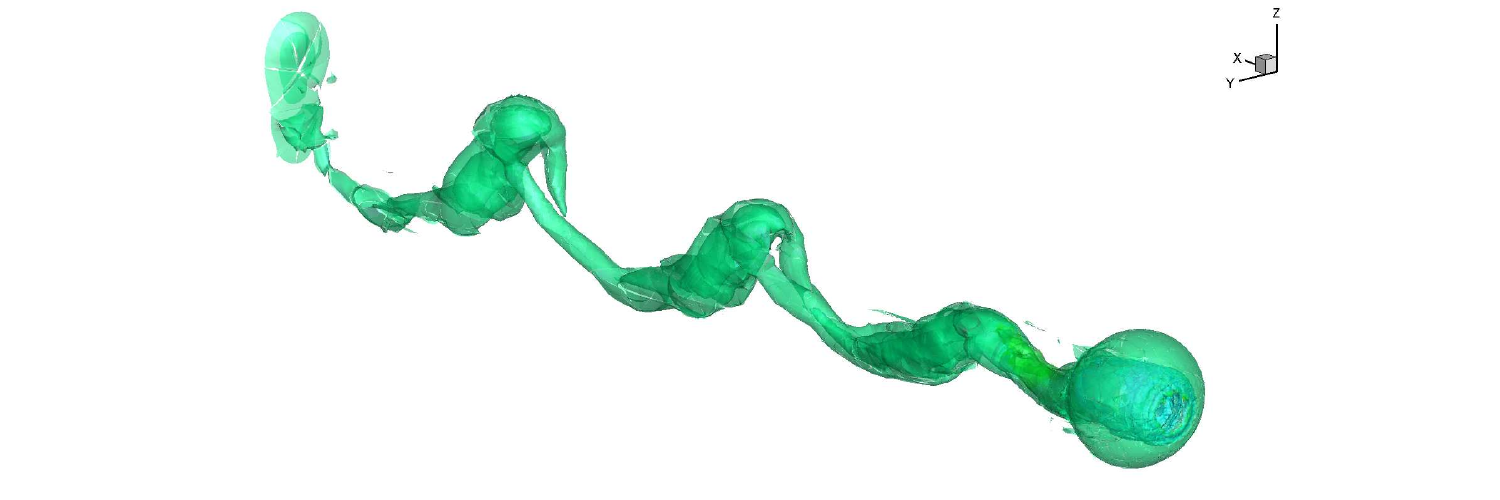}
    \caption{Isosurface plot of the span-wise velocities $v$, $w$ and density $\rho$ at $t=200$.}
    \label{fig.Sph.2}
		\end{center}
\end{figure}
\begin{figure}[ht]
    \begin{center}
		\includegraphics[width=0.49\textwidth]{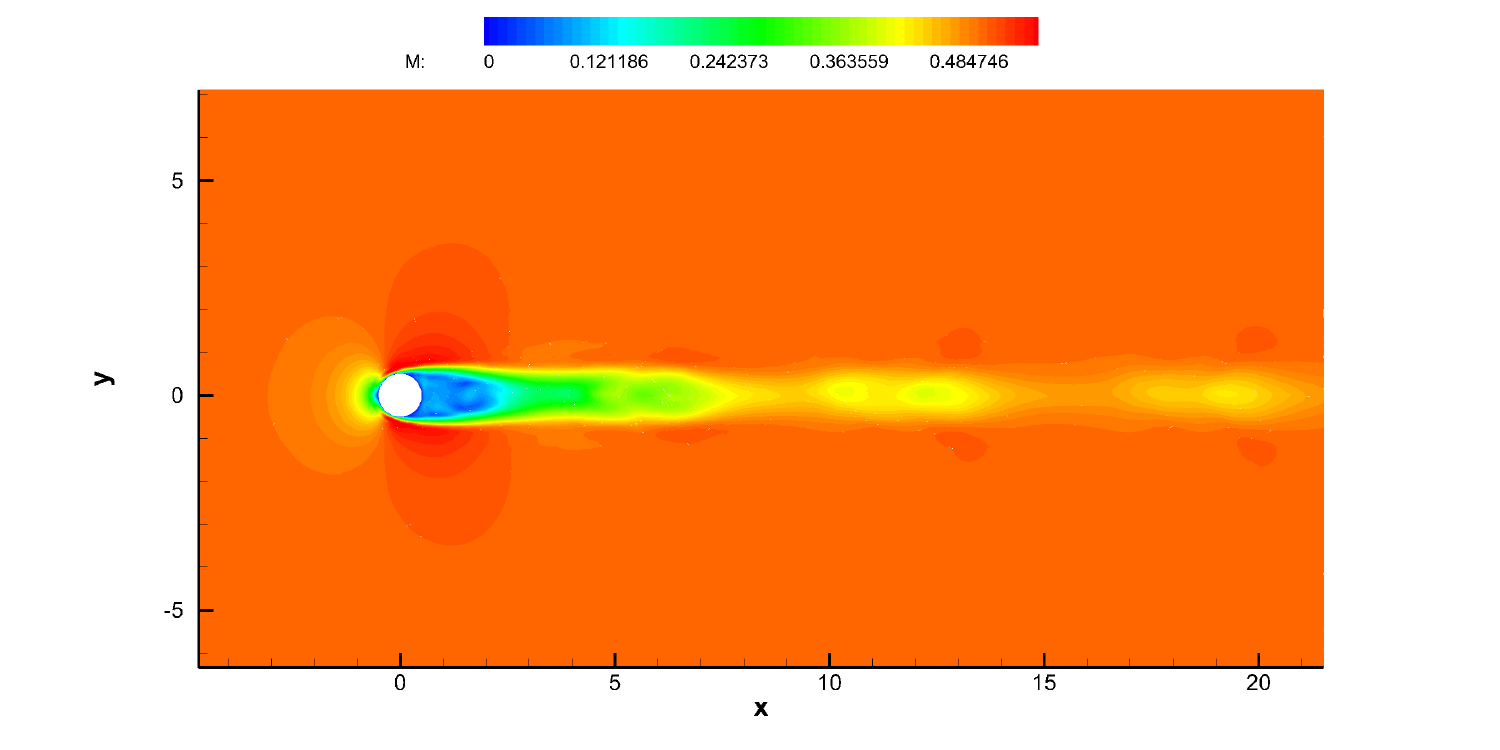}
    \includegraphics[width=0.49\textwidth]{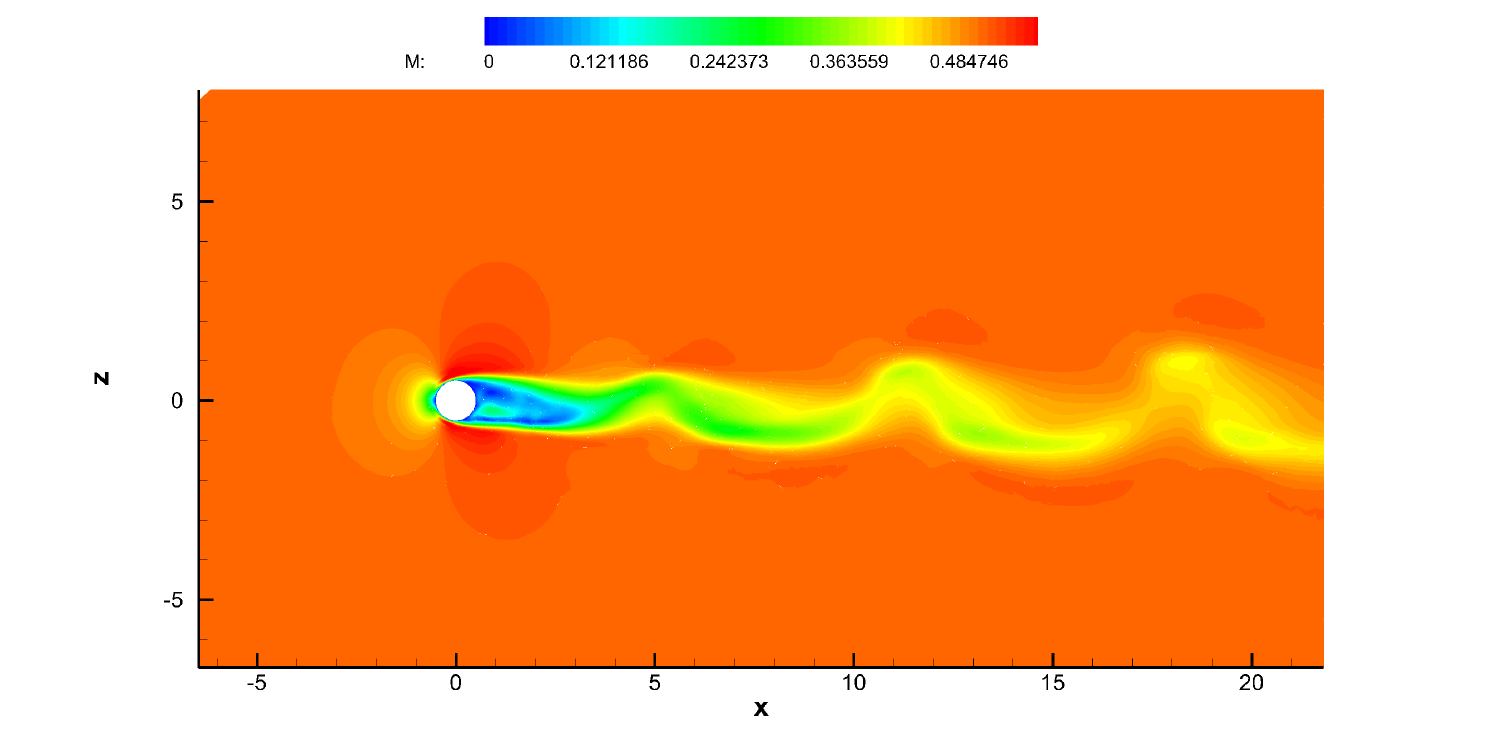}
    \caption{Mach number contours at $t=200$ in the plane $x-y$ (left) and $x-z$ (right).}
    \label{fig.Sph.3}
		\end{center}
\end{figure}


The isosurfaces of the span-wise velocities $v$ and $w$ as well as a density contour plot are reported in Figure \ref{fig.Sph.2} at time $t=200$ and show a 
very similar pattern compared to the one observed in experiments \cite{Sakamoto1990}. A plot of the Mach number magnitude is reported in Figure \ref{fig.Sph.3} 
on the two main planes $x-y$ and $x-z$ and show a very similar pattern compared to the ones observed in \cite{ADERNSE}. 
The frequency of the wake signal is obtained by computing the Fourier spectrum of $v$ and $w$ in the observation point $(5,0,0)$ and then taking the dominant frequency. 
The resulting Strouhal number for this simulation is $St=0.139$, which is close the the one obtained in \cite{ADERNSE} ($St=0.137$) and the one reported by 
Johnson and Patel in \cite{johnson1999}, who obtained a Strouhal number of $St=0.137$.

\section{Conclusions}
\label{sec.concl}
A novel semi-implicit space-time discontinuous Galerkin method for the solution of the two and three-dimensional compressible Navier-Stokes equations on unstructured staggered 
curved meshes was derived and discussed. The resulting formulation is in principle arbitrary high order accurate in space and time. The nonlinear terms become 
explicit thanks to the use of an outer Picard iteration. Experimentally, we recover the 
same good properties of the main system for the pressure as in the incompressible case for piecewise constant polynomials in time ($p_\gamma=0$). 
The general structure derived in this paper allows to treat also the viscous stress tensor and the heat flux implicitly, with a particular benefit in the case 
of the artificial viscosity method, where an explicit discretization of the artificial viscosity terms would lead to very small time steps. The proposed numerical 
method was tested on a large set of two and three-dimensional test problems, ranging from very low over moderate to very high Mach numbers. In all cases, the same 
discretization was used, hence the proposed pressure-based DG scheme is a genuine all Mach number flow solver. 

Future research will concern the application of the presented algorithm to geophysical and meteorological flows, e.g. for the simulation of meteorological flows 
in the alpine environment at the regional scale, where the geometry is rather complex due to the topography. Another important class of applications 
might be flow problems involving natural convection, where the typical Mach numbers are small or very small, but where density variations and heat transfer are 
very important. In this case the high order semi-implicit discretization will allow much larger time steps compared to a purely explicit DG method. 
The use of unstructured meshes should also be very well-suited for the discretization of very complex geometries as they appear for example in biological 
applications,  such as blood flow in the human cardiovascular system or air flow in the human respiratory system. 
The method should be in principle also well suited for industrial all Mach number flows in complex geometries. 

Finally, another potential future topic of research may be the extension and application of the new staggered semi-implicit space-time DG method to the 
magneto-hydrodynamics equations (MHD) and to the Baer-Nunziato model of compressible multi-phase flows.

\section*{Acknowledgments}
The research presented in this paper was funded by the European Research Council (ERC) under the European Union's Seventh Framework
Programme (FP7/2007-2013) within the research project \textit{STiMulUs}, ERC Grant agreement no. 278267.

The authors also acknowledge the Leibniz Rechenzentrum (LRZ) in Munich, Germany, for awarding access to the \textit{SuperMUC} 
supercomputer, as well as the support of the HLRS in Stuttgart, Germany, for awarding access to the \textit{Hazel Hen} supercomputer.

\bibliographystyle{elsarticle-num}
\bibliography{SIDG}

\end{document}